\documentclass[a4paper,10pt]{article}
\usepackage[latin1]{inputenc}
\usepackage[T1]{fontenc}
\usepackage[english]{babel}
\usepackage{amssymb, amsmath,amscd,amsxtra,amsfonts, euscript,latexsym, dsfont, mathrsfs, longtable,  textcomp, marvosym}%
\usepackage{graphics,graphicx}
\usepackage{epsf}
\usepackage{enumerate}
\usepackage{theorem}
\usepackage{fullpage}
\usepackage{url}

\usepackage{mathpazo}
\usepackage{hyperref}
\hypersetup{
    backref=true, 
    pagebackref=true,
    plainpages=false,
    bookmarks=true,         
    unicode=false,          
    pdftoolbar=true,        
    pdfmenubar=true,        
    pdffitwindow=true,      
    pdftitle={My title},    
    pdfauthor={Author},     
    pdfsubject={Subject},   
    pdfcreator={Creator},   
    pdfproducer={Producer}, 
    pdfkeywords={keywords}, 
    pdfnewwindow=true,      
    colorlinks=true,       
    linkcolor=blue,          
    citecolor=blue,        
    filecolor=blue,      
    urlcolor=magenta,          
    urlbordercolor=0 1 1
}
\newcommand{\dimo}{\noindent{\bf Proof: }}

\newcommand{\qed}{\thinspace\null\nobreak\hfill\hbox{\vbox{\kern-.2pt\hrule
height.2pt depth.2pt\kern-.2pt\kern-.2pt \hbox to2.5mm{\kern-.2pt\vrule
width.4pt \kern-.2pt\raise2.5mm\vbox to.2pt{}\lower0pt\vtop to.2pt{}\hfil
\kern-.2pt \vrule width.4pt\kern-.2pt}\kern-.2pt\kern-.2pt\hrule
height.2pt depth.2pt \kern-.2pt}}\par\medbreak}

\usepackage[normalem]{ulem}

\def\mdiv{\operatorname{div}}
\def\mrot{\operatorname{curl}}

\def\ds{\displaystyle}
\def\oF{\overline{F}}
\def\oW{\overline{W}}

\def\R{\mathbb{R}}

\def\BV{{\mathrm{BV}}}

\def\wpq#1#2{{{\mathrm W}^{#1,#2}}}
\def\lp#1{{\mathrm L}^{#1}}
\def\mdiv{\operatorname{div}}

\def\cont{{\mathrm{C}}}

\def\cW{{\mathcal W}}

\def\hausp#1{{\cal H}^{\mathrm{#1}}}

\def\haushn{\hausp{N-1}}
\def\one#1{\mathds{1}_{#1}}

\newtheorem{thm}{Theorem}[section]

\newtheorem{prop}[thm]{Proposition}

\newtheorem{algo}{Algorithm}
\newtheorem{claim}[thm]{Claim}
\newtheorem{oss}[thm]{Remark}

\newcommand{\boss}{\begin{oss}\rm }
\newcommand{\eoss}{\end{oss}}
\renewcommand{\sectionmark}[1]%
{\markright{\MakeUppercase{\thesection.\#1}}}
\makeatletter
\let\@fnsymbol\@arabic
\makeatother

\title{Phase-field approximations \\of the Willmore functional and flow}
\author{E. Bretin\thanks{Université de Lyon, CNRS UMR 5208, INSA de Lyon, Institut Camille Jordan, 20, avenue Albert Einstein, F-69621 Villeurbanne Cedex, France. Email: {\tt bretin@cmap.polytechnique.fr}} \and S. Masnou\thanks{Universit\'e de Lyon, CNRS UMR 5208, Universit\'e Lyon 1, Institut Camille Jordan, 43 boulevard du 11 novembre 1918, F-69622 Villeurbanne-Cedex, France. Email: {\tt masnou@math.univ-lyon1.fr}} \and E. Oudet\thanks{Lab. Jean Kuntzmann, Université Joseph Fourier, Tour IRMA, BP 53, 51, rue des Mathématiques, F-38041 Grenoble Cedex 9, France.  Email: {\tt edouard.oudet@imag.fr}}}

\begin{document}
\maketitle

\begin{abstract}
We discuss in this paper phase-field approximations of the Willmore functional and the associated $\lp{2}$-flow. After recollecting known results on the approximation of the Willmore energy and its $\lp{1}$ relaxation, we derive the expression of the flows associated with various approximations, and we show their behavior by formal arguments based on matched asymptotic expansions. We introduce an accurate numerical scheme, whose local convergence can be proved, to describe with more details the behavior of two flows, the classical and the flow associated with an approximation model due to Mugnai. We propose a series of numerical simulations in 2D and 3D to illustrate their behavior in both smooth and singular situations.
\end{abstract}
\section{Introduction}
Phase-field approximations of the Willmore functional have raised quite a lot of interest in recent years, both from the theoretical and the numerical viewpoints. In particular, attention has been given to understanding the continuous and numerical approximations of both smooth and singular sets with finite relaxed Willmore energy. Various approximation models have been proposed so far, whose properties are known only partially. Our main motivation in this paper is a better understanding of these models, and more precisely:
\begin{enumerate}
\item Exhibiting algebraic differences/similarities between the various approximations;
\item Deriving the $\lp{2}$-flows associated with these models;
\item Studying the asymptotic behavior of the flows, at least in smooth situations;
\item Simulating numerically these flows, and observing whether and how singularities may appear.
\end{enumerate}
We focus on four models due, respectively, to De Giorgi, Bellettini, and Paolini~\cite{degiorgi-a,belpao93}, Bellettini~\cite{Bellettini1997}, Mugnai~\cite{mugnai2010}, and Esedoglu, R{\"a}tz, and R{\"o}ger~\cite{Esedoglu_12}. The paper is organized as follows: Section~\ref{sec:intro} is an introductory section where we collect known results on the diffuse approximation of the perimeter, the diffuse approximation of the Willmore energy, and the critical issue of approximating singular sets with finite relaxed Willmore energy. We also recall the definitions of the above mentioned approximations. In Section~\ref{sec:remarks}, we make new observations on the differences between these different diffuse energies. Section~\ref{sec:flow} is devoted to the derivation of the $\lp{2}$-flows associated with, respectively, Bellettini's, Mugnai's, and Esedoglu-R{\"a}tz-R{\"o}ger's models (actually a variant of the latter), and, for every flow, we use the formal method of matched asymptotic expansions to derive the asymptotic velocity of the limit interface as the diffuse approximation becomes asymptotically sharp. We show in particular that, in dimensions $2$ and $3$ for all flows, and in any dimension for some of them, they correspond asymptotically to the continuous Willmore flow as long as the interface is smooth. In Section~\ref{sec:num}, we focus on the numerical simulation of De Giorgi-Bellettini-Paolini's flow (which we shall refer to as the {\it classical} flow) and Mugnai's flow, and we propose a fixed-point algorithm whose local convergence can be proved. We illustrate with various numerical examples the behavior of both flows in space dimensions $2$ and $3$, both in smooth and singular situations. We show in particular that our scheme can capture with good accuracy well-known singular configurations yielded by the classical flow, and that these configurations evolve as if the parametric Willmore flow were used. We also illustrate with several simulations that, in contrast, Mugnai's flow prevents the creation of singularities.
\section{What is known?}\label{sec:intro}
\subsection{Genesis : the van der Waals-Cahn-Hilliard interface model and the diffuse approximation of perimeter}
In his 1893 paper on the thermodynamic theory of capillary (see an English translation, with interesting comments, in~\cite{rowlinson}), van der Waals studied the free energy of a liquid-gas interface. Arguing that the density of molecules at the interface can be modelled as a continuous function of space $u$, he used thermodynamic and variational arguments to derive an expression of the free energy, in a small volume $V$ enclosing the interface, as $\int_V(f_0(u)+\lambda|\nabla u|^2)dx$, where $f_0(u)$ denotes the energy of a homogeneous phase at density $u$ and $\lambda$ is the capillarity coefficient. The same expression was derived by Cahn and Hilliard in 1958 in their paper~\cite{CAHN:1958} on the interface energy, to a first approximation, of a binary alloy with $u$ denoting the mole fraction of one component. Cahn and Hilliard argued that both terms in the energy have opposite contributions: if the transition layer's size increases, then the gradient term diminishes, but this is possible only by introducing more material of nonequilibrium composition, and thus at the expense of increasing $\int_Vf_0(u)dx$. Rescaling the energy, and changing the notations in the obvious way, yields the general form
\begin{equation}
F_\varepsilon(u)=\int_V(\frac\varepsilon 2|\nabla u|^2+\frac{W(u)}{\varepsilon})dx\label{vdW}
\end{equation}
In the original papers of van der Waals, Cahn and Hilliard, $f_0$ was a smooth double-well function, yet with a slope between the local minima. For simplicity, since it does not modify the mathematical analysis, $W$ will denote in the sequel a smooth double-well function with no slope (we will take in general $W(s)=\frac 1 2s^2(1-s)^2$).
\par Two equations are usually associated with the van der Waals-Cahn-Hilliard energy, and will be used in this paper: the Allen-Cahn and the Cahn-Hilliard equations. The evolution Allen-Cahn equation is the $\lp{2}$-gradient descent associated with \eqref{vdW} and is written
$$u_t=\frac{W'(u)}\varepsilon-\varepsilon\Delta u.$$
We shall also refer to the stationary Allen-Cahn equation 
$$W'(u)-\varepsilon^2\Delta u=0.$$
The Cahn-Hillard evolution equation is derived in a different manner: from a mathematical viewpoint, it is the $H^{-1}$-gradient flow associated with the van der Waals-Cahn-Hilliard energy~\cite{Cahn19941045,MR1772733}. The physical derivation of the equation is also instructive~\cite{cahnspin,CAHN:1958}: since $\nabla_u F_\varepsilon(u)=-\varepsilon\Delta u+\frac{W'(u)}{\varepsilon}$ quantifies how the energy changes when molecules change position, it coincides with the chemical potential $\mu$. Fick's first law states that the flux of particles is proportional to the gradient of $\mu$, i.e. $J=-\alpha\nabla \mu$. Finally, the conservation law $u_t+\mdiv{J}=0$ yields the Cahn-Hilliard evolution equation
$$u_t=\alpha\Delta(-\varepsilon\Delta u+\frac{W'(u)}{\varepsilon}).$$
To summarize, the Allen-Cahn equation describes the motion of phase boundaries driven by surface tension, whereas the Cahn-Hilliard equation is a conservation law that characterizes the motion induced by the chemical potential, which is the gradient of the surface tension.
\par Let us now recall the asymptotic behavior of $F_\varepsilon(u)$, as $\varepsilon\to 0^+$, that has been exhibited by Modica and Mortola~\cite{Modica1977a} following a conjecture of De Giorgi. We first fix some notations. Let $\Omega\subset\R^n$ be open, bounded and with Lipschitz boundary. $W\in\cont^{3}(\R,\R^+)$ is a double-well potential with two equal minima (in the sequel we will work, unless specified, with $W(s)=\frac 1 2s^2(1-s)^2$).
Modica and Mortola have shown that the $\Gamma$-limit in $\lp{1}(\Omega)$ of the family of functionals
$$P_\varepsilon(u)=\left\{\begin{array}{ll}\ds\int_\Omega \left(\frac\varepsilon 2|\nabla u|^2+\frac{W(u)}\varepsilon \right) dx&\mbox{if }u\in\wpq{1}{2}(\Omega)\\
+\infty&\mbox{otherwise in }\lp{1}(\Omega)\end{array}\right.$$
is $c_0 P(u)$ where
$$P(u)=\left\{\begin{array}{ll}\ds |Du|(\Omega)&\mbox{if }u\in\BV(\Omega,\{0,1\})\\
+\infty&\mbox{otherwise in }\lp{1}(\Omega)\end{array}\right.$$
and $c_0=\int_{0}^1\sqrt{2W(s)}ds$. In particular, if  $E\subset\R^n$ has finite perimeter in $\Omega$ and $u:=\one{E}\in\BV(\Omega,\{0,1\})$, then one can build a sequence of functions $(u_\varepsilon)\in\wpq{1}{2}(\Omega)$ such that $u_\varepsilon\to u\in\lp{1}(\Omega)$ and $P_\varepsilon(u_\varepsilon)\to c_0|Du|(\Omega)=c_0P(E,\Omega)$ with $P(E,\Omega)$ the perimeter of $E$ in $\Omega$. 
\par To prove it, it is enough by density to restrict to smooth sets. Being $E$ smooth, a good approximating sequence is given by $u_\varepsilon=q(\frac{d(x)}\varepsilon)$ (actually a variant of this expression, but we shall skip the details for the moment) where $d$ is the signed distance function at $\partial E$, i.e. $d(x)=-d(x,\partial E)$ if $x\in E$ and $d(x,\partial E)$ else, and $q(t)=\frac{1-\tanh(t)}2$ is the unique decreasing minimizer of 
\begin{equation}
\int_\R(\frac{|\varphi'(t)|^2}2+W(\varphi(t)))dt\label{transvers}
\end{equation}
under the assumptions $\lim_{t\to-\infty}\varphi(t)=1$,  $\lim_{t\to \infty}\varphi(t)=0$ and $\varphi(0)=\frac 1 2$. In other words, the approximation of $u=\one{E}$ is done by a suitable rescaling of the level lines of the distance function to $\partial E$. Such rescaling is optimal, in the sense that it minimizes the transversal energy \eqref{transvers} and forces the concentration as $\varepsilon\to 0^+$.
\par Observe now that 
$$\int_\Omega \left(\frac\varepsilon 2|\nabla u|^2+\frac{W(u)}\varepsilon \right) dx\geq \int_\Omega \sqrt{2}|\nabla u|\sqrt{W(u)}$$
and the equality holds if $\frac\varepsilon 2|\nabla u|^2=\frac{W(u)}\varepsilon$. Therefore, by lower semicontinuity arguments, the quality of an approximation depends on the so-called {\it discrepancy measure}
$$\xi_\varepsilon= (\frac\varepsilon 2|\nabla u|^2-\frac{W(u)}\varepsilon){\cal L}^2,$$
that will play an even more essential role for some diffuse approximations of the Willmore functional.

\subsection{De Giorgi-Bellettini-Paolini's approximation of the Willmore energy}
Based on a conjecture of De Giorgi~\cite{degiorgi-a}, several authors~\cite{belpao93,tone02,belmug04,moser,Roger_schatzle_2006,naga-tone} have investigated the diffuse approximation of the Willmore functional, which is for a set $E\subset\R^N$ with smooth boundary in~$\Omega$:
$$W(E,\Omega)= \frac 1 2\int_{\partial E\cap\Omega} |H_{\partial E}(x)|^2~d\haushn$$
where $H_{\partial E}(x)$ is the classical mean curvature vector at $x\in\partial E$. The approximation functionals are defined as
$$\cW_\varepsilon(u)=\left\{\begin{array}{ll}\ds \frac 1{2\varepsilon}\int_\Omega \left( \varepsilon\Delta u-\frac{W'(u)}{\varepsilon} \right)^2dx&\mbox{if }u\in\lp{1}(\Omega)\cap\wpq{2}{2}(\Omega)\\
+\infty&\mbox{otherwise in }\lp{1}(\Omega)\end{array}\right.$$
Introduced by Bellettini and Paolini in~\cite{belpao93}, they differ from the original De Giorgi's conjecture in the sense that the perimeter is not explicitly encoded in the expression. They have however the advantage to be directly related to the Cahn-Hilliard equation, whose good properties~\cite{chen-96} play a key role in the approximation. In the sequel, we shall refer to these functionals as the {\it classical} approximation model.
\par The reason why $\varepsilon\Delta u-\frac{W'(u)}{\varepsilon}$ is related to the mean curvature can be simply understood at a formal level: it suffices to observe that the mean curvature of a smooth surface is associated with the first variation of its area, and that $-\varepsilon\Delta u+\frac{W'(u)}{\varepsilon}$ is the $\lp{2}$ gradient of $\frac\varepsilon 2|\nabla u|^2+\frac{W(u)}\varepsilon$ that appears in the approximation of the surface area.
\par The results on the asymptotic behavior of $\cW_\varepsilon$ as $\varepsilon\to 0^+$ have started with the proof by Bellettini and Paolini~\cite{belpao93} of a $\Gamma-\limsup$ property, i.e. the Willmore energy of a smooth hypersurface $E$ is the limit of $\cW_\varepsilon(u_\varepsilon)$, up to a multiplicative constant, where $u_\varepsilon$ is defined exactly as for the approximation of the perimeter.
\par The $\Gamma-\liminf$ property is much harder to prove. The contributions on this point~\cite{tone02,belmug04,moser,Roger_schatzle_2006,naga-tone} culminated with the proof by R{\"o}ger and Sch{\"a}tzle~\cite{Roger_schatzle_2006} in space dimensions $N=2,3$ and, independently, by Nagase and Tonegawa~\cite{naga-tone} in dimension $N=2$, that the result holds true for smooth sets. More precisely, given $u=\one{E}$ with $E\in\cont^{2}(\Omega)$, and $u_\varepsilon$ converging to $u$ in $\lp{1}$ with a uniform control of $P_\varepsilon(u_\varepsilon)$, then
$$c_0 W(E,\Omega)\leq\liminf_{\varepsilon\to 0^+} \cW_\varepsilon(u_\varepsilon).$$
The proof is based on a careful control of the discrepancy measure $\xi_\varepsilon= (\frac\varepsilon 2|\nabla u|^2-\frac{W(u)}\varepsilon){\cal L}^2$ that guarantees good concentration properties, i.e. the varifolds $v_\varepsilon=|\nabla u_\varepsilon|{\cal L}^2\otimes \delta_{\nabla u_\varepsilon^\perp}$ (whose mass is related naturally to the variations of the approximating functions $u_\varepsilon$) concentrate to a limit integer varifold that has generalized mean curvature in $\lp{2}$ and that is supported on a supset of $\partial E$. Then, in R{\"o}ger and Sch{\"a}tzle's proof, a lower semicontinuity argument and the locality of integer varifolds' mean curvature yields the result. It holds in dimensions $N=2,3$ at most due to dimensional requirements for Sobolev embeddings and for the control of singular terms used in the proof. The result in higher dimension is still open. 

\par What about unsmooth sets? Can the approximation results be extended to the relaxed Willmore functional? The answer is negative in general, as discussed below.
\subsection{The approximation does not hold in general for unsmooth limit sets}\label{sec:approx}
\begin{figure}[!ht]
\begin{center}
  {\begin{minipage}[c]{8cm}
    \includegraphics[width=7cm]{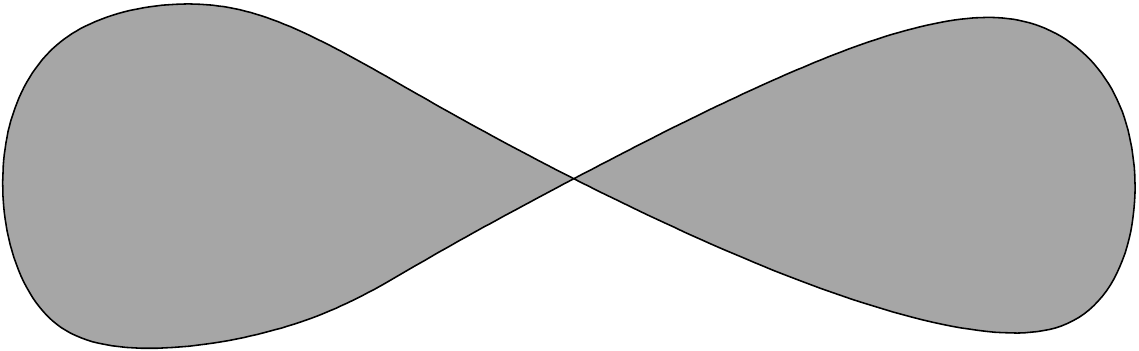}
    \end{minipage}}\vspace*{0.2cm}\\
\includegraphics[width=3.5cm]{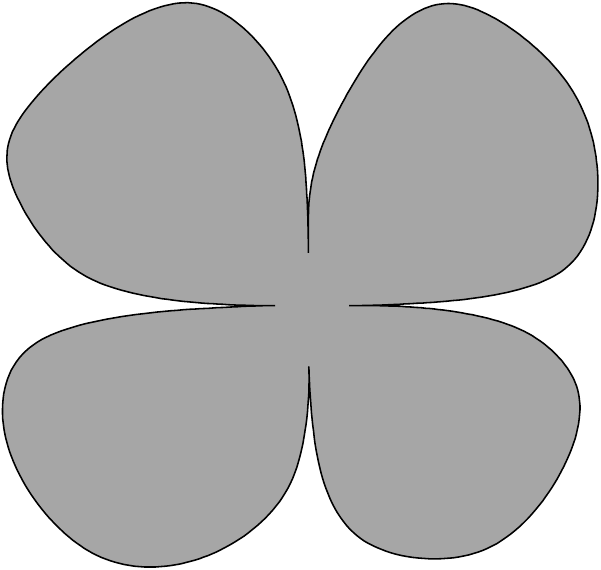}\qquad\includegraphics[width=3.5cm]{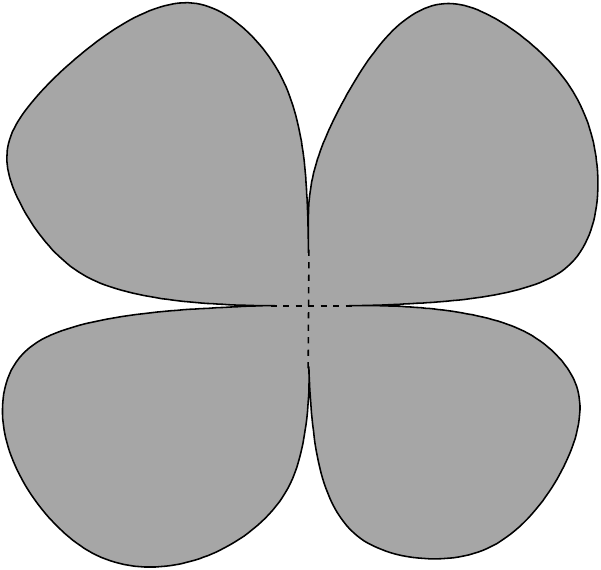}\qquad\includegraphics[width=3.5cm]{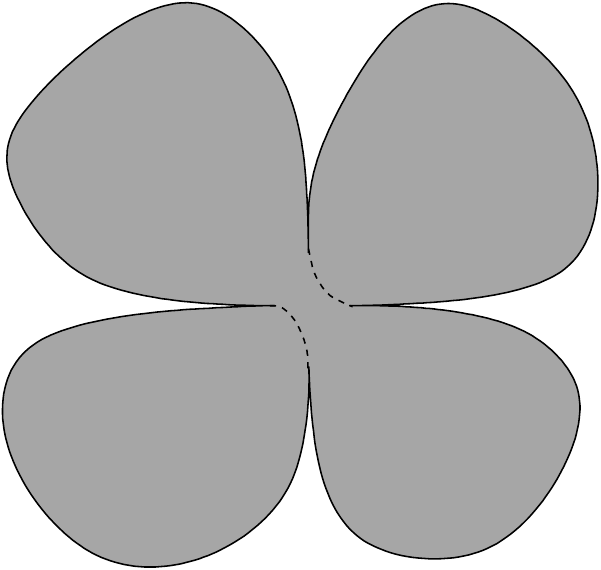}
  \caption{Top: a set $E_1$ such that $\Gamma-\lim \cW_\varepsilon(E_1)<\infty$ and $\oW(E_1)=+\infty$. Bottom, from left to right, a set $E_2$, the limit configuration whose energy coincides with $\Gamma-\lim \cW_\varepsilon(E_2)$, a configuration whose energy coincides with $\oW(E_2)$.}\label{fig:gloups}
\end{center}
  \end{figure}

Define for any set $E$ of finite perimeter in $\Omega$ its relaxed Willmore functional
$$\oW(E,\Omega)=\inf\{\liminf W(E_h,\Omega),\;\partial E_h\in\Omega\in\cont^2,\;E_h\to E\mbox{ in }\lp{1}(\Omega)\}.$$
The properties of this relaxation are fully known in dimension $2$~\cite{BDP,BM1,BM,MN1} and partially known in higher dimension~\cite{AM,LM,MN2}. It is natural to ask whether the $\Gamma$-convergence of $\cW_\varepsilon$ to $W$ can be extended to $\oW$. Unfortunately, this is not the case as it follows from the following observations (for simplicity we denote $\Gamma-\lim \cW_\varepsilon(E)=\Gamma-\lim \cW_\varepsilon(\one{E})$) that are illustrated in Figure~\ref{fig:gloups}:
\begin{enumerate}
\item there exists a bounded set $E_1\subset\R^2$ of finite perimeter such that 
$$\Gamma-\lim \cW_\varepsilon(E_1)<\infty\qquad\mbox{and}\qquad\oW(E_1)=+\infty$$
\item there exists a bounded set $E_2\subset\R^2$ of finite perimeter such that 
$$\Gamma-\lim \cW_\varepsilon(E_2)<\oW(E_2)<+\infty$$
\end{enumerate}
\par   The reason why $\oW(E_1)=+\infty$ is a result by Bellettini, Dal Maso and Paolini~\cite{BDP} according to which a non oriented tangent must exist everywhere on the boundary. Besides, $\oW(E_2)<+\infty$ because, still by a result of Bellettini, Dal Maso and Paolini, the boundary is smooth out of evenly many cusps. Let us now explain why, in both cases,  $\Gamma-\lim \cW_\varepsilon(E_{1,2})<+\infty$.
The reason for this is the existence of smooth solutions with singular nodal sets for the Allen-Cahn equation
$$\Delta u-W'(u)=0$$
According to Dang, Fife and Peletier~\cite{Dang_Fife_Peletier_92}, there exists for such equation in $\R^2$ a unique saddle solution $u$ with values in $(-1,1)$. By saddle solution, it is meant that $u(x,y)>0$ in quadrants I and III, and $u(x,y)<0$ in quadrants II and IV, in particular $u(x,y)=0$ on the nodal set $xy=0$. Considering $u_\varepsilon(x)=u(\varepsilon x)$, we immediately get that
$$\varepsilon^2\Delta u_\varepsilon-W'(u_\varepsilon)=0$$
thus the second term in $\cW_\varepsilon(u_\varepsilon)$ {\it vanishes}, and the first term being obviously bounded, it follows from the lower semicontinuity of the $\Gamma$-limit that $\Gamma-\lim \cW_\varepsilon(E_{1,2})<+\infty$.
Furthermore, the approximation of $E_2$ can be made so as to create a cross in the limit, as in bottom-middle figure. The limit energy is therefore lower than the energy obtained by pairwise connection without crossing of the cusps (bottom-right figure). Thus,  $\Gamma-\lim \cW_\varepsilon(E_2)<\oW(E_2)<+\infty$.

For the reader not familiar with varifolds, it must be emphasized that this is not in contradiction with the results described in the previous section, and more precisely with the fact that the discrepancy measure guarantees the concentration of the diffuse varifolds at a limit integer varifold with generalized curvature in $\lp{2}$. Indeed, the boundary curves of $E_1$ and of the bottom-middle set can be canonically associated with a varifold having $\lp{2}$ generalized curvature because, by compensation between the tangents associated with each branch meeting at the cross, there is no singularity.

We end this section with the question that follows naturally from the discussion above: is it possible to find a diffuse approximation that $\Gamma$-converges to $\oW$ (up to a multiplicative constant) whenever $\oW(E)<+\infty$ ?

\subsection{Diffuse approximations of the relaxed Willmore functional}
\subsubsection{Bellettini's approximation in dimension $N\geq 2$}
In~\cite{Bellettini1997}, G. Bellettini proposed a diffuse model for approximating the relaxations of geometric functionals of the form $\int_{\partial E}(1+f(x,\nabla d_E,\nabla^2 d_E))d\haushn$ where $E$ is smooth and $d_E$ is the signed distance function from $\partial E$. Such functionals include the Willmore energy since, on $\partial E$, $H=(\Delta d_E)\nabla d_E=\operatorname{tr}(\nabla^2d_E)\nabla d_E$ thus $|H|^2=|\operatorname{tr}(\nabla^2d_E)|^2$. Particularizing Bellettini's approximation model to this case yields the smooth functionals
$$\cW^{\mbox{\tiny Be}}_\varepsilon(u)=\left\{\begin{array}{ll}\ds\frac 1 2\int_{\Omega\setminus \{|\nabla u|=0\}}(|\mdiv\frac{\nabla u}{|\nabla u|}|^2)(\frac\varepsilon 2|\nabla u|^2+\frac{W(u)}\varepsilon)dx&\mbox{if }u\in\cont^\infty(\Omega)\\
+\infty&\mbox{otherwise in }\lp{1}(\Omega)\end{array}\right.$$
Then, according to Bellettini~\cite[Thms 4.2,4.3]{Bellettini1997}, in any space dimension,
$$(\Gamma-\lim_{\varepsilon\to 0} P_{\varepsilon} +  \cW^{\mbox{\tiny Be}}_\varepsilon)(E)=c_0 (P(E) + \oW(E))\qquad\mbox{for every $E$ of finite perimeter such that $\oW(E)<+\infty$}$$
The constructive part of the proof is based, as usual, on using approximating functions of the form $u_\varepsilon=q(\frac{d_E}{\varepsilon})$. As for the lower semicontinuity part, it is facilitated by the explicit appearance of the mean curvature in the expression. Recall indeed that, $u$ being smooth, for almost every $t$, $H_u(x):=(\mdiv\frac{\nabla u}{|\nabla u|})\frac{\nabla u}{|\nabla u|}(x)$ is the mean curvature at a point $x$ of the isolevel $\{y,\,u(y)=t\}$.
Let $(u_\varepsilon)$ be a sequence of smooth functions that approximate $u=\one{E}$ in $\lp{1}(\R^N)$ and has uniformly bounded total variation. Then, by the coarea formula,
$$\cW^{\mbox{\tiny Be}}_\varepsilon(u_\varepsilon)\geq \frac 1 2\int_{\{|\nabla u_h|\not=0\}}|\nabla u_\varepsilon|\sqrt{2W(u_\varepsilon)}|H_{u_\varepsilon}|^2dx= \frac 1 2\int_{0}^1\sqrt{2W(t)}\int_{\{u_\varepsilon=t\}\cap \{|\nabla u_\varepsilon|\not=0\}}|H_{u_\varepsilon}|^2d\haushn\,dt.$$
The last inequality is important: it guarantees a control of the Willmore energy of the isolevel surfaces of $u_\varepsilon$. This is a major difference with the classical approximation, for which such control does not hold.
\par Then, it suffices to observe that, by the Cavalieri formula and for a suitable subsequence, $|\{u_\varepsilon\geq t\}\Delta \{u\geq t\}|\to 0$ for almost every $t$. In addition, $\{u\geq t\}=E$ for almost every $t$, and by the lower semicontinuity of the relaxation $\oF(E)\leq\liminf_{\varepsilon\to 0}F(u_\varepsilon)$. Fatou's Lemma finally gives
$$\liminf_{\varepsilon\to 0}\cW^{\mbox{\tiny Be}}_\varepsilon(u_\varepsilon)\geq\oF(E) \int_{0}^1\sqrt{2W(t)}dt=c_0\oF(E).$$
Bellettini's approximation has however a drawback that will be explained with more details later: when one computes the flow associated with the functional, the 4th order term is nonlinear, which raises difficulties at the numerical level since it cannot be treated implicitly.
\subsubsection{Mugnai's approximation in dimension $N=2$}
 In the regular case and in dimensions 2,3, it follows from the results of Bellettini and Mugnai~\cite{belmu10} that, up to a uniform control of the perimeter, the $\Gamma$-limit of the functionals defined by
$$\cW^{\mbox{\tiny Mu}}_\varepsilon(u) =\left\{\begin{array}{ll} \ds\frac 1{2\varepsilon}\int_{\Omega} \left| \varepsilon\nabla^2 u-\frac{W'(u)}\varepsilon \nu_u\otimes\nu_u \right|^2 dx&\mbox{if }u\in\cont^2(\Omega)\\
+\infty&\mbox{otherwise in }\lp{1}(\Omega),\end{array}\right.$$
where $\nu_u=\frac{\nabla u}{|\nabla u|}$ when $|\nabla u| \not=0$, and $\nu_u=${\it constant unit vector} on $\{|\nabla u|=0\}$,
coincides with $$c_0\int_{\Omega\cap\partial E}|A_{\partial E}(x)|^2dx$$ for every smooth $E$, with $A_{\partial E}(x)$ the second fundamental form of $\partial E$ at $x$. Again, this approximation allows a control of the mean curvature of the isolevel surfaces of an approximating sequence $u_\varepsilon$, thus prevents from the creation of saddle solutions to the Allen-Cahn equation since, by~\cite[Lemma 5.3]{belmu10} and ~\cite[Lemma 5.2]{mugnai2010}, 
$$ |\nabla u||\mdiv\frac{\nabla u}{|\nabla u|}|\leq \frac 1\varepsilon| \varepsilon\nabla^2 u-\frac{W'(u)}\varepsilon \nu_u\otimes\nu_u|$$
In dimension $2$, the second fundamental form along a curve coincides with the curvature. Therefore, by identifying the limit varifold obtained when $u_\varepsilon$ converges to $u=\one{E}$, and using the representation results of~\cite{BM}, Mugnai was able to prove in~\cite{mugnai2010} that, in dimension $2$, the $\Gamma$-limit of $\cW^{\mbox{\tiny Mu}}_\varepsilon$ (with uniform control of the perimeter) coincides with $\oF(E)$ for any $E$ with finite perimeter.

\subsubsection{Esedoglu-Rätz-Röger's approximation in dimension $N\geq 2$}\label{sec:err}
The model of Esedoglu, R{\"a}tz, and R{\"o}ger in~\cite{Esedoglu_12} is a modification of the classical energy that aims to preserve the ``parallelity'' of the level lines of the approximating functions, and avoids the formation of saddle points, by constraining the level lines' mean curvature using a term {\it \`a la} Bellettini. More precisely, one can calculate that
$$\varepsilon\Delta u-\frac{W'(u)}{\varepsilon}=\varepsilon|\nabla u|\mdiv{\frac{\nabla u}{|\nabla u|}}-\nabla \xi_\varepsilon\cdot\frac{\nabla u}{|\nabla u|^2}.$$
with $\xi_\varepsilon= \left(\frac\varepsilon 2|\nabla u|^2-\frac{W(u)}\varepsilon \right)$ the discrepancy function (with a small abuse of notation, we use the same notation for the discrepancy measure and its density).
\par 
Therefore, $\varepsilon\Delta u-\frac{W'(u)}{\varepsilon}$ approximates correctly the mean curvature (up to a multiplicative constant) if the projection of $\nabla\xi_\varepsilon$ on the orthogonal direction $\nabla u$ is small. Equivalently, it can be required that	
 $$\varepsilon\Delta u-\frac{W'(u)}{\varepsilon}-\varepsilon|\nabla u|\mdiv{\frac{\nabla u}{|\nabla u|}}$$
be small, therefore a natural profile-forcing approximation model is (with $\alpha\geq 0$ a parameter):
$$\cW^{\mbox{\tiny EsR\"aR\"o}}_\varepsilon(u)=\left\{\begin{array}{ll}\ds\frac{1}{2\varepsilon}\int_\Omega \left( \varepsilon\Delta u-\frac{W'(u)}{\varepsilon} \right)^2dx+\\
 \quad\quad\ds\frac 1{2\varepsilon^{1+\alpha}}\int_\Omega(\varepsilon\Delta u-\frac{W'(u)}{\varepsilon}-\varepsilon|\nabla u|\mdiv{\frac{\nabla u}{|\nabla u|}})^2dx&\mbox{if }u\in\cont^\infty(\Omega)\\
+\infty&\mbox{otherwise in }\lp{1}(\Omega)\end{array}\right.$$
To simplify the theoretical analysis, the model proposed by Esedoglu, R{\"a}tz, and R{\"o}ger is slightly different. It uses the fact that, if a phase field $u_\varepsilon$ resembles $q(\frac d\varepsilon)$, one has $\varepsilon|\nabla u|\sim\sqrt{2W(u)}$, which leads Esedolu, R{\"a}tz, and R{\"o}ger to penalize
 $$\varepsilon\Delta u-\frac{W'(u)}{\varepsilon}-(\varepsilon|\nabla u|(2W(u))^\frac 1 2)^{\frac 1 2}\mdiv{\frac{\nabla u}{|\nabla u|}}.$$
Finally, they propose the following approximating functional
 $$\widehat{\cW^{\mbox{\tiny EsR\"aR\"o}}_\varepsilon}(u)=\left\{\begin{array}{ll}\ds\frac{1}{2\varepsilon}\int_\Omega \left( \varepsilon\Delta u-\frac{W'(u)}{\varepsilon} \right)^2dx+\\
 \quad\quad\ds\frac 1{2\varepsilon^{1+\alpha}}\int_\Omega(\varepsilon\Delta u-\frac{W'(u)}{\varepsilon}-(\varepsilon|\nabla u|\sqrt{2W(u)})^{\frac 1 2}\mdiv{\frac{\nabla u}{|\nabla u|}})^2dx&\mbox{if }u\in\cont^\infty(\Omega)\\
+\infty&\mbox{otherwise in }\lp{1}(\Omega)\end{array}\right.$$
This energy controls the mean curvature of the level lines of an approximating function since (see ~\cite{Esedoglu_12})
$$\widehat{ \cW^{\mbox{\tiny EsR\"aR\"o}}_\varepsilon}(u)\geq\frac{\varepsilon^{-\alpha}}{2+2\varepsilon^{-\alpha}}\int_0^1\sqrt{2W(t)}\int_{\{u=t\}\cap\{\nabla u\not=0\}}(\mdiv{\frac{\nabla u}{|\nabla u|}})^2d\haushn\,dt,$$
which, once again, excludes Allen-Cahn solutions. With the control above, the authors prove with the same argument as Bellettini~\cite{Bellettini1997} that, for any $\alpha>0$, 
$$\Gamma-\lim_{\varepsilon\to 0} P_\varepsilon +  \widehat{\cW^{\mbox{\tiny EsR\"aR\"o}}_\varepsilon}=c_0 \left(P + \oW\right)\qquad\mbox{in }\lp{1}(\Omega).$$
With $\alpha=0$ the $\Gamma$-convergence result does not hold anymore, but instead, with a uniform control of the perimeter,
$$\Gamma-\lim_{\varepsilon\to 0} \widehat{\cW^{\mbox{\tiny EsR\"aR\"o}}_\varepsilon}\geq \frac {c_0} 2\oW.$$
which still guarantees a control of $\oW$.

For the sake of numerical simplicity, another version is tackled numerically in~\cite{Esedoglu_12}, based again on the approximation  $\varepsilon|\nabla u|\sim\sqrt{2W(u)}$:
 $$\widehat{\widehat{ \cW^{\mbox{\tiny EsR\"aR\"o}}_\varepsilon}}(u)=\left\{\begin{array}{ll}\ds\frac{1}{2\varepsilon}\int_\Omega \left( \varepsilon\Delta u-\frac{W'(u)}{\varepsilon} \right)^2dx+\\
 \quad\quad\ds\frac 1{2\varepsilon^{1+\alpha}}\int_\Omega(\varepsilon\Delta u-\frac{W'(u)}{\varepsilon}-\sqrt{2W(u)}\mdiv{\frac{\nabla u}{|\nabla u|}})^2dx&\mbox{if }u\in\cont^\infty(\Omega)\\
+\infty&\mbox{otherwise in }\lp{1}(\Omega)\end{array}\right.$$

We will focus in the sequel on $\cW^{\mbox{\tiny EsR\"aR\"o}}_\varepsilon$, whose flow will be derived, as well as its asymptotic behavior as $\varepsilon$ goes to $0$.

\subsection{Few remarks on the connections between the different approximations}\label{sec:remarks}
\subsubsection{From Mugnai's model to Esedoglu-Rätz-Röger's}
We saw previously that the  phase-field approximations $\cW^{\mbox{\tiny Be}}_\varepsilon$ and $\cW^{\mbox{\tiny EsR\"aR\"o}}_\varepsilon$ $\Gamma$-converge, up to a uniform control of perimeter,  to $c_0 \overline{W}$ in any dimension, and the same holds true in dimension $2$ for $\cW^{\mbox{\tiny Mu}}_\varepsilon$. We will now emphasize the connections between these approximations. More precisely, we will see that Mugnai's approximation $\cW^{\mbox{\tiny Mu}}_\varepsilon$  can be viewed as the sum of a geometric-type 
approximation of the Willmore energy plus a profile penalization term of the same kind as in Esedoglu, R{\"a}tz, R{\"o}ger's model (or, more precisely, the initial model $\cW^{\mbox{\tiny EsR\"aR\"o}}_\varepsilon$).
Indeed we have, denoting $\nu=\frac{\nabla u}{|\nabla u|}$ when $|\nabla u| \not=0$, and $\nu=${\it constant unit vector} on $\{|\nabla u|=0\}$,
\begin{eqnarray*}
W^{\mbox{\tiny Mu}}_\varepsilon(u) &=&  \ds\frac{1}{2\varepsilon}\int_{\Omega} \left| \varepsilon\nabla^2 u-\frac{W'(u)}\varepsilon \nu\otimes\nu \right|^2 dx  \\
                &=&  \ds \frac{1}{2\varepsilon} \int_{\Omega\setminus\{|\nabla u|=0\}} \left( \varepsilon\nabla^2 u : \nu\otimes\nu -\frac{W'(u)}\varepsilon \right)^2 dx+  \ds \int_{\Omega}\frac \varepsilon 2\left( |\nabla^2 u |^2 - (\nabla^2 u :  \nu\otimes\nu )^2 \right)dx. 
\end{eqnarray*}
where, being $A,B$ two matrices, we denote as $A:B = \sum_{i,j} A_{ij} B_{ij}$ the usual matrix scalar product. Using $ A : e_1 \otimes e_2 = <A e_2,e_1 >$,  we observe that  $\nabla^2 u : \nu\otimes\nu=\Delta u-|\nabla u|\mdiv{\frac{\nabla u}{|\nabla u|}}$, therefore the first term of $\cW^{\mbox{\tiny Mu}}_\varepsilon$ coincides with the second term of $\cW^{\mbox{\tiny EsR\"aR\"o}}_\varepsilon$ for $\alpha=0$. The second term of $\cW^{\mbox{\tiny Mu}}_\varepsilon$ can be splitted as 
\begin{eqnarray*}
   \ds\int_{\Omega} \frac \varepsilon 2\left( |\nabla^2 u |^2 - (\nabla^2 u :  \nu\otimes\nu )^2 \right)dx  & =&  \ds \frac 1 2  \int_{\Omega} \left| \nabla \left( \frac{\nabla u}{|\nabla u|} \right) \right|^2 \left( \varepsilon |\nabla u|^2 \right)dx,\\
 &~& \quad \quad  + \ds \int_{\Omega} \frac \varepsilon 2\left( |\nabla^2u~\nu|^2 - |\nabla^2 u  : \nu\otimes \nu|^2 \right)dx.\\
\end{eqnarray*}
Note that
$$  \int_{\Omega} \varepsilon \left( |\nabla^2u~\nu|^2 - \left(\nabla^2 u : \nu\otimes \nu\right)^2 \right)dx \geq 0,$$
is positive and vanishes for all functions $u$ of the general form $u = \eta\left(d(x)\right)$ with $\eta$ smooth. It is therefore a 
soft profile-penalization term that forces the approximating function to be a profile, yet not necessarily the optimal profile $q$. As for the term $\ds\int_{\Omega} \left| \nabla \left( \frac{\nabla u}{|\nabla u|} \right) \right|^2 \left( \varepsilon |\nabla u|^2 \right)dx$, it is purely geometric and constrains the approximating function's level lines mean curvature. It would therefore be worth addressing the $\Gamma$-convergence of the new functional
$$\cW^{\mbox{\tiny New}}_\varepsilon(u)=\frac 1{2\varepsilon}\int_{\Omega} (\varepsilon\Delta u-\frac 1 \varepsilon W'(u))^2dx+\frac 1{2\varepsilon^\alpha}\int \varepsilon \left( |\nabla^2u~\nu|^2 - \left(\nabla u ^2 : \nu\otimes \nu\right)^2 \right)dx.$$
The reason why such approximation would be interesting is that, if it indeeds $\Gamma$-converges, the associated flow would not be influenced by the asymptotic behavior of the penalization term, since it vanishes for approximating functions that are profiles. More precisely, the Willmore flow could be captured at low order of $\varepsilon$, and not at the numerically challenging order $\varepsilon^3$ as for the Esedoglu-R\"atz-R\"oger model with $\alpha=0$ or $1$.

\subsubsection{Towards a modification of Mugnai's energy that forces the $\Gamma$-convergence in dimension $\geq 3$}
Obviously, we cannot expect that  Mugnai's energy  $\Gamma$-converges  to the  Willmore energy  in dimension greater than $2$, since for $E$ smooth
$$|A_{\partial E}|^2=|H_{\partial E}|^2-\sum_{i\neq j } \kappa_i \kappa_j$$
where $\kappa_1, \kappa_2 \dots \kappa_{N-1}$ are the principal curvatures. 
This identity suggests however that a suitable correction could force the $\Gamma$-convergence, i.e. by subtracting to $\cW^{\mbox{\tiny Mu}}_{\varepsilon}$ an approximation of 
$$ J(E,\Omega) =  \int_{\partial E \cap \Omega} \sum_{i\neq j } \kappa_i \kappa_j ~d\haushn.$$

Recalling our assumptions that $d<0$ in $E$, and as an easy consequence of Lemma 14.17 in~\cite{GT} (see also~\cite{Ambrosio2000}), we obtain in a small tubular neighborhood of $\partial E$:
\begin{eqnarray*}
 \mdiv \left(\Delta d(x) \nabla d(x) \right) &=& (\Delta d(x))^2 + \nabla \Delta d(x) \cdot \nabla d(x)\\
                                             &=&   \left(\sum_{i}  \frac{\kappa_i(\pi(x))}{1 + d(x) \kappa_i(\pi(x))} \right)^2  -  \sum_i \frac{\kappa_i(\pi(x))^2}{(1 + d(x) \kappa_i(\pi(x)))^2} \\
                                             & \simeq & \sum_{i\neq j } \kappa_i \kappa_j \quad \text{on} \quad \partial E,
\end{eqnarray*}
where $\pi(x)$ is the projection of $x$ on $\Gamma$. Thus, a possible approximation of $c_0 J(E,\Omega)$ is  
$$  J^1_{\varepsilon}(u) =  - \frac{2}{\varepsilon} \int_{\Omega} \left( \varepsilon \Delta u - \frac{1}{\varepsilon} W'(u) \right) \frac{W'(u)}{\varepsilon} dx.$$
Indeed, with $u = q(d/\varepsilon)$ and with a suitable truncation of $q$ so that $q'(d/\varepsilon)$ vanishes on $\partial\Omega$ (which is always possible if $E\subset\!\subset\Omega$), integrating by parts yields:
\begin{eqnarray*}
  J^1_{\varepsilon}(u) &=&   -\frac{2}{\varepsilon^2}  \int_{\Omega}  \Delta d q'\left( \frac{d}{\varepsilon}\right) q''\left( \frac{d}{\varepsilon}\right) dx = -\frac{1}{\varepsilon}  \int_{\Omega}  \Delta d \nabla \left( q'\left( \frac{d}{\varepsilon} \right)^2 \right). \nabla d ~  dx \\
           &=& \frac{1}{\varepsilon} \int_{\Omega} \mdiv \left(\Delta d \nabla d \right) q'\left( \frac{d}{\varepsilon}  \right)^2 ~ dx \simeq {c_0} \int_{\partial E \cap \Omega} \sum_{i\neq j } \kappa_i \kappa_j ~d\haushn.
\end{eqnarray*}
Remark also that since a profile function $u = q(d/\varepsilon)$ satisfies $\frac{1}{\varepsilon} W'(u) = \varepsilon\nabla^2 u : {\mathbb N}(u)$ where 
${\mathbb N}(u)= \nu \otimes \nu = \frac{\nabla u}{|\nabla u|} \otimes\frac{\nabla u}{|\nabla u|} $, the following energies
$$
\begin{cases}
 J^2_{\varepsilon}(u) &=  - \frac{2}{\varepsilon^2} \int_{\Omega} \left( \varepsilon \Delta u - \varepsilon\nabla^2 u : {\mathbb N}(u)\right) W'(u) dx \\
 J^3_{\varepsilon}(u) &=  - 2\int_{\Omega} \left( \varepsilon \Delta u - \frac{1}{\varepsilon} W'(u) \right) \nabla^2 u : {\mathbb N}(u)dx\\
 J^4_{\varepsilon}(u) &=  - 2\int_{\Omega} \left( \varepsilon \Delta u - \varepsilon \nabla^2 u : {\mathbb N}(u)\right) \nabla^2 u : {\mathbb N}(u)dx
\end{cases}
$$ 
approximate also $c_0J(E,\Omega)$. In particular, as a modified version of Mugnai's energy,  we can consider 
$$ \widetilde{\cW^{\mbox{\tiny Mu}}_{\varepsilon}} = \cW^{\mbox{\tiny Mu}}_{\varepsilon} +  \frac 1 2(J^1_{\varepsilon}(u) +  J^3_{\varepsilon}(u) - J^4_{\varepsilon}(u))$$
We have indeed
$$ J^1_{\varepsilon}(u) + J^3_{\varepsilon}(u) -  J^2_{\varepsilon}(u) - J^4_{\varepsilon}(u) = \frac{2}{\varepsilon} \int_{\Omega} \left( \varepsilon\nabla^2 u : \nu\otimes\nu -\frac{W'(u)}\varepsilon \right)^2 dx.$$
and
$$ \cW^{\mbox{\tiny Mu}}_{\varepsilon} =  \cW_{\varepsilon}(u)  -  \frac 1 2 J^2_{\varepsilon}(u)$$
since $|\nabla^2u|^2=(\Delta u)^2$, ${\mathbb N}(u):{\mathbb N}(u)=1$,
$$|\varepsilon\nabla^2u-\varepsilon^{-1}W'(u){\mathbb N}(u)|^2=\varepsilon^2|\nabla^2u|^2-2W'(u)\nabla^2u:{\mathbb N}(u)+\varepsilon^{-2}W'(u)^2{\mathbb N}(u):{\mathbb N}(u)$$
and
$$(\varepsilon\Delta u-\varepsilon^{-1}W'(u))^2=\varepsilon^2(\Delta u)^2-2W'(u)\Delta u+\varepsilon^{-2}(W'(u))^2.$$
Therefore
$$  \widetilde{\cW^{\mbox{\tiny Mu}}_{\varepsilon}}  = \cW_{\varepsilon}(u) +   \ds \frac{1}{\varepsilon} \int_{\Omega} \left( \varepsilon\nabla^2 u : \nu\otimes\nu -\frac{W'(u)}\varepsilon \right)^2 dx,$$
which resembles Esedoglu-Rätz-Röger's approximation with $\alpha = 0$ since, in both approximations, the second term forces $u$ to be a ``profile'' function, and vanishes at the limit. In view of the approximation result of Esedoglu, Rätz, and Röger, it is reasonable to expect that $\widetilde{\cW^{\mbox{\tiny Mu}}_{\varepsilon}}$ $\Gamma$-converges to the relaxed Willmore energy in any dimension.
\section{The Willmore flow and its approximation by the evolution of a diffuse interface}\label{sec:flow}
This section is devoted to the approximation of the Willmore flow by $\lp{2}$-gradient flows associated with the approximating 
energies introduced above. In particular, we shall derive explicitly each approximating gradient flow and, 
using the matched asymptotic expansion method~\cite{fife,pego,belpaoqo,Loreti_march}, we will show that, at least formally and for smooth interfaces, there is convergence to the Willmore flow, at least in dimensions 2 and 3 for all flows, and in any dimension for some of them. The general question {\it ``if a sequence of functionals $\Gamma$-converges to a limit functional, is there also convergence of the associated flows?''} is rather natural, since $\Gamma$-convergence implies convergence of minimizers, up to the extraction of a subsequence. However, the question is difficult and remains open for the Willmore functional. Our results below give formal indications that the convergence holds. Serfaty discussed in~\cite{Serfaty} a general theorem on the $\Gamma$-convergence of gradient flows, provided that the generalized gradient of the associated functional can be controled (see in particular the discussion on the Cahn-Hilliard flow). Such control is so far out or reach for the Willmore functional. 

\subsection{On the Willmore flow}
Let $E(t)$, $0 \leq t \leq T$, denote the evolution by the Willmore flow of smooth domains, i.e. the outer normal velocity $V(t)$ is given at $x\in \partial E(t)$ by 
$$ V = \Delta_SH - \frac{1}{2} H^3 + H\|A\|^2,$$
where $\Delta_S$ is the Laplace-Beltrami operator on $\partial E(t)$, $H$ the scalar mean curvature, $A$ the second fundamental form, and $\|A\|^2$ is the sum of the squared coefficients of $A$.

In the plane, the Willmore flow coincides with the flow of curves associated with the Bernoulli-Euler elastica energy, i.e., denoting by $\kappa$ the scalar curvature
$$V = \Delta_S \kappa + \frac{1}{2}\kappa^3.$$
The long time existence of  a single curve evolving by this flow  is established in \cite{Existence_curve}, and any curve with fixed length converges to an elastica.

\par In higher dimension,  Kuwert and Sch{\"a}tzle give in \cite{Existence_flow_sphere,KuScha04} a long time existence proof of the Willmore flow and 
the convergence to a round sphere for sufficiently small initial energy. Singularities may appear for larger initial energies, as indicated by numerical simulations~\cite{Mayer01anumerical}.

\subsection{Approximating the Willmore flow with the classical De Giorgi-Bellettini-Paolini approach}
The $\lp{2}$-gradient flow of the approximating energy
$$ \cW_\varepsilon(u)= \frac{1}{2\varepsilon}\int_\Omega \left( \varepsilon\Delta u-\frac{W'(u)}{\varepsilon} \right)^2dx,$$
is equivalent to the evolution equation
$$ \partial_t u = -\Delta \left( \Delta u -\frac{1}{\varepsilon^2} W'(u) \right) + \frac{1}{\varepsilon^2} W''(u) \left( \Delta u -\frac{1}{\varepsilon^2} W'(u)  \right), $$
that can be rewritten as the phase field system
\begin{equation}
 \begin{cases}
     \varepsilon^2 \partial_t u = \Delta \mu -  \frac{1}{\varepsilon^2} W''(u) \mu \\
     \mu =  W'(u) - \varepsilon^2 \Delta u.
    \end{cases}\label{flow:classic}\end{equation}

\paragraph{Existence and well-posedness}
The well-posedness  of the phase field model \eqref{flow:classic} at fixed parameter $\varepsilon$ 
has been studied in \cite{colli1} with a volume constraint fixing the average of $u$, and in \cite{colli2} with both volume and area constraints.
\paragraph{Convergence to the Willmore flow}
Loreti and March showed in~\cite{Loreti_march}, by using the formal method of matched asymptotic expansions, that if $\partial E$ is smooth and evolves by Willmore flow, it can be approximated by level lines of the solution $u_\varepsilon$ of the phase field system \eqref{flow:classic} as $\varepsilon$ goes to $0$. In addition, $u_{\varepsilon}$ and $\mu_{\varepsilon}$ are expected to take the form
$$
\left\{\begin{array}{l}
  u_{\varepsilon}(x,t) =  q \left( \frac{d(x,E(t))}{\varepsilon}\right) + \varepsilon^2 \left( \|A\|^2 - \frac{1}{2}H^2 \right) \eta_1\left( \frac{d(x,E(t))}{\varepsilon}\right)+O(\varepsilon^3)\\
\mu_{\varepsilon}(x,t)=- \varepsilon  H q' \left( \frac{d(x,E(t))}{\varepsilon}\right) + \varepsilon^2  H^2 \eta_2 \left( \frac{d(x,E(t))}{\varepsilon}\right)+O(\varepsilon^3)
\end{array}\right.,$$
where $\eta_1$ and $\eta_2$  are two functions depending only of the double well potential $W$, and defined as the solutions of 
$$
\begin{cases}
 \eta_1''(s) - W''(q(s))\eta_1(s) = s q'(s), \quad \text{with} \quad \lim_{s \to \pm \infty} \eta_1(s) = 0, \\
 \eta_2''(s) - W''(q(s))\eta_2(s) = q''(s), \quad \text{with} \quad \lim_{s \to \pm\infty} \eta_2(s) = 0. \\
\end{cases}
$$

 An important point is that 
the second-order term in the  asymptotic expansion of $u_{\varepsilon}$ has an influence on the limit law as $\varepsilon$ goes to zero~\cite{Loreti_march}. This is a major difference with the Allen-Cahn equation, for which the velocity law follows from the expansion at zero and first orders only~\cite{belpaoqo}. As a consequence, addressing numerically the Willmore flow is more delicate and requires using a high accuracy approximation in space to guarantee a sufficiently good approximation of the expansion of $u_{\varepsilon}$.

\subsection{Approximating the Willmore flow with Bellettini's model}
We focus now on the approximation model
$$ \cW^{\mbox{\tiny Be}}_\varepsilon(u)=  \frac{1}{2}\int_\Omega \mdiv \left( \frac{\nabla u}{|\nabla u|} \right)^2\left( \frac{\varepsilon}{2} |\nabla u|^2 + \frac{1}{\varepsilon} W(u) \right)dx$$
We will prove in the next section that its $\lp{2}$-gradient flow is equivalent to the evolution equation
\begin{equation}
\quad \partial_t u =  \frac{ K(u)^2 }{2} \left( \Delta u  - \frac{1}{\varepsilon^2} W'(u)\right) + \frac{1}{2} \nabla[K(u)^2].\nabla u - \frac{1}{\varepsilon}  \mdiv \left( P^u \frac{\nabla  \left[ K(u) h_{\varepsilon}(u) \right]}{|\nabla u|} \right),\label{flow:bel}
\end{equation}
where  $P^u = I_d -\frac{\nabla u}{|\nabla u|}\otimes \frac{\nabla u}{|\nabla u|}=I_d-{\mathbb N}(u)$, 
$h_{\varepsilon}(u) = \left( \frac{\varepsilon}{2} |\nabla u|^2 + \frac{1}{\varepsilon} W(u) \right)$ and  $K(u) = \mdiv \left( \frac{\nabla u}{|\nabla u|} \right)$.
 ~\\ ~\\
Existence and well-posedness of this equation are open questions. 
Numerical simulations performed with this flow are shown in \cite{Esedoglu_12}.  Note that the fourth-order nonlinear term makes numerics harder. 

\par Using the formal method of matched asymptotic expansions, we show below that  the phase field model~\eqref{flow:bel} converges in any dimension, at least formally, to the Willmore flow 
as $\varepsilon$ goes to $0$. More precisely, we observe an asymptotic expansion of $u_{\varepsilon}$ of the form
$$ u_{\varepsilon}(x,t) =   q \left( \frac{d(x,E(t))}{\varepsilon}\right) + O(\varepsilon^2), $$
where the second-order term does not have any influence 
on the limit velocity law as $\varepsilon$ goes to zero, in contrast with the classical approximation of the previous section.
\subsubsection{Derivation of the $\lp{2}$-gradient flow of $\cW^{\mbox{\tiny Be}}_\varepsilon(u)$} ~\\
\begin{prop}
The $\lp{2}$-gradient flow of Bellettini's model is equivalent to
\begin{eqnarray*}
\partial_t u &=&  \frac{ K(u)^2 }{2} \left( \Delta u  - \frac{1}{\varepsilon^2} W'(u)\right) + \frac{1}{2} \nabla[K(u)^2].\nabla u - \frac{1}{\varepsilon}  \mdiv \left( P^u \frac{\nabla  \left[ K(u) h_{\varepsilon}(u) \right]}{|\nabla u|} \right),
\end{eqnarray*}
where  $P^u = I_d -\frac{\nabla u}{|\nabla u|}\otimes \frac{\nabla u}{|\nabla u|}$, $h_{\varepsilon}(u) = \left( \frac{\varepsilon}{2} |\nabla u|^2 + \frac{1}{\varepsilon} W(u) \right)$ and  $K(u) = \mdiv \left( \frac{\nabla u}{|\nabla u|} \right)$.
\end{prop}
\dimo
The differential of $K$ at $u$ satisfies 
$$ K'(u)(w) = \lim_{t\to 0}\frac{K(u+tw)-K(u)}t=\mdiv\left(  \frac{\nabla w}{ |\nabla u|} - \frac{ \nabla u.\nabla w \nabla u }{|\nabla u|^3} \right),$$
therefore
\begin{eqnarray*}
 (\cW^{\mbox{\tiny Be}}_\varepsilon(u))'(w) &=& \ds\int_{\Omega} \left[ K(u) h_{\varepsilon}(u) \right]  \mdiv\left(  \frac{\nabla w}{ |\nabla u|} - \frac{ \nabla u.\nabla w \nabla u }{|\nabla u|^3} \right) dx \\
                          &+& \frac{1}{2} \ds \int_{\Omega}  K(u)^2 \left( \varepsilon \nabla u \nabla w + \frac{1}{\varepsilon} W'(u) w \right) dx 
\end{eqnarray*}
It follows that the $\lp{2}$-gradient of $\cW^{\mbox{\tiny Be}}_\varepsilon$ reads as
\begin{eqnarray*}
 \nabla \cW^{\mbox{\tiny Be}}_\varepsilon(u) &=&    \mdiv \left( \frac{\nabla  \left[ K(u) h_{\varepsilon}(u) \right]}{|\nabla u|} \right) -  \mdiv \left( \frac{\nabla  \left[ K(u) h_{\varepsilon}(u) \right]}{|\nabla u|} \cdot \frac{\nabla u}{|\nabla u|} \frac{\nabla u}{|\nabla u|} \right)  \\
                          &~& \quad \quad -  \frac{1}{2} \left(\varepsilon \mdiv \left( K(u)^2 \nabla u \right) - \frac{1}{\varepsilon} K(u)^2 W'(u)\right).
\end{eqnarray*}
whence the $\lp{2}$-gradient flow of $\cW^{\mbox{\tiny Be}}_\varepsilon(u)$ follows.\qed

\subsubsection{Asymptotic analysis} ~\\
In this section, we compute the formal expansions of the solution $u_{\varepsilon}(x,t)$ to the phase field model \eqref{flow:bel}.

\paragraph{Preliminaries}

We assume without loss of generality that the isolevel set $\Gamma(t)=\{u_{\varepsilon}=\frac 1 2\}$ is a smooth $n-1$ dimensional boundary $\Gamma(t)=\partial E(t) = \partial\{ x \in \R^d ; u_{\varepsilon}(x,t) \geq 1/2 \}$. 
We follow the method of matched asymptotic expansions proposed in  \cite{fife,pego,belpaoqo,Loreti_march}. 
We assume that the so-called outer expansion of $u_{\varepsilon}$, i.e. the expansion far from the front $\Gamma$, is of the form 
$$u_{\varepsilon}(x,t) = u_0(x,t) + \varepsilon u_1(x,t) + \varepsilon^2 u_2(x,t) + O(\varepsilon^3)$$
In a small neighborhood of $\Gamma$, we define the stretched normal distance to the front,
$$ z = \frac{d(x,t)}{\varepsilon},$$
were $d(x,t)$ denotes the signed distance to $E(t)$ such that $d(x,t)<0$ in $E(t)$. 
We then focus on inner expansions of $u_{\varepsilon}(x,t)$, i.e. expansions close to the front, of the form 

$$
 u_{\varepsilon}(x,t) = U(z,x,t) = U_0(z,x,t) + \varepsilon U_1(z,x,t) + \varepsilon^2 U_2(z,x,t) +O(\varepsilon^3)
$$
Let us define a unit normal $ m$ to $\Gamma$ and the normal velocity $V$ to the front as  
$$ V = -\partial_t d(x,t) , \quad m = \nabla d(x,t),\qquad x\in\Gamma,$$
where $\nabla$ refers to spatial derivation only (the same holds for further derivation operators used in the sequel). Following~\cite{pego,Loreti_march} we assume that $U(z,x,t)$ does not change when $x$ varies normal to $\Gamma$ with $z$ held fixed, or equivalently
$\nabla_x U . m = 0$. This amounts to requiring that the blow-up with respect to the parameter $\varepsilon$ is coherent with the flow. 

\begin{claim}
In a suitable regime provided by the method of matched asymptotic expansions, the normal velocity of the $\frac 1 2$-front $\Gamma(t)=\partial E(t)$ associated with a solution $u_{\varepsilon}(x,t) $ to Bellettini's phase field model~\eqref{flow:bel}
is the Willmore velocity
$$ V =  \Delta_{\Gamma} H  +   \|A\|^2 H  -\frac{H^3}{2},$$
and
$$ u_{\varepsilon}(x,t) =   q ( \frac{d(x,E(t))}{\varepsilon}) + O(\varepsilon^2)$$
\end{claim}

Following~\cite{pego,Loreti_march}, it is easily seen that
$$
\begin{cases}
 \nabla u = \nabla_x U + \varepsilon^{-1} m \partial_z U \\
 \Delta u = \Delta_x U + \varepsilon^{-1} \Delta d \partial_z U + \varepsilon^{-2} \partial^2_{zz} U \\
 \partial_t u = \partial_t U - \varepsilon^{-1} V \partial_z U. 
\end{cases}
$$
Recall also that in a sufficiently small neighborhood of $\Gamma$, according to Lemma 14.17 in~\cite{GT} (see also~\cite{Ambrosio2000}), we have
$$ \Delta d (x,t) = \sum_{i=1}^{n-1} \frac{\kappa_i(\pi(x))}{1 + \kappa_i(\pi(x)) d(x,t)} = \sum_{i=1}^{n-1} \frac{\kappa_i(\pi(x))}{1 + \kappa_i(\pi(x)) \varepsilon z }$$
where $\pi(x)$ is the projection of $x$ on $\Gamma$, and $\kappa_i$ are the principal curvatures on $\Gamma$. \\ 
In particular this implies that 
$$ \Delta d (x,t) = H - \varepsilon z \|A\|^2 + O(\varepsilon^2), $$
where $H$ and $\|A\|^2$ denote, respectively, the mean curvature and the squared $2$-norm of the second fundamental form on $\Gamma$ at $\pi(x)$.

\paragraph{Outer solution:}
We now compute the solution $u_{\varepsilon}$ in the outer region. 
By equation~\eqref{flow:bel},  $u_0$ satisfies $W'(u_0) = 0$ and 
$$u_0(x,t) = 
\begin{cases}
 1 &  \text{if } x \in E(t) \\
 0 & \text{otherwise}
\end{cases}
$$
We also see that $u_1 = 0$ is a possible solution at the first order. 
\paragraph{Matching condition :} The  inner and outer expansions are related by the following matching condition
$$ u_0(x,t) + \varepsilon u_1(x,t) + \cdots = U_0(z,x,t) + \varepsilon U_1(z,x,t) + \cdots$$
with $x$ near the front $\Gamma$ and $\varepsilon z$ between $O(\varepsilon)$ and $\circ(1)$. With the notation 
$$ u_i^{\pm}(x,t) =  \lim_{s \to 0 \pm} u_i(x + s m,t),$$
one has that
$$
\begin{cases}
 u_0^{\pm}(x,t) = \lim_{z \to \pm \infty} U_0(z,x,t) \\
 \lim_{z \to \pm \infty} u_1^{\pm}(x,t) + z m\cdot \nabla u_0^{\pm}(x,t)   = \lim_{z \to\pm \infty} U_1(z,x,t) 
\end{cases}
$$
In particular, for the phase field model~\eqref{flow:bel},  it follows that 
$$   \lim_{z \to + \infty} U_0(z,x,t) = 0, \lim_{z \to - \infty} U_0(z,x,t) = 1 \quad \text{and}  \lim_{z \to \pm \infty} U_1(z,x,t) = 0,$$
\paragraph{Inner solution:}
Note that 
$$  \frac{\nabla u}{|\nabla u|} = \frac{m - \varepsilon \nabla_x U / \partial_z U}{ \sqrt{1 + \varepsilon^2 |\nabla_x U|^2 / (\partial_z U)^2}},$$
therefore, using the orthogonality  condition $\nabla_x U . m = 0$:
$$ \begin{cases}
    K(u) = \Delta d + O(\varepsilon) \\
    h_{\varepsilon}(u) = \frac{1}{\varepsilon} \left[ \frac{1}{2}(\partial_{z} U)^2 + W(U)  \right] + O(1) \\
    \frac{1}{2} \nabla \left[ K(u)^2 \right]. \nabla u = \frac{1}{\varepsilon} \left( \Delta d \nabla (\Delta d).\nabla d \right) \partial_z U + O(1) \\
    \frac{1}{\varepsilon}  \mdiv \left( P^u \frac{\nabla  \left[ K(u) h_{\varepsilon}(u) \right]}{|\nabla u|}  \right) = \frac{1}{\varepsilon} \mdiv \left( \nabla(\Delta d) - \nabla (\Delta d)\cdot\nabla d \nabla d \right) \left( \frac{\frac{1}{2} (\partial_z U)^2 + W(U)}{|\partial_z U|} \right) + O(1)
   \end{cases}
 $$
Recall also that  in a sufficiently small neighborhood of $\Gamma$, 
$$ \Delta d (x,t) = \sum_{i=1}^{n-1} \frac{\kappa_i(\pi(x))}{1 + \kappa_i(\pi(x)) \varepsilon z },$$
thus
$$  
\begin{cases}
 \left( \Delta d \nabla (\Delta d).\nabla d \right) = - \|A\|^2 H  + O(\varepsilon) \\
 \mdiv \left( \nabla(\Delta d) - \nabla (\Delta d)\cdot\nabla d \nabla d \right) = \Delta_{\Gamma} H + O(\varepsilon).
\end{cases}
$$
Then, the first order in $\varepsilon^{-2}$ of Equation~\eqref{flow:bel} reads 
$$ \frac{H^2}{2} \left( \partial^2_{zz} U_0 - W'(U_0) \right) = 0.$$
Adding the boundary condition obtained from the matching condition, and using $U_0(0,x,t) = 1/2$ leads to 
$$U_0  = q(z).$$
Moreover, the second order in $\varepsilon^{-1}$ of~\eqref{flow:bel} shows that 
\begin{eqnarray*}
 -V \partial_z U_0 =  \frac{H^2}{2}\left( \partial^2_{zz} U_1 - W''(U_0)U_1 \right) + \frac{H^3}{2} \partial_{z} U_0 - \|A\|^2 H  \partial_{z} U_0 
- \Delta_{\Gamma} H  \left( \frac{\frac{1}{2} (\partial_z U_0)^2 + W(U_0)}{|\partial_z U_0|} \right).   
\end{eqnarray*}
As $U_0(z,x,t) = q(z)$  and $q' = -\sqrt{2 W(q)}$, we obtain  
\begin{eqnarray*}
 -V q'  =  \frac{H^2}{2}\left( \partial^2_{zz} U_1 - W''(q)U_1 \right) + \left( \frac{H^3}{2} - \|A\|^2 H -  \Delta_{\Gamma} H \right)  q'.  
\end{eqnarray*}
Then, multiplying by $q'$ and integrating over $\R$, it follows that 
$$ V =  \Delta_{\Gamma} H  +   \|A\|^2 H  -\frac{H^3}{2},$$
thus the sharp interface limit $\partial E(t)$ as $\varepsilon$ goes to zero evolves, at least formally, as the Willmore flow. In addition, we have $U_1 = 0$, therefore
$$ u_{\varepsilon}(x,t) =   q ( \frac{d(x,E(t))}{\varepsilon}) + O(\varepsilon^2)$$
and the second-order term does not appear in the expression of $V$. This explains the numerical stability, despite the use of an explicit Euler scheme, observed by Esedoglu, R\"atz and R\"oger in~\cite{Esedoglu_12}.
\subsection{Approximating the Willmore flow with Mugnai's model}

The aim of this section is the derivation and the study of the $\lp{2}$-gradient flow associated with Mugnai's energy 
$$ \cW^{\mbox{\tiny Mu}}_{\varepsilon}(u) =   \frac{1}{2 \varepsilon} \ds \int_{\Omega} \left|\varepsilon D^2 u - \frac{1}{\varepsilon} W'(u) \frac{\nabla u}{|\nabla u|} \otimes  \frac{\nabla u}{|\nabla u|}   \right|^2 dx.$$
We will show that the flow is equivalent to the phase field system
\begin{equation}\begin{cases}
     \varepsilon^2 \partial_t u =  \Delta \mu -  \frac{1}{\varepsilon^2} W''(u) \mu  + W'(u){\cal B}(u)   \\
     \mu =   \frac{1}{\varepsilon^2} W'(u) - \Delta u, 
    \end{cases}\label{flow:mugnai}
\end{equation}
where 
$$ {\cal B}(u) =    \mdiv \left( \mdiv\left(\frac{\nabla u}{|\nabla u|}\right)\frac{\nabla u}{|\nabla u|} \right) -  \mdiv \left(  \nabla \left( \frac{\nabla u}{|\nabla u|} \right)\frac{\nabla u}{|\nabla u|}\right).$$
Note that this system coincides with the classical one, up to the addition of a penalty term ${\cal L}(u) = W'(u) {\cal B}(u)$. \\

The well-posedness  of the phase field model \eqref{flow:mugnai} at fixed parameter $\varepsilon$ is open, and requires presumably a regularization of the term ${\cal B}(u)$ as done numerically in the next section. 
\begin{claim}
In a suitable regime provided by the method of matched asymptotic expansions, the normal velocity of the $\frac 1 2$-front $\Gamma(t)=\partial E(t)$ associated with a solution $(u_{\varepsilon}, \mu_{\varepsilon})$ to Mugnai's phase field model~\eqref{flow:mugnai}
is
$$ V = \Delta_{\Gamma} H + \sum_{i}\kappa_i^3 - \frac{1}{2} \|A\|^2 H.$$
and
$$\left\{\begin{array}{l}
  u_{\varepsilon}(x,t) =  q \left( \frac{d(x,E(t))}{\varepsilon}\right) + \varepsilon^2 \frac{\|A\|^2}{2} \eta_1\left( \frac{d(x,E(t)}{\varepsilon}\right)+O(\varepsilon^3) \\
 \mu_{\varepsilon}(x,t)= - \varepsilon H q'\left( \frac{d(x,E(t)}{\varepsilon}\right)+  \|A\|^2 \varepsilon^2 \eta_2 \left( \frac{d(x,E(t)}{\varepsilon}\right)  O(\varepsilon^3), 
 \end{array}\right.$$
where $\eta_1$ and $\eta_2$ are profile functions.
\end{claim}
\boss
The front velocity associated with Mugnai's phase field model coincides, up to a multiplicative constant, with the velocity of the $\lp{2}$-flow of the squared second fundamental form energy $\int_{\Gamma}\|A\|^2d\haushn$. Indeed, according to~\cite[Section 5.3]{amb-mant-96}, the latter is
$$\tilde V=2\Delta_\Gamma H+2H\|A\|^2-H^3+6\sum_{i<j<\ell}\kappa_i\kappa_j\kappa_\ell$$
Observing that $H\|A\|^2 = \sum_i \kappa_i^3  + \sum_{i\not=j} \kappa_i \kappa_j^2$ and
$$H^3 = \sum \kappa_i^3 + 3 \sum_{i\not=j} \kappa_i \kappa_j^2 + 6 \sum_{i<j<\ell} \kappa_i \kappa_j \kappa_\ell,$$
one has 
$$H\|A\|^2 - \frac 1 2 H^3 + 3\sum_{i<j<\ell}\kappa_i\kappa_j\kappa_\ell = \frac 1 2(\sum_i \kappa_i^3  -   \sum_{i\not=j} \kappa_i\kappa_j^2).$$
Since
$$\sum_{i}\kappa_i^3- \frac 1 2 H\|A\|^2  = \frac  1 2(\sum \kappa_i^3 - \sum_{i\not=j} \kappa_i \kappa_j^2)$$
we finally get that $\tilde V=2V$.
\eoss
\boss 
It is easily seen that, in dimensions $2$ and $3$, Mugnai's flow coincides with the Willmore flow. It is obvious in dimension $2$, whereas in dimension $3$ one has
$$\sum\kappa_i^3- \frac{1}{2} H\|A\|^2=\kappa_1^3+\kappa_2^3-\frac 1 2(k_1+k_2)(k_1^2+k_2^2)=\|A\|^2 H  -\frac{H^3}{2}.$$
Another explanation involves Gauss-Bonnet Theorem.  In Mugnai's model, the energy associated with the squared $2$-norm of the second fundamental form prevents from topological changes. By Gauss-Bonnet Theorem, this energy coincides with the Willmore energy up to a topological additive constant, and thus both associated flows coincide.
\eoss
\subsubsection{Derivation of the gradient flow}
\begin{prop}
The $\lp{2}$-gradient flow of Mugnai's model is equivalent to
$$\left\{\begin{array}{lll}
 \varepsilon^2 \partial_t u &=&  \Delta \mu -  \frac{1}{\varepsilon^2} W''(u) \mu  + W'(u){\cal B}(u)   \\
    \mu &=&   \frac{1}{\varepsilon^2} W'(u) - \Delta u, 
    \end{array}\right.$$
where 
$$ {\cal B}(u) =    \mdiv \left( \mdiv\left(\frac{\nabla u}{|\nabla u|}\right)\frac{\nabla u}{|\nabla u|} \right) -  \mdiv \left(  \nabla \left( \frac{\nabla u}{|\nabla u|} \right)\frac{\nabla u}{|\nabla u|}\right).$$
\end{prop}
\dimo
Let
 $$ {\mathbb V}(u) = \varepsilon D^2 u - \frac{1}{\varepsilon} W'(u) \frac{\nabla u}{|\nabla u|} \otimes  \frac{\nabla u}{|\nabla u|}.$$
The differential of ${\mathbb V}$ in the direction $w$ is
 \begin{eqnarray*}
   {\mathbb V}'(u)(w) &=& \lim_{t \to 0} ({\mathbb V}(u+tw) - {\mathbb V}(u))/t \\ 
            &=&  \varepsilon D^2 w - \frac{1}{\varepsilon} W''(u) w \frac{\nabla u}{|\nabla u|} \otimes  \frac{\nabla u}{|\nabla u|} -  \frac{1}{\varepsilon} W'(u) \left( \frac{\nabla u\otimes \nabla w + \nabla w\otimes \nabla u }{|\nabla u|^2}  \right) \\
            &~& \quad \quad \quad +  \frac{2}{\varepsilon} W'(u) \left( \frac{\nabla u \otimes \nabla u }{|\nabla u|^4} <\nabla u,\nabla w> \right) 
 \end{eqnarray*}
Denoting ${\mathbb N}(u)= \frac{\nabla u \otimes \nabla u}{|\nabla u|^2}$, we have 
\begin{multline*}
  \varepsilon\nabla \cW^{\mbox{\tiny Mu}}_{\varepsilon}(u) = \varepsilon D^2 : {\mathbb V}(u) - \frac{1}{\varepsilon} W''(u) {\mathbb N}(u):{\mathbb V}(u)  \\+ \frac{2}{\varepsilon} \mdiv\left( W'(u) \frac{ {\mathbb V}(u) \nabla u}{|\nabla u|^2}\right)  - \frac{2}{\varepsilon} \mdiv\left( W'(u) \left( {\mathbb V}(u) : {\mathbb N}(u)\right) \frac{\nabla u }{|\nabla u|^2}  \right),
\end{multline*}
where $D^2 : {\mathbb V}(u) =\nabla\otimes\nabla :{\mathbb V}(u)=\sum_{ij}\partial^2_{ij}{\mathbb V}_{ij}(u)=\mdiv(\mdiv {\mathbb V}(u))$ with the same abuse of notation as when one writes $\mdiv w=\nabla\cdot w$. 
The gradient of $\cW^{\mbox{\tiny Mu}}_{\varepsilon}$ can be also expressed as  
\begin{multline*}
 \varepsilon\nabla \cW^{\mbox{\tiny Mu}}_{\varepsilon}(u) = \varepsilon D^2 : {\mathbb V}(u) - \frac{1}{\varepsilon} W''(u) {\mathbb N}(u): {\mathbb V}(u) +\\
 \frac{2}{\varepsilon}  \mdiv \left( W'(u) \left( \frac{{\mathbb V}(u) \nabla u}{|\nabla u |^2} - < {\mathbb V}(u) \nabla u/|\nabla u|^2,\nabla u /|\nabla u|> \frac{\nabla u}{|\nabla u|} \right) \right)\end{multline*}
We now give an explicit expression of each previous term.

\paragraph{Evaluation of $\varepsilon D^2 :  {\mathbb V}(u)$} ~\\
For any operator $\Lambda$  and any real-valued function $u\mapsto\rho(u)$, we have
$$ D^2 : (\rho(u)\Lambda(u)) =  \Lambda(u) : D^2 \rho(u) + 2 < \nabla \rho(u), \mdiv(\Lambda(u))> + \rho(u) D^2 : \Lambda(u).$$
In particular, applying to $\Lambda(u)={\mathbb V}(u)$
 \begin{eqnarray*}
  \varepsilon D^2:{\mathbb V}(u) &=& \varepsilon^2 D^2 : D^2 u - D^2 : \left( W'(u) {\mathbb N}(u)\right) \\
                 &=& \varepsilon^2 \Delta^2 u - \left(  W'(u) D^2 : {\mathbb N}(u)+ 2 W''(u) <\nabla u, \mdiv({\mathbb N}(u))> +  (D^2(W'(u))): {\mathbb N}(u) \right). 
 \end{eqnarray*}
The last term reads as follows 
 \begin{eqnarray*}
  D^2(W'(u)): {\mathbb N}(u)&=& \left( W^{(3)}(u) \nabla u \otimes \nabla u + W''(u) \nabla^2 u  \right) : {\mathbb N}(u)\\
                  &=&  W^{(3)}(u) |\nabla u|^2 + W''(u)\frac{<D^2 u \nabla u,\nabla u >}{|\nabla u|^2} \\
                  &=& \Delta \left( W'(u)\right) - W''(u) \left( \Delta u  - \frac{<D^2 u \nabla u,\nabla u >}{|\nabla u|^2}\right),
 \end{eqnarray*}
where we used
$$  \Delta W'(u) = W''(u) \Delta u + W^{(3)}(u) |\nabla u|^2.$$
Recalling that for all vector fields $w_1, w_2$,
 $$ \mdiv(w_1 \otimes w_2) = \mdiv(w_2) w_1 + (\nabla w_1) w_2,$$
and applying to the estimation of $\mdiv({\mathbb N}(u))$, one gets that
 $$ \mdiv({\mathbb N}(u)) = \mdiv\left(\frac{\nabla u}{|\nabla u|}\right)\frac{\nabla u}{|\nabla u|} + \nabla \left( \frac{\nabla u}{|\nabla u|} \right)\frac{\nabla u}{|\nabla u|}.$$
Note that 
\begin{eqnarray*}
  \left[\nabla \left( \frac{\nabla u}{|\nabla u|} \right)\frac{\nabla u}{|\nabla u|}\right]\cdot \nabla u = \left[ \frac{D^2 u \nabla u}{|\nabla u|^2} - < \frac{D^2 u \nabla u}{|\nabla u|^2}, \frac{\nabla u}{|\nabla u|} >  \frac{\nabla u}{|\nabla u|}\right] \cdot \nabla u = 0. 
\end{eqnarray*}
Therefore
 \begin{eqnarray*}
  2 W''(u) <\nabla u, \mdiv({\mathbb N}(u))> = 2 W''(u) |\nabla u|  \mdiv\left(\frac{\nabla u}{|\nabla u|}\right) =  2 W''(u) \left(\Delta u - \frac{< D^2 u \nabla u, \nabla u>}{|\nabla u|^2} \right) \\
 \end{eqnarray*}
Lastly,
 \begin{eqnarray*}
   W'(u) D^2 : {\mathbb N}(u)&=& W'(u) \mdiv \left( \mdiv \left( {\mathbb N}(u)\right) \right) \\
                   &=& W'(u) \mdiv \left( \mdiv\left(\frac{\nabla u}{|\nabla u|}\right)\frac{\nabla u}{|\nabla u|} + \nabla \left( \frac{\nabla u}{|\nabla u|} \right)\frac{\nabla u}{|\nabla u|} \right) \\
                   &=& W'(u) \left[   \mdiv \left( \mdiv\left(\frac{\nabla u}{|\nabla u|}\right)\frac{\nabla u}{|\nabla u|} \right)  + \mdiv \left(  \nabla \left( \frac{\nabla u}{|\nabla u|} \right)\frac{\nabla u}{|\nabla u|}\right) \right]
 \end{eqnarray*}
therefore
 \begin{eqnarray*}
   \varepsilon D^2 : {\mathbb V}(u) &=&  \varepsilon^2 \Delta^2 u - \Delta W'(u) -  W''(u)\left(\Delta u - \frac{< D^2 u \nabla u, \nabla u>}{|\nabla u|^2} \right) \\
                  &-&  W'(u)  \left[   \mdiv \left( \mdiv\left(\frac{\nabla u}{|\nabla u|}\right)\frac{\nabla u}{|\nabla u|} \right)  + \mdiv \left(  \nabla \left( \frac{\nabla u}{|\nabla u|} \right)\frac{\nabla u}{|\nabla u|}\right) \right]. 
  \end{eqnarray*}
 ~\\
 \paragraph{Sum of the first two terms of $\varepsilon\nabla \cW^{\mbox{\tiny Mu}}_{\varepsilon}(u)$}~\\
Let $ I_1= \ds\varepsilon D^2 : {\mathbb V}(u) -\frac{1}{\varepsilon} W''(u) {\mathbb N}(u): {\mathbb V}(u)$. Remark that 
 \begin{eqnarray*}
  \frac{1}{\varepsilon} W''(u) {\mathbb N}(u): {\mathbb V}(u) &=& \frac{1}{\varepsilon} W''(u) {\mathbb N}(u): \left(\varepsilon D^2 u -\frac{1}{\varepsilon} W'(u) {\mathbb N}(u)\right) \\
                                    &=& W''(u) \left(\frac{< D^2 u \nabla u , \nabla u  >}{|\nabla u|^2}  - \frac{1}{\varepsilon^2}W'(u) \right) \\
                                    &=& W''(u) \left(\Delta u - \frac{1}{\varepsilon^2} W''(u) \right) - W''(u) \left( \Delta u - \frac{< \nabla^2 u \nabla u, \nabla u>}{|\nabla u|^2}\right)   
 \end{eqnarray*}
Combining with the previous estimation of $\varepsilon D^2 : {\mathbb V}(u)$, we obtain
 \begin{eqnarray*}
   I_1 &=& \varepsilon D^2 : {\mathbb V}(u)  - \frac{1}{\varepsilon} W''(u) {\mathbb N}(u): {\mathbb V}(u) \\
 &=& \varepsilon \Delta \left[\varepsilon \Delta u - \frac{1}{\varepsilon} W'(u) \right]  - \frac{1}{\varepsilon} W''(u) \left[\varepsilon \Delta u - \frac{1}{\varepsilon} W'(u) \right] \\ 
                                        &~& \quad \quad \quad - W'(u) \left[   \mdiv \left( \mdiv\left(\frac{\nabla u}{|\nabla u|}\right)\frac{\nabla u}{|\nabla u|} \right)  + \mdiv \left(  \nabla \left( \frac{\nabla u}{|\nabla u|} \right)\frac{\nabla u}{|\nabla u|}\right) \right]. 
 \end{eqnarray*}

\medskip
\paragraph{Estimation of the divergence term}~\\

Let $I_2 = \ds\frac{2}{\varepsilon}  \mdiv  W'(u) \left( \frac{{\mathbb V}(u) \nabla u}{|\nabla u |^2} -   W'(u)  ({\mathbb V}(u) : {\mathbb N}) \frac{\nabla u}{|\nabla u|^2} \right) $. On the one hand, with
 $$ {\mathbb V}(u) \frac{\nabla u}{|\nabla u|^2} = \varepsilon D^2 u \frac{\nabla u}{|\nabla u|^2} - \frac{1}{\varepsilon}W'(u) \frac{\nabla u}{|\nabla u|^2},$$
we see that 
 \begin{eqnarray*}
  \mdiv \left(W'(u) {\mathbb V}(u) \frac{\nabla u}{|\nabla u|^2} \right) &=& W'(u) \left[  \varepsilon \mdiv \left( \frac{D^2 u \nabla u}{|\nabla u|^2} \right)   - \frac{1}{\varepsilon}  \mdiv\left(  W'(u) \frac{\nabla u}{|\nabla u|^2} \right) \right] \\
                                                          &+& W''(u) \left(  \varepsilon \frac{< D^2 u \nabla u , \nabla u  >}{|\nabla u|^2} - \frac{1}{\varepsilon} W'(u) \right)
 \end{eqnarray*}
On the other hand, 
 $$ W'(u)  ({\mathbb V}(u) : {\mathbb N}) \frac{\nabla u}{|\nabla u|^2} = W'(u) \left(\varepsilon \frac{< D^2 u \nabla u , \nabla u  >}{|\nabla u|^2} \frac{\nabla u}{|\nabla u|^2 }- \frac{1}{\varepsilon} W'(u)  \frac{\nabla u}{|\nabla u|^2} \right)   $$
and
 \begin{eqnarray*}
  \mdiv \left( W'(u)  ({\mathbb V}(u) : {\mathbb N}) \frac{\nabla u}{|\nabla u|^2} \right)  &=& W''(u) \left[ \varepsilon \frac{< D^2 u \nabla u , \nabla u  >}{|\nabla u|^2} - \frac{1}{\varepsilon} W'(u)\right]  \\
                                                                    &+& W'(u) \mdiv \left[ \varepsilon \frac{< D^2 u \nabla u , \nabla u  >}{|\nabla u|^2} \frac{\nabla u}{|\nabla u|^2 }- \frac{1}{\varepsilon} W'(u)  \frac{\nabla u}{|\nabla u|^2}\right].
 \end{eqnarray*}
Finally,
 \begin{eqnarray*}
 I_2 &=&   2 W'(u) \mdiv \left( \frac{D^2 u \nabla u}{|\nabla u|^2} - < \frac{D^2 u \nabla u}{|\nabla u|}, \frac{\nabla u}{|\nabla u|} >  \frac{\nabla u}{|\nabla u|^2}   \right) \\
    &=& 2 W'(u)~\mdiv \left( \nabla \left( \frac{\nabla u}{|\nabla u|} \right) \frac{\nabla u}{|\nabla u|}  \right)
 \end{eqnarray*}
 \paragraph{Evaluation of the energy gradient}~\\
 \begin{eqnarray*}
   \varepsilon\nabla \cW^{\mbox{\tiny Mu}}_{\varepsilon}(u)&=&  I_1 + I_2 \\
                         &=& \varepsilon \Delta \mu  - \frac{1}{\varepsilon}W''(u) \mu - W'(u) \left[   \mdiv \left( \mdiv\left(\frac{\nabla u}{|\nabla u|}\right)\frac{\nabla u}{|\nabla u|} \right)  - \mdiv \left(  \nabla \left( \frac{\nabla u}{|\nabla u|} \right)\frac{\nabla u}{|\nabla u|}\right) \right], 
 \end{eqnarray*}
where
 $$\mu = \varepsilon   \Delta u  - \frac{1}{\varepsilon} W'(u).$$
 whence the $\lp{2}$-gradient flow associated with Mugnai's model follows.\qed

\subsubsection{Formal asymptotic expansions}
We apply in this section the formal method of matched asymptotic expansions to the solution $(u_{\varepsilon},\mu_{\varepsilon})$ of~\eqref{flow:mugnai}
  $$  \quad \begin{cases}
     \varepsilon^2 \partial_t u = \Delta \mu -  \frac{1}{\varepsilon^2} W''(u) \mu + W'(u) {\cal B}(u) \\
     \mu =  W'(u) - \varepsilon^2 \Delta u.
    \end{cases}.$$
Again, we assume without loss of generality that the isolevel set $\Gamma(t)=\{u_{\varepsilon}=\frac 1 2\}$ is a smooth $n-1$ dimensional boundary $\Gamma(t)=\partial E(t) = \partial\{ x \in \R^d ; u_{\varepsilon}(x,t) \geq 1/2 \}$. In addition, we assume that there exist outer expansions of $u_{\varepsilon}$ and $\mu_{\varepsilon}$ far from the front $\Gamma$ of the form 
$$ \begin{cases}
    u_{\varepsilon}(x,t) = u_0(x,t) + \varepsilon u_1(x,t) + \varepsilon^2 u_2(x,t) + O(\varepsilon^3) \\
  \mu_{\varepsilon}(x,t) = \mu_0(x,t) + \varepsilon \mu_1(x,t) + \varepsilon^2 \mu_2(x,t) +O(\varepsilon^3)   \end{cases}
$$
Considering the stretched variable $z = \frac{d(x,t)}{\varepsilon}$ on  a small neighborhood of $\Gamma$, 
we also look for inner expansions of $u_{\varepsilon}(x,t)$ and $\mu_{\varepsilon}(x,t)$ of the form 
$$
\begin{cases}
 u_{\varepsilon}(x,t) = U(z,x,t) = U_0(z,x,t) + \varepsilon U_1(z,x,t) + \varepsilon^2 U_2(z,x,t) +O(\varepsilon^3) \\
 \mu_{\varepsilon}(x,t) = W(z,x,t) =  W_0(z,x,t) + \varepsilon W_1(z,x,t) + \varepsilon^2 W_2(z,x,t) +O(\varepsilon^3) 
\end{cases}
$$
As before,  we define a unit normal $ m$ to $\Gamma$ and the normal velocity $V$ to the front as  
$$ V = -\partial_t d(x,t) , \quad m = \nabla d(x,t),\qquad x\in\Gamma.$$
Let us now expand $u_\varepsilon$ and $\mu_\varepsilon$.
\paragraph{Outer expansion:} Analogously to \cite{Loreti_march}, we obtain
$$u_0(x,t) = 
\begin{cases}
 1 & \text{ if } x \in E(t) \\
 0 & \text{otherwise}
\end{cases}, \quad \text{and} \quad   u_1 = u_2 = u_3 = \mu_0 = \mu_1= \mu_2 = 0.  
$$

\paragraph{Matching conditions:} 
The matching conditions (see \cite{Loreti_march} for more details) imply in particular that  
$$   \lim_{z \to + \infty} U_0(z,x,t) = 0, \lim_{z \to - \infty} U_0(z,x,t) = 1, \quad \lim_{z \to \pm \infty} U_1(z,x,t) = 0  \text{ and }  \quad \lim_{z \to \pm \infty} U_2(z,x,t) = 0$$
and 
$$   \lim_{z \to \pm \infty} W_0(z,x,t) = 0, \quad  \lim_{z \to \pm \infty} W_1(z,x,t) = 0  \text{ and }  \quad \lim_{z \to \pm \infty} W_2(z,x,t) = 0$$

\paragraph{Penalization term ${\cal B}(u)$}:  
With
$$  \frac{\nabla u}{|\nabla u|} = \frac{m - \varepsilon \nabla_x U / \partial_z U}{ \sqrt{1 + \varepsilon^2 |\nabla_x U|^2 / (\partial_z U)^2}},$$
and using $\nabla_x U . m = 0$, it follows that
\begin{eqnarray*}
 {\cal B}(u) &=& \left[ \mdiv \left( \mdiv\left(\frac{\nabla u}{|\nabla u|}\right)\frac{\nabla u}{|\nabla u|} \right)  - \mdiv \left(  \nabla \left( \frac{\nabla u}{|\nabla u|} \right)\frac{\nabla u}{|\nabla u|}\right) \right] \\
      &=& \mdiv \left( \Delta d \nabla d \right) - \mdiv \left( \nabla^2 d \nabla d \right) + O(\varepsilon)=  \left(\sum_{i}  \frac{\kappa_i(\pi(x))}{1 +  z \varepsilon  \kappa_i(\pi(x)} \right)^2  -  \sum_i \frac{\kappa_i(\pi(x))^2}{(1 + z \varepsilon \kappa_i(\pi(x)))^2} + O(\varepsilon) \\
      &=&  \left( H^2 - \|A\|^2 \right) +  O(\varepsilon).  
\end{eqnarray*}

\paragraph{Inner expansion:} We can derive the asymptotics of the second equation of the system \eqref{flow:mugnai}
$$  \mu =  W'(u) - \varepsilon^2 \Delta u $$
as follows   
$$
\begin{cases}
 W_0  = W'(U_0) - \partial^2_z U_0 \\
 W_1 =  W''(U_0)U_1 - \partial^2_z U_1 - \kappa \partial_z U_0 \\
 W_2 =  W''(U_0)U_2 - \partial^2_z U_2   + \frac{1}{2} W^{(3)}(U_0) U^2_1 -  H \partial_z U_1  + z \|A\|^2 \partial_z U_0   - \Delta_x U_0,
\end{cases}
$$
As for the first equation 
$$\varepsilon^2 \partial_t u = \Delta \mu -  \frac{1}{\varepsilon^2} W''(u) \mu + W'(u) {\cal B}(u)$$ 
therefore
$$
\begin{cases}
 0 = \partial^2_z W_0 - W''(U_0) W_0, \\
 0 = \partial^2_z W_1 + H \partial_z W_0 - \left(W''(U_0) W_1 +  W^{(3)}(U_0) U_1 W_0 \right),\\
\end{cases}.
$$
and 
\begin{eqnarray*}
 0 &=& \partial^2_z W_2 + H \partial_z \tilde{\mu}_1 - \|A\|^2 z  \partial_z W_0 + \Delta_x W_0  + W'(U_0)  \left( H^2 - \|A\|^2 \right) \\
   &-& \left(W''(U_0) W_2 +   W^{(3)}(U_0) U_1 W_1 +  W^{(3)}(U_0) U_2 W_0 + \frac{1}{2} W^{(4)}(U_0) U_1^2 W_0  \right) 
\end{eqnarray*}

\noindent {\it \underline{First order}}: ~\\ The two following equations  
$$ \partial^2_z W_{0} - W''(U_0) W_0  =0, \quad \text{and} \quad     W_0 =  W'(U_0)-\partial^2_z U_0,$$
associated  with the boundary conditions 
$$\lim_{z \to - \infty} U_0(z,x,t) = 1, \quad   \lim_{z \to + \infty} U_0(z,x,t) = 0, \quad \text{and} \quad  \lim_{z \to \pm \infty} W_0(z,x,t) = 0, $$
admit as solution pair
$$  U_0(z,x) = q(z), \quad \text{and} \quad W_0 =0.$$ 

\noindent {{\it \underline{Second order}}}:~\\ The second order  gives
$$ \partial^2_z W_{1} - W''(q) W_1  = 0 , \quad   W_1 =  W''(q)U_1 - \partial^2_z U_1 - H q'(z),$$
which has the solution 
$$ U_1 = 0,  \quad \text{and } W_1 = - H q'(z).$$

\noindent {{\it \underline{Third order}}}: ~\\
Using  $U_0 = q$, $W_0 = U_1 = 0$ and $W_1 = -Hq'$, the first  equation can be rewritten as 
\begin{eqnarray*}
  \partial^2_z W_2 - W''(q) W_2  &=& -  H \partial_z W_1 -  \left( H^2 - \|A\|^2 \right)W'(q) \\
                                   &=&  H^2 q''(z) -    \left( H^2 - \|A\|^2 \right) q''(z)  = \|A\|^2  q''(z),
\end{eqnarray*}
and implies that 
$$ W_2 = \|A\|^2 \eta_2(z) + c(x,t) q'(z),$$
where $\eta_2$ is defined as the solution of 
$$ \eta_2''(z) - W''(q(z))\eta_2(z) = q''(z), \quad \text{with} \quad \lim_{z \to \pm \infty} \eta_2(z) = 0.$$
Remark that $\eta_2$ can be expressed as 
$$\eta_2(z) = \frac{1}{2} z q'(z).$$

Note that the second equation also reads as 
$$ \partial^2_z U_2 - W''(q) U_2 =   z \|A\|^2 q'(z) - W_2 =  \frac{1}{2} z \|A\|^2 q'(z) - c(x,t) q'(z).$$
In particular, multiplying by $q'$ and integrating over $\R$ in $z$ shows that $c(x,t) = 0$. 
We then deduce that  
$$ U_2 = \frac{1}{2} \|A\|^2 \eta_1(z),$$
where $\eta_1$ is defined as the solution of 
$$ \eta_1''(s) - W''(q(s))\eta_1(s) = s q'(s), \quad \text{with} \quad \lim_{s \to \pm \infty} \eta_1(s) = 0.$$ 
In conclusion, we have
$$ W_2 =  \frac{1}{2} \|A\|^2 z q'(z) \quad \text{and} \quad  U_2 =  \frac{1}{2} \|A\|^2 \eta_1(z).$$

\noindent {{\it \underline{Fourth order and estimation of the velocity $V$}}}: \\
 We can now explicit the   term of order $1$ in $\varepsilon$ of ${\cal B}(u)$. Indeed we have 
 $ U(z,x,t) = q(z) + O(\varepsilon^2)$, and as  $\nabla_x q(z) = 0$,  we have
\begin{eqnarray*}
 {\cal B}(u) &=& \left[ \mdiv \left( \mdiv\left(\frac{\nabla u}{|\nabla u|}\right)\frac{\nabla u}{|\nabla u|} \right)  - \mdiv \left(  \nabla \left( \frac{\nabla u}{|\nabla u|} \right)\frac{\nabla u}{|\nabla u|}\right) \right] \\
      &=& \mdiv \left( \Delta d \nabla d \right) - \mdiv \left( \nabla^2 d \nabla d \right) + O(\varepsilon^2)=  \left(\sum_{i}  \frac{\kappa_i(\pi(x))}{1 -  z \varepsilon  \kappa_i(\pi(x))} \right)^2  -  \sum_i \frac{\kappa_i(\pi(x))^2}{(1 - z \varepsilon \kappa_i(\pi(x)))^2} + O(\varepsilon^2) \\
      &=&  \left( H^2 - \|A\|^2 \right) -  \varepsilon 2 z  \left( H \|A\|^2 - \Theta^3 \right) + O(\varepsilon^2),  
\end{eqnarray*}
where $\Theta^3 = \sum_{i} k_i(\pi(x))^3$. \\

The fourth order of the first equation now reads
\begin{eqnarray*}
  - V q' &=&  \left[ \partial^2_{z} W_3  -  W''(q) W_3 \right] - W^{(3)}(q) U_2 W_1  + \left( H \partial_z W_2 - \|A\|^2 z \partial_z W_1 \right) + \Delta_x W_1 -  z W'(q) 2 \left(H \|A\|^2 - \Theta^3 \right)   \\ 
       &=&   \left[ \partial^2_{z} W_3  -  W''(q) W_3 \right] + \frac{1}{2} W^{(3)}(q) H \|A\|^2 \eta_1 q' + \frac{1}{2} \|A\|^2 H \left(3 z q'' + q'\right)- \left( \Delta_{\Gamma} H \right) q' - 2 \left( H \|A\|^2 - \Theta^3 \right)z q'' \\
       &=&  \left[ \partial^2_{z} W_3  -  W''(q) W_3 \right]  + \frac{1}{2} W^{(3)}(q) H \|A\|^2 \eta_1 q' + \left(-\frac{1}{2}\|A\|^2 H + 2 B^3 \right)z q'' + \frac{1}{2} \|A\|^2 Hq'  - \left( \Delta_{\Gamma} H \right) q'
\end{eqnarray*}
Multiplying by $q'$ and integrating over $\R$ leads to 
$$ V = - \frac{1}{S} \left[ \left( \frac{1}{2} \|A\|^2 H S +  \left(-\frac{1}{2}\|A\|^2 H + 2 \Theta^3 \right) \int_{\R} z q'' q' dz + \frac{1}{2} \|A\|^2 H \int_{\R} W^{(3)}(q) \eta_1 (q')^2 dz  \right) - \Delta_{\Gamma} H S \right],$$ 
where $S = \int_{\R} q'(z)^2 dz$. \\
Remark also that 
$$ \int_{\R} z q'' q' dz = \frac{1}{2} \int_{\R} z ((q')^2)' dz = - \frac{1}{2} \int_{\R} (q')^2 dz = -\frac{1}{2} S.$$
Moreover, recall that $\eta_1$ satisfies 
$$ \begin{cases}
    \eta_1^{''} - W''(q)\eta_1 = z q' \\
    \eta_1^{'''} - W''(q) \eta_1' -  W^{(3)}(q)q' \eta_1 = (z q')',
   \end{cases}
$$
then we have
\begin{eqnarray*}
 \int_{\R}  W^{(3)}(q) \eta_1 (q')^2 dz &=& \int_{\R} \left( \eta_1^{'''} - W''(q)\eta_1' \right) q' dz - \int_{\R} (z q')' q' dz = -\frac{1}{2}S,
\end{eqnarray*}
and we conclude that  
$$ V = \Delta_{\Gamma} H + \Theta^3 - \frac{1}{2} \|A\|^2 H.$$

\subsection{Approximating the Willmore flow with Esedoglu-R\"atz-R\"oger's energy}

We now consider the following variant of the Esedoglu-R\"atz-R\"oger's energy, which we introduced in Section~\ref{sec:err}
$$\cW^{\mbox{\tiny EsR\"aR\"o}}_\varepsilon(u)= \frac{1}{2\varepsilon}\int_\Omega \left( \varepsilon\Delta u-\frac{W'(u)}{\varepsilon} \right)^2dx 
+  \beta J_{\varepsilon}(u)$$
where the penalization term $J_{\varepsilon}(u)$ reads 
$$
 J_{\varepsilon}(u) =  \frac 1{\varepsilon^{1+\alpha}}\int_\Omega \left( \varepsilon \nabla^2 u : {\mathbb N}(u) -\frac{W'(u)}{\varepsilon} \right)^2dx, \quad \text{ and } \quad {\mathbb N}(u)= \frac{\nabla u}{|\nabla u|} \otimes \frac{\nabla u}{|\nabla u|}. 
$$
The aim of this section is to derive and study the PDE obtained as the $\lp{2}$-gradient flow of $\cW^{\mbox{\tiny EsR\"aR\"o}}_\varepsilon(u)$.
We will show that the flow is equivalent to the phase field system
\begin{equation}
\begin{cases}
     \varepsilon^2 \partial_t u = \Delta \mu -  \frac{1}{\varepsilon^2} W''(u) \mu - \beta \widetilde{L}(u) \\
     \mu =  W'(u) - \varepsilon^2 \Delta u. \\
     \xi_{\varepsilon} = \varepsilon \nabla^2 u : {\mathbb N}(u) -\frac{W'(u)}{\varepsilon}  \\
     \widetilde{L}(u) =  2\varepsilon^{1 - \alpha} \left[   \left( {\mathbb N}(u): \nabla^2 \xi_{\varepsilon} - \frac{1}{\varepsilon^2} W''(u)  \xi_{\varepsilon}  \right) +  2\left(  \mdiv\left( \frac{\nabla u}{|\nabla u|} \right)\frac{\nabla u}{|\nabla u|} \right) \cdot \nabla \xi_{\varepsilon} + {\cal B}(u)\xi_{\varepsilon} \right]\\
    {\cal B}(u) =    \mdiv \left( \mdiv\left(\frac{\nabla u}{|\nabla u|}\right)\frac{\nabla u}{|\nabla u|} \right) -  \mdiv \left(  \nabla \left( \frac{\nabla u}{|\nabla u|} \right)\frac{\nabla u}{|\nabla u|}\right)
\end{cases}\label{flow:err}\end{equation}

Note that this system coincides with the classical one, up to the addition of a penalty term $-\beta \widetilde{L}(u)$. \\

The well-posedness  of the phase field model \eqref{flow:err} at fixed parameter $\varepsilon$ is open, and requires presumably a regularization 
as done numerically in \cite{Esedoglu_12}.\\

By formal arguments involving matched asymptotic expansions again, we will show that this approximating flow is expected to converge, as $\varepsilon$ goes to zero, to the Willmore flow in dimension $N\geq 2$, 
at least whenever $\alpha=0$ or $\alpha=1$. More precisely we will show the
\begin{claim}
In a suitable regime provided by the method of matched asymptotic expansions, the normal velocity of the $\frac 1 2$-front $\Gamma(t)=\partial E(t)$ associated with a solution $(u^\alpha_{\varepsilon}, \mu^\alpha_{\varepsilon},\xi^\alpha_\varepsilon)$ to Esedoglu-R\"atz-R\"oger's phase field model~\eqref{flow:err} in both cases $\alpha=0$ and $\alpha=1$ is the Willmore velocity
$$V =  \Delta_{\Gamma} H + \|A\|^2 H - \frac{H^3}{2}.$$
In addition, for $\alpha=0$:
$$\left\{\begin{array}{l}
  u^0_{\varepsilon}(x,t)=q \left( \frac{d(x,E(t))}{\varepsilon}\right) + \varepsilon^2 \frac{\|A\|^2-H^2}{1+2\beta} \eta_1\left( \frac{d(x,E(t)}{\varepsilon}\right)+O(\varepsilon^3)\\
\mu^0_{\varepsilon}(x,t)=- \varepsilon  H q' \left( \frac{d(x,E(t)}{\varepsilon}\right) + \varepsilon^2   \left(  H^2 - 2 \beta \left[ \frac{2 \|A\|^2 - H^2}{1 + 2 \beta} \right] \right) \eta_2\left( \frac{d(x,E(t)}{\varepsilon}\right)+O(\varepsilon^3),\\
\xi^0_\varepsilon(x,t)=\varepsilon\left( \frac{2 \|A\|^2 - H^2}{1 + 2 \beta} \right) \eta_2\left( \frac{d(x,E(t)}{\varepsilon}\right)+O(\varepsilon^2)\end{array}\right.$$
where $\eta_2(z)=zq'(z)$ is a profile function. For $\alpha=1$:
$$\left\{\begin{array}{l}
  u^1_{\varepsilon}(x,t)=q \left( \frac{d(x,E(t))}{\varepsilon}\right) + O(\varepsilon^3)\\
\mu^1_{\varepsilon}(x,t)=- \varepsilon  H q' \left( \frac{d(x,E(t)}{\varepsilon}\right) + 2\varepsilon^2  \|A\|^2 \eta_2\left( \frac{d(x,E(t)}{\varepsilon}\right)+O(\varepsilon^3),\\
\xi^1_\varepsilon(x,t)=\varepsilon^2\frac{(2 \|A\|^2  - H^2)}{ 4 \beta } \eta_2\left( \frac{d(x,E(t)}{\varepsilon}\right)+O(\varepsilon^3)\end{array}\right.$$
\end{claim}
\boss
The previous claim gives indications on the design of a numerical scheme for simulating the Esedoglu-R\"atz-R\"oger's flow in the cases $\alpha=0,1$. Clearly, the flow acts at the second order for $u$ in the case $\alpha=0$, and not less than at the third order (at least) whenever $\alpha=1$. This implies that capturing with accuracy the motion of the interface should be much more delicate when $\alpha=1$.
\eoss

\subsubsection{Derivation of the $\lp{2}$-gradient flow of $\cW^{\mbox{\tiny EsR\"aR\"o}}_\varepsilon(u)$}
\begin{prop}
The $\lp{2}$-gradient flow of Esedoglu-R\"atz-R\"oger's 's model is equivalent to
$$\left\{\begin{array}{lll}
     \varepsilon^2 \partial_t u &=& \Delta \mu -  \frac{1}{\varepsilon^2} W''(u) \mu - \beta \widetilde{L}(u) \\
     \mu &=&  W'(u) - \varepsilon^2 \Delta u. \\
     \xi_{\varepsilon} &=& \varepsilon \nabla^2 u : {\mathbb N}(u) -\frac{W'(u)}{\varepsilon}  \\
     \widetilde{L}(u) &=&  2\varepsilon^{1 - \alpha} \left[   \left( {\mathbb N}(u): \nabla^2 \xi_{\varepsilon} - \frac{1}{\varepsilon^2} W''(u)  \xi_{\varepsilon}  \right) +  2\left(  \mdiv\left( \frac{\nabla u}{|\nabla u|} \right)\frac{\nabla u}{|\nabla u|} \right) \cdot \nabla \xi_{\varepsilon} + {\cal B}(u)\xi_{\varepsilon} \right]\\
    {\cal B}(u) &=&    \mdiv \left( \mdiv\left(\frac{\nabla u}{|\nabla u|}\right)\frac{\nabla u}{|\nabla u|} \right) -  \mdiv \left(  \nabla \left( \frac{\nabla u}{|\nabla u|} \right)\frac{\nabla u}{|\nabla u|}\right)
\end{array}\right.$$
\end{prop}
\dimo
Since 
$$ \xi_{\varepsilon}(u) = \varepsilon \nabla^2 u : {\mathbb N}(u)  -\frac{W'(u)}{\varepsilon} = \varepsilon  \left( \nabla^2 u \frac{\nabla u}{|\nabla u|} \right) \cdot \frac{\nabla u}{|\nabla u|} - \frac{1}{\varepsilon} W'(u), $$
one has that
$$  \xi'_{\varepsilon}(u)(w) =  \varepsilon \left( \nabla^2 w : {\mathbb N}(u) + 2 \frac{\nabla^2 u : \nabla w \otimes \nabla u  }{|\nabla u|^2} - 2 \nabla^2 u : {\mathbb N}(u) \frac{\nabla u \cdot \nabla w}{|\nabla u|^2} \right) - \frac{1}{\varepsilon} W''(u) w.$$
The gradient of $J_{\varepsilon}(u)$ follows, recalling that $P^u = I_d -{\mathbb N}(u)$:
\begin{eqnarray*}
 \nabla  J_{\varepsilon}(u) &=& \frac{2}{\varepsilon^{\alpha}} \left[ \nabla^2 : [{\mathbb N}(u)\xi_{\varepsilon}] - \frac{1}{\varepsilon^2} W''(u) \xi_{\varepsilon} - 2 \mdiv \left( \frac{\xi_{\varepsilon}(u) \nabla^2 u \nabla u }{|\nabla u|^2} \right) 
+ 2 \mdiv \left( \frac{\nabla^2 u : {\mathbb N}(u)\xi_{\varepsilon} \nabla u}{|\nabla u|^2} \right) \right] \\
                        &=&  \frac{2}{\varepsilon^{\alpha}} \left[ \nabla^2 : [{\mathbb N}(u)\xi_{\varepsilon}] - \frac{1}{\varepsilon^2} W''(u) \xi_{\varepsilon} - 2 \mdiv \left( \frac{\xi_{\varepsilon}(u) P^u \nabla^2 u \nabla u }{|\nabla u|^2} \right) \right]
\end{eqnarray*}
More precisely, using 
$$ \nabla^2 : [{\mathbb N}(u)\xi_{\varepsilon}] = (\nabla^2 : {\mathbb N}(u)) \xi_{\varepsilon} + 2 \mdiv({\mathbb N}(u)) \cdot \nabla \xi_{\varepsilon} + {\mathbb N}(u): \nabla^2 \xi_{\varepsilon},$$ 
\begin{eqnarray*}
   \mdiv({\mathbb N}(u)) &=&  \left( \mdiv\left( \frac{\nabla u}{|\nabla u|} \right) \frac{\nabla u}{|\nabla u|} +  \nabla \left[ \frac{\nabla u}{|\nabla u|} \right] \frac{\nabla u}{|\nabla u|}\right) \\
                                           &=&  \left(  \mdiv\left( \frac{\nabla u}{|\nabla u|} \right)\frac{\nabla u}{|\nabla u|}  +   \frac{P^u \nabla^2 u \nabla u }{|\nabla u|^2}  \right)
\end{eqnarray*}
and
$$
(\nabla^2 : {\mathbb N}(u)) = \mdiv ( \mdiv {\mathbb N}(u)) =  \mdiv \left(  \mdiv\left( \frac{\nabla u}{|\nabla u|} \right)\frac{\nabla u}{|\nabla u|}  +   \frac{P^u \nabla^2 u \nabla u }{|\nabla u|^2}  \right),
$$
one gets
\begin{eqnarray*}
  \nabla  J_{\varepsilon}(u) &=& \frac{2}{\varepsilon^{\alpha}} \left[  \left( {\mathbb N}(u): \nabla^2 \xi_{\varepsilon} - \frac{1}{\varepsilon^2} W''(u)  \xi_{\varepsilon}  \right) +  2\left(  \mdiv\left( \frac{\nabla u}{|\nabla u|} \right)\frac{\nabla u}{|\nabla u|} \right) \cdot \nabla \xi_{\varepsilon} + {\cal B}(u)\xi_{\varepsilon} \right].
\end{eqnarray*}
from which the $\lp{2}$-gradient flow of Esedoglu-R\"atz-R\"oger's 's model follows.\qed

\subsubsection{Asymptotic analysis} ~\\
In this section, we perform the formal method of matched asymptotic expansions 
to the system solution $(u_{\varepsilon},\mu_{\varepsilon}, \xi_{\varepsilon})$ of 
 $$  \begin{cases}
     \varepsilon^2 \partial_t u = \Delta \mu -  \frac{1}{\varepsilon^2} W''(u) \mu - \beta \widetilde{L}(u) \\
     \mu =  W'(u) - \varepsilon^2 \Delta u. \\
     \xi_{\varepsilon} = \varepsilon \nabla^2 u : {\mathbb N}(u) -\frac{W'(u)}{\varepsilon}  \\
     \widetilde{L}(u) =  2\varepsilon^{1-\alpha} \left[   \left( {\mathbb N}(u): \nabla^2 \xi_{\varepsilon} - \frac{1}{\varepsilon^2} W''(u)  \xi_{\varepsilon}  \right) +  2\left(  \mdiv\left( \frac{\nabla u}{|\nabla u|} \right)\frac{\nabla u}{|\nabla u|} \right) \cdot \nabla \xi_{\varepsilon} + {\cal B}(u)\xi_{\varepsilon} \right],
    \end{cases}$$
in both cases $\alpha=0$ and $\alpha=1$. \\

As previously, we assume that the $1/2$-isolevel set of  $u_{\varepsilon}$ is a smooth $(n-1)-$dimensional interfaces $\Gamma(t)$ 
defined  as the boundary of a set $E(t) = \left\{ x \in \R^d ; u_{\varepsilon}(x,t) \geq 1/2 \right\}$. 

We assume that there exist outer expansions of $u_{\varepsilon}$, $\mu_{\varepsilon}$ and $\xi_{\varepsilon}$ far from the front $\Gamma$ of the form 
$$ \begin{cases}
    u_{\varepsilon}(x,t) = u_0(x,t) + \varepsilon u_1(x,t) + \varepsilon^2 u_2(x,t) + O(\varepsilon^3)  \\
  \mu_{\varepsilon}(x,t) = \mu_0(x,t) + \varepsilon \mu_1(x,t) + \varepsilon^2 \mu_2(x,t) + O(\varepsilon^3) \\
  \xi_{\varepsilon}(x,t) = \frac{1}{\varepsilon} \xi_{-1}(x,t) +  \xi_{0}(x,t) + \varepsilon \xi_{1}(x,t) + \varepsilon^2 \xi_{2}(x,t) + O(\varepsilon^3) 
   \end{cases}
$$
Considering the stretched variable $z = \frac{d(x,t)}{\varepsilon}$ on  a small neighborhood of $\Gamma$,
we also look for inner expansions of $u_{\varepsilon}(x,t)$, $\mu_{\varepsilon}(x,t)$ and $\xi_{\varepsilon}(x,t)$ of the form 
$$\begin{cases}
   u_{\varepsilon}(x,t) &= U(z,x,t) = U_0(z,x,t) + \varepsilon U_1(z,x,t) + \varepsilon^2 U_2(z,x,t) +O(\varepsilon^3)  \\
   \mu_{\varepsilon}(x,t) &= W(z,x,t) = W_0(z,x,t) + \varepsilon W_1(z,x,t) + \varepsilon^2 W_2(z,x,t) +O(\varepsilon^3)  \\
    \xi_{\varepsilon}(x,t) &= \varPhi(z,x,t) =  \varepsilon^{-1}\varPhi_{-1}(z,x,t) +  \varPhi_0(z,x,t) + \varepsilon \varPhi_1(z,x,t)+\varepsilon^2 \varPhi_2(z,x,t)+O(\varepsilon^3)  
  \end{cases}
$$
In particular, remark that the third equation of \eqref{flow:err} yields
$$ \varPhi(z,x,t) =  \frac{1}{\varepsilon} \left(  \frac{ \partial^2_{zz} U}{\left(1  + \varepsilon^2 \frac{ |\nabla_x U|^2  }{ (\partial_z U)^2 }\right)}  - W'(U)\right) + \varepsilon \left( \frac{\partial_z \left( |\nabla_x U|^2 \right)}{\partial_z U} \right)\left(1  + \varepsilon^2 \frac{ |\nabla_x U|^2  }{ (\partial_z U)^2 }\right)^{-1}.$$

As before, it can be observed for the outer expansions that
$$u_0(x,t) = 
\begin{cases}
 1 & \text{if } x \in E(t) \\
 0 & \text{otherwise}
\end{cases}, \quad \text{and} \quad   u_1 = u_2 = u_3 = \mu_0 = \mu_1= \mu_2 = \xi_{-1} = \xi_0 = \xi_2 =  0.  
$$
The matching conditions imply the following boundary conditions on the inner expansions: 
$$
\begin{cases}
  \lim_{z \to + \infty} U_0(z,x,t) = 0 \\
\lim_{z \to - \infty} U_0(z,x,t) = 1
\end{cases}, \lim_{z \to \pm \infty} U_i(z,x,t) = 0 \text{ for }  i \in  \{ 1,2\},  
$$
and
$$
\lim_{z \to \pm \infty} W_i(z,x,t) = 0, \text{ for }  i \in  \{0,1,2\}, \quad \text{and} \quad \lim_{z \to \pm \infty} \varPhi_i(z,x,t) = 0, \text{ for } i \in \{-1,0,1,2\}.
$$
~\\

\paragraph{Inner expansion with $\alpha = 0$ :} ~\\
This paragraph is devoted to the derivation of the expression of the inner expansion in the  special case $\alpha = 0$.\\

\noindent {\it \underline{First order:}}
\par We have the following system 
$$ \begin{cases}
   0 &= \partial^2_z W_{0} - W''(U_0) W_0  - 2 \beta \left[ \partial_{zz} \varPhi_{-1} - W''(U_0) \varPhi_{-1} \right] \\
   W_0 &=  W'(U_0)-\partial^2_z U_0  \\
   \varPhi_{-1} &= \partial^2_z U_{0} - W'(U_0) 
   \end{cases},
$$
which admits the solution triplet
$$ U_0 = q(z), \quad W_0 = 0, \quad \text{and} \quad  \varPhi_{-1} = 0.
$$

\noindent {\it \underline{Second order}}:
\par At second order, we obtain
$$ \begin{cases}
   0 &= \partial^2_z W_{1} - W''(q) W_1  - 2 \beta \left[ \partial_{zz} \varPhi_{0} - W''(q) \varPhi_{0} \right] \\
   W_1 &=  W''(q) U_1 - \partial^2_z U_1 - H q'  \\
   \varPhi_{0} &= \partial^2_z U_{1} - W''(q) U_{1} 
   \end{cases},
$$
whose solution is given by 
$$ U_1 = 0, \quad W_1 = - H q', \quad \text{and} \quad  \varPhi_{0} = 0.
$$

\noindent {\it \underline{Third order}}:
\par At third order,
$$ \begin{cases}
   0 &= \partial^2_z W_{2} - W''(U_0) W_2  +  H \partial_z W_1  - 2 \beta \left[ \partial_{zz} \varPhi_{1} - W''(q) \varPhi_{1} \right] \\
   W_2 &=  W''(q) U_2 - \partial^2_z U_2 + \|A\|^2 z q'  \\
   \varPhi_{1} &= \partial^2_z U_{2} - W''(q) U_{2} 
   \end{cases},
$$
and we are now looking for a system of solutions of the form 
$$ W_2 = c_W(x,t) \eta_2(z), \quad U_2 = c_U(x,t) \eta_1(z), \quad \text{and} \quad \varPhi_1 = c_{\varPhi}(x,t) \eta_2(z),$$
where the two profiles $\eta_1$ and $\eta_2$ are solutions, respectively, of 
$$ \eta_1'' - W''(q)\eta_1 = z q' \quad \text{and} \quad \eta_2'' - W''(q)\eta_2 = q''.$$
Furthermore, the first equation gives
\begin{eqnarray*}
  0 = c_W \left( \eta_2'' - W''(q)\eta_2 \right) - H^2 q'' - 2c_{\varPhi}\left( \eta_2'' - W''(q)\eta_2 \right)  =  \left(  c_W  - H^2   - 2 \beta c_{\varPhi} \right)  q'',
\end{eqnarray*}
the second equation implies that
\begin{eqnarray*}
 \frac{1}{2} c_W z q' = - c_U \left(  \eta_1'' - W''(q)\eta_1 \right) + \|A\|^2 z q' = (- c_U + \|A\|^2)z q'.  
\end{eqnarray*}
and the third equation shows that
\begin{eqnarray*}
 \frac{1}{2} c_{\varPhi} z q'  = c_U \left(  \eta_1'' - W''(q)\eta_1 \right) = (c_U)z q'.  
\end{eqnarray*}
This provides a linear system  
$$ c_W  - \beta  2 c_{\varPhi} = H^2, \quad c_W  + 2 c_U = 2 \|A\|^2 \quad \text{and} \quad  c_{\varPhi} = 2 c_U,
$$
which admits as solutions
$$ c_W = H^2 - 2 \beta \left[ \frac{2 \|A\|^2 - H^2}{1 + 2 \beta} \right], \quad c_U = \frac{\|A\|^2 - H^2/2}{1 + 2 \beta}, \quad \text{and} \quad c_{\varPhi} = \frac{2 \|A\|^2 - H^2}{1 + 2 \beta}.$$
Therefore,
$$ \begin{cases}
     W_2 =  \left(  H^2 - 2 \beta \left[ \frac{2 \|A\|^2 - H^2}{1 + 2 \beta} \right] \right) \eta_2, \\
     U_2 = \left( \frac{\|A\|^2 - H^2/2}{1 + 2 \beta} \right) \eta_1, \\
     \varPhi_1 =   \left( \frac{2 \|A\|^2 - H^2}{1 + 2 \beta} \right) \eta_2.
   \end{cases}
$$

\noindent {\it \underline{Fourth order and estimation of the velocity $V$}}:
\par The fourth order reads as follows
\begin{eqnarray*}
  - V q' &=&  \left[ \partial^2_{z} W_3  -  W''(q) W_3 \right] - W^{(3)}(q) U_2 W_1  + \left( H \partial_z W_2 - \|A\|^2 z \partial_z W_1 \right) + \Delta_x W_1  \\
          &~& \quad \quad \quad \quad  - 2 \beta \left( \left[ \partial_{zz} \varPhi_{2} - W''(q) \varPhi_{2} \right]  + 2 H \partial_z \varPhi_{1}    \right) \\
        &=& \left[ \partial^2_{z} ( W_3 - 2 \beta \varPhi_{2})   -  W''(q) (W_3 - 2 \beta \varPhi_{2}) \right] + \left( \frac{H \|A\|^2 - H^3/2}{1 + 2 \beta} \right)  W^{(3)}(q) \eta_1 q' \\
        &~& \quad \quad \quad  - \Delta_{\Gamma} H q' + \left( H^3/2   + \|A\|^2 H - H \beta  \left( \frac{2 \|A\|^2 - H^2}{1 + 2 \beta} \right)  \right)z q'' + \left( H^3/2 - H \beta  \left( \frac{2 \|A\|^2 - H^2}{1 + 2 \beta} \right)   \right) q'
\end{eqnarray*}
Recalling that
$$ 
\begin{cases}
 \int_{\R} (q'(z))^2 dz = S \\
 \int_{\R} z q'' q' dz = -\frac{1}{2} S \\
 \int_{\R}  W^{(3)}(q) \eta_1 (q')^2 dz = -\frac{1}{2}S,
\end{cases}
$$
multiplying the last equation by $q'$ and integrating over $\R$ leads to 
$$  V =  \Delta_{\Gamma} H  + \frac{\|A\|^2 H}{2} - \frac{H^3}{4}  +  \frac{ H \|A\|^2 - H^3/2}{2(1 + 2 \beta)} (1 + 2 \beta ) = \Delta_{\Gamma} H +    \|A\|^2 H - \frac{1}{2} H^3
$$
~\\

\paragraph{Inner expansion with $\alpha = 1$ :} ~\\
We are now looking for the inner expansion of the PDE system~\eqref{flow:err} in the  special case $\alpha = 1$ : \\

\noindent{\it \underline{First order}}: ~\\
We have the following system 
$$ \begin{cases}
   0 &= - 2 \beta \left[ \partial_{zz} \varPhi_{-1} - W''(U_0) \varPhi_{-1} \right] \\
   W_0 &=  W'(U_0)-\partial^2_z U_0  \\
   \varPhi_{-1} &= \partial^2_z U_{0} - W'(U_0) 
   \end{cases},
$$
whose solution is given by 
$$ U_0 = q(z), \quad W_0 = 0, \quad \text{and} \quad  \varPhi_{-1} = 0.
$$

\noindent{\it \underline{Second order}}: ~\\
At second order
$$ \begin{cases}
   0 &=   \partial^2_{zz} W_{0} - W''(q) W_0  - 2 \beta \left[ \partial^2_{zz} \varPhi_{0} - W''(q) \varPhi_{0} \right] \\
   W_1 &=  W''(q) U_1 - \partial^2_{zz} U_1 - H q'  \\
   \varPhi_{0} &= \partial^2_{zz} U_{1} - W''(q) U_{1} 
   \end{cases},
$$
which admits as solution
$$ U_1 = 0, \quad W_1 = - H q', \quad \text{and} \quad  \varPhi_{0} = 0.
$$

\noindent{\it \underline{Third order}}: ~\\
At third order
$$ \begin{cases}
   0 &= \partial^2_{zz} W_{1} - W''(q) W_1  - 2 \beta \left[ \partial_{zz} \varPhi_{1} - W''(q) \varPhi_{1} \right] \\
   W_2 &=  W''(q) U_2 - \partial^2_{zz} U_2 + \|A\|^2 z q'  \\
   \varPhi_{1} &= \partial^2_z U_{2} - W''(q) U_{2} 
   \end{cases},
$$
whose solution triplet is
$$  U_2 = 0, \quad W_2 = \|A\|^2 z q', \quad \text{and} \quad  \varPhi_{1} = 0.$$

\noindent{\it \underline{Fourth order}}:\\
From
$$   0 = \partial^2_z W_{2} - W''(q) W_2  +  H \partial_z W_1  - 2 \beta \left[ \partial_{zz} \varPhi_{2} - W''(q) \varPhi_{2} \right],$$
we deduce that 
$$  \left[ \partial_{zz} \varPhi_{2} - W''(q) \varPhi_{2} \right] =  \frac{(2 \|A\|^2 - H^2)}{2 \beta } q'',$$
and then 
$$\varPhi_{2} = \frac{(2 \|A\|^2 - H^2)}{2 \beta } \eta_2 =  \frac{(2 \|A\|^2  - H^2)}{ 4 \beta } z q'.$$

\noindent{\it \underline{Last order and estimation of the velocity $V$}}: ~\\
We have
\begin{eqnarray*}
  - V q' &=&  \left[ \partial^2_{z} W_3  -  W''(q) W_3 \right] - W^{(3)}(q) U_2 W_1  + \left( H \partial_z W_2 - \|A\|^2 z \partial_z W_1 \right) + \Delta_x W_1  \\
          &~& \quad \quad \quad \quad  - 2 \beta \left( \left[ \partial_{zz} \varPhi_{3} - W''(q) \varPhi_{3} \right]  + 2 H \partial_z \varPhi_{2} \right) \\
        &=& \left[ \partial^2_{z} ( W_3 - 2 \beta \varPhi_{3})   -  W''(q) (W_3 - 2 \beta \varPhi_{3}) \right]  - \Delta_{\Gamma} H q' \\
        &~& \quad \quad \quad  + \left(2 \|A\|^2 H  -  H \beta  \left( \frac{2 \|A\|^2 - H^2}{\beta} \right)  \right)z q'' + \left( \|A\|^2 H  - H \beta  \left( \frac{2 \|A\|^2 - H^2}{ \beta} \right)   \right) q' \\
        &=& \left[ \partial^2_{z} ( W_3 - 2 \beta \varPhi_{3})   -  W''(q) (W_3 - 2 \beta \varPhi_{3}) \right]  - \Delta_{\Gamma} H q'  +  H^3 z q'' + \left( - \|A\|^2 H  + H^3 \right) q'. \\
\end{eqnarray*}

As previously, this shows that the velocity $V$ of the interfaces equals the Willmore velocity
$$V =  \Delta_{\Gamma} H + \|A\|^2 H - \frac{H^3}{2}.$$

\boss
The analysis of the asymptotic behavior for $\alpha$ non integer is more delicate because it requires studying non integer orders of $\varepsilon$ and it is far from being clear how integer and non integer scales may combine. As for integer values of $\alpha>1$, a careful study at higher orders of $\varepsilon$ should be possible but is out of the scope of the present paper.
\eoss

\section{2D and 3D numerical simulations for the classical and Mugnai's diffuse flows}\label{sec:num}
There is an important literature on numerical methods for the approximation of interfaces evolving by a geometric law. 
They can be roughly classified into three categories: parametric methods  \cite{Deckelnick1999,Rusu_elastic_flow,Deckelnick2005,Barrett2008_2,Barrett:2008,Polthier-Hab,HsuKusnerSullivan92}, level-set formulations  \cite{OsherSethian,Osherbook1,Osherbook2,Evans_spruck,chen_giga_goto}, phase-field approaches \cite{Modica1977,Chen1992,bellettini_1996,Paolini_anisotropic}. See for instance \cite{Deckelnick2005} for a complete review (with a particular emphasis on the mean curvature flow, but fourth-order flows are also addressed) and a comparison between the different strategies.
In the context of fourth order geometric evolution equations, in particular the Willmore flow,
parametric approaches have been proposed in \cite{Dziuk:2008,Barrett:2007:PFE:1225314.1225638,BarrettGN08}  for curves and surfaces
using a semi-implicit finite element method. In \cite{OlRu10}, a fully implicit approach
 via a variational formulation is also analyzed for the approximation of anisotropic Willmore flow.  
The level set methods have been applied for the first time in \cite{Droske04alevel}.
Concerning the phase field approach,  semi-implicit schemes including standard finite element differences, finite elements,
and Fourier spectral methods are developed in \cite{Du_Wang_2004,Du_wang:2006,Esedoglu_12} and analyzed in \cite{Du_convergenceof}.
A fully implicit scheme coupled with a finite element method has been more recently introduced in \cite{FrRuWi11} via a variational formulation.
An adaptation to fourth order geometric evolution equations of the Bence-Merriman-Osher algorithm \cite{BenceMerrimanOsher} is also proposed in
\cite{BMO_willmore}. Let us finally mention the discrete methods involving surface triangulations and discrete curvature operators~\cite{HsuKusnerSullivan92,BobenkoSchroeder05,Wardetzky-et-al-07}.\\

In this paper, we will consider a quite different and new scheme to solve both the classical and Mugnai's phase field systems~\eqref{flow:classic} and~\eqref{flow:mugnai}. The simulations can be compared with those obtained by Esedoglu, R\"atz and R\"oger in~\cite{Esedoglu_12} for their phase field system~\eqref{flow:err} (actually a variant of it, see Section~\ref{sec:err}), and for Bellettini's phase field system~\eqref{flow:bel}.

Here, we use an implicit scheme to ensure the decreasing of the
diffuse Willmore energy, and a Fourier spectral method in order to get high accuracy approximation in space. 
At each step time, it is necessary to solve a nonlinear equation. A Newton algorithm like in  \cite{FrRuWi11} 
appears to be very efficient in practice, but not in accordance with a Fourier spectral discretization, so we opted for a fixed point approach. \\

\subsection{New numerical schemes for the approximation of classical and Mugnai's flows}
\subsubsection{Classical diffuse approximation flow} 
We introduce a new scheme to approximate numerically some solutions of the phase field system 
$$
 \begin{cases}
      \partial_t u =  \ds\frac{1}{\alpha \varepsilon^2 }\Delta \mu -  \frac{1}{\alpha \varepsilon^4} W''(u) \mu  \\
     \mu =  \alpha W'(u) - \alpha \varepsilon^2 \Delta u,  
    \end{cases} 
 $$
where $\alpha$ is a positive constant. The particular phase field system~\eqref{flow:classic} corresponding to the classical diffuse Willmore flow is obtained for $\alpha=1$. 

We compute the solution for any time $t \in [0,T]$ in a box $\Omega = [0,1]^N$ with periodic boundary conditions. 
We use a Euler semi-implicit discretization in time:  
  $$
 \begin{cases}
  u^{n+1} = \delta_t \left[ \frac{1}{\alpha \varepsilon^2} \Delta \mu^{n+1} -  \frac{1}{\alpha \varepsilon^4} W''(u^{n+1}) \mu^{n+1} \right] + u^{n} \\
  \mu^{n+1} =   \alpha W'(u^{n+1}) - \alpha \varepsilon^2 \Delta u^{n+1},
 \end{cases}
 $$
where $\delta_t$ is the time step, $u^{n}$ and $\mu^n$ are the approximations of the solutions $u$ and $\mu$, respectively, evaluated at time $t_n = n~\delta_t$.  The system can be written as
$$\left\{\begin{array}{l}
u^{n+1} -\frac{\delta_t}{\alpha \varepsilon^2}\Delta \mu^{n+1}=E\\
\mu^{n+1}+\alpha\varepsilon^2\Delta u^{n+1}=F,\end{array}\right.$$
with $E=u^{n}-  \frac{\delta_t}{\alpha \varepsilon^4} W''(u^{n+1}) \mu^{n+1}$, $F= \alpha W'(u^{n+1})$. Therefore,
$$\left\{\begin{array}{l}
u^{n+1}+\delta_t\Delta^2u^{n+1}=E+\frac{\delta_t}{\alpha \varepsilon^2}\Delta F\\
\mu^{n+1}+\delta_t\Delta^2\mu^{n+1}=F-\alpha\varepsilon^2\Delta E.\end{array}\right.$$
Thus, $(u^{n+1},\mu^{n+1})$ is the solution of the nonlinear equation 
\begin{equation}\left (
    \begin{array}{c}
      u^{n+1} \\
      \mu^{n+1}
    \end{array}
 \right)  = \phi\left( \begin{array}{c}
      u^{n+1} \\
      \mu^{n+1}
    \end{array} \right),\label{fixed}\end{equation}
where
 $$
  \phi\left( \begin{array}{c}
      u^{n+1} \\
      \mu^{n+1}
    \end{array} \right)
  = 
 \left(I_d + \delta_t \Delta^2 \right)^{-1} 
 \left( \begin{array}{cc}
      I_d &  \frac{\delta_t}{\alpha \varepsilon^2} \Delta \\
      - \alpha \varepsilon^2 \Delta & I_d
    \end{array}  \right) \left (
    \begin{array}{c}
      u^n -  \frac{\delta_t}{\alpha \varepsilon^4} W''(u^{n+1}) \mu^{n+1}   \\
       \alpha W'(u^{n+1})
    \end{array}
 \right)
 $$
A natural way to approximate the solution $(u^{n+1},\mu^{n+1})$ to~\eqref{fixed} is a fixed point iterative method.

\par The space discretization is built with Fourier series. It has the advantage of preserving a high order approximation in space while allowing a fast and simple processing of the homogeneous operator
$$ {\mathbb G}=\left(I_d + \delta_t \Delta^2 \right)^{-1}\begin{pmatrix}
      I_d &  \frac{\delta_t}{\alpha \varepsilon^2} \Delta \\
      - \alpha \varepsilon^2 \Delta & I_d\end{pmatrix}=\begin{pmatrix}
      I_d &  -\frac{\delta_t}{\alpha \varepsilon^2} \Delta \\
      \alpha \varepsilon^2 \Delta & I_d\end{pmatrix}^{-1}.$$  
In practice, the solutions $u(x, t_n)$ and $\mu(x,t_n)$ at time $t_n = n \delta_t$ are approximated by the truncated Fourier series :
$$ 
u^n_{\cal P}(x) = \sum_{\|p\|_{\infty} \leq {\cal P}} u^n_p e^{2 i\pi x\cdot p}, \quad \text{and} \quad \mu^n_{\cal P}(x) = \sum_{\|p\|_{\infty} \leq {\cal P}} \mu^n_p e^{2 i\pi x\cdot p}
$$
where $\|p\|_{\infty} = \max_{1 \leq i \leq N} |p_i|$, ${\cal P}$ is the maximal number of Fourier modes in each direction, and the coefficients $u^n_p$, $\mu^n_p$ are derived from a prior Fourier decomposition of $\begin{pmatrix}
      u^n -  \frac{\delta_t}{\alpha \varepsilon^4} W''(u^{n+1}) \mu^{n+1}   \\
       \alpha W'(u^{n+1})
    \end{pmatrix}$ combined with an application in the Fourier domain of the operator ${\mathbb G}$. More precisely, the fixed point algorithm that we propose reads as follows:
\begin{algo}\label{algo1}
\begin{itemize}
 \item[] Initialization : $u_0^{n+1} = u^{n}$, $\mu_0^{n+1} = \mu_n$ 
  \item[] { While $\| u^{n+1}_{k+1} - u^{n+1}_{k} \| +  \| \mu^{n+1}_{k+1} - \mu^{n+1}_{k} \|> 10^{-8}$}, perform the loop on $k$:
 \begin{itemize}
      \item[1)] Compute 
$$h^n_{\cal P} = u^n -  \frac{\delta_t}{\alpha \varepsilon^4} W''(u_k^{n+1}) \mu_k^{n+1}, \quad \text{and} \quad \tilde{h}^n_{\cal P} =  \alpha W'(u_k^{n+1})$$
      \item[2)] Using the Fast Fourier Transform, compute the  truncated Fourier series of $h^n_{\cal P}$ and $\tilde h^n_{\cal P}$ : 
$$ h^n_{\cal P}(x) = \sum_{\|p\|_{\infty} \leq {\cal P}} h^n_p ~ e^{2 i\pi x\cdot p}, \quad \text{and} \quad  \tilde{h}^n_{\cal P}(x) = \sum_{\|p\|_{\infty} \leq {\cal P}} \tilde{h}^n_p~ e^{2 i\pi x\cdot p}$$ 
      \item[3)] Compute 
$$ u^{n+1}_{k+1}(x) = \sum_{\|p\|_{\infty} \leq {\cal P}} (u_{k+1})^n_p ~ e^{2 i\pi x\cdot p}, \quad \text{and} \quad \mu^{n+1}_{k+1}(x) = \sum_{\|p\|_{\infty} \leq {\cal P}} (\mu_{k+1})^n_p ~ e^{2 i\pi x\cdot p},$$
where
$$  \begin{cases}
(u_{k+1})^n_p &= \frac{1}{1 + \delta_t (4 \pi^2 |p|)^2 } \left(  h^n_p - \frac{\delta_t}{\alpha \varepsilon^2} 4 \pi^2 |p|^2 \tilde{h}^n_p \right) \\
 (\mu_{k+1})^n_p &=  \frac{1}{1 + \delta_t (4 \pi^2 |p|)^2 }\left(  \tilde{h}^n_p + \alpha \varepsilon^2 4 \pi^2 |p|^2 h^n_p   \right)
\end{cases} $$ 
\end{itemize}
 \item[]  { End}
\end{itemize}
\end{algo}
Note that the semi-implicit scheme 
$$  \begin{cases}
  u^{n+1} = \delta_t \left[  \frac{1}{\varepsilon^2 \alpha}\Delta \mu^{n+1} -  \frac{1}{\varepsilon^4 \alpha} W''(u^{n}) \mu^{n} \right] + u^{n} \\
  \mu^{n+1} =   \alpha W'(u^{n}) - \alpha \varepsilon^2 \Delta u^{n+1},
 \end{cases}$$
implies the following scheme on $u^{n}$
$$   u^{n+1} = \left(I_d  + \delta_t \Delta^2 \right)^{-1} \left[ u^n + \frac{\delta_t}{\varepsilon^2} \Delta W'(u^n) + \frac{\delta_t}{\varepsilon^2}W''(u^n) \left( \Delta u^n - \frac{1}{\varepsilon^2} W'(u^{n-1}) \right) \right],$$
which is expected to be stable under a condition of the form 
$$ \delta_t \leq C \min \left\{ \varepsilon^2 \delta x^2, \varepsilon^4 \right\}$$
where  $\delta x = 1/{(2 \cal P)}$ and $C$ is a constant  depending only the double-well potential $W$. \\
In practice, using fixed point iterations instead of a semi-implicit scheme appears more accurate numerically. This can be justified with the following proposition:
\begin{prop}\label{prop:convalgo1}
Algorithm~\ref{algo1} converges locally under the assumptions
\begin{equation}
\max \left\{ [\alpha M_2]^2+2[\frac{\delta_t}{ \varepsilon^4} M_3 ( M_1 +  N^{3/2} \pi^2 \frac{\varepsilon^2}{\delta_x^{5/2}})]^2,2[\frac{\delta_t}{\alpha \varepsilon^4} M_2]^2  \right\}\|<1\label{eq:cond}
\end{equation}
where $M_i= \sup_{s\in [0,1]} |W^{(i)}(s)|$.
\end{prop}
\dimo We look for the conditions such that
$$ \| D \phi(u^{n+1},\mu^{n+1}) (\delta_u, \delta_\mu) \|^2 < \| (\delta_u, h_{\mu}) \|^2,$$
where the differential of $\phi$ is such that
$$ D \phi(u^{n+1},\mu^{n+1}) (\delta_u, \delta_\mu)  =  \left( \begin{array}{cc}
      I_d & - \frac{\delta_t}{\alpha \varepsilon^2} \Delta \\
      \alpha \varepsilon^2 \Delta & I_d
    \end{array}  \right)^{-1} \left (
    \begin{array}{c}
      -  \frac{\delta_t}{\alpha \varepsilon^4} \left( W^{(3)}(u^{n+1}) \mu^{n+1} \delta_u +  W^{(2)}(u^{n+1}) \delta_\mu \right)    \\
       \alpha W^{(2)}(u^{n+1}) \delta_u     \end{array} \right).
   $$
Note that the eigenvalues of the operator  
$$  \left( \begin{array}{cc}
      I_d & - \frac{\delta_t}{\alpha \varepsilon^2} \Delta \\
      \alpha \varepsilon^2 \Delta & I_d
    \end{array}  \right) $$
are
$$ \lambda_p  = 1 \pm 4 \pi^2 i \sqrt{\delta_t} |p|^2, \quad \text{ for } \|p\|_{\infty} \in [0,\cal P].$$
In particular, this implies that 
$$ \left\|  \left( \begin{array}{cc}
      I_d & - \frac{\delta_t}{\alpha \varepsilon^2} \Delta \\
      \alpha \varepsilon^2 \Delta & I_d
    \end{array}  \right)^{-1}  \right\| \leq 1.$$
Moreover, remark also that 
\begin{eqnarray*}
 | \mu^{n+1}| &=& | \alpha W'(u^{n+1}) - \alpha \varepsilon^2 \Delta u^{n+1}| \leq  \alpha \left( M_1 +  N^{3/2} \pi^2 \frac{\varepsilon^2}{\delta_x^{5/2}} \right).
\end{eqnarray*}
It follows that
$$
\begin{cases}
  |\alpha W^{(2)}(u^{n+1}) \delta_u| \leq \alpha M_2 |\delta_u| \\
  \left| W^{(3)}(u^{n+1}) \mu^{n+1} \delta_u +  W^{(2)}(u^{n+1}) \delta_\mu \right| \leq \left( \alpha M_3( M_1 +  N^{3/2} \pi^2 \frac{\varepsilon^2}{\delta_x^{5/2}}) |\delta_{u}| + M_2 |\delta_\mu|  \right) 
\end{cases}
$$
and then
\begin{multline*}
\left\| \left(   \begin{array}{c}
      -  \frac{\delta_t}{\alpha \varepsilon^4} \left( W^{(3)}(u^{n+1}) \mu^{n+1} \delta_u +  W^{(2)}(u^{n+1}) \delta_\mu \right)    \\
       \alpha W^{(2)}(u^{n+1}) \delta_u     \end{array}  \right) \right\|^2 \\
\leq \max \left\{ [\alpha M_2]^2+2[\frac{\delta_t}{ \varepsilon^4} M_3 ( M_1 +  N^{3/2} \pi^2 \frac{\varepsilon^2}{\delta_x^{5/2}})]^2,2[\frac{\delta_t}{\alpha \varepsilon^4} M_2]^2  \right\}\|  (\delta_u, \delta_\mu) \|^2,\end{multline*}
which concludes the convergence proof of the fixed-point iteration procedure if conditions~\eqref{eq:cond} above are fulfilled.\qed

\subsubsection{Mugnai's flow}\label{sec:mugflow}
We now use a similar scheme for the following generalization of Mugnai's phase field system:
$$
 \begin{cases}
      \partial_t u = \ds\frac{1}{\varepsilon^2 \alpha}\Delta \mu -  \frac{1}{\varepsilon^4 \alpha} W''(u) \mu + \widetilde{\cal B}(u) \\
     \mu =  \ds\alpha W'(u) -  \alpha \varepsilon^2 \Delta u,  
    \end{cases} 
$$
with $\widetilde{\cal B}(u)=\frac{W'(u){\cal B}(u)}{\alpha\varepsilon^4}$. The exact Mugnai's flow~\eqref{flow:mugnai} corresponds to the choice $\alpha=1$ (up to a time rescaling). We use again a semi-implicit scheme:
  $$
 \begin{cases}
  u^{n+1} = \delta_t \left[ \frac{1}{\varepsilon^2 \alpha }\Delta \mu^{n+1} -  \frac{1}{ \alpha \varepsilon^4} W''(u^{n+1}) \mu^{n+1} + \widetilde{\cal B}(u^{n}) \right] + u^{n} \\
  \mu^{n+1} =   \alpha  W'(u^{n+1}) - \varepsilon^2 \alpha \Delta u^{n+1},
 \end{cases}
 $$
where the penalization term  $\widetilde{\cal B}(\cdot)$ is treated explicitly. We use a fixed point iteration to approximate the solution pair $(u^{n+1},\mu^{n+1})$ to the system:
 $$
(u^{n+1},\mu^{n+1})=  \tilde\phi(u^{n+1},\mu^{n+1})=
 \left(I_d + \delta_t \Delta^2 \right)^{-1} 
 \left( \begin{array}{cc}
      I_d &  \frac{\delta_t}{\alpha \varepsilon^2} \Delta \\
      - \alpha  \varepsilon^2 \Delta & I_d
    \end{array}  \right) \left (
    \begin{array}{c}
      u^n -  \frac{\delta_t}{ \alpha \varepsilon^4} W''(u^{n+1}) \mu^{n+1} + \delta_t \widetilde{\cal B}(u^{n})\\
      \alpha W'(u^{n+1})
    \end{array}
 \right).
 $$
For it is highly singular, the penalization term $\widetilde{\cal B}(\cdot)$ needs to be regularized to avoid numerical errors. Observing that
 \begin{eqnarray*}
 \widetilde{\cal B}(u)=  W'(u) \left[ \left( \left| \nabla \left(\frac{\nabla u}{|\nabla u|}\right)\right|^2 -  \left| \mdiv\left(\frac{\nabla u}{|\nabla u|}\right)\right|^2 \right) - \mrot\left(\mrot\left( \frac{\nabla u}{|\nabla u|} \right) \right)\cdot \frac{\nabla u}{|\nabla u|} \right] 
 \end{eqnarray*}
we consider the regularized penalization term 
$$\widetilde{\cal B}_{\sigma}(u) = W'(u) \left[ \left( \left| \nabla \nu_{u,\sigma} \right|^2 -  \left| \mdiv \nu_{u,\sigma} \right|^2 \right) - \mrot\left(\mrot\left(  \nu_{u,\sigma}  \right) \right)\cdot \nu_{u,\sigma}\right]$$
where  $\nu_{u,\sigma} =  \frac{\nabla u}{ \sqrt{|\nabla u|^2 + \sigma^2}}$ with $\sigma$ a small regularization parameter. In particular, the positivity  
$$ \left( \left| \nabla \nu_{u,\sigma} \right|^2 -  \left| \mdiv \nu_{u,\sigma} \right|^2 \right) \geq 0,$$
is ensured, which is in accordance with the continuous case. In practice, finite differences are used for the numerical evaluation of $\widetilde{\cal B}_{\sigma}(u)$. Finally, we propose the following algorithm:
\begin{algo}\label{algo2}
\begin{itemize}
 \item[] Initialization : $u_0^{n+1} = u^{n}$, $\mu_0^{n+1} = \mu_n$ 
  \item[] { While $\| u^{n+1}_{k+1} - u^{n+1}_{k} \| +  \| \mu^{n+1}_{k+1} - \mu^{n+1}_{k} \|> 10^{-8}$, perform the loop on $k$:}
 \begin{itemize}
      \item[1)] Using finite differences to evaluate $\widetilde{\cal B}_{\sigma}(u_k^{n+1})$, compute 
$$h^n_{\cal P} = u^n -  \frac{\delta_t}{\alpha \varepsilon^4} W''(u_k^{n+1}) \mu_k^{n+1}+\delta_t \widetilde{\cal B}_{\sigma}(u^{n}), \quad \text{and} \quad \tilde{h}^n_{\cal P} =  \alpha W'(u_k^{n+1})$$
      \item[2)] Using the Fast Fourier Transform, compute the  truncated Fourier series of $h^n_{\cal P}$ and $\tilde h^n_{\cal P}$ : 
$$ h^n_{\cal P}(x) = \sum_{\|p\|_{\infty} \leq {\cal P}} h^n_p ~ e^{2 i\pi x\cdot p}, \quad \text{and} \quad  \tilde{h}^n_{\cal P}(x) = \sum_{\|p\|_{\infty} \leq {\cal P}} \tilde{h}^n_p~ e^{2 i\pi x\cdot p}$$ 
      \item[3)] Compute 
$$ u^{n+1}_{k+1}(x) = \sum_{\|p\|_{\infty} \leq {\cal P}} (u_{k+1})^n_p ~ e^{2 i\pi x\cdot p}, \quad \text{and} \quad \mu^{n+1}_{k+1}(x) = \sum_{\|p\|_{\infty} \leq {\cal P}} (\mu_{k+1})^n_p ~ e^{2 i\pi x\cdot p},$$
where
$$  \begin{cases}
(u_{k+1})^n_p &= \frac{1}{1 + \delta_t (4 \pi^2 |p|)^2 } \left(  h^n_p - \frac{\delta_t}{\alpha \varepsilon^2} 4 \pi^2 |p|^2 \tilde{h}^n_p \right) \\
 (\mu_{k+1})^n_p &=  \frac{1}{1 + \delta_t (4 \pi^2 |p|)^2 }\left(  \tilde{h}^n_p + \alpha \varepsilon^2 4 \pi^2 |p|^2 h^n_p   \right)
\end{cases} $$ 
\end{itemize}
 \item[]  { End}
\end{itemize}
\end{algo}
\begin{prop}
Algorithm~\ref{algo2} converges locally under the assumptions
\begin{equation}
\max \left\{ [\alpha M_2]^2+2[\frac{\delta_t}{ \varepsilon^4} M_3 ( M_1 +  N^{3/2} \pi^2 \frac{\varepsilon^2}{\delta_x^{5/2}})]^2,2[\frac{\delta_t}{\alpha \varepsilon^4} M_2]^2  \right\}\|<1\label{eq:cond2}
\end{equation}
where $M_i= \sup_{s\in [0,1]} |W^{(i)}(s)|$.
\end{prop}
\dimo Being the penalty term $\widetilde{\cal B}_{\sigma}(u)$ treated explicitly, it does not appear in the expression of the differential of $\tilde\phi$, thus the same proof as for Proposition~\ref{prop:convalgo1} can be used.\qed

\subsection{Numerical simulations of the classical flow}

The following simulations have been realized using Matlab. The isolevel sets 
$\Gamma(t) = \{ x:\, u(x,t) =  \frac 1 2\}$
are computed and drawn using the Matlab functions \verb|contour| in 2D and \verb|isosurface| in 3D. 
We  use the double-well potential $W(s) = \frac{1}{2} s^2 (1-s)^2$ and consider the PDE system
$$
 \begin{cases}
      \partial_t u = \Delta \mu -  \frac{1}{\varepsilon^2} W''(u) \mu  \\
     \mu =  \frac{1}{\varepsilon^2} W'(u) - \Delta u. \\
    \end{cases} 
 $$
with initial conditions $u(x,0)$ and $\mu(x,0)$ of the form
$$ \begin{cases}
    u(x,0) &= \gamma \left( \frac{d(x,E)}{\varepsilon} \right) \\
    \mu(x,0) &= - \frac{1}{\varepsilon}\Delta d(x,E)  \gamma' \left( \frac{d(x,E)}{\varepsilon} \right) 
   \end{cases}
$$
\paragraph{Evolution of a circle}

\begin{figure}[!ht] 
\centering
\includegraphics[width=6cm]{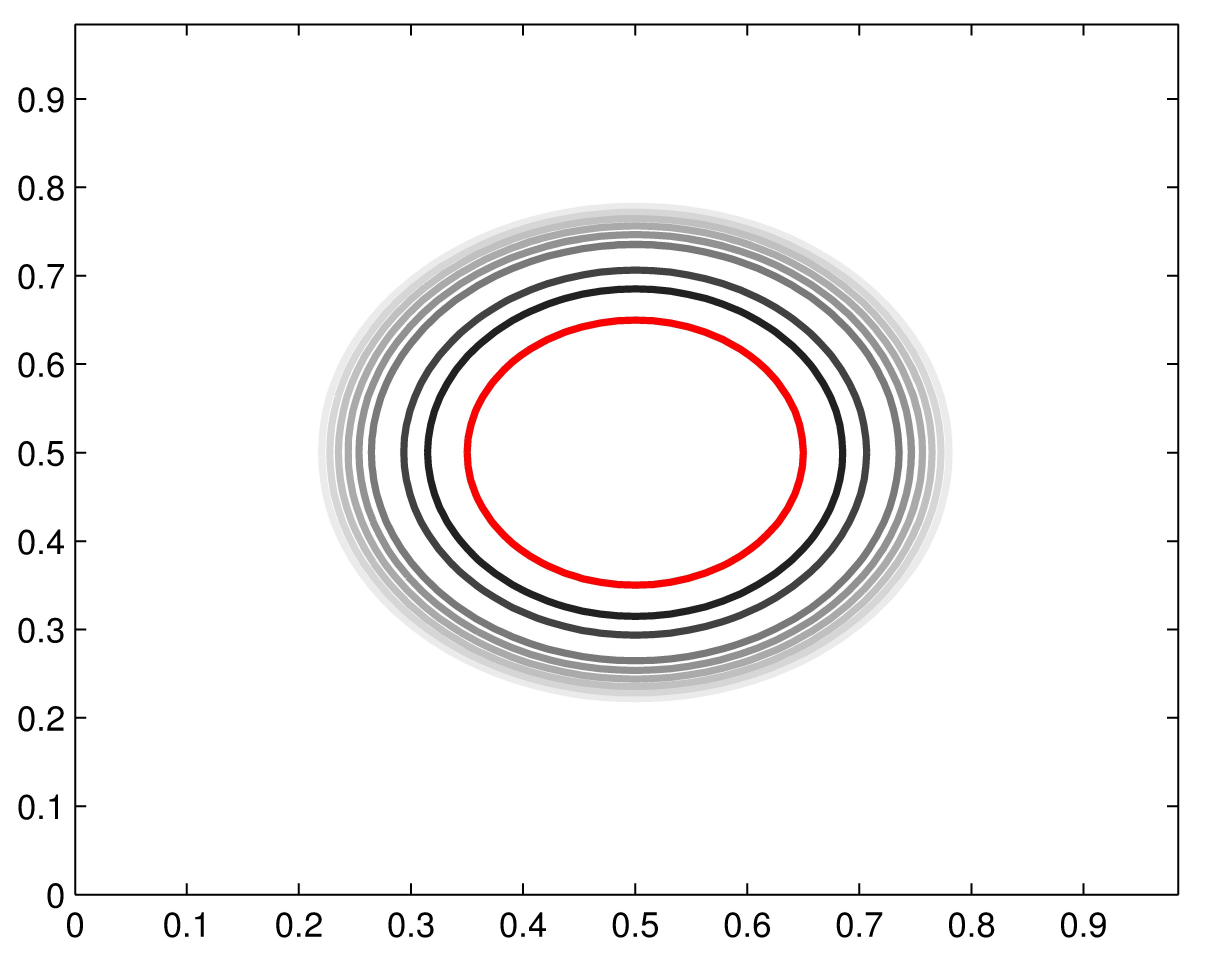}
\includegraphics[width=6cm]{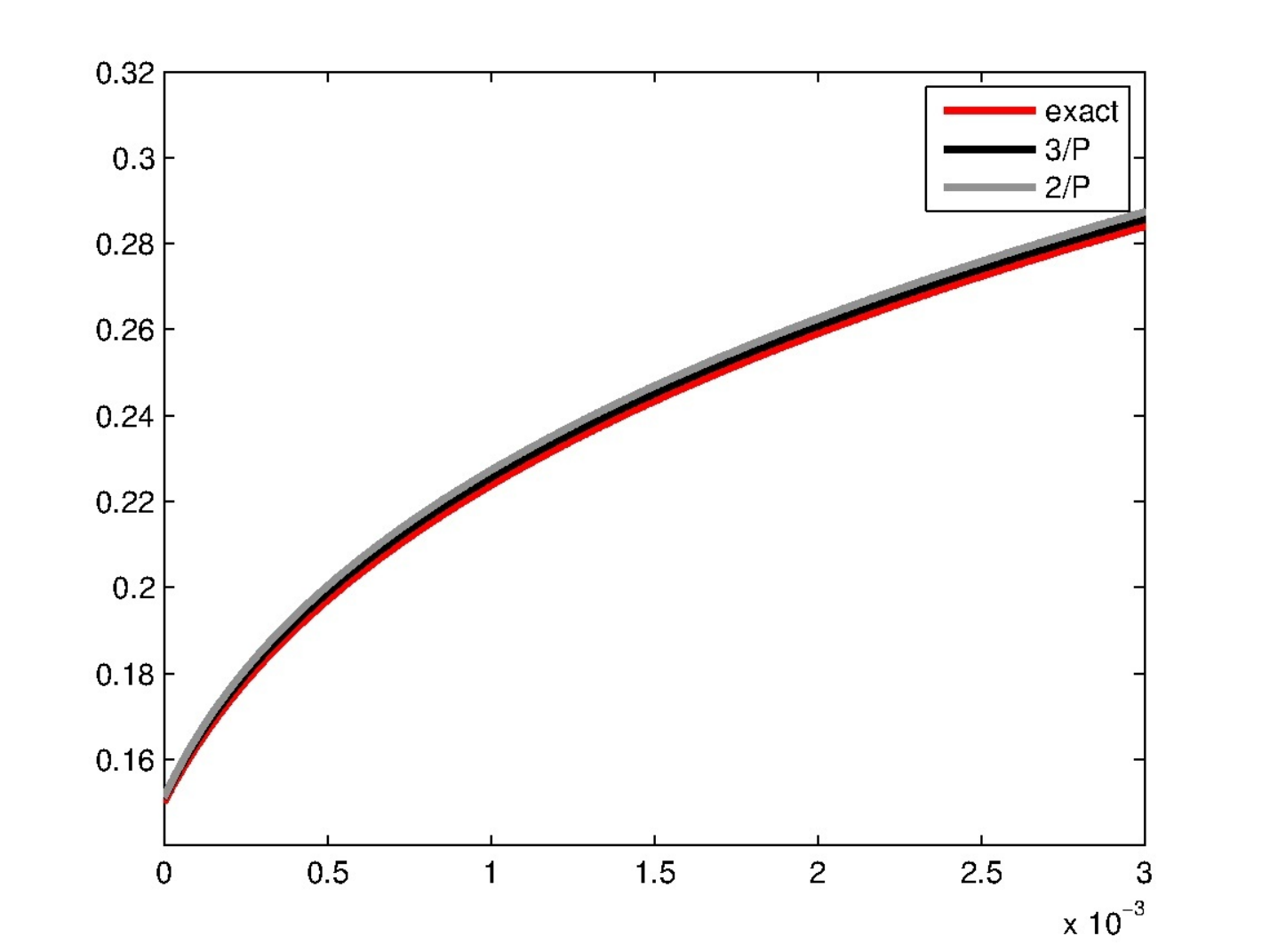}
\caption{Left : sampling of $\Gamma(t)$ at different times $t$ ; Right : The graphs of $t \to R_{\varepsilon}(t)$ for $\varepsilon=\frac 2 P$ and $\varepsilon=\frac 3 P$, compared with the exact solution.}
\label{fig:error_circle_classic}
\end{figure}

The first test plotted in Figure~\ref{fig:error_circle_classic}  illustrates the good behavior of our scheme with respect to the exact solution.
The initial set  $E$ is a circle of radius $R_0 = 0.15$. The continuous Willmore flow preserves the circle yet increases the radius according the law
$$
R(t) = \left(R_0^4 + 2 t \right)^{1/4}.
$$
The left picture in Figure~\ref{fig:error_circle_classic} represents the interfaces $\Gamma(t) $ at different times $t$ 
obtained with the following numerical parameters: ${\cal P} = 2^7$, $\varepsilon = 2/{\cal P}$ and  $\delta_t = \frac{\varepsilon^2}{2{\cal P}^2}$.  
The right picture in Figure~\ref{fig:error_circle_classic} depicts the error between the numerical radius $R_{\varepsilon}(t)$ and the theoretical radius $R(t)$, at different times and for two different values of $\varepsilon$ (the other parameters are kept unchanged). It is reasonable to believe that this experiment illustrates the numerical convergence of $R_{\varepsilon}(t)$ to $R(t)$ as $\varepsilon$ goes to zero.

\paragraph{Evolution of two disjoint circles and formation of singularities}

One of our motivations in this study is to understand and observe the behavior of the diffuse Willmore solution in the situations where singularities appear. As it was discussed in Section~\ref{sec:approx}, this may happen for instance with the classical approximation flow. We consider as initial set $\Omega_0$ the union of two disjoint circles of radius $R = 0.15$.
Each circle should evolve as a circle with increasing radius, up to the contact occurs. To the best of our knowledge, the theoretical Willmore flow is not clearly defined after this critical collision time. Therefore, the asymptotic limit of the solution $t\mapsto u_{\varepsilon}(\cdot,t)$ as $\varepsilon$ goes to $0$ could be a  good candidate for the definition of a weak Willmore flow. However, different behaviors of $t\mapsto u_{\varepsilon}(\cdot,t)$  have been observed in the literature. For instance, the two circle merge  in \cite{FrRuWi11} whereas a crossing of interfaces appears at collision time in \cite{Esedoglu_12}. \\ 

\begin{figure}[!ht] 
\centering
\includegraphics[width=5cm]{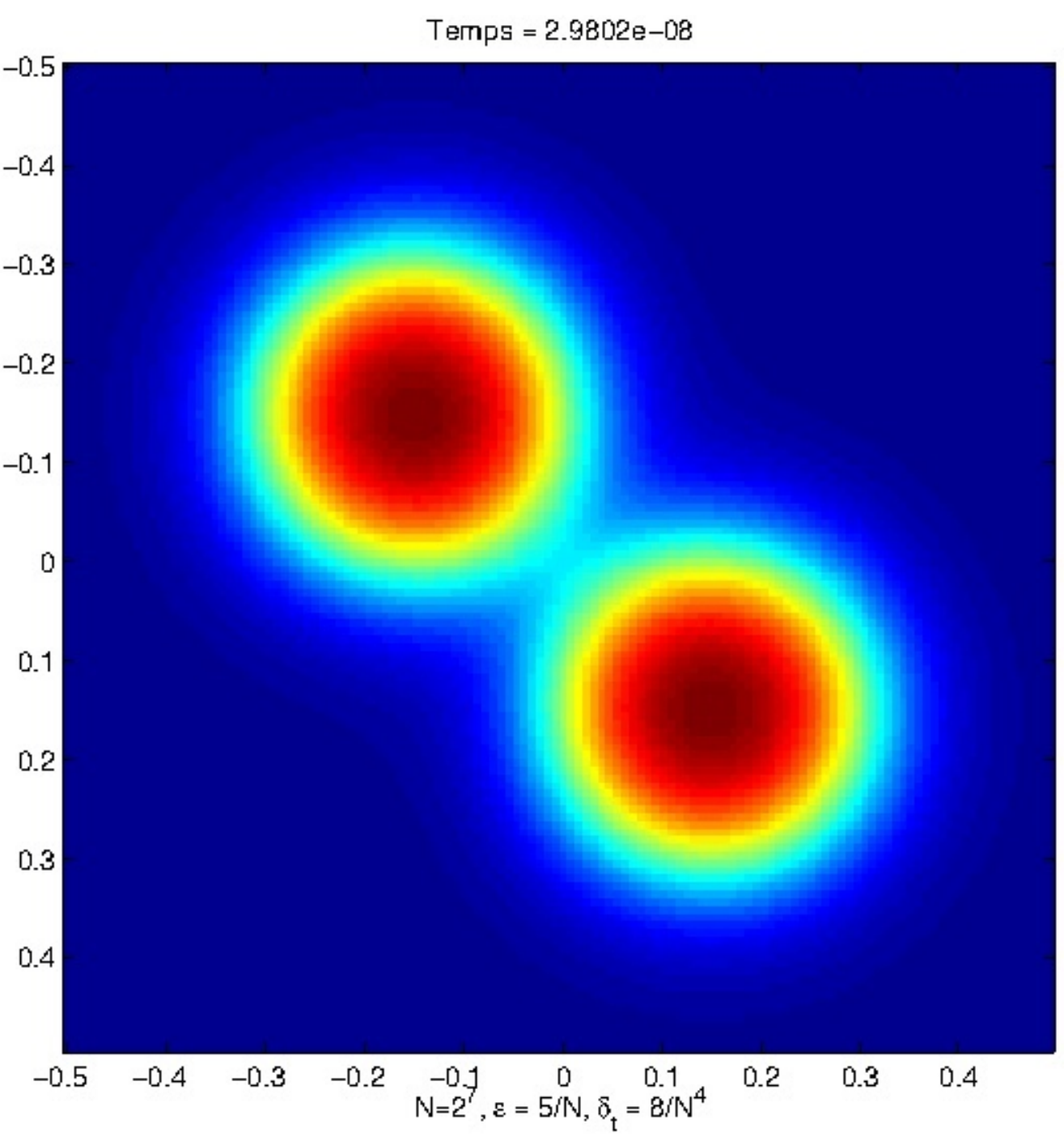}
\includegraphics[width=5cm]{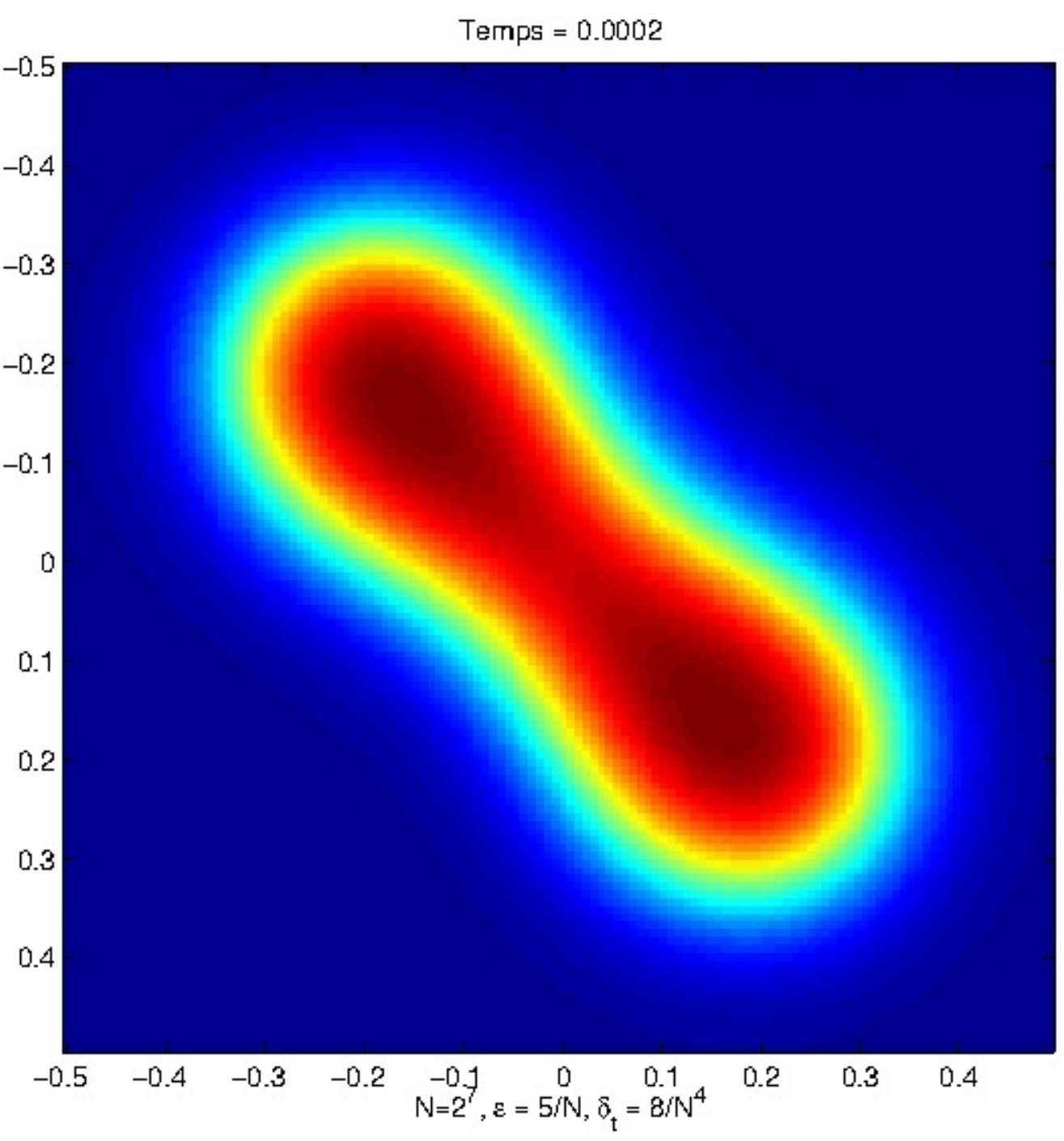}
\includegraphics[width=5cm]{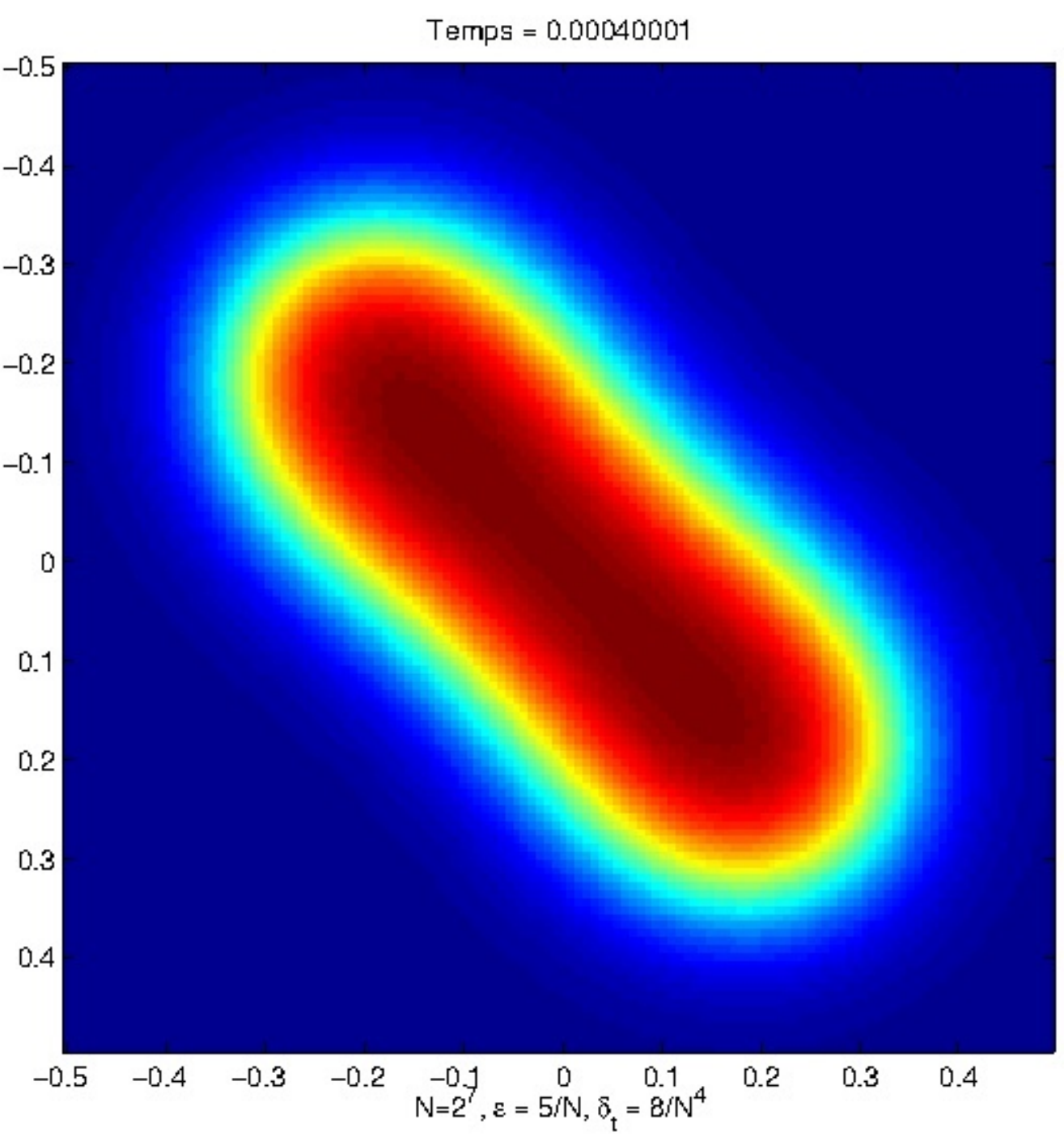} \\
\includegraphics[width=5cm]{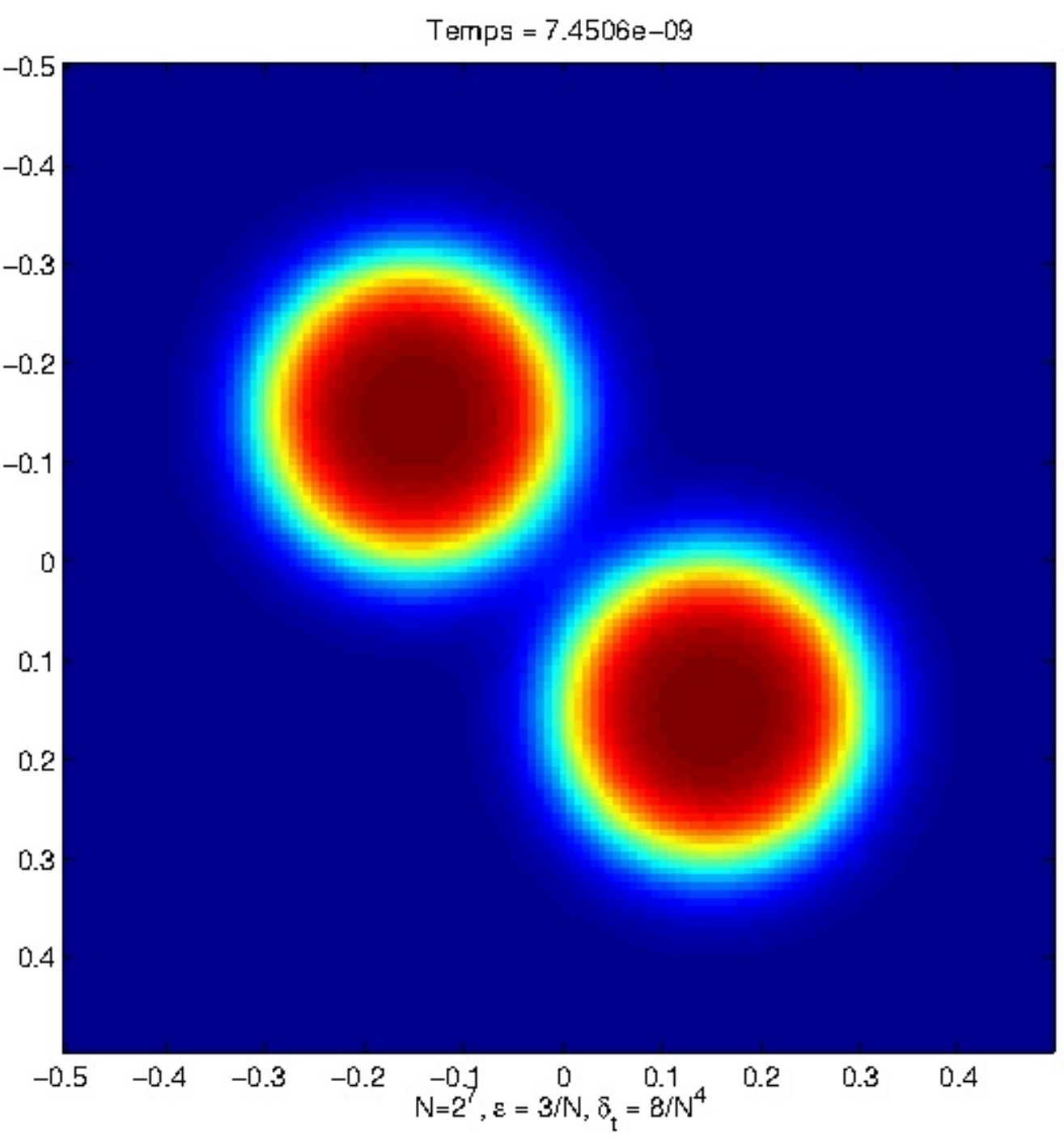}
\includegraphics[width=5cm]{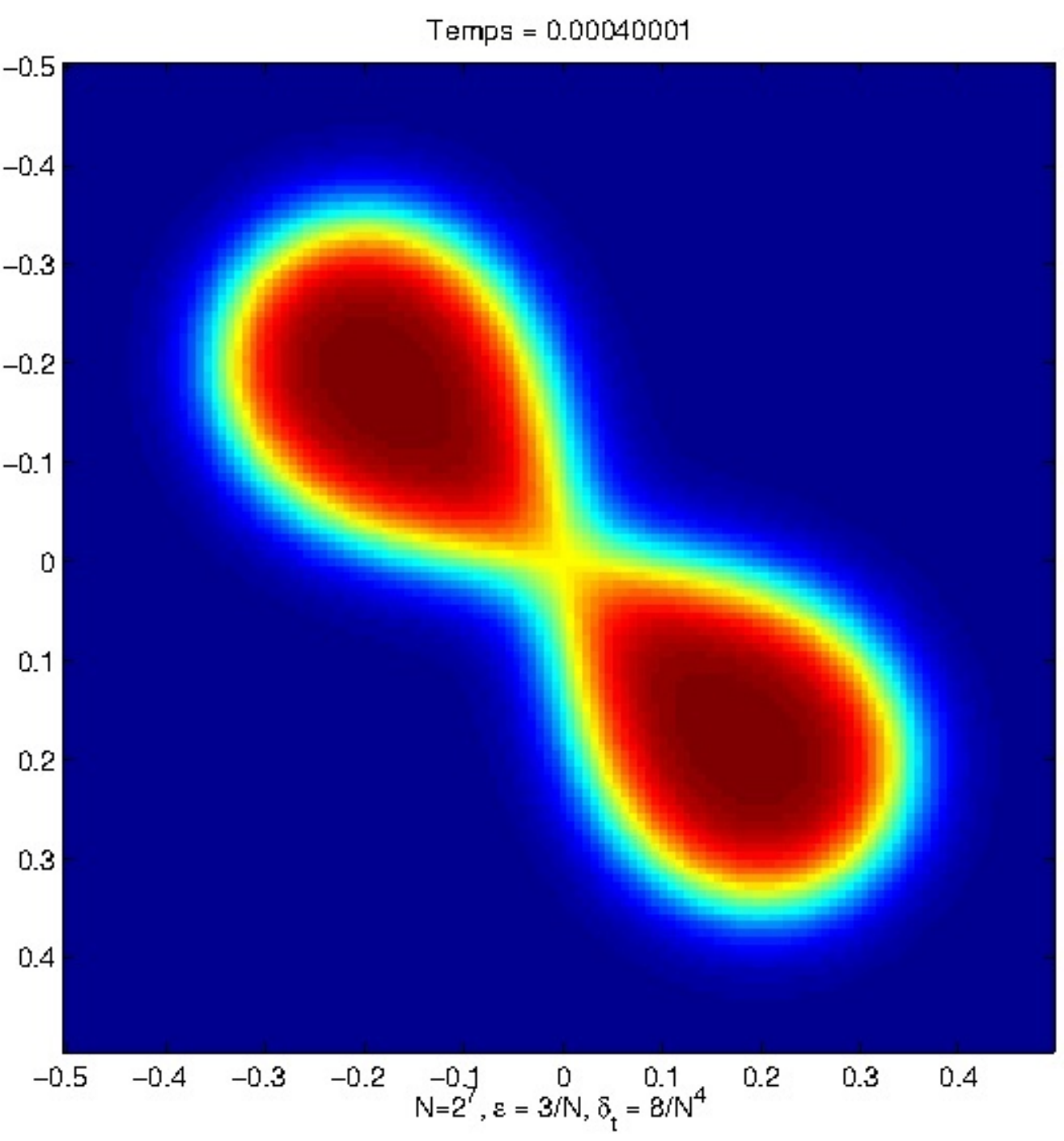}
\includegraphics[width=5cm]{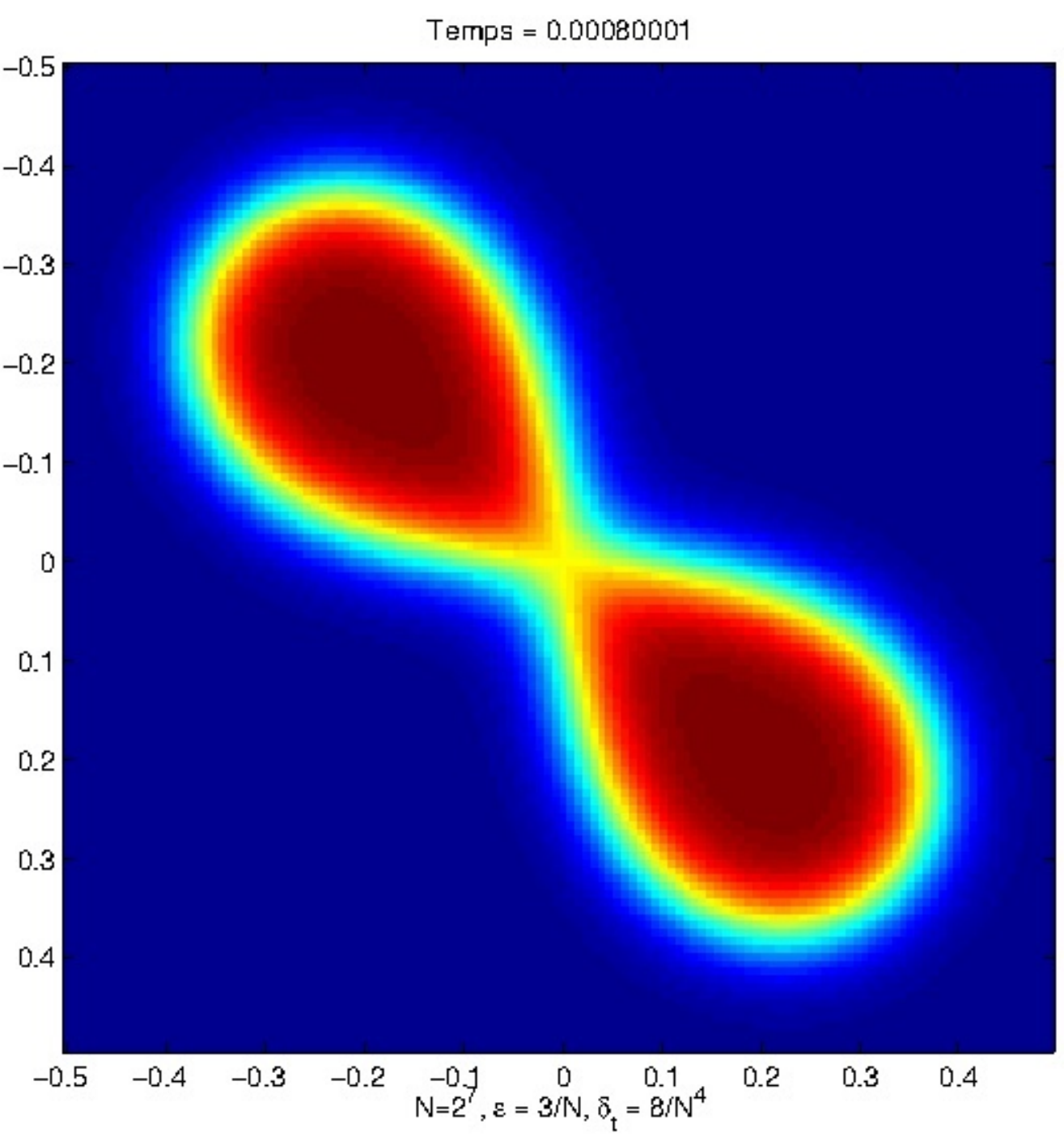} \\
\includegraphics[width=5cm]{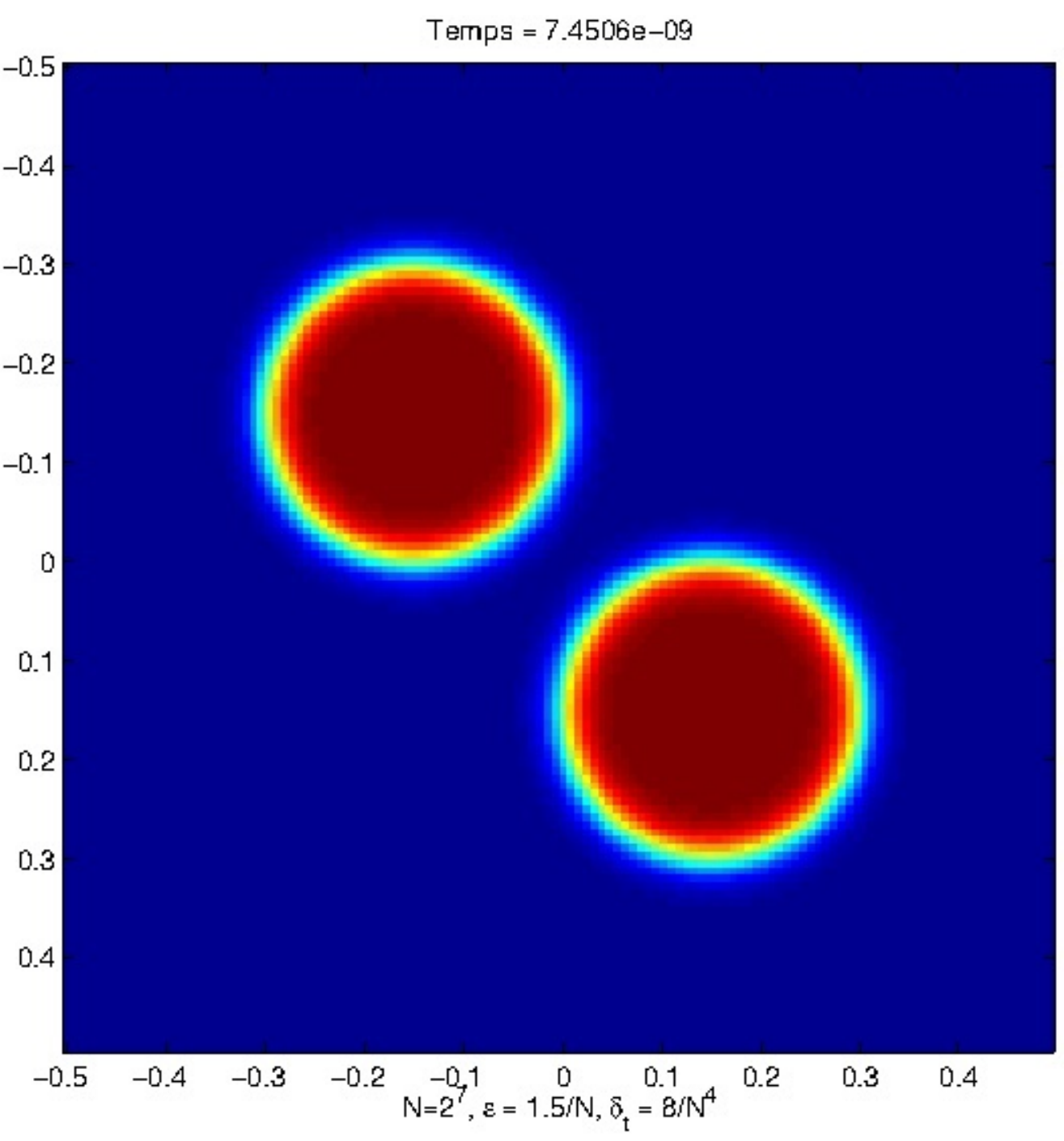}
\includegraphics[width=5cm]{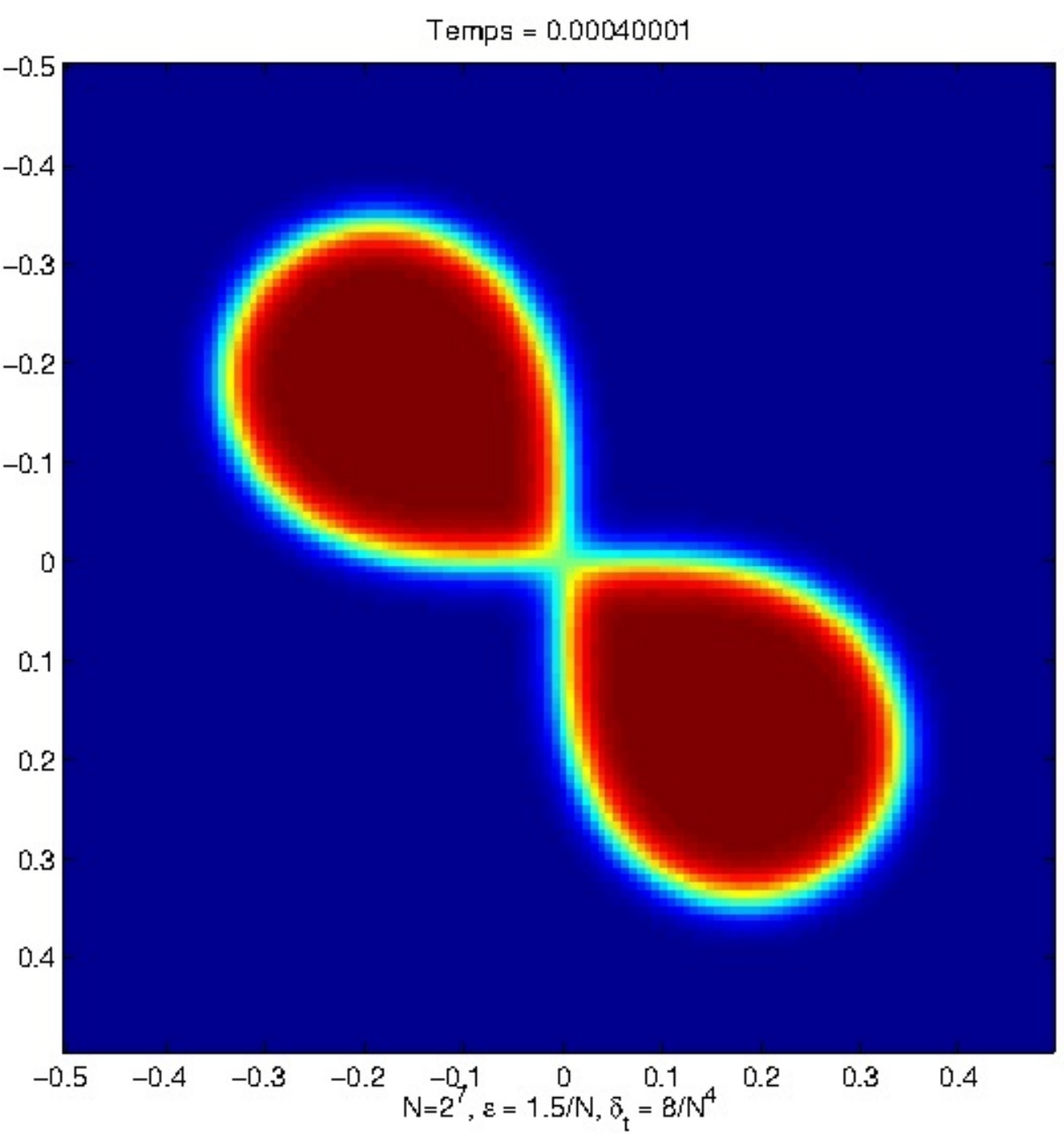}
\includegraphics[width=5cm]{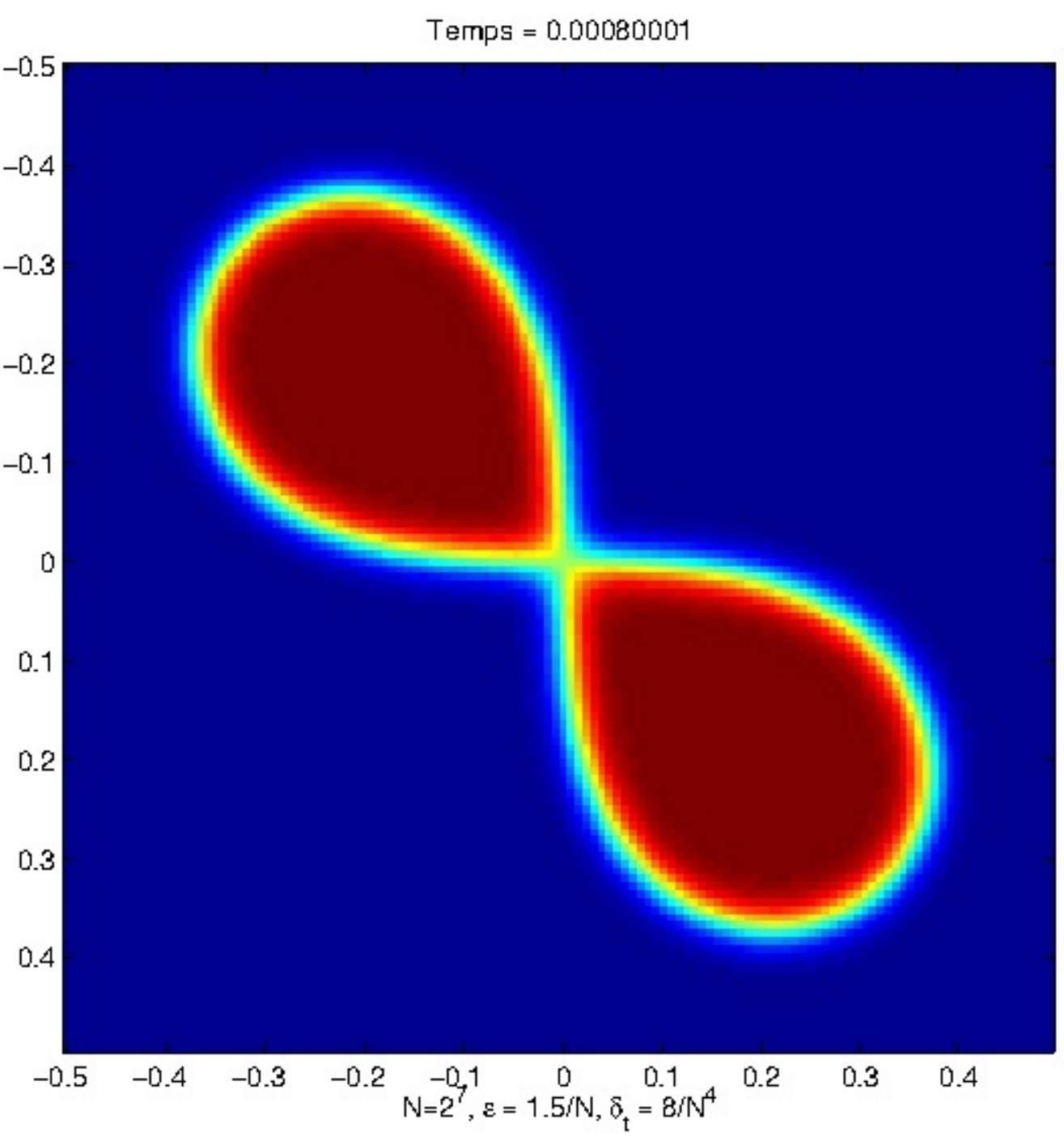} \\
\caption{Evolution by the classical approximation of the Willmore flow of two disjoint circles, for various values of $\varepsilon$. First line: $\varepsilon = 5/{\cal P}$ ; second line: $\varepsilon = 3/{\cal P}$ ; third line: $\varepsilon = 1.5/{\cal P}$; The curve $\Gamma(t)$ is observed at times: $t=0$ (left), $t = 0.0004$ (middle), $t = 0.0008$ (right).}
\label{fig:evolution_two_circle}
\end{figure}

We plot on Figure \ref{fig:evolution_two_circle} the graph of  $t\mapsto u_{\varepsilon}(\cdot,t)$ computed for different values of $\varepsilon$. We choose for the other parameters: $P = 2^7$, $\delta_t = 1/{\cal P}^{-4}$. 
In the first experiment obtained with $\varepsilon = 5/{\cal P}$,  the two circles merge. In contrast, a crossing of interfaces appear
for the cases $\varepsilon = 3/{\cal P}$ and $\varepsilon = 1.5$. More precisely, we can distinguish three different periods in the last two experiments: in the first period, both circles evolve independently one from the other. The second period begins when the distance between the two circles 
is about the size of the diffuse interfaces and the formation of a crossing is observed. This corresponds to a solution of the Allen-Cahn equation with unsmooth nodal set. After contact, the interfaces continues to evolve while the crossing seems to be numerically stable and does not influence the interface evolution.
More precisely, the interface $\Gamma(t)$ seems to converge to a growing eight, which is one of the closed planar elasticae described in Langer and Singer's work~\cite{Langer_Singer}.  

\begin{figure}[!ht] 
\centering
\includegraphics[width=6cm]{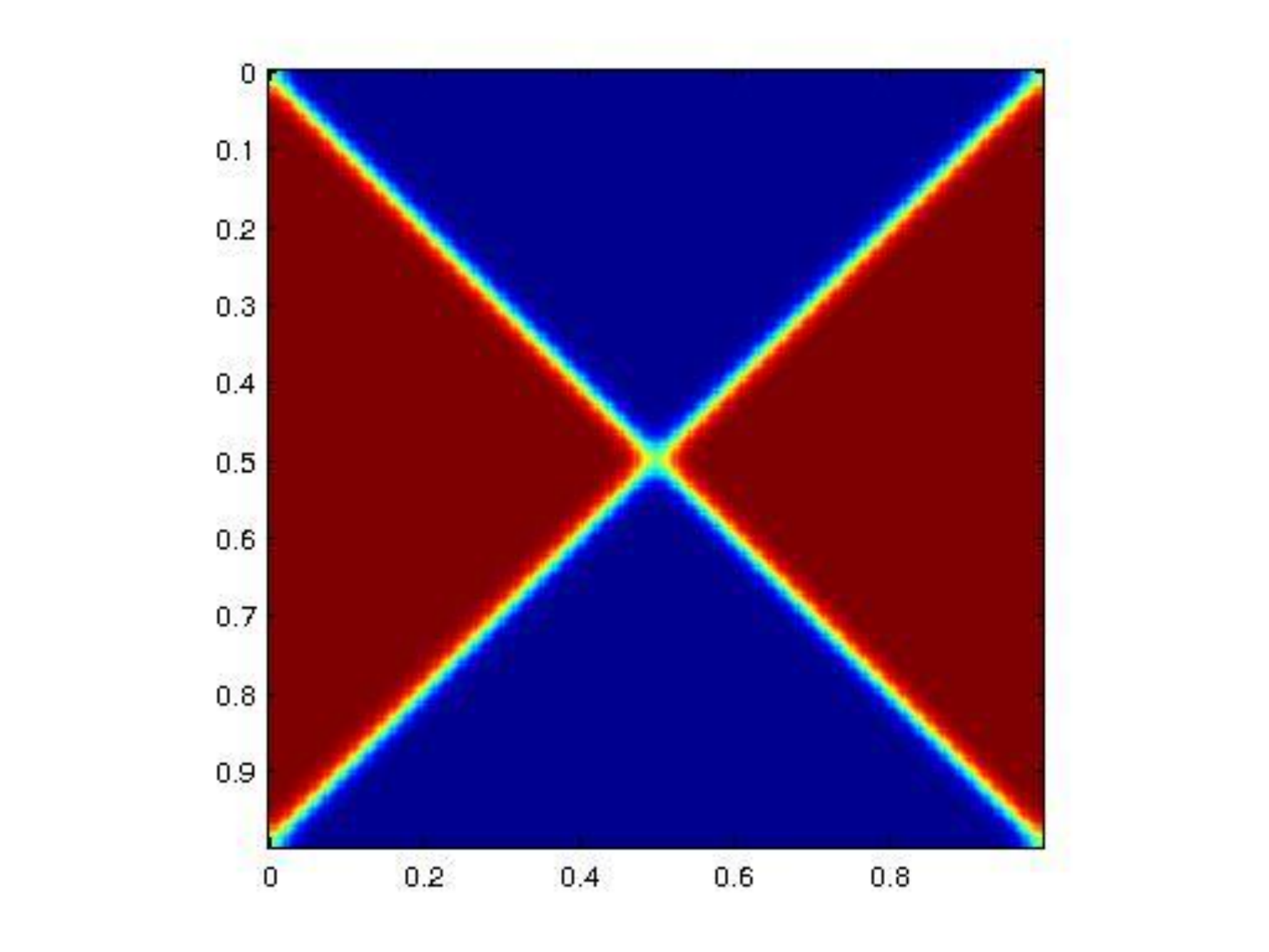}\hspace*{-0.8cm}
\includegraphics[width=6cm]{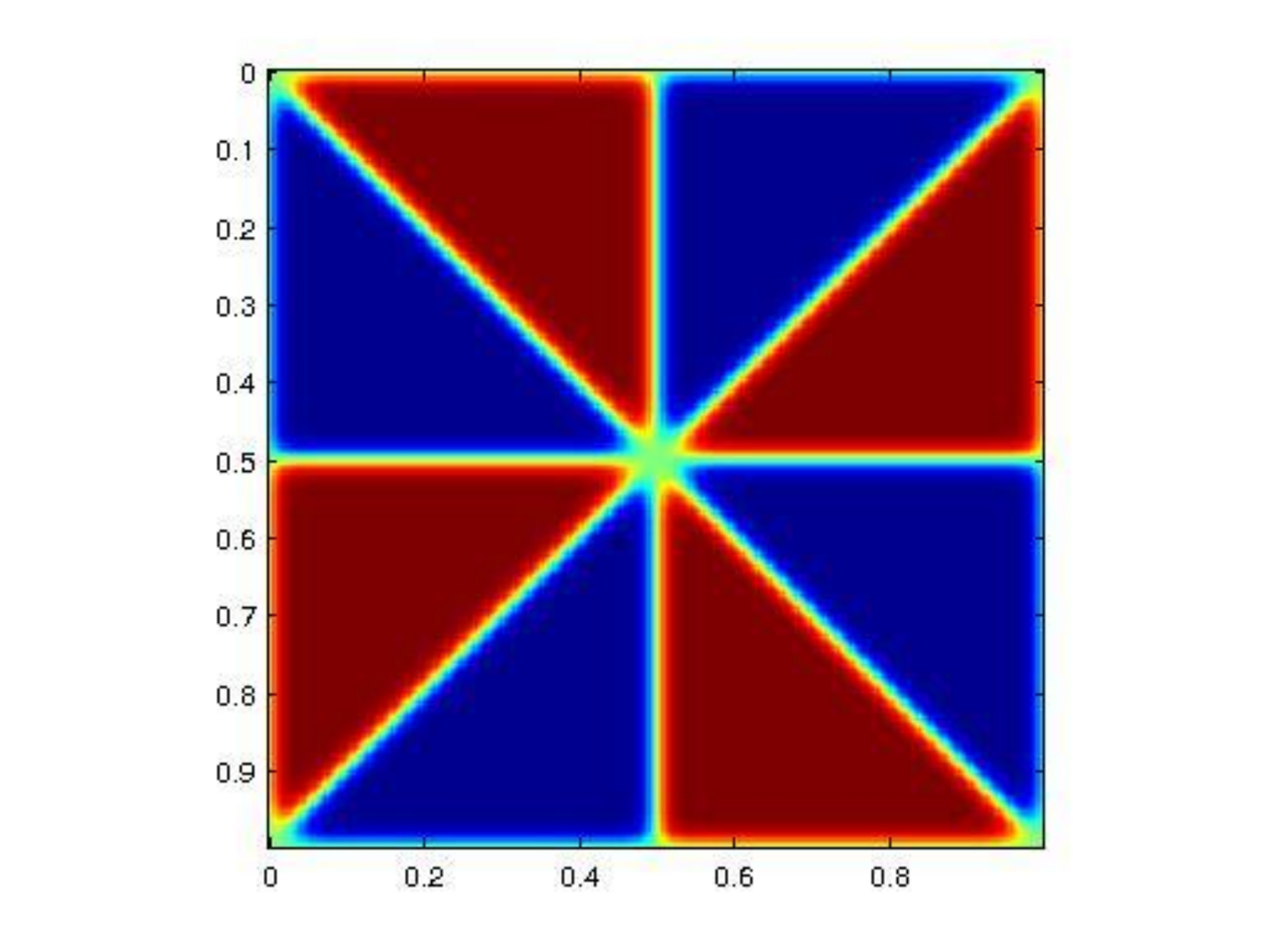}\hspace*{-0.8cm}
\includegraphics[width=6cm]{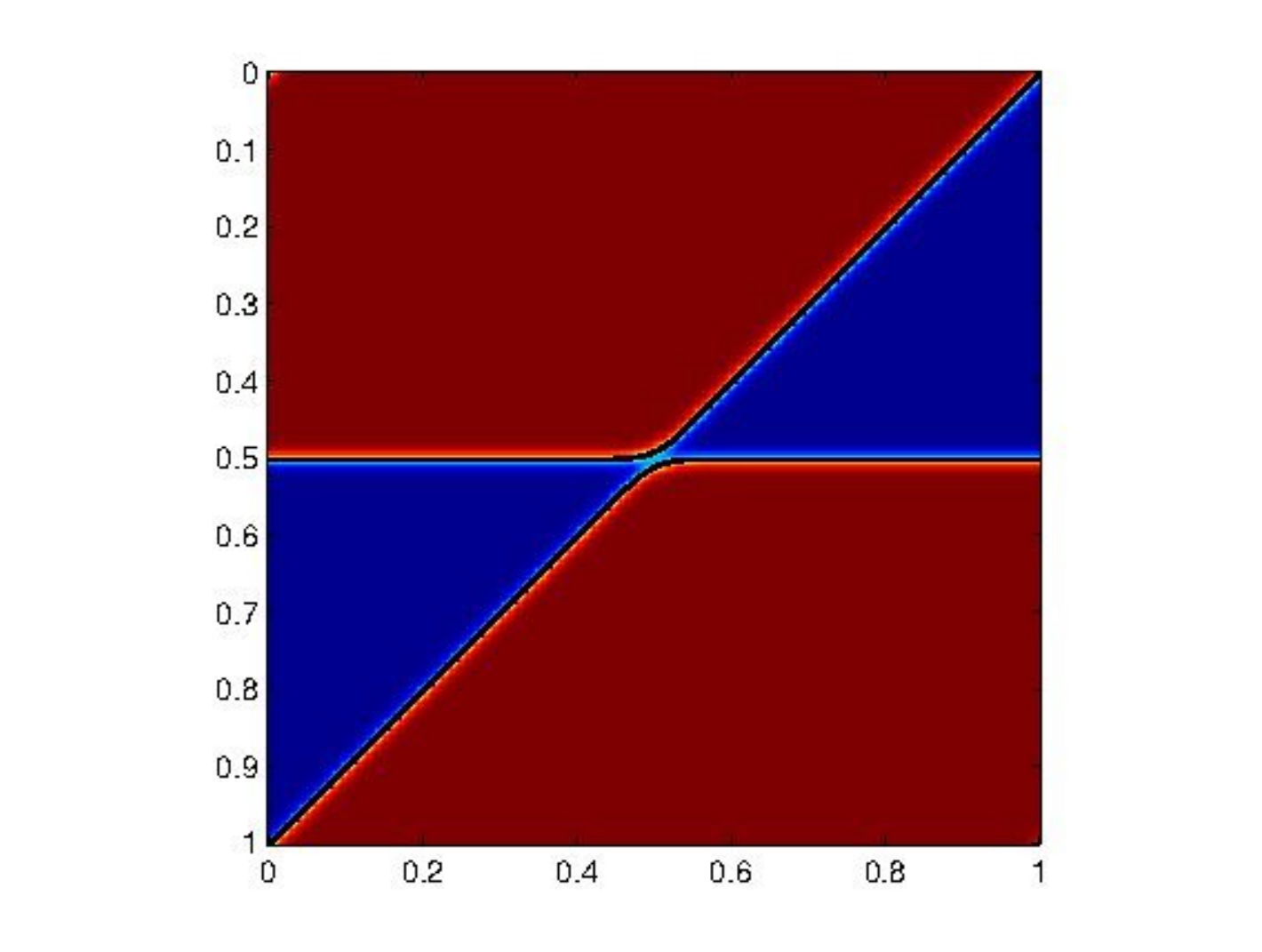}
\caption{Left : two examples of saddle-shaped Allen-Cahn solutions with $4$ and $8$ ends. Right: Allen-Cahn solution without dihedral symmetry (thus without saddle point, the $1/2$-isolevel line is shown in black).}
\label{fig:saddle_shape_solution}
\end{figure}

\paragraph{Numerical examples of saddle-shaped solutions of Allen-Cahn equation}
We already mentioned in Section~\ref{sec:approx} the existence result due to Dang, Fife, and Peletier~\cite{Dang_Fife_Peletier_92} of an entire solution in the plane to the Allen-Cahn equation whose nodal set coincides with $\{(x,y),\,xy=0\}$ (it can be generalized to every even dimension~\cite{cabre-terra}). By restricting to a sector and using consecutive reflections, it is possible to build solutions whose nodal set has an arbitrary number of branches with the property of dihedral symmetry, i.e. of equal angle between two consecutive branches~\cite{gui}. Actually, by a result of Hartman and Wintner~\cite{hart-wintn}, saddle-shaped solutions must satisfy the equal angle property. We illustrate numerically in Figure~\ref{fig:saddle_shape_solution}, left, two solutions with different branching degrees. These examples are classical, and have been previously obtained by various authors~\cite{Esedoglu_p,Lowengrub,Esedoglu_12}.
\par To be complete, let us mention that $2k$-ended solutions, i.e. solutions whose nodal set coincides outside any compact set with the union of $2k$ straight lines which cross at the origin, do exist without the dihedral symmetry requirement~\cite{delpino-etal}. In the particular case of $4$-ended solutions, the result can even be proved for arbitrary angles between the lines~\cite{Kowalczyk_pacard}. Of course, by Hartman and Wintner's result, the nodal set itself cannot self-intersect at the origin if the dihedral symmetry does not hold, but remains smooth instead. This is illustrated in Figure~\ref{fig:saddle_shape_solution}, right, where the $1/2$-isolevel line has been represented in black.
 
\begin{figure}[!ht]  
\begin{center}
\includegraphics[width=6cm]{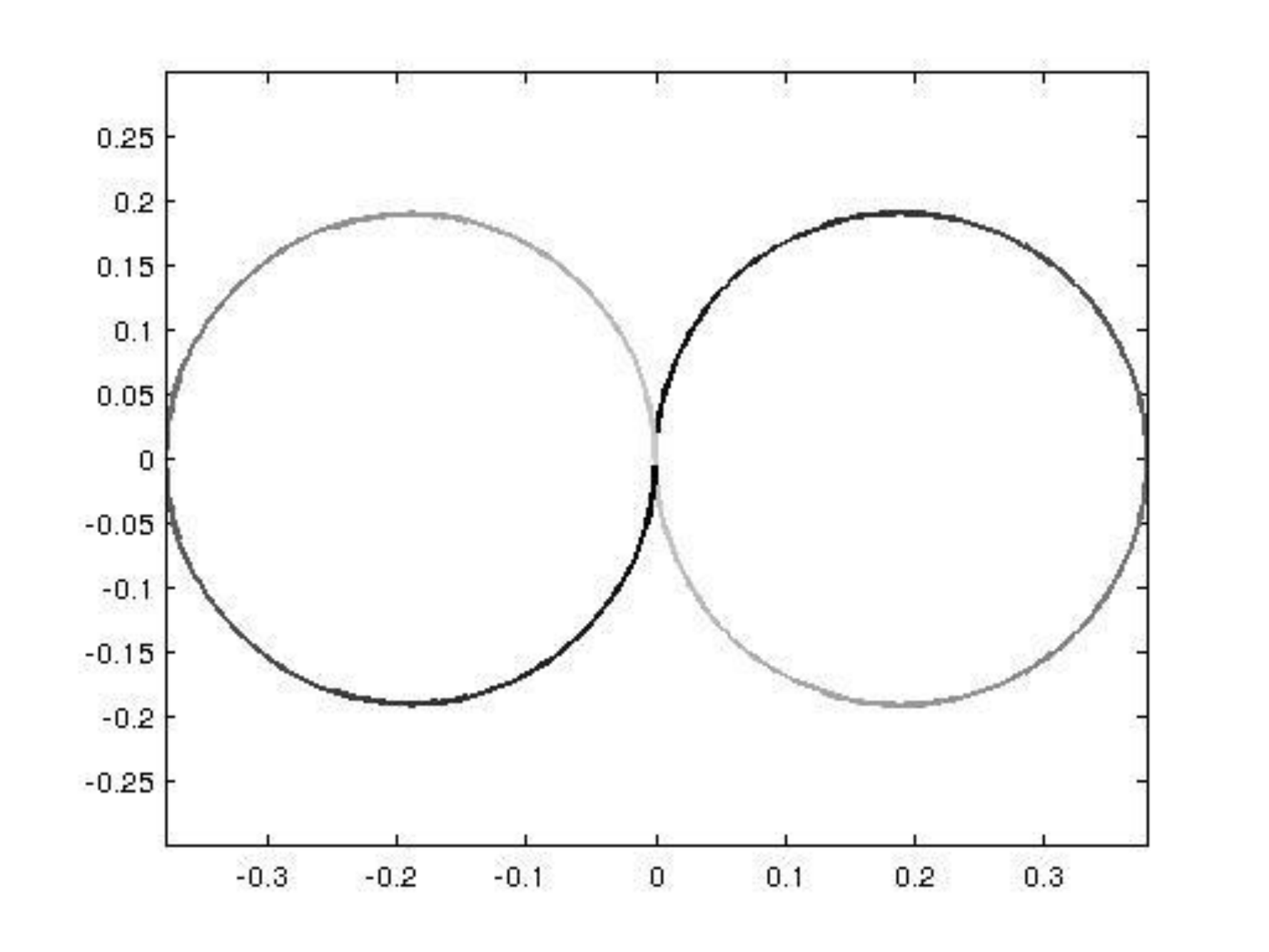}
\includegraphics[width=6cm]{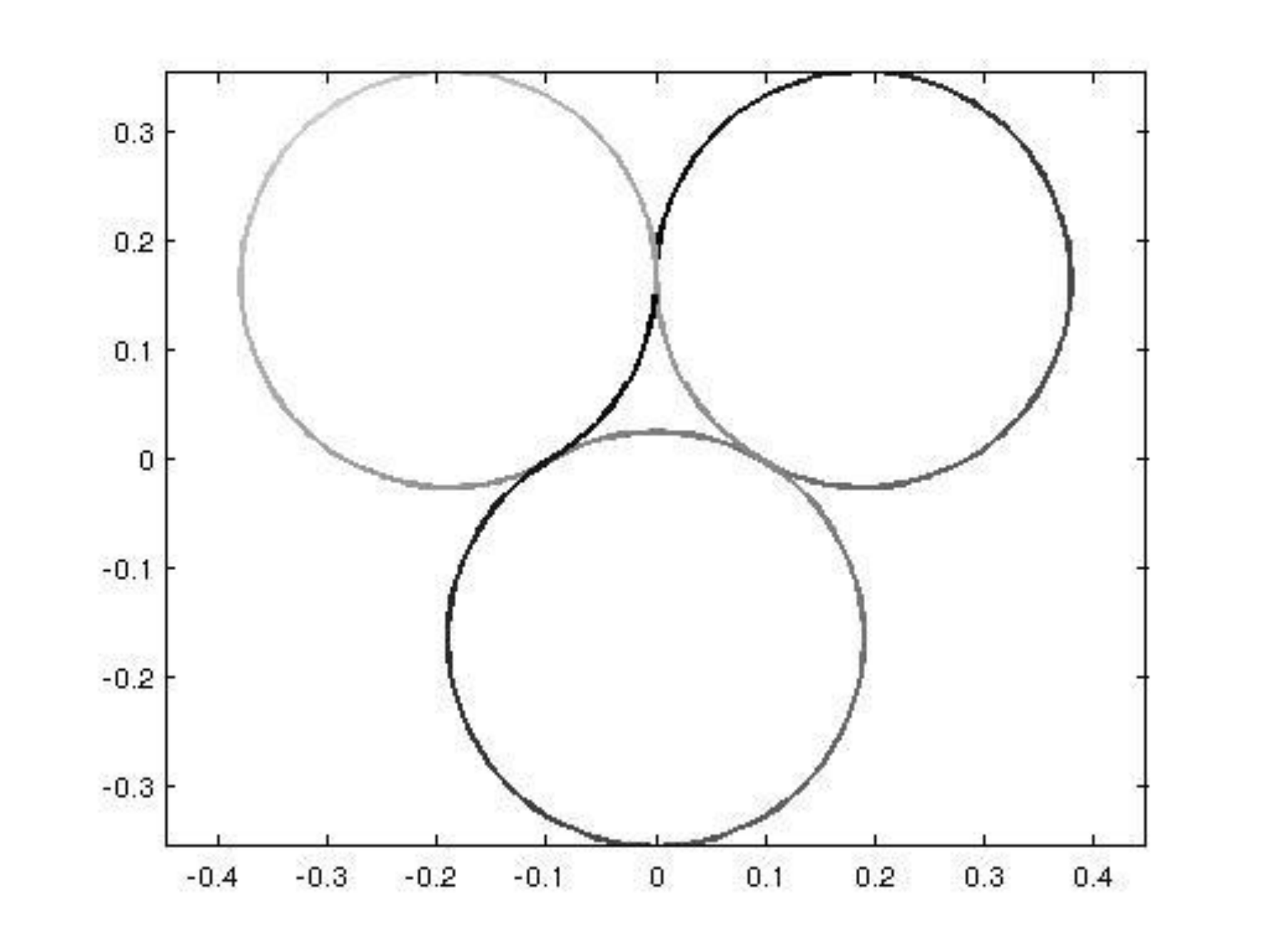}
\caption{Two different choices of a smooth parametric initial curve $\Gamma(0)$ forming either two or three circles.}
\label{fig:saddle_shape_solution_init}
\end{center}
\end{figure}

\begin{figure}[!ht] 
\begin{center}
\includegraphics[width=5cm]{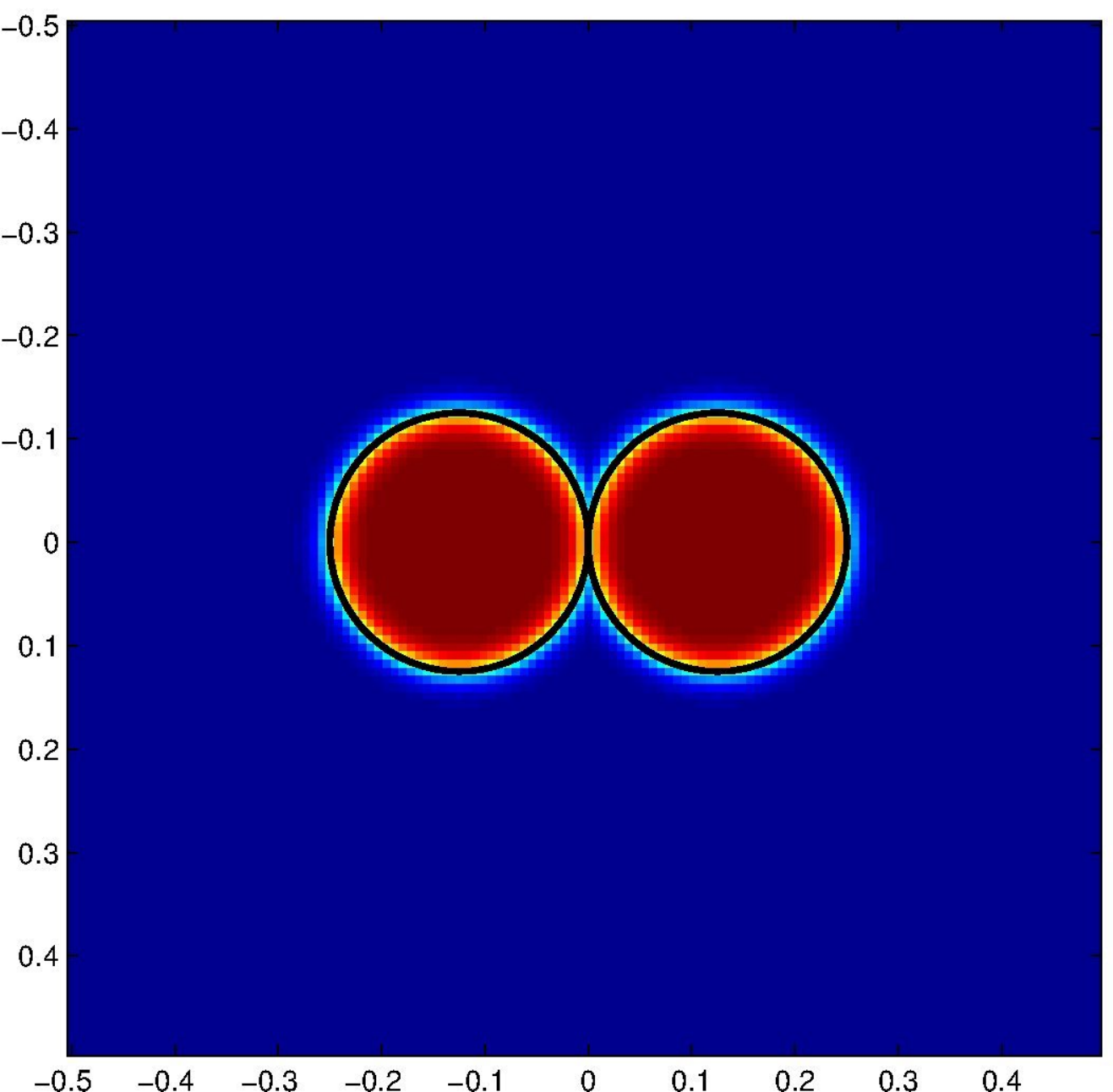}
\includegraphics[width=5cm]{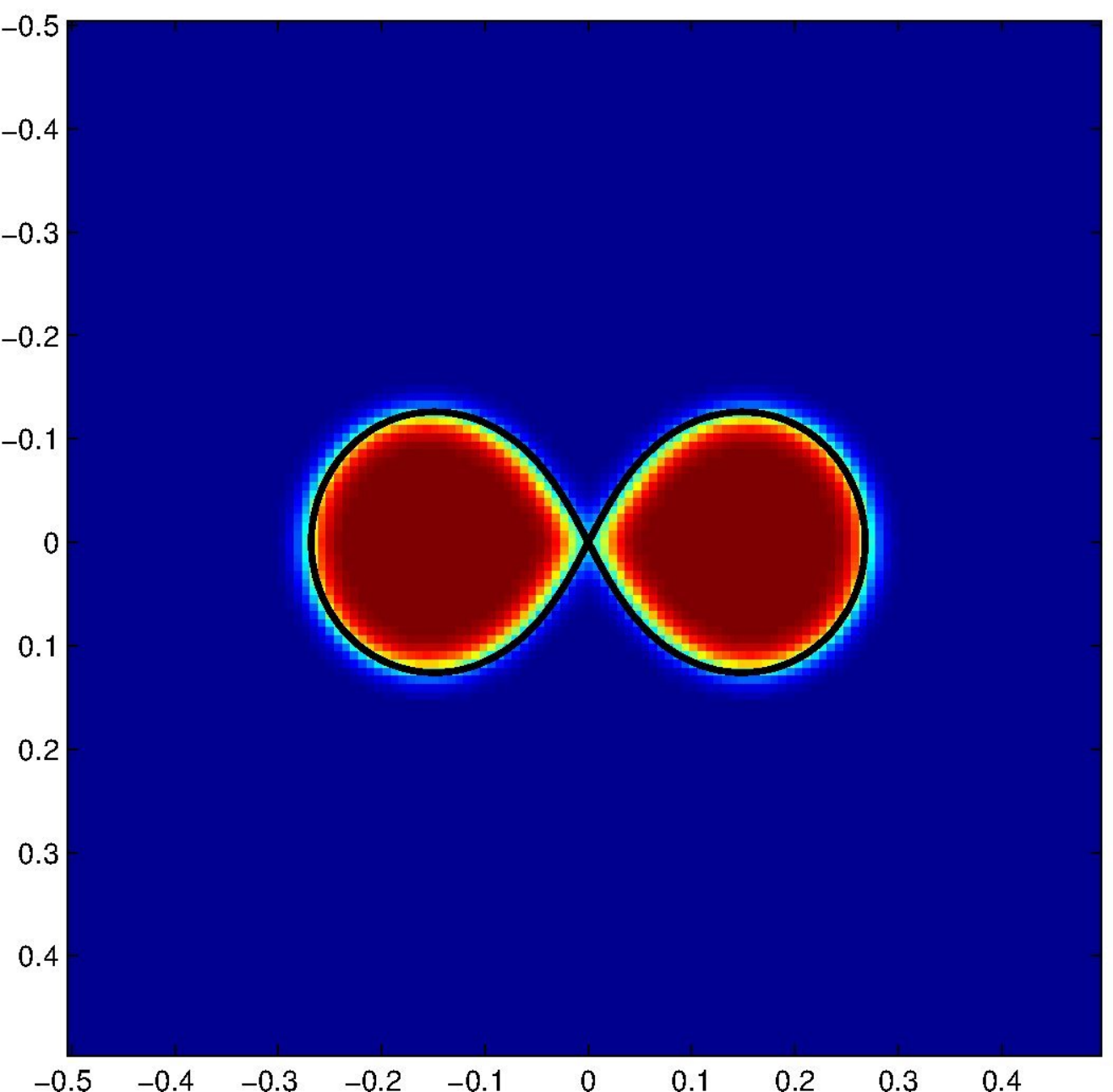}
\includegraphics[width=5cm]{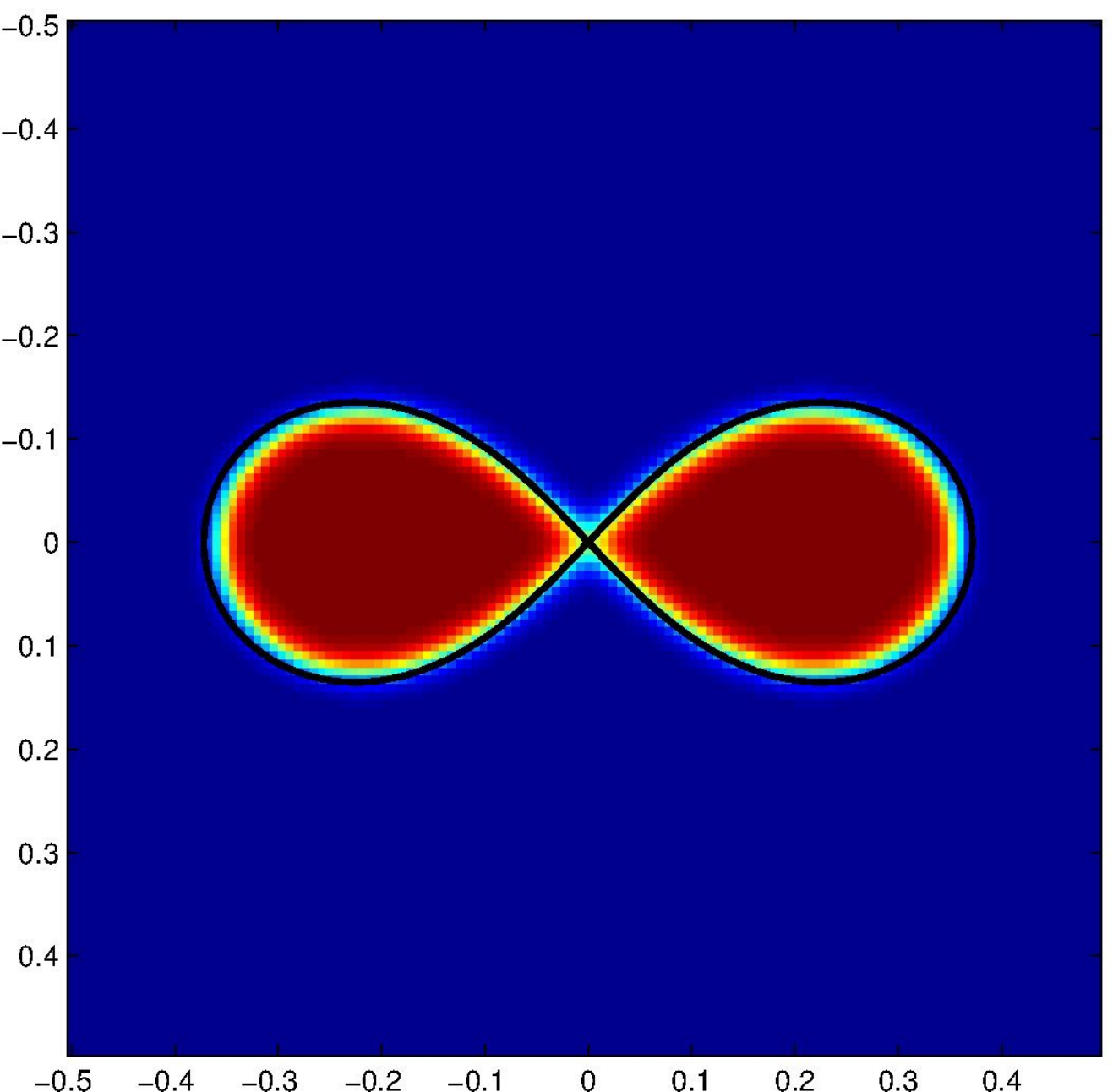} \\
\includegraphics[width=5cm]{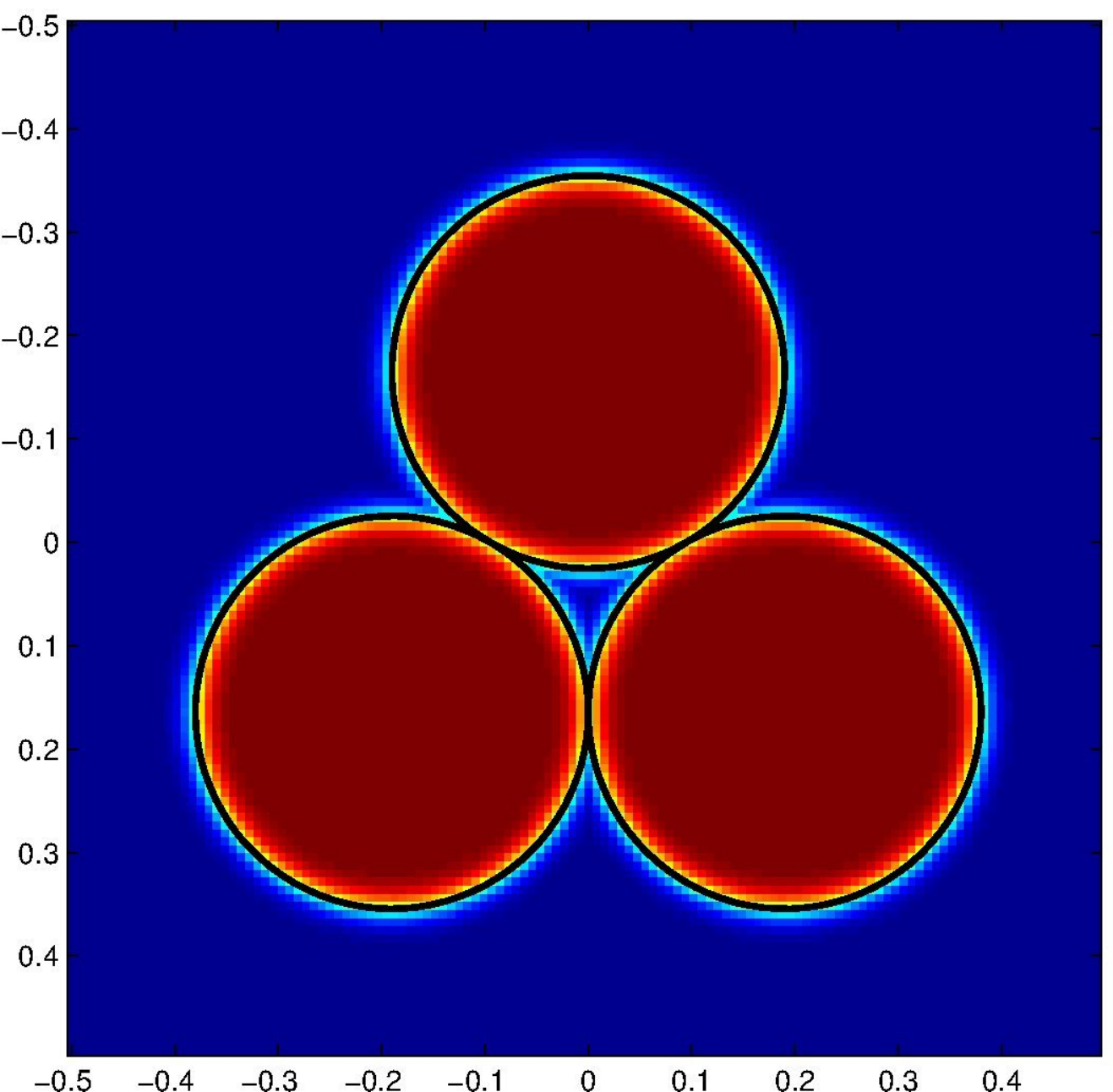}
\includegraphics[width=5cm]{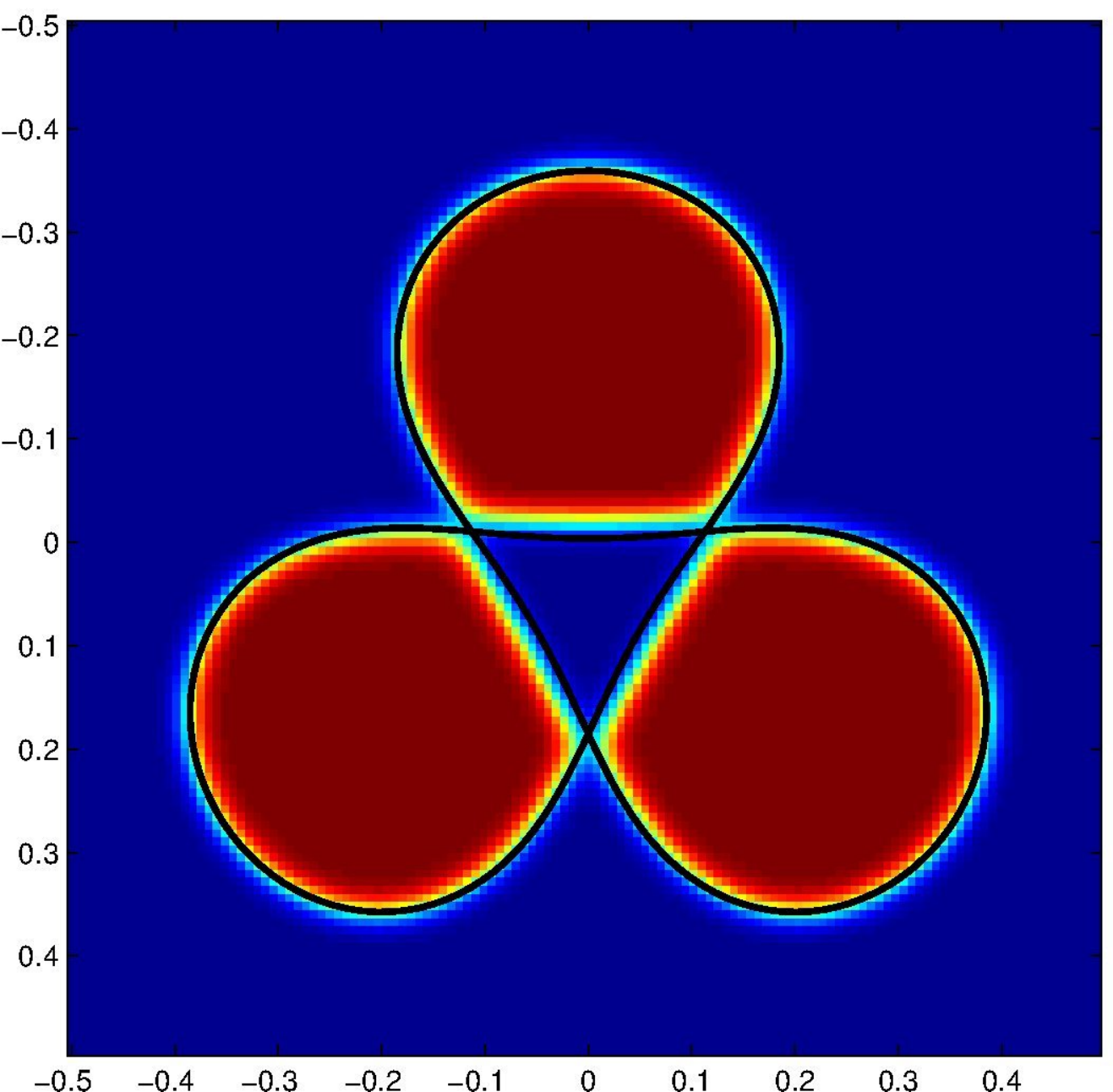}
\includegraphics[width=5cm]{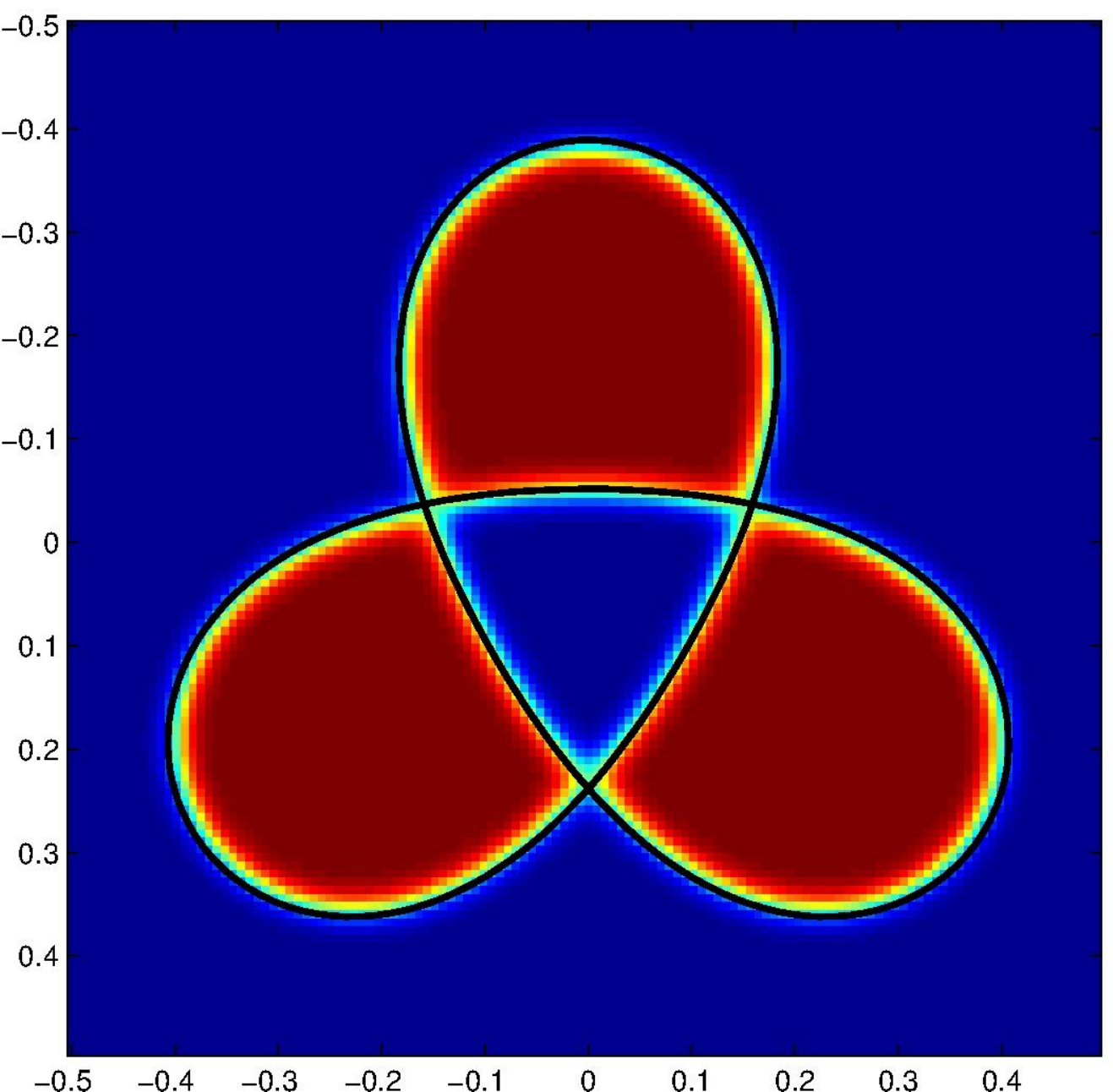} \\
\caption{Comparison between the classical phase field flow and the parametric Willmore flow (black line) starting from the initial curves of Figure~\ref{fig:saddle_shape_solution_init}; Left : $t=0$ ; Middle : $t = 5. 10^{-5}$ ; Right : $t = 5. 10^{-4}$. Both flows yield the same numerical solution.}
\label{fig:evolution_two_circles_pfp}
\end{center}
\end{figure}

\paragraph{Comparison between phase field and parametric approaches}
We observed previously that the evolution of two disjoint circles after contact and creation of a crossing is similar to the evolution of a eight-like single curve. To highlight this point, we tested on the same configurations both the evolution provided by the classical diffuse flow and the evolution  with respect to a discrete parametric Willmore flow. 
In particular, we consider two different initial conditions corresponding, respectively, 
to the union of two or three contiguous circles. The parametric choice  of $\Gamma(0)$ is illustrated 
on Figure \ref{fig:saddle_shape_solution_init} and corresponds to using a single smooth $C^{1,1}$ curve that covers two or three circles, respectively.  The discrete parametric Willmore flow is computed with a finite element method as proposed and analyzed  by Dziuk in \cite{Dziuk:2008}. 
The phase field simulations are done with the set of parameters: ${\cal P} = 2^7$, $\varepsilon = 1/{\cal P}$ and $\delta_t =  \varepsilon {\cal P}^{-2}/10$. Numerical results are shown on Figure \ref{fig:evolution_two_circles_pfp}. As expected, 
these two different approaches give very similar results.  This suggests that the interface obtained by  a phase field approximation converges, 
 after apparition of a singularity, to an interface which evolves as a regular parametric Willmore flow. This is actually very much in favor of a {\it varifold interpretation}, at least on this example, of both flows. What really cares is the support, and its geometry, and not the fact that it is seen either as an isolevel set or as a parametrized set.
 
\begin{figure}[!ht] 
\centering
\includegraphics[width=5cm]{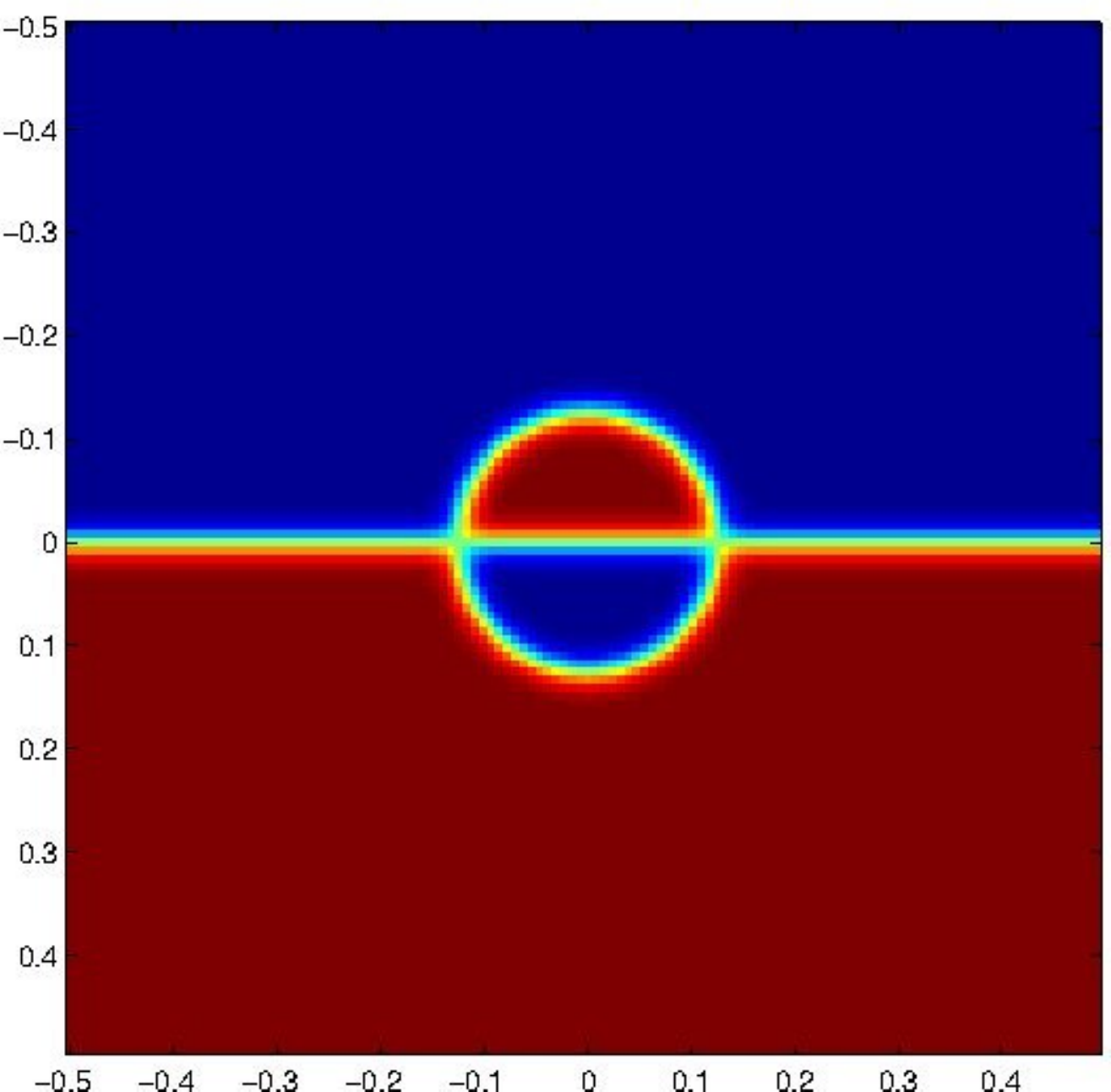}
\includegraphics[width=5cm]{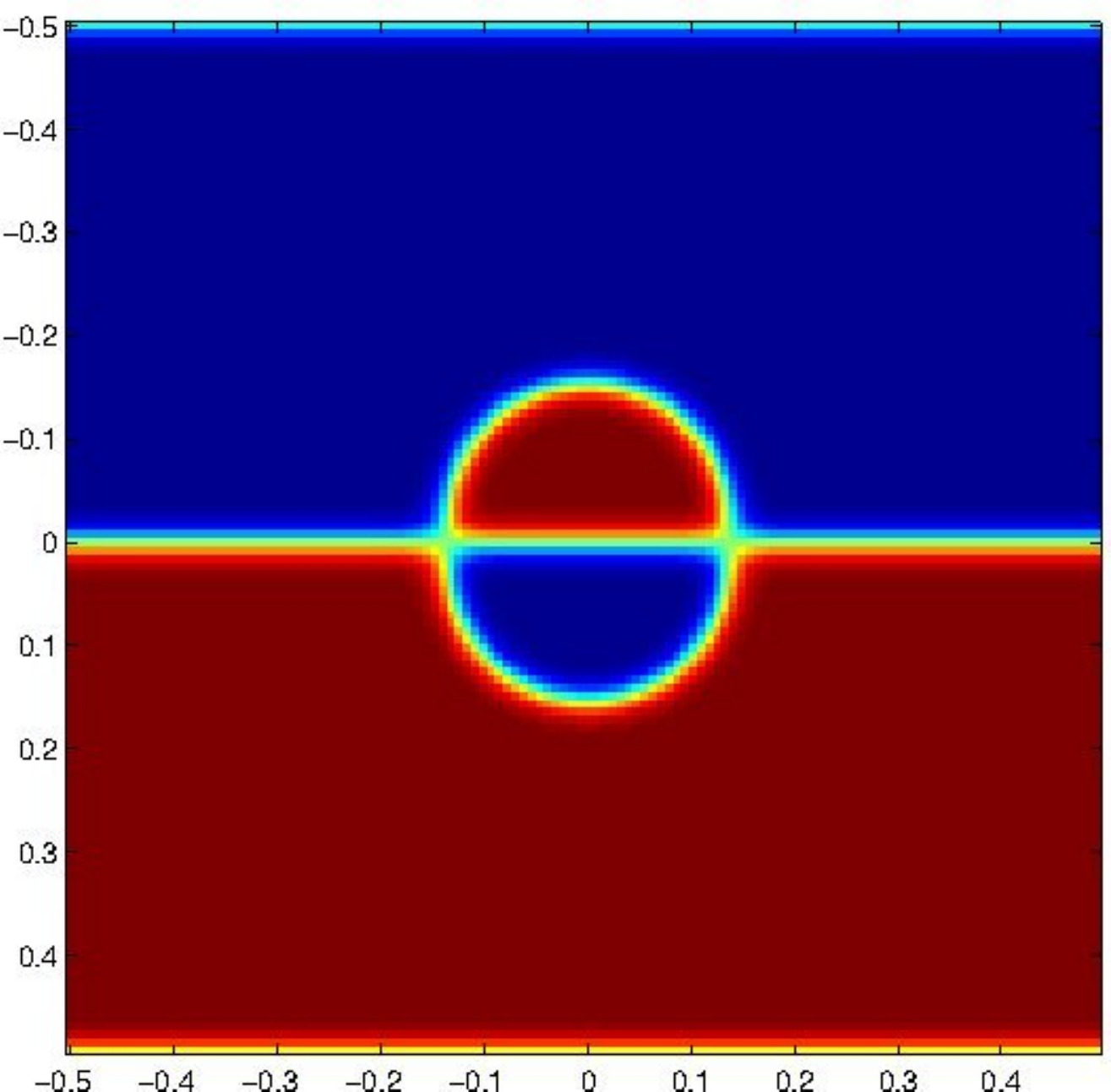}
\includegraphics[width=5cm]{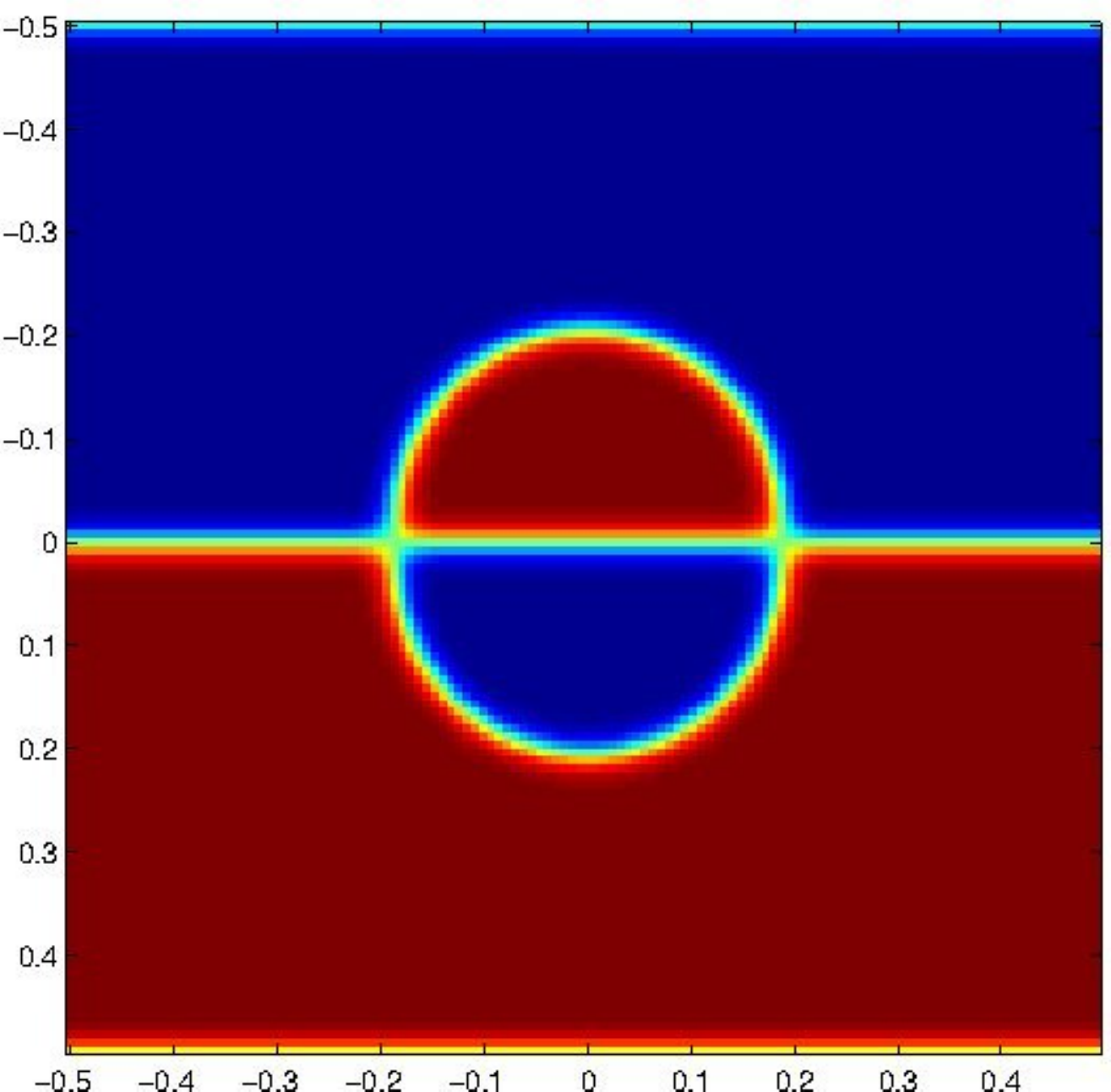} \\
\includegraphics[width=5cm]{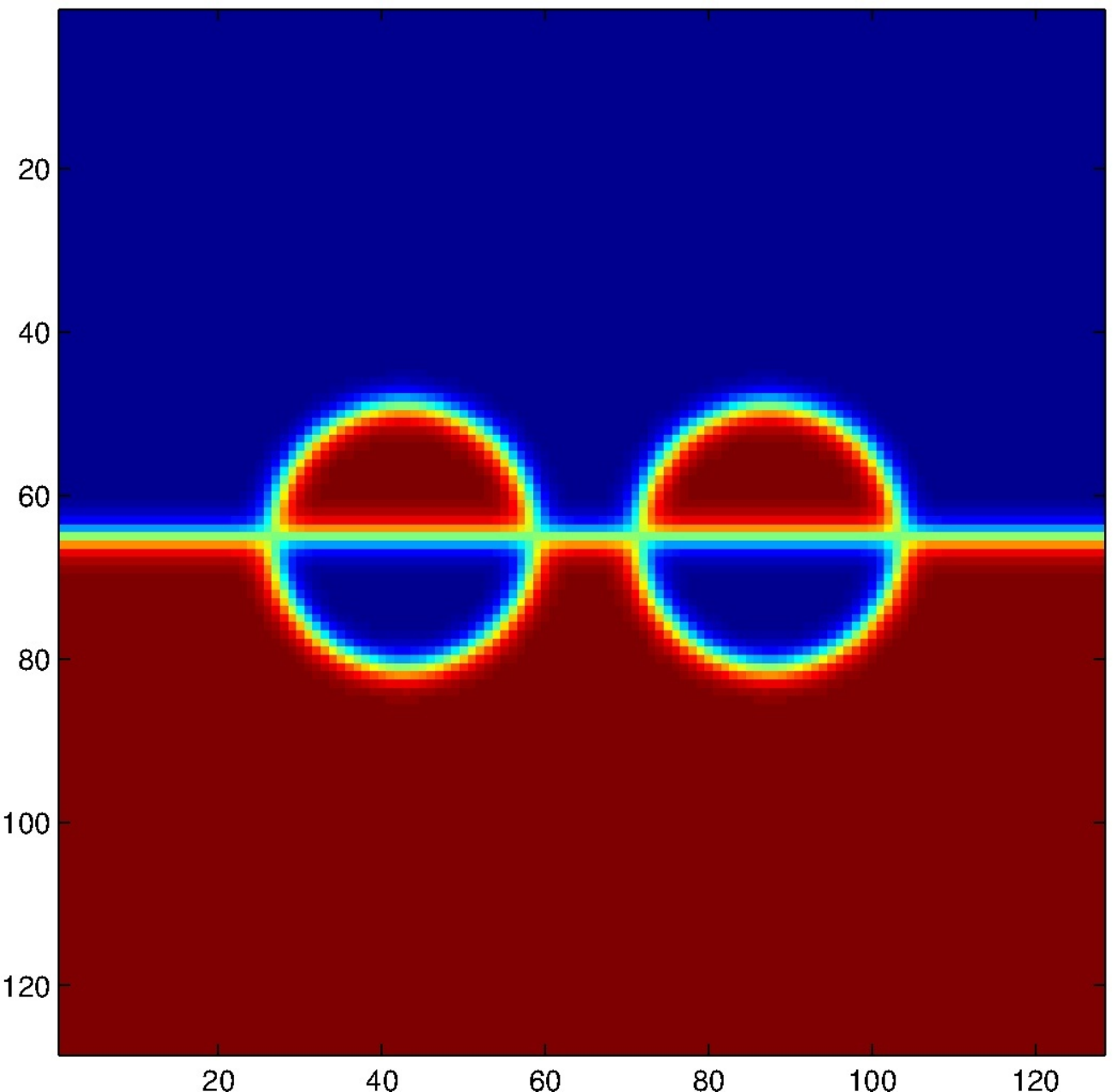}
\includegraphics[width=5cm]{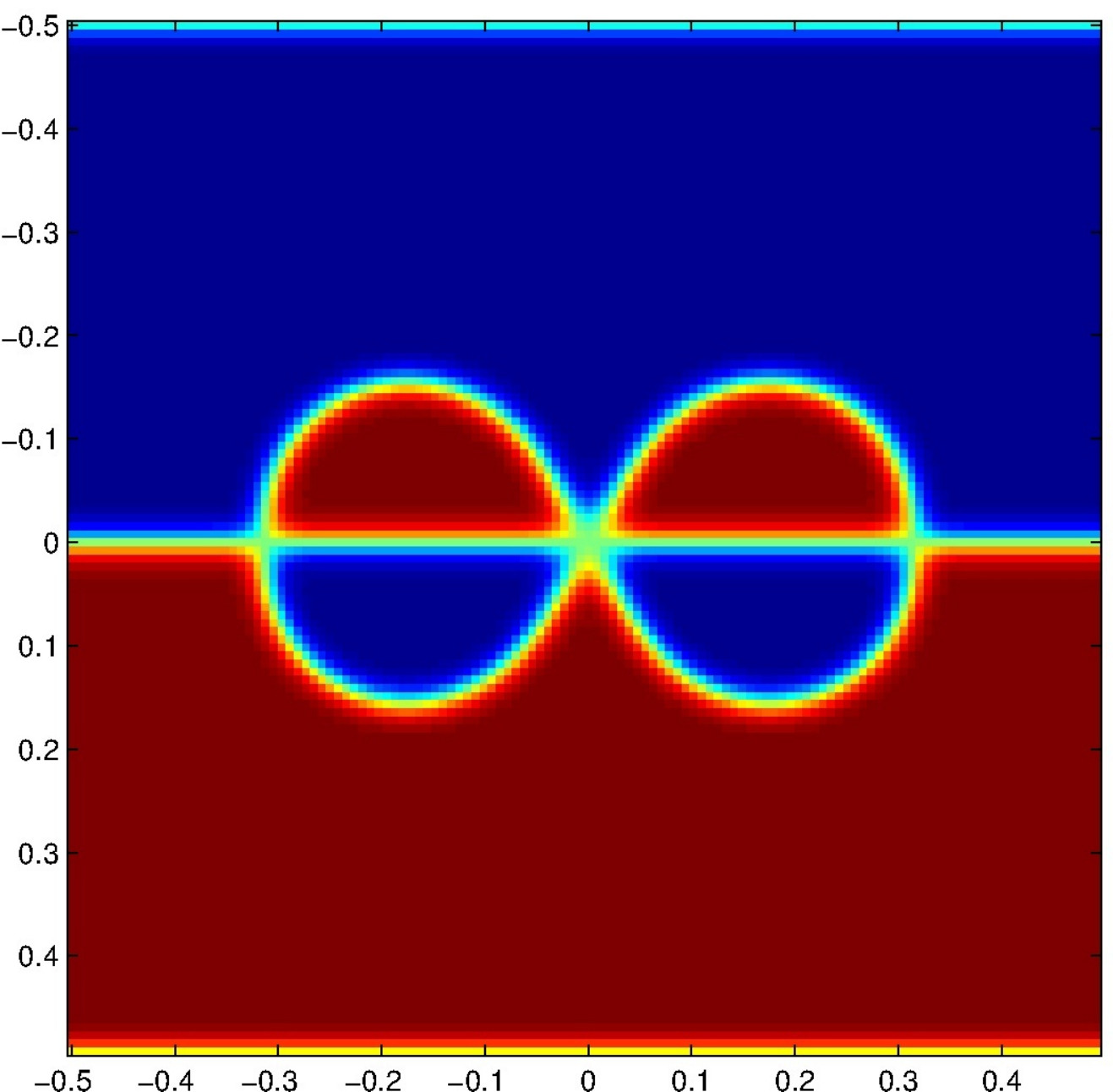}
\includegraphics[width=5cm]{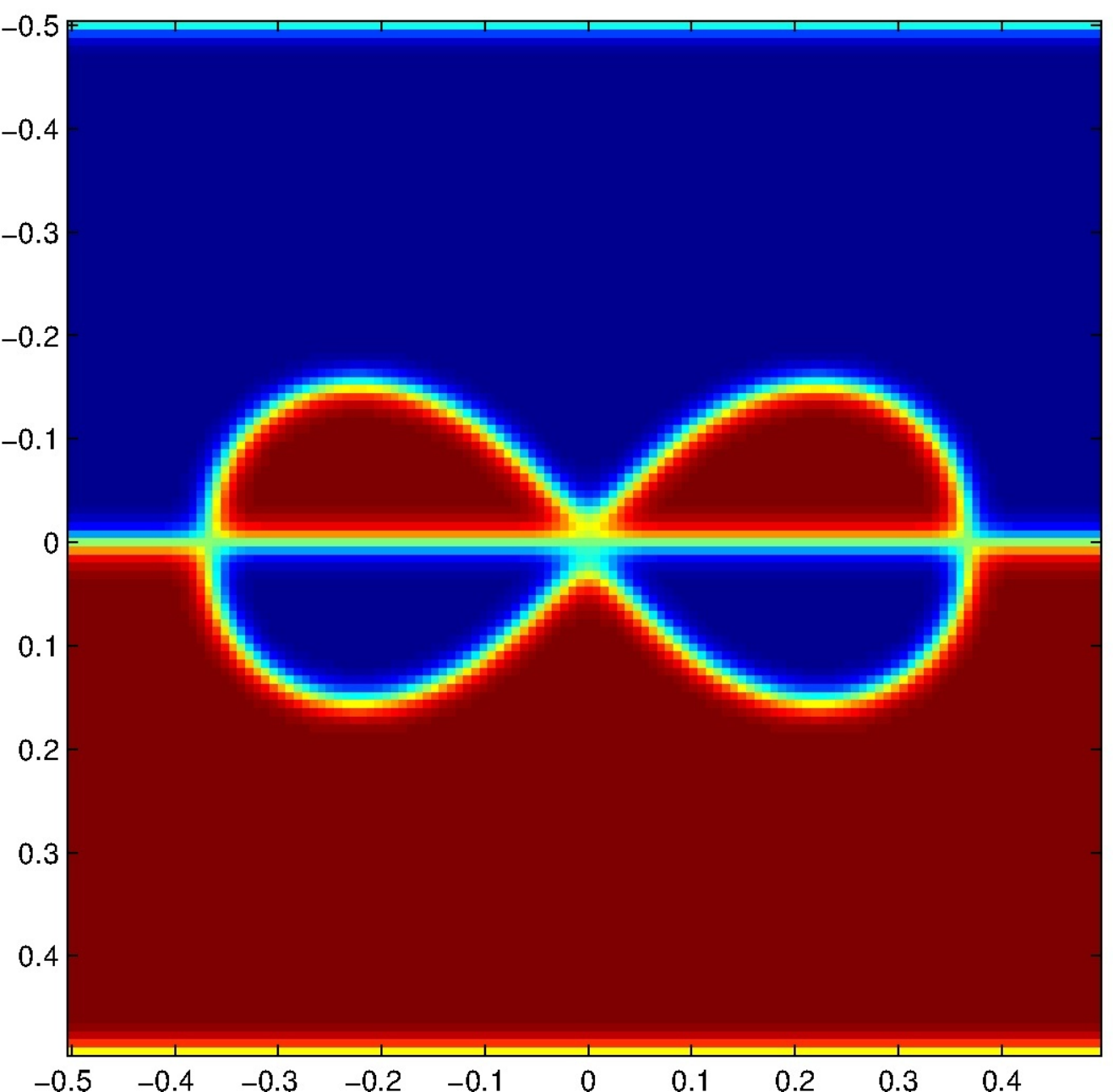} \\
\caption{Two examples where the evolution of either one or two circles is not altered by an additional separating line.}
\label{fig:other_experiment_2D}
\end{figure}

~\\
The experiments on Figure  \ref{fig:other_experiment_2D} are in the same spirit. On the first line, we illustrate the evolution of two phases forming a circle cut by a straight line. Note that the circle seems to evolve independently of the line.  The second situation is quite similar with two disjoint circles 
cut by a line, and the same conclusion holds.

\begin{figure}[!ht] 
\centering
\includegraphics[width=5cm]{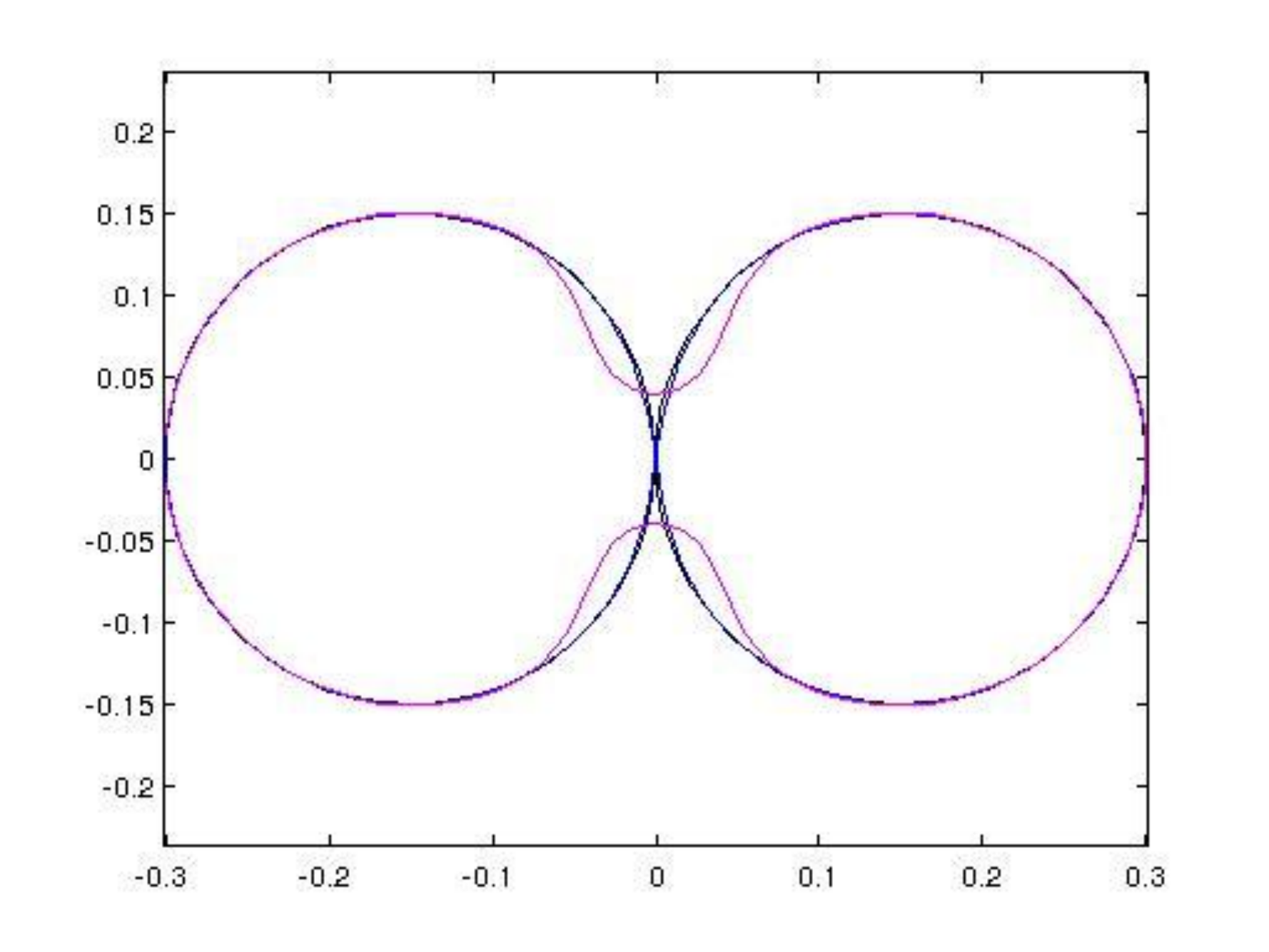}
\includegraphics[width=5cm]{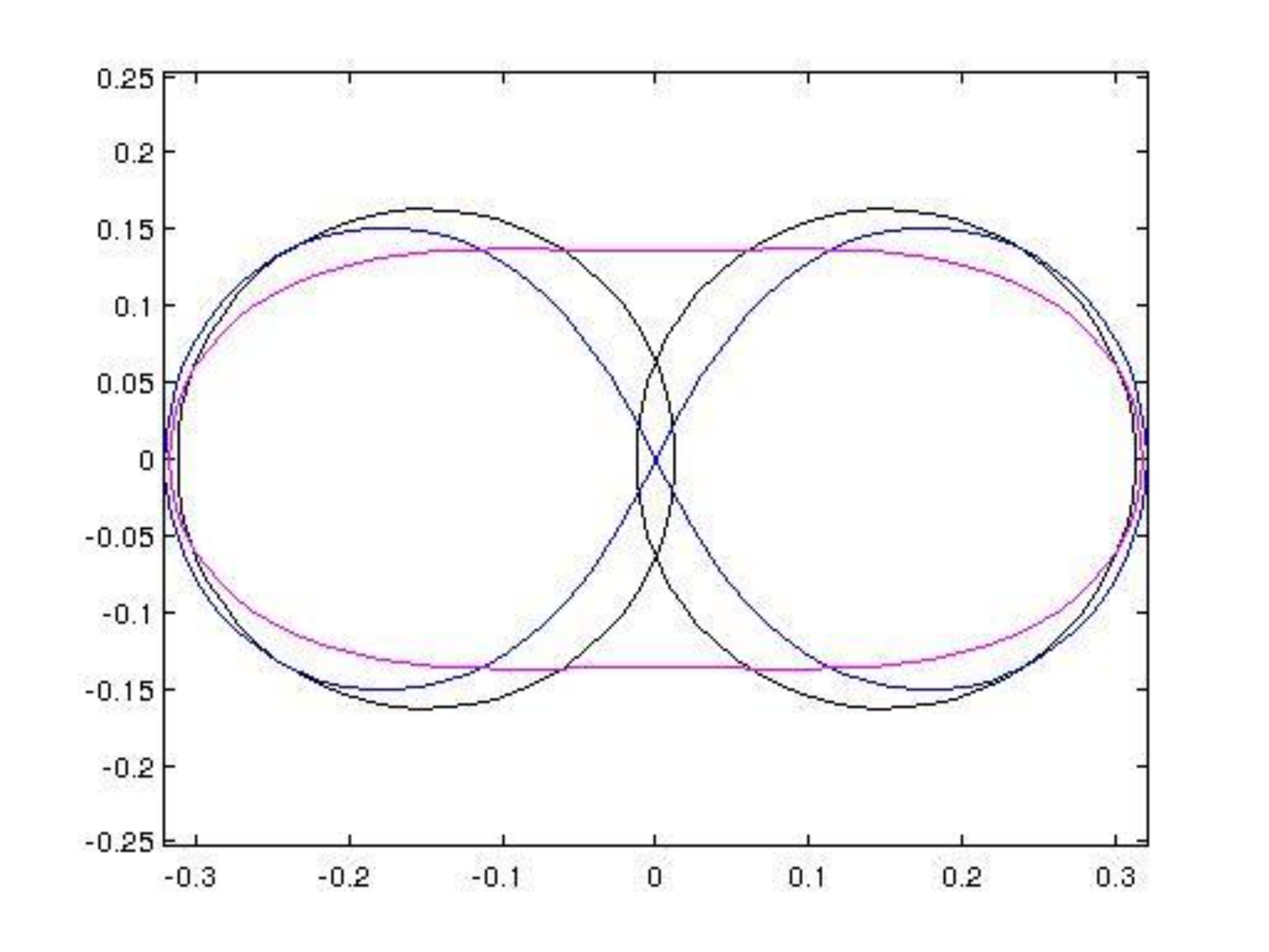}
\includegraphics[width=5cm]{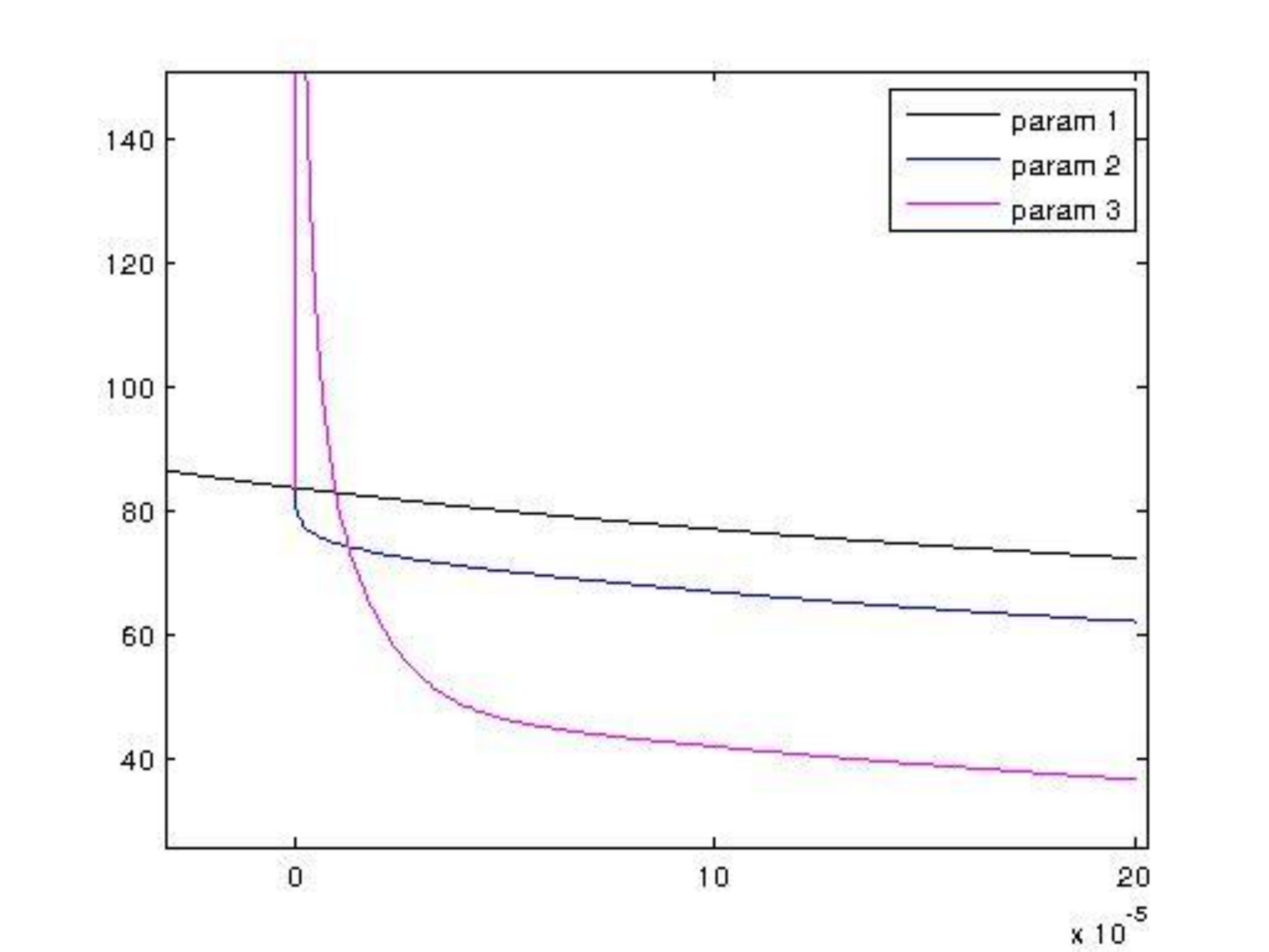} \\
\caption{Parametric evolution of two contiguous circles associated to three different initial parametrizations; 
Left : Interfaces at  $t=10^{-5}$ ; Middle : Interfaces at $t = 5.10^{-5}$ ; 
 Right : Evolution of the Willmore energy of each curve.}
\label{fig:evolution_two_circles_comp}
\end{figure}

~\\
\paragraph{Evolution by the classical diffuse Willmore flow of two contiguous circles}~\\
We now compare in Figure \ref{fig:evolution_two_circles_comp}  the evolution by a discrete parametric Willmore flow of the different interfaces obtained from three different initial parameterizations of two contiguous circles. In the first parametrization (in blue), both circles are parametrized independently. The second parameterization (in magenta) corresponds to the covering of the two circles by a unique smooth parametric curve that self-crosses at the origin. The third parameterization (in red) corresponds to the singular curve that does not cross the horizontal axis at the origin (thus forming a double cusp point). 
The first two pictures in Figure \ref{fig:evolution_two_circles_comp} show the three interfaces obtained at different times.  The third one depicts the evolution of the Willmore energy associated to each  evolving interface. It is interesting to compare the second and the third parameterization. The energy of the second parameterization (i.e. passing from two circles to the eight-type curve after contact) decreases smoothly, therefore the parameterization seems to relate naturally to a continuous evolution. In contrast, the energy of the third parameterization explodes at contact, and then decreases strongly to become the lowest (after time $t>10^{-5}$) with respect to the other parameterizations. This experiment illustrates clearly the bifurcation at contact, and justifies why different configurations have been observed in the literature. For instance, there is no crossing observed in~\cite{FrRuWi11} because the authors used a large time step $\delta_t$ and therefore ignored the contact. However, after a while, the energy of the non-crossing configuration is indeed the best. It may be argued that a continuous flow should go up to contact, and therefore a numerical flow that is accurate enough to capture the singularity should be the best. On the other hand, once the crossing configuration has been chosen, there is no way to have an energy as low as the non crossing configuration's energy.

\paragraph{Experiments in space dimension $3$}
For the 3D simulations presented hereafter, we used the parameters : ${\cal P} = 2^7$, $\varepsilon = 1.5/{\cal P}$ and 
$\delta_t = 1/10~{\cal P}^{-2}~\varepsilon^2$. \\

\begin{figure}[!ht] 
\centering
\includegraphics[width=5cm]{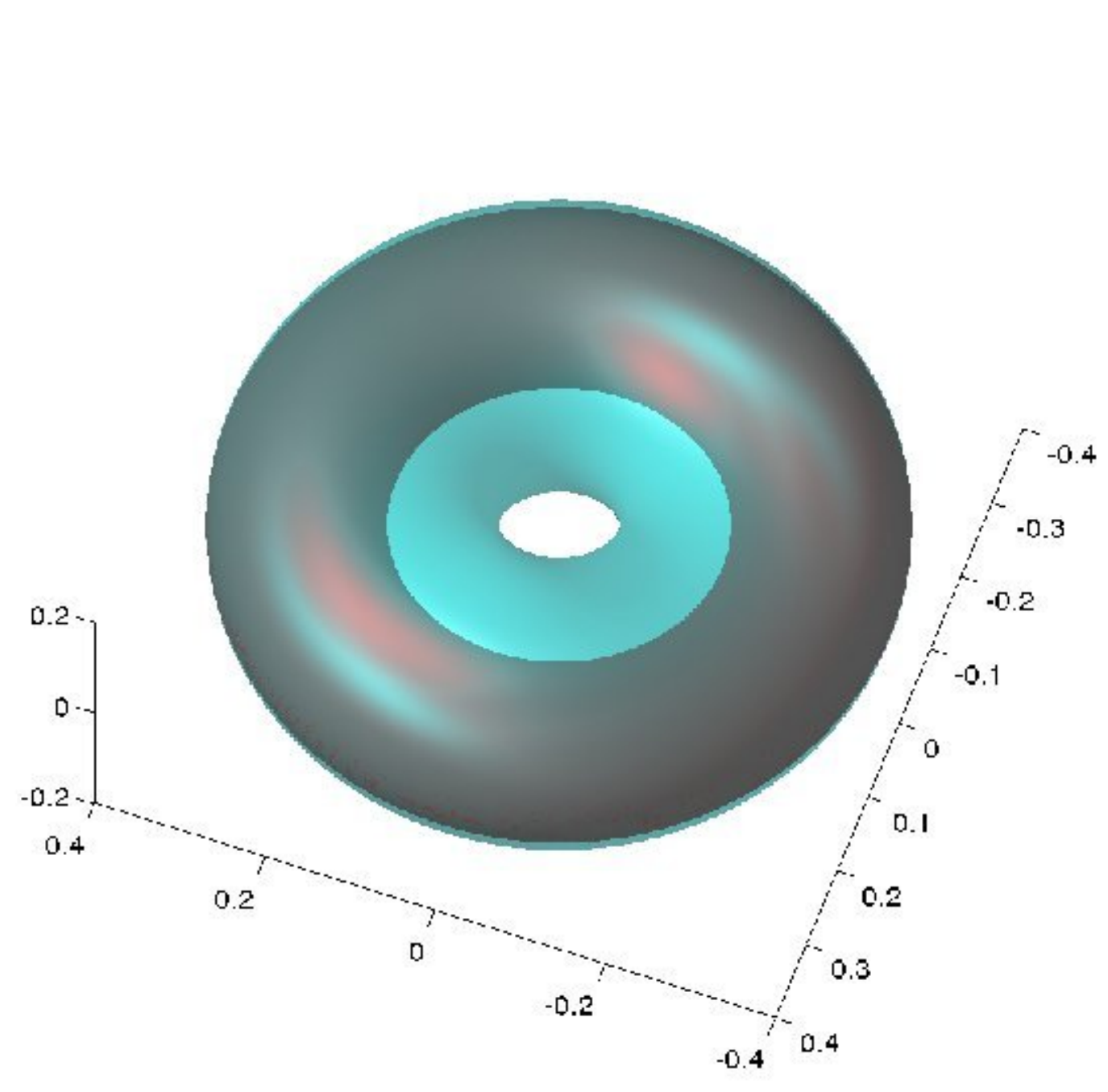}
\includegraphics[width=5cm]{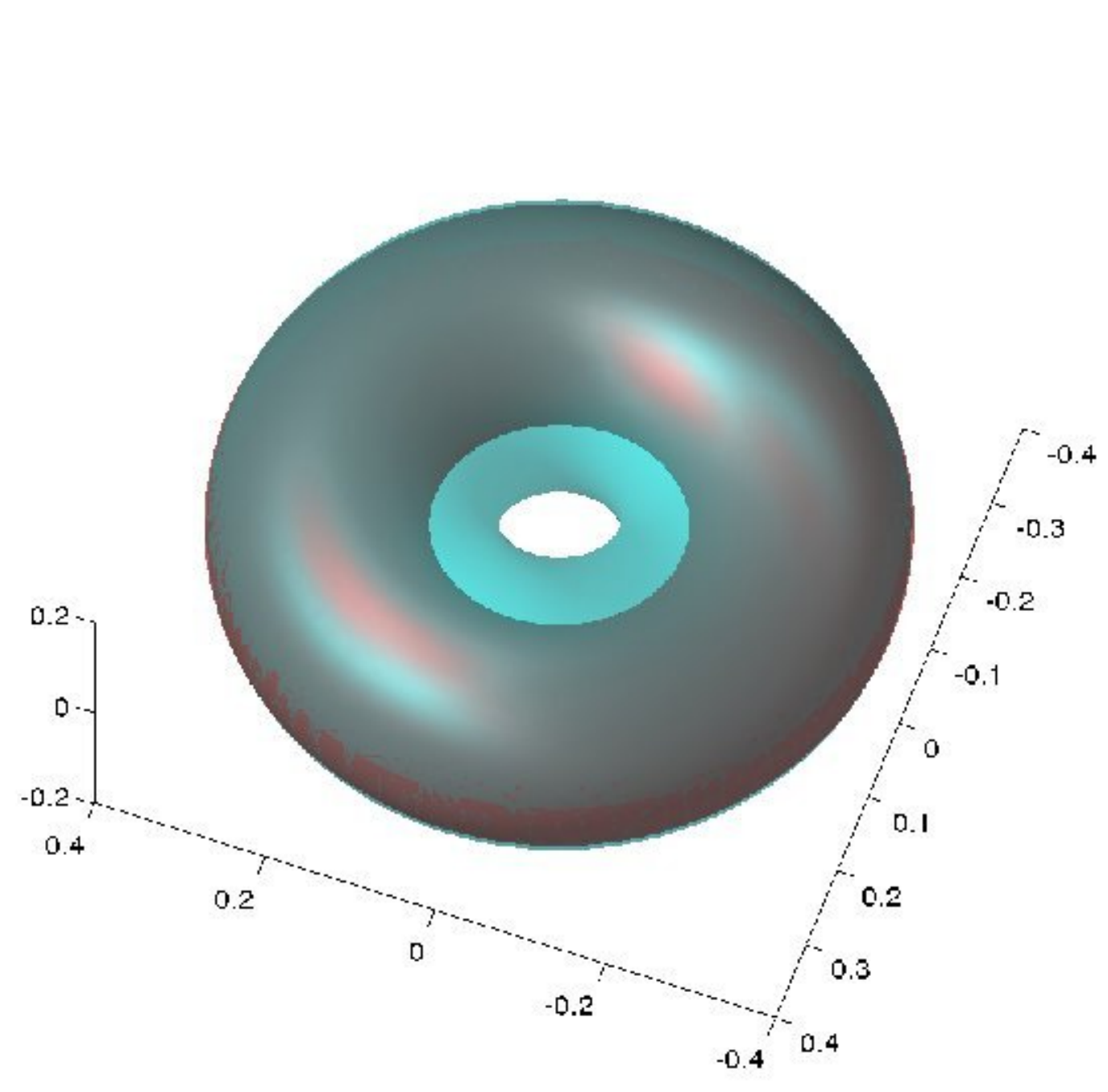}
\includegraphics[width=5cm]{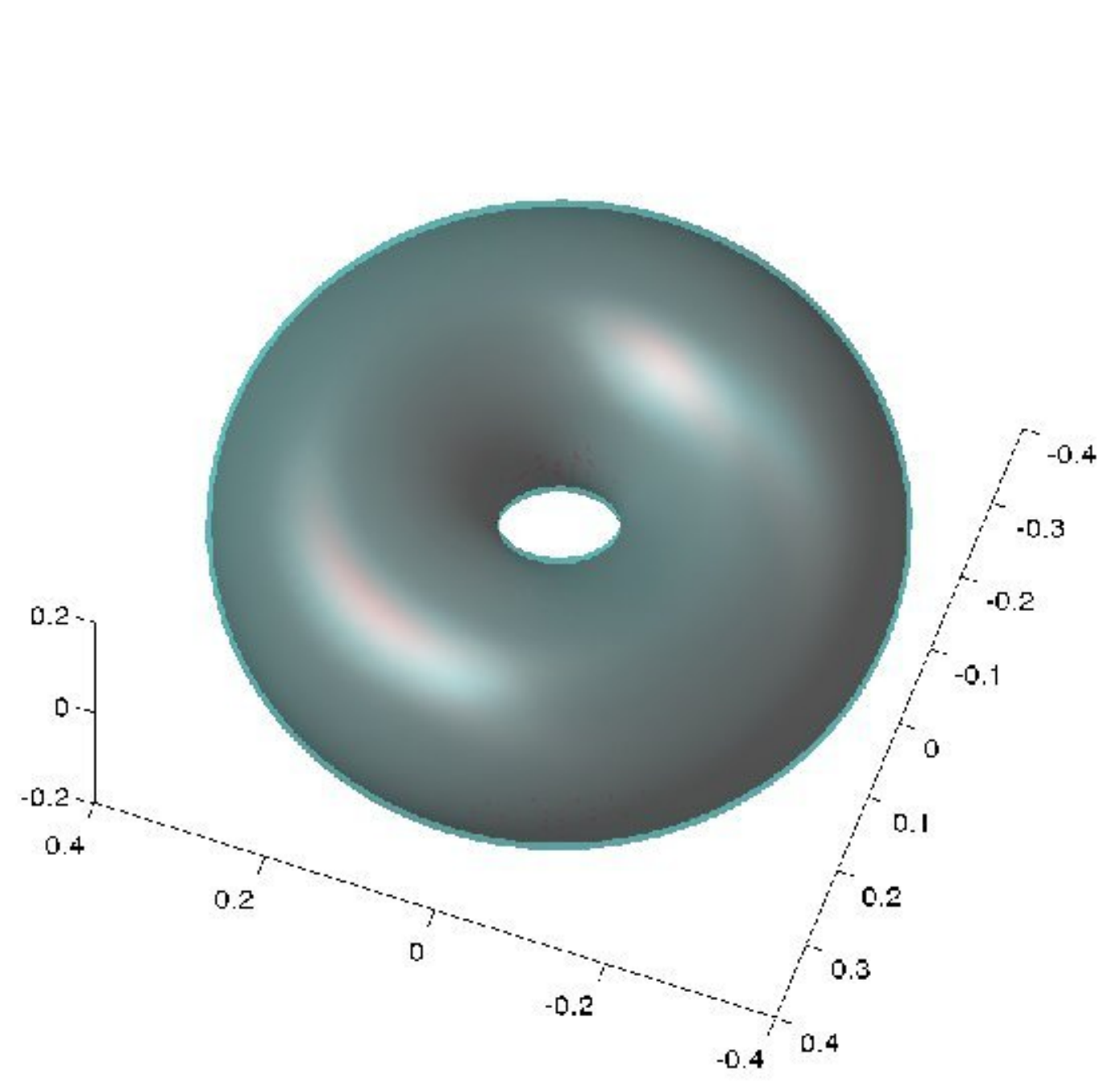}\vspace*{-0.2cm} \\
\includegraphics[width=5cm]{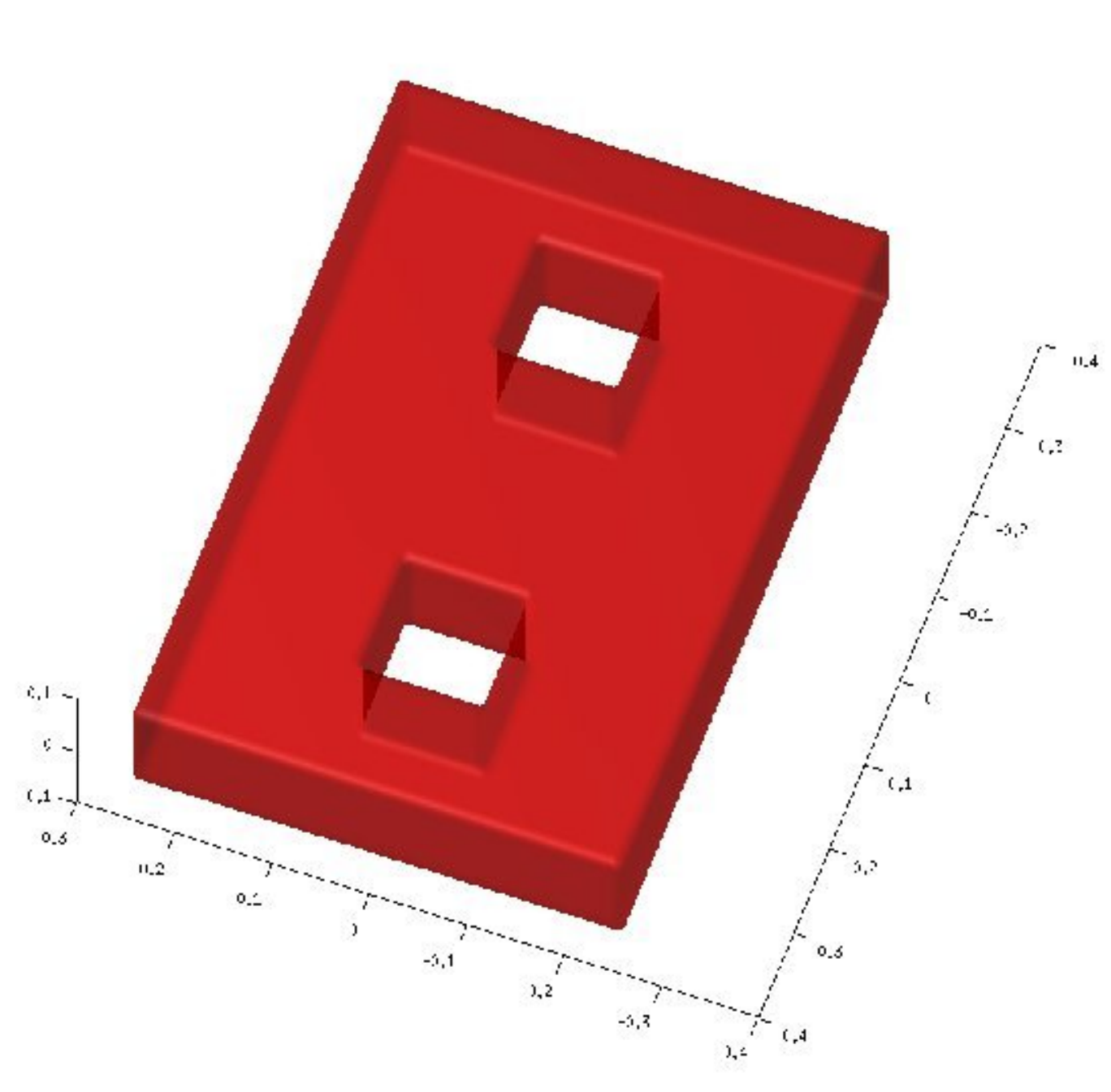}
\includegraphics[width=5cm]{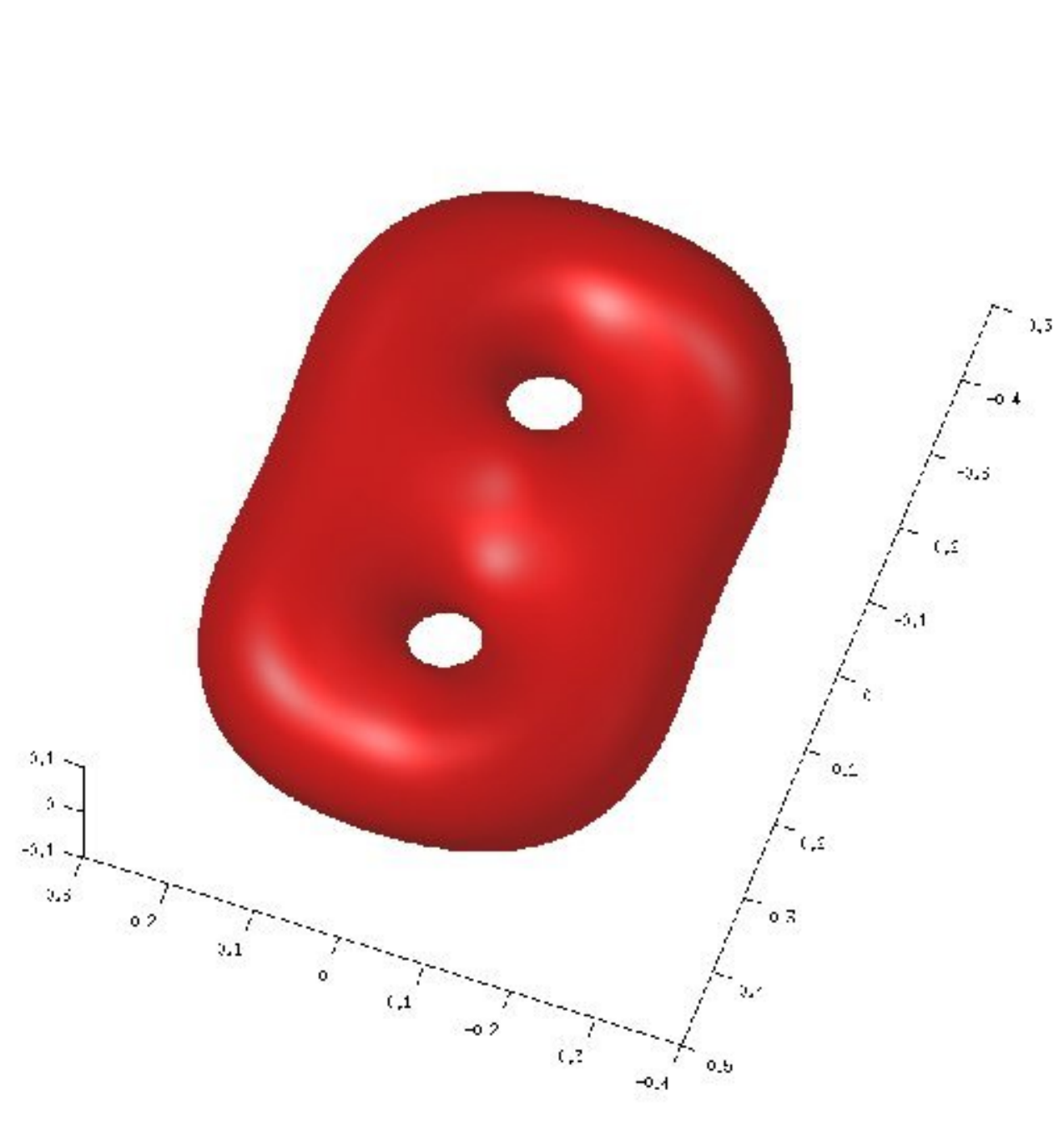}
\includegraphics[width=5cm]{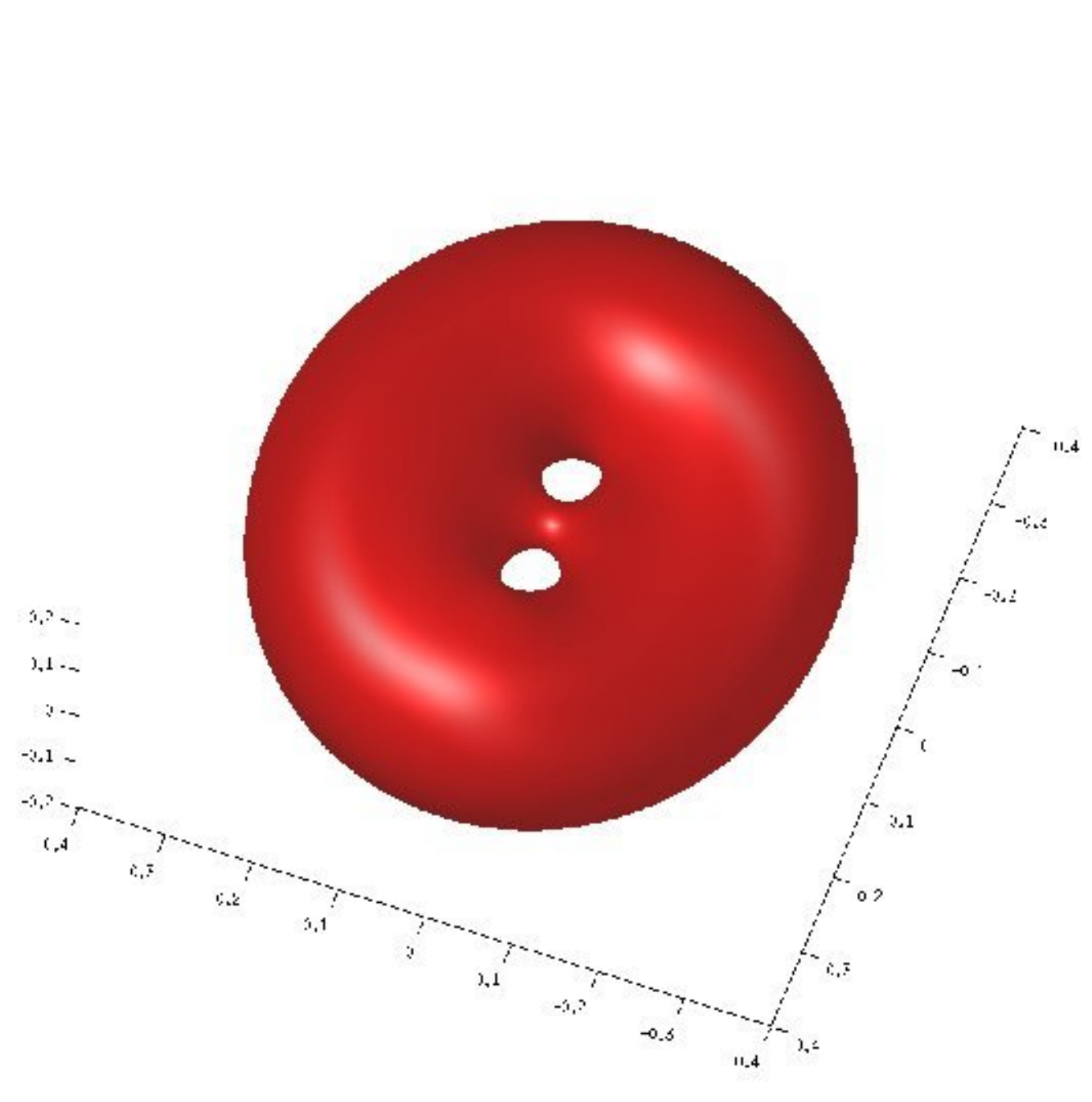}\vspace*{-0.2cm} \\
\includegraphics[width=5cm]{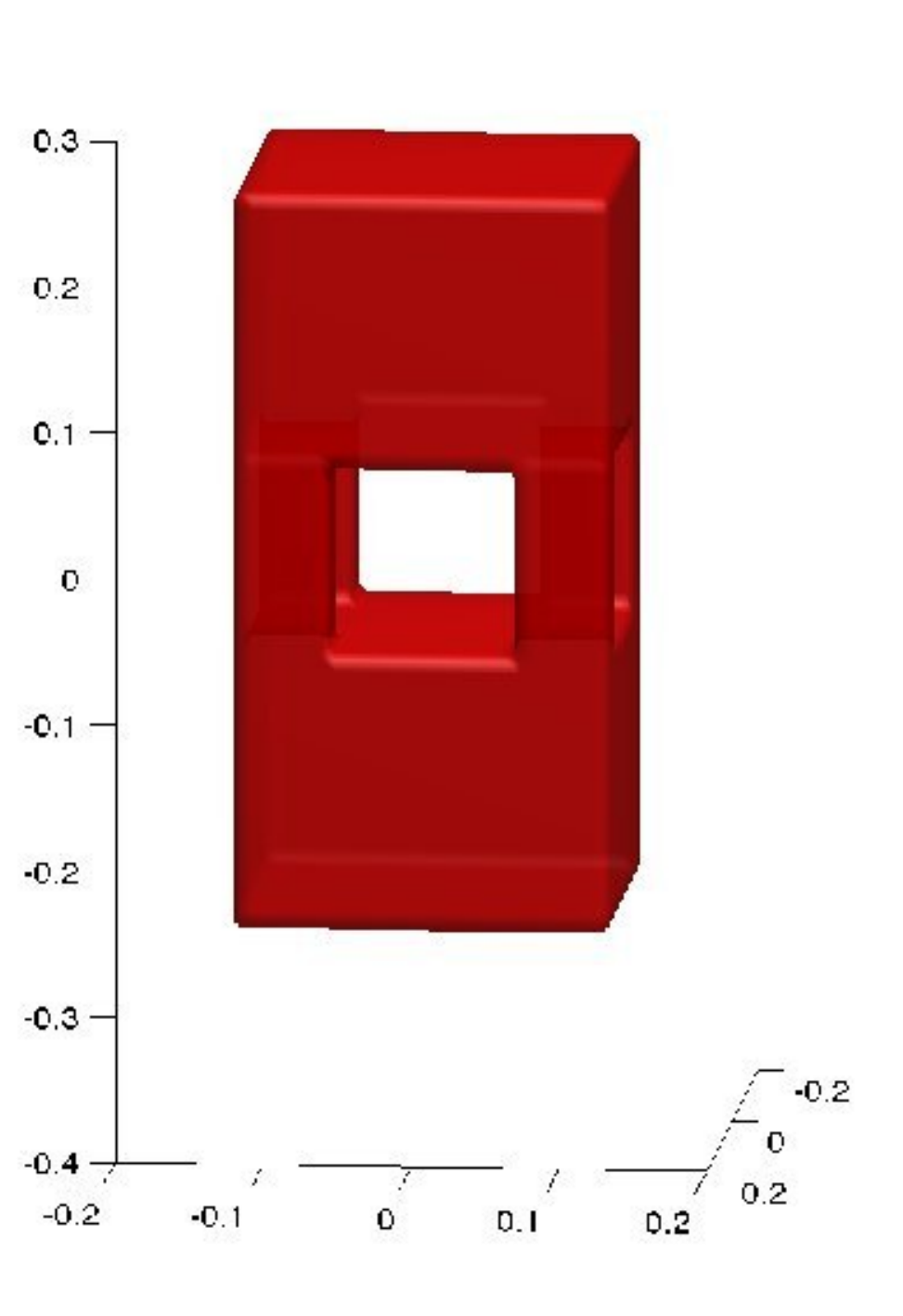}
\includegraphics[width=5cm]{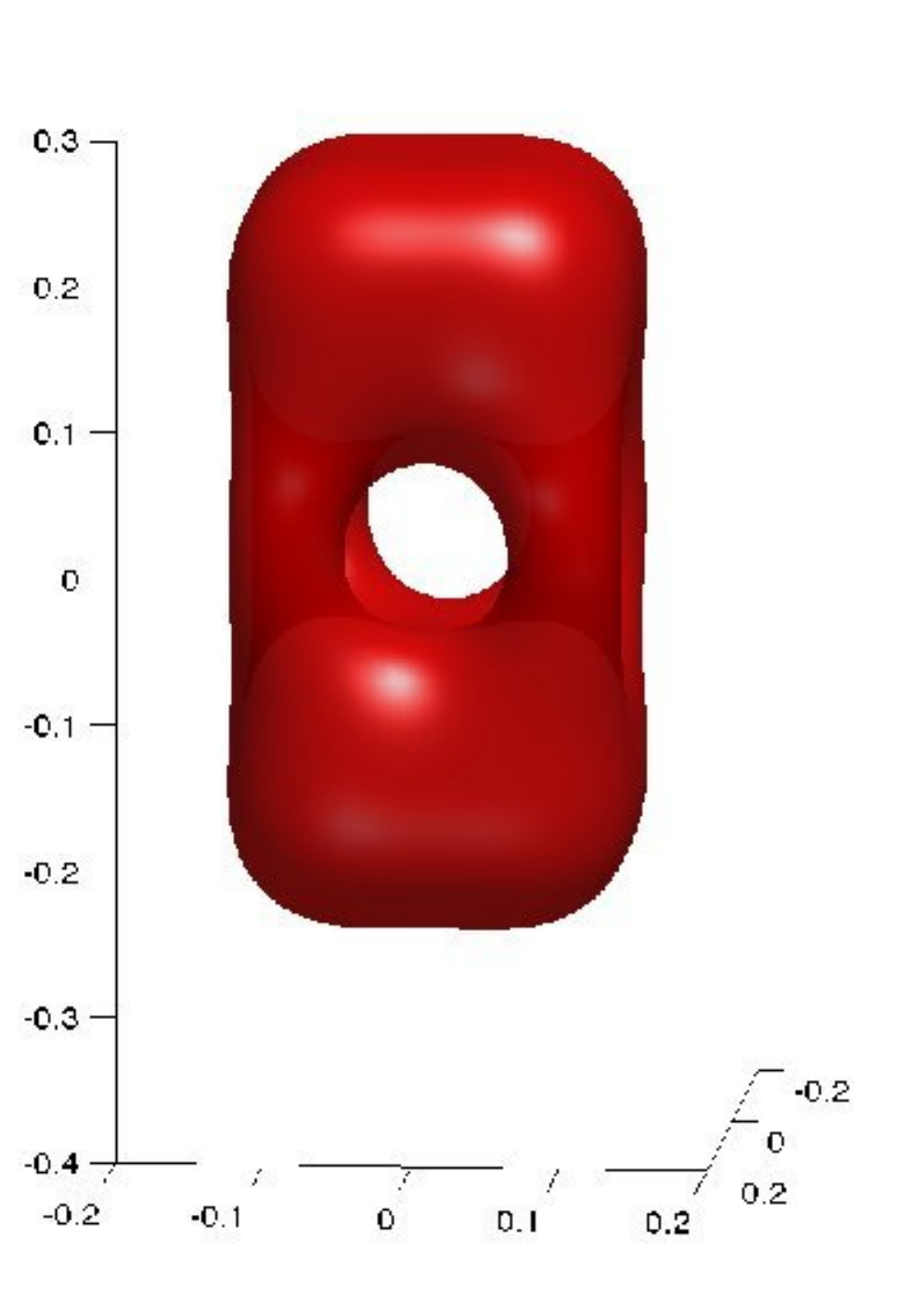}
\includegraphics[width=5cm]{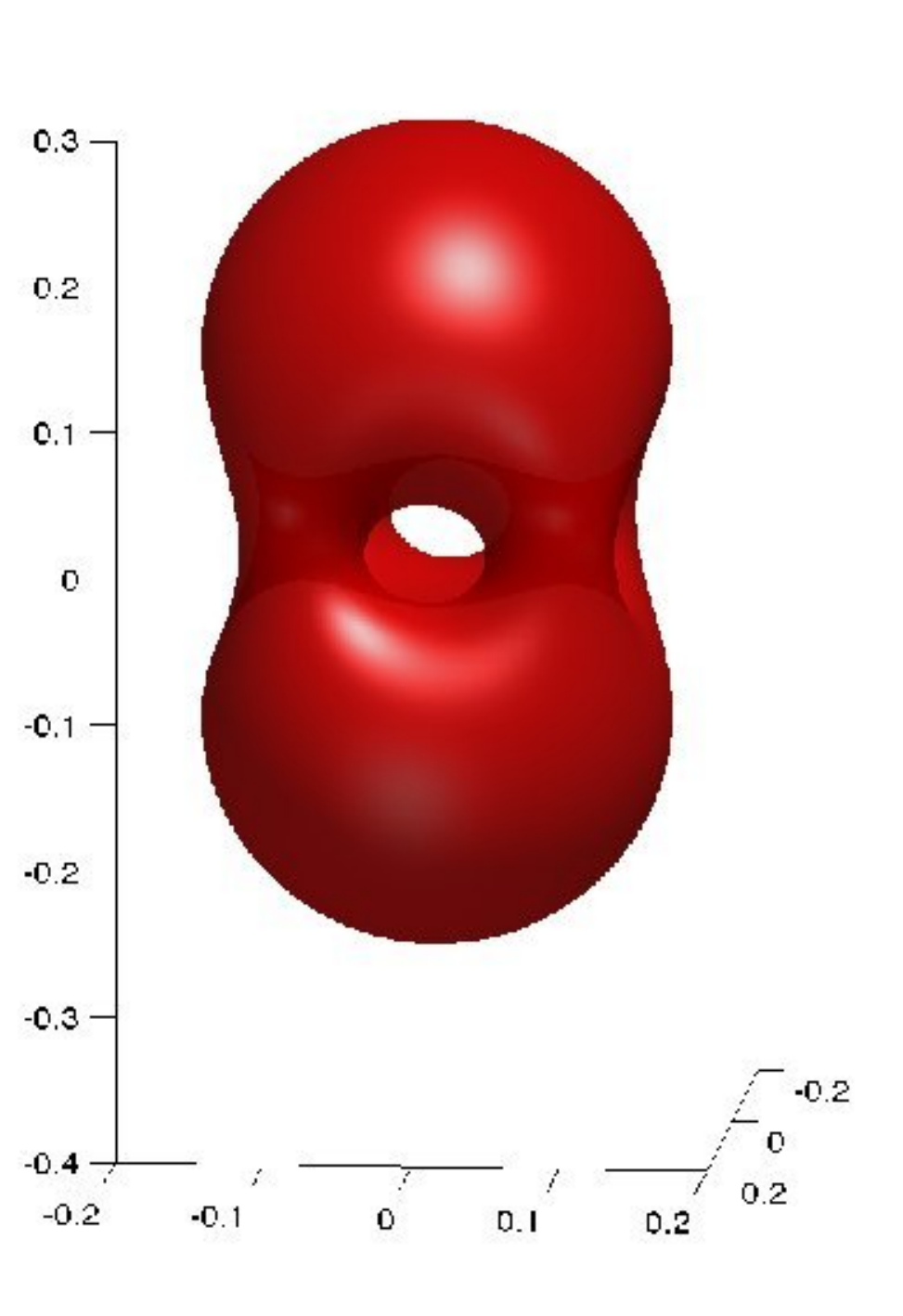}
\caption{Smooth evolution of $\Gamma(t)$ by the classical diffuse Willmore flow in $3$ D. First line: convergence to a Clifford torus (in blue). Second and third lines: convergence to Lawson-Kusner's surfaces of genus 2 and 4, respectively.}
\label{fig:evolution_tore}
\end{figure}

The first simulation illustrates the evolution of a torus. According to the Willmore conjecture, which seems to have been proved in~\cite{Willmore_conjecture}, the torus that minimizes the Willmore energy
is Clifford's, whose ratio between both radii equals $\sqrt{2}$. 
In the first line of Figure \ref{fig:evolution_tore}, we plot  for different values of $t$, a Clifford torus (in blue) 
and the interface $\Gamma(t)$ (in red) obtained numerically. As expected, the interface $\Gamma(t)$ converges to the Clifford torus. 
The second line of Figure~\ref{fig:evolution_tore}  shows the evolution of a parallelepiped with two holes. The interface converges to a Lawson-Kusner surface of genus $2$, that is conjectured to minimize the Willmore energy among surfaces with genus $2$~\cite{Kusner,HsuKusnerSullivan92}. The same experiment is done for a genus $4$ surface on the last line, and there is again convergence to a Lawson-Kusner surface. 
We believe that these simulations illustrate the good quality of our numerical scheme and its ability to recover some critical points for the Willmore energy.  \\

\begin{figure}[!ht] 
\centering
\includegraphics[width=5cm]{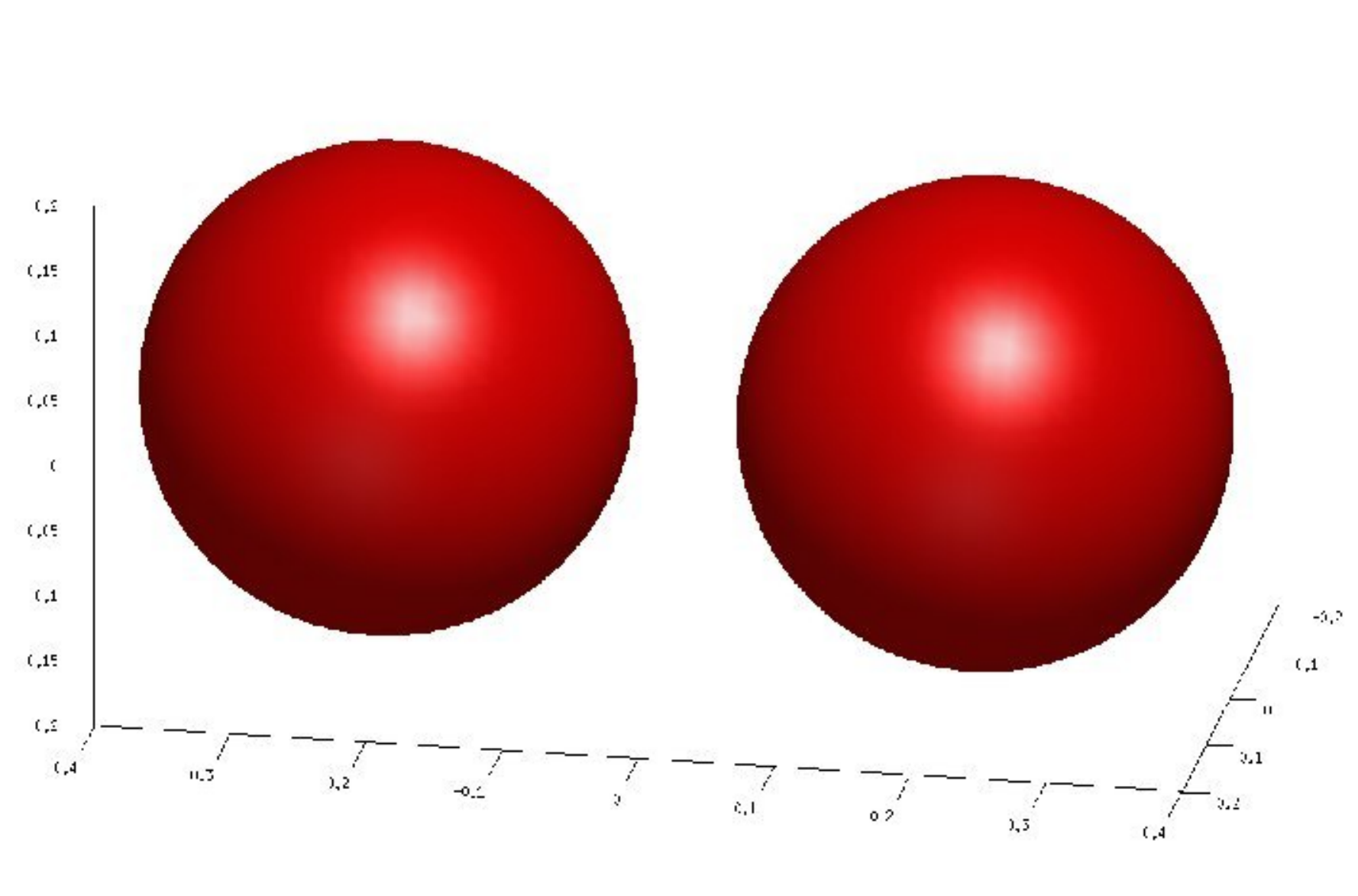}
\includegraphics[width=5cm]{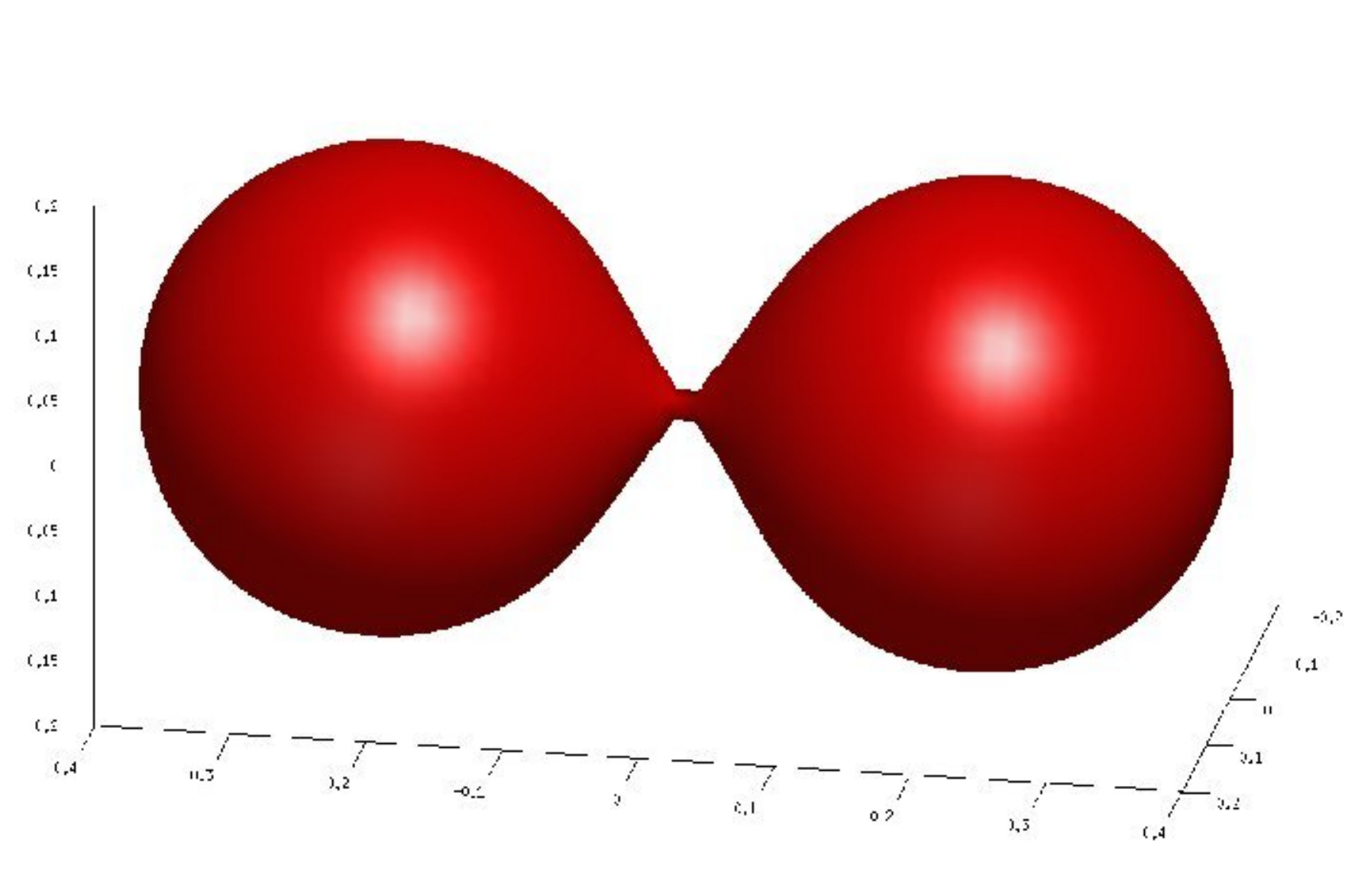}
\includegraphics[width=5cm]{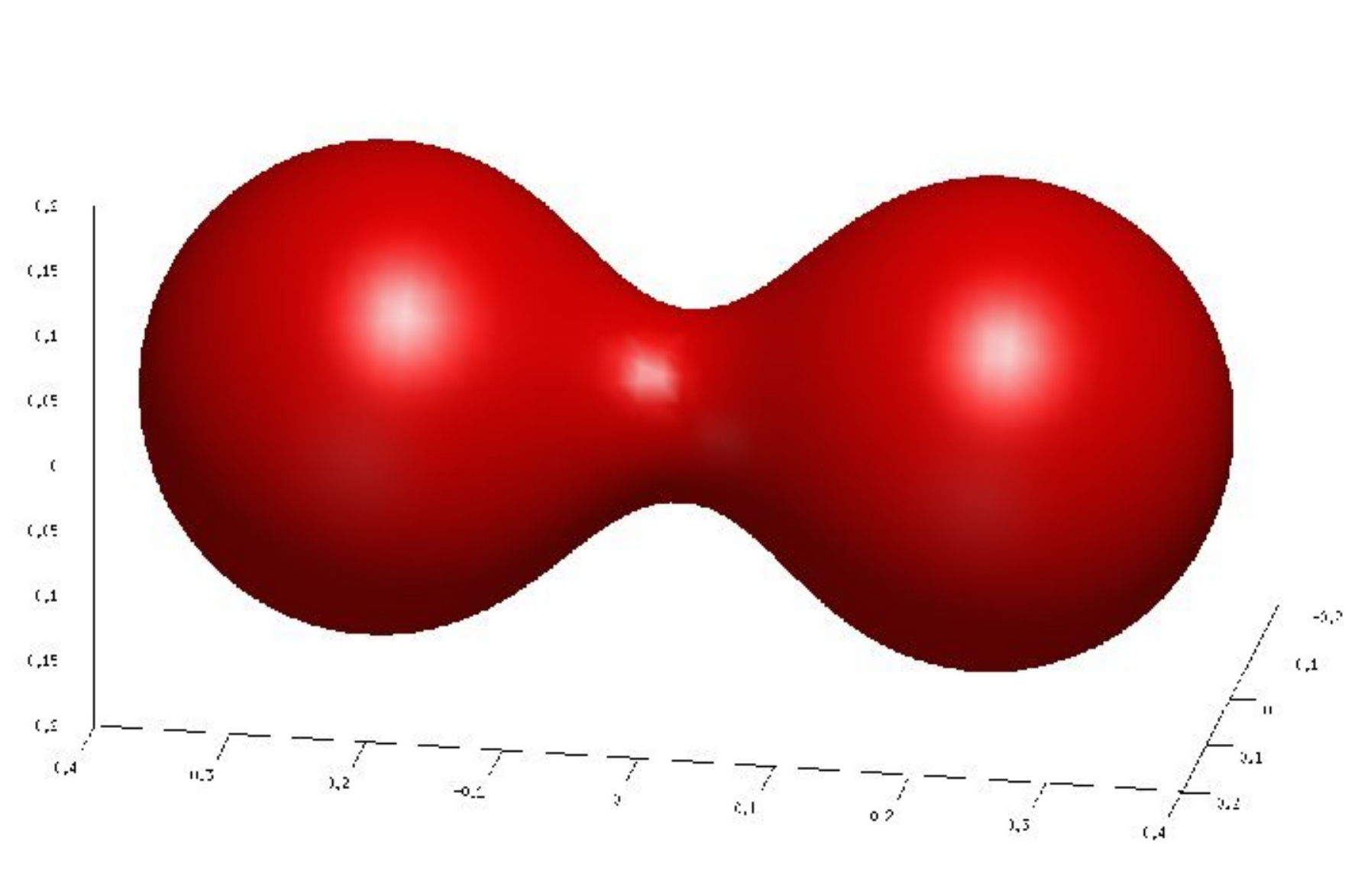} \\
\includegraphics[width=5cm]{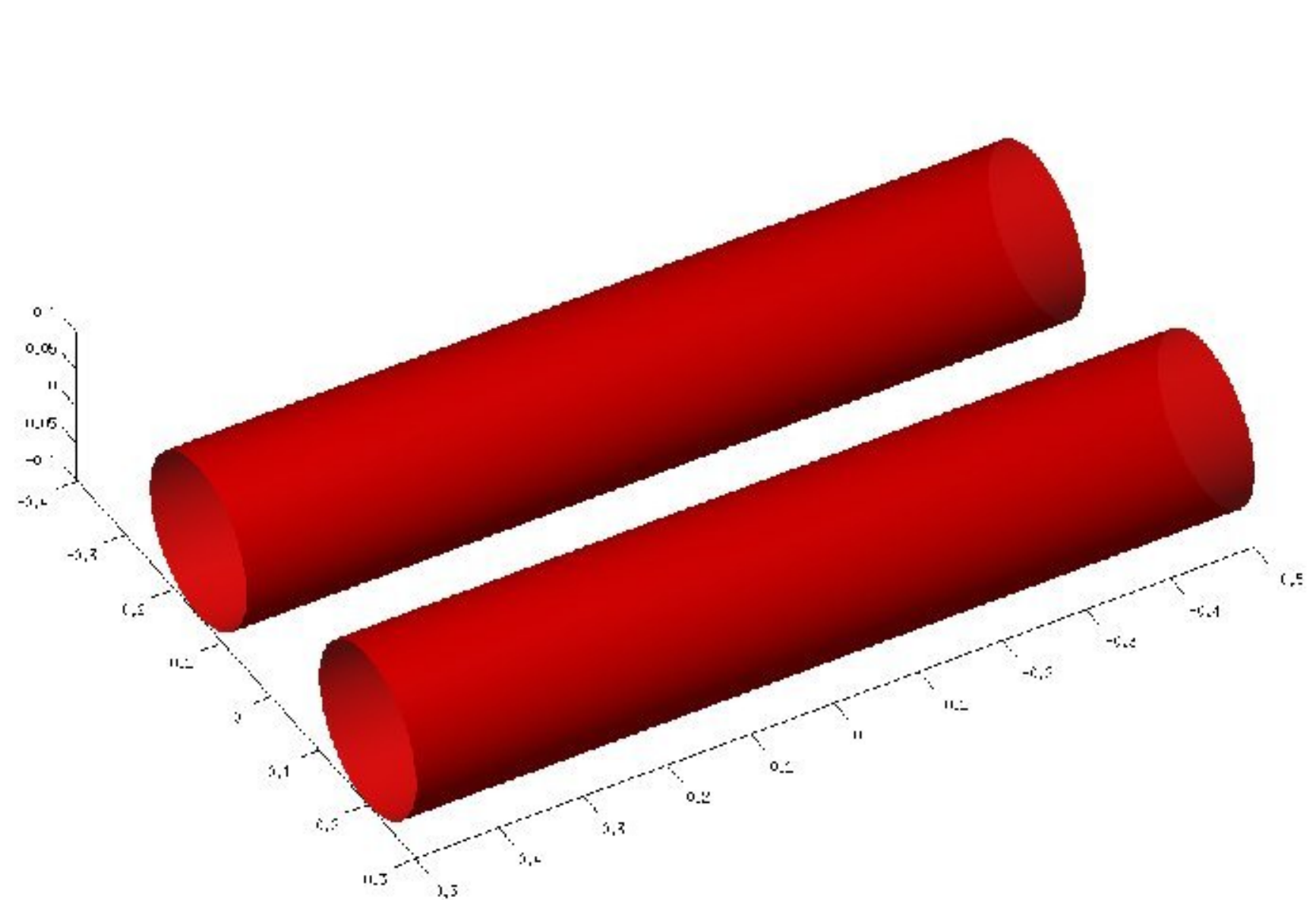}
\includegraphics[width=5cm]{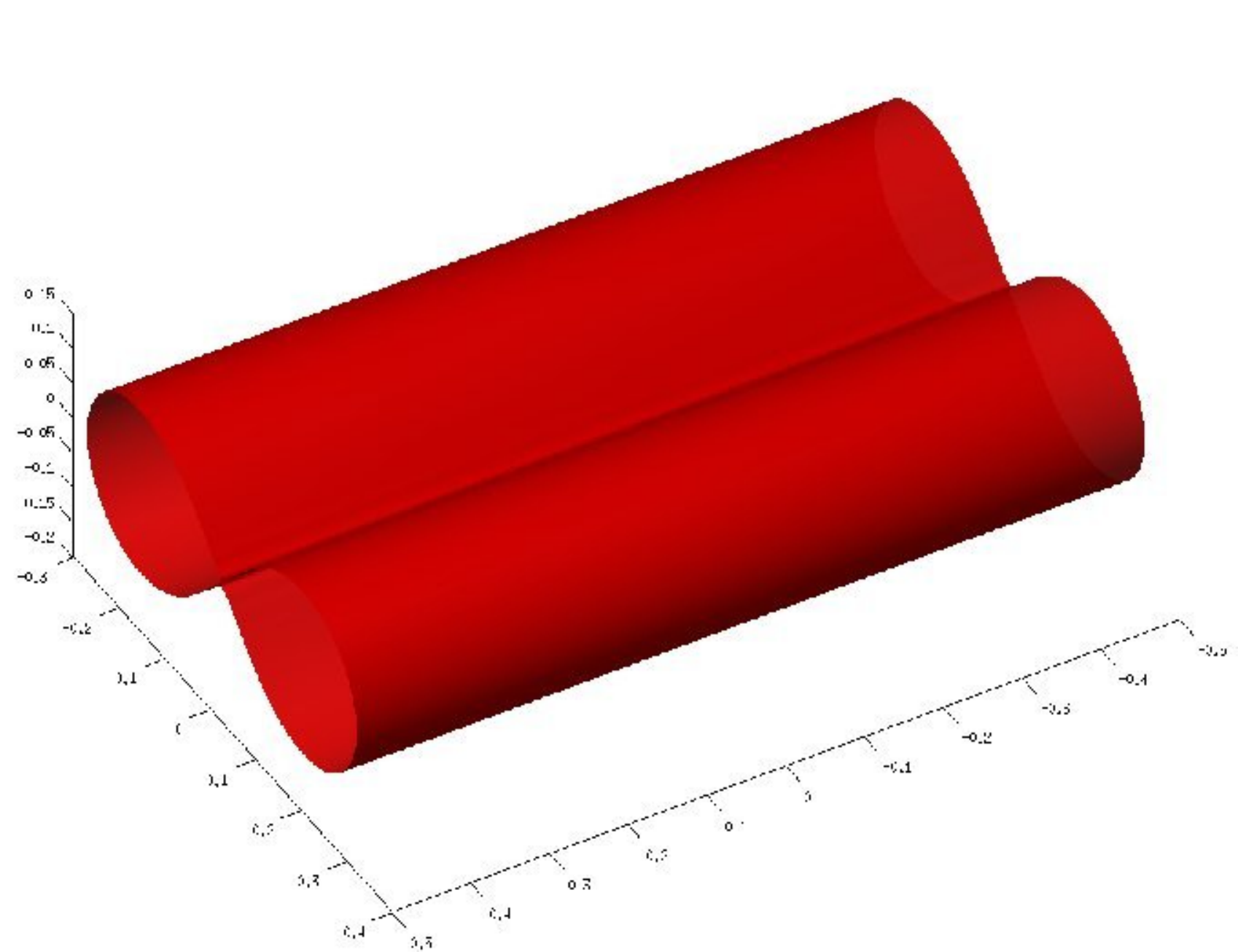}
\includegraphics[width=5cm]{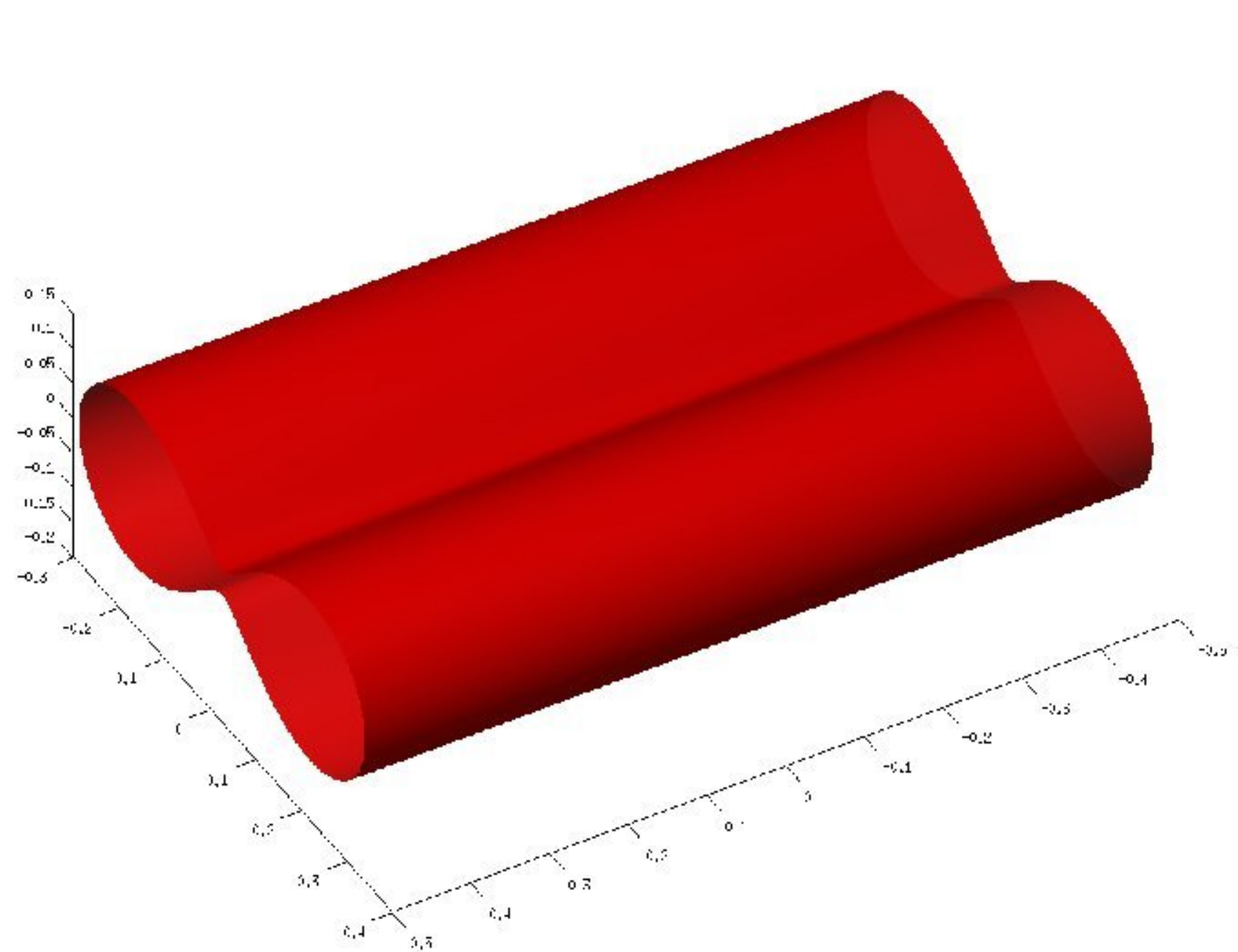} \\
\includegraphics[width=5cm]{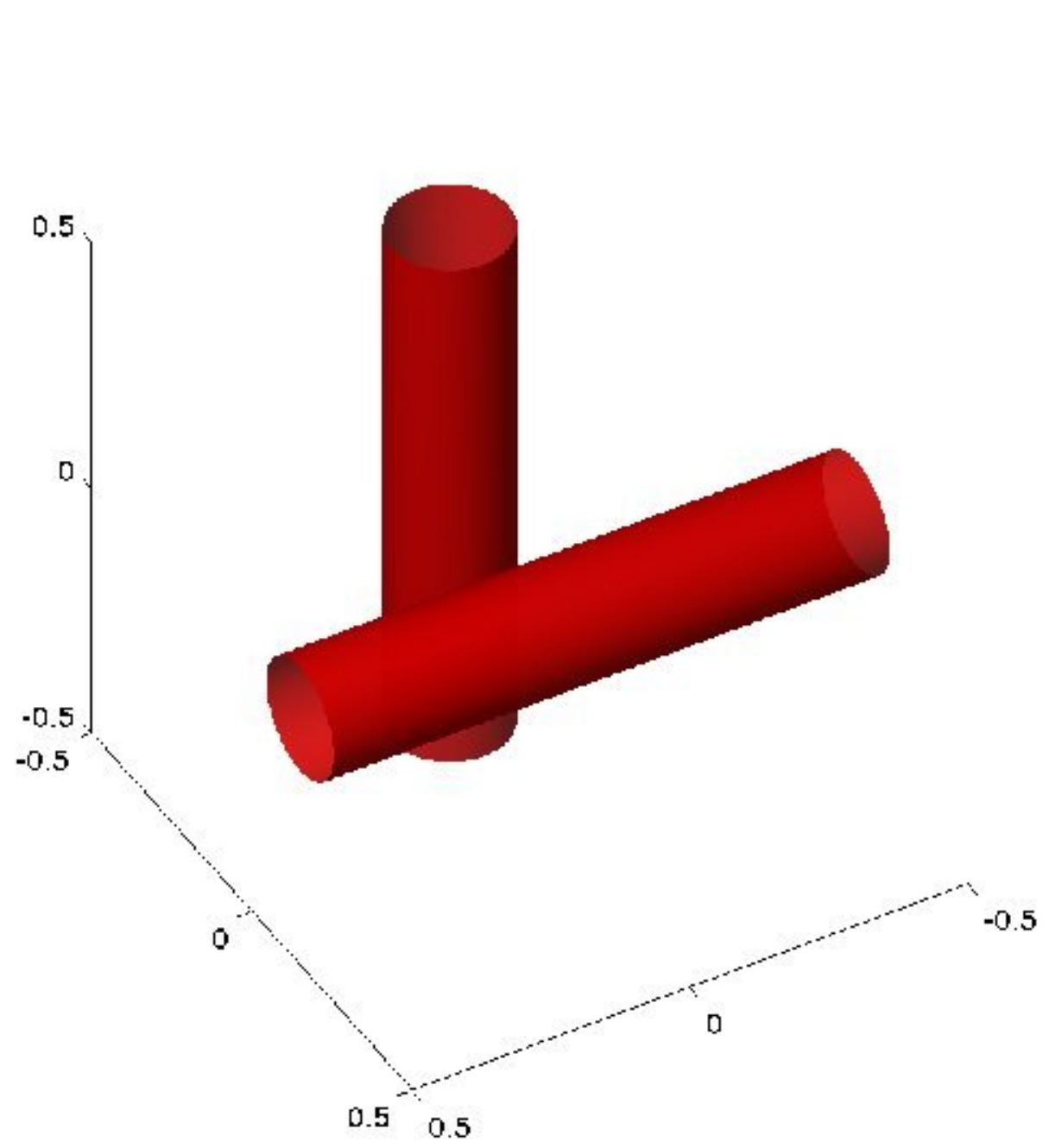}
\includegraphics[width=5cm]{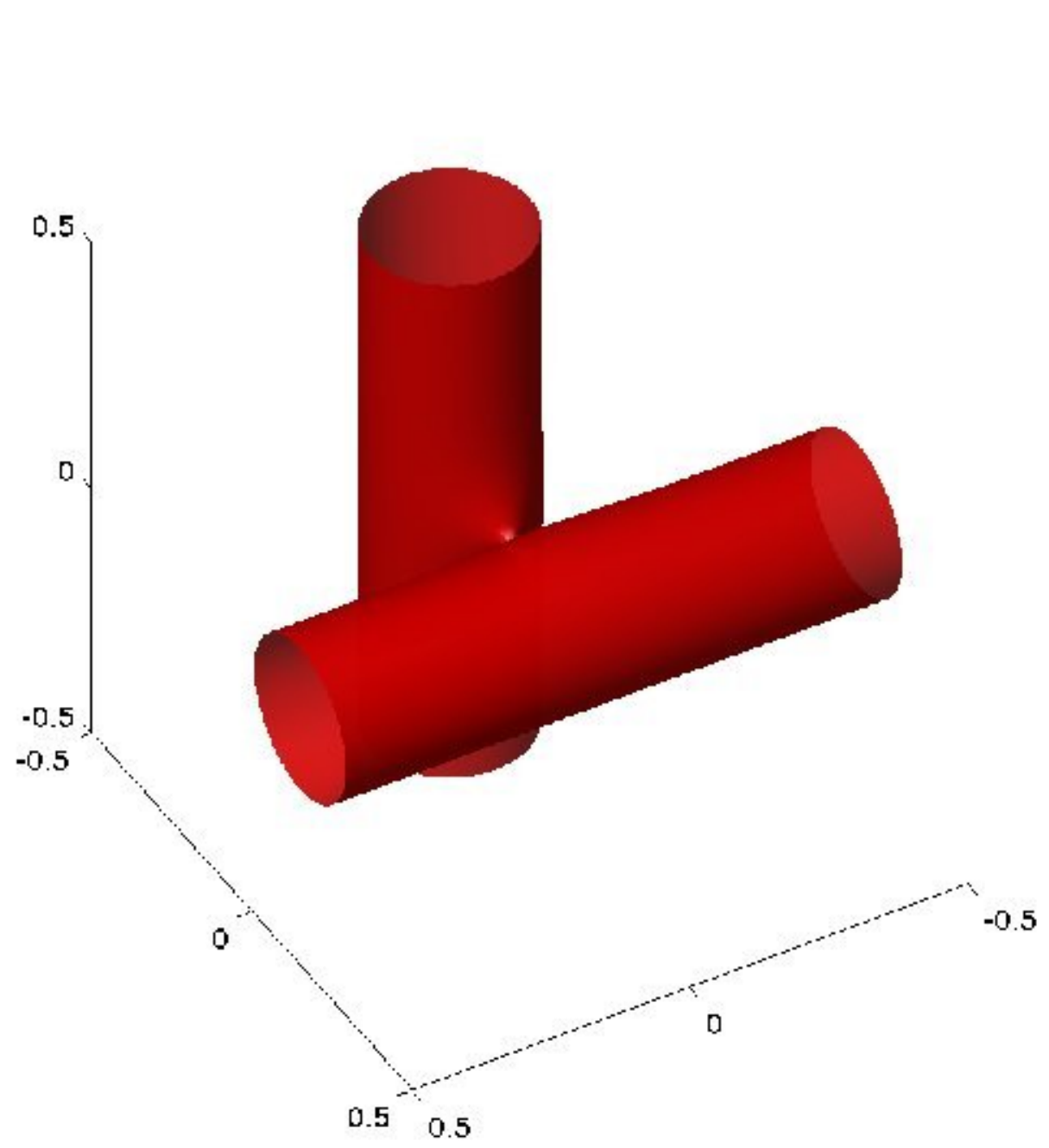}
\includegraphics[width=5cm]{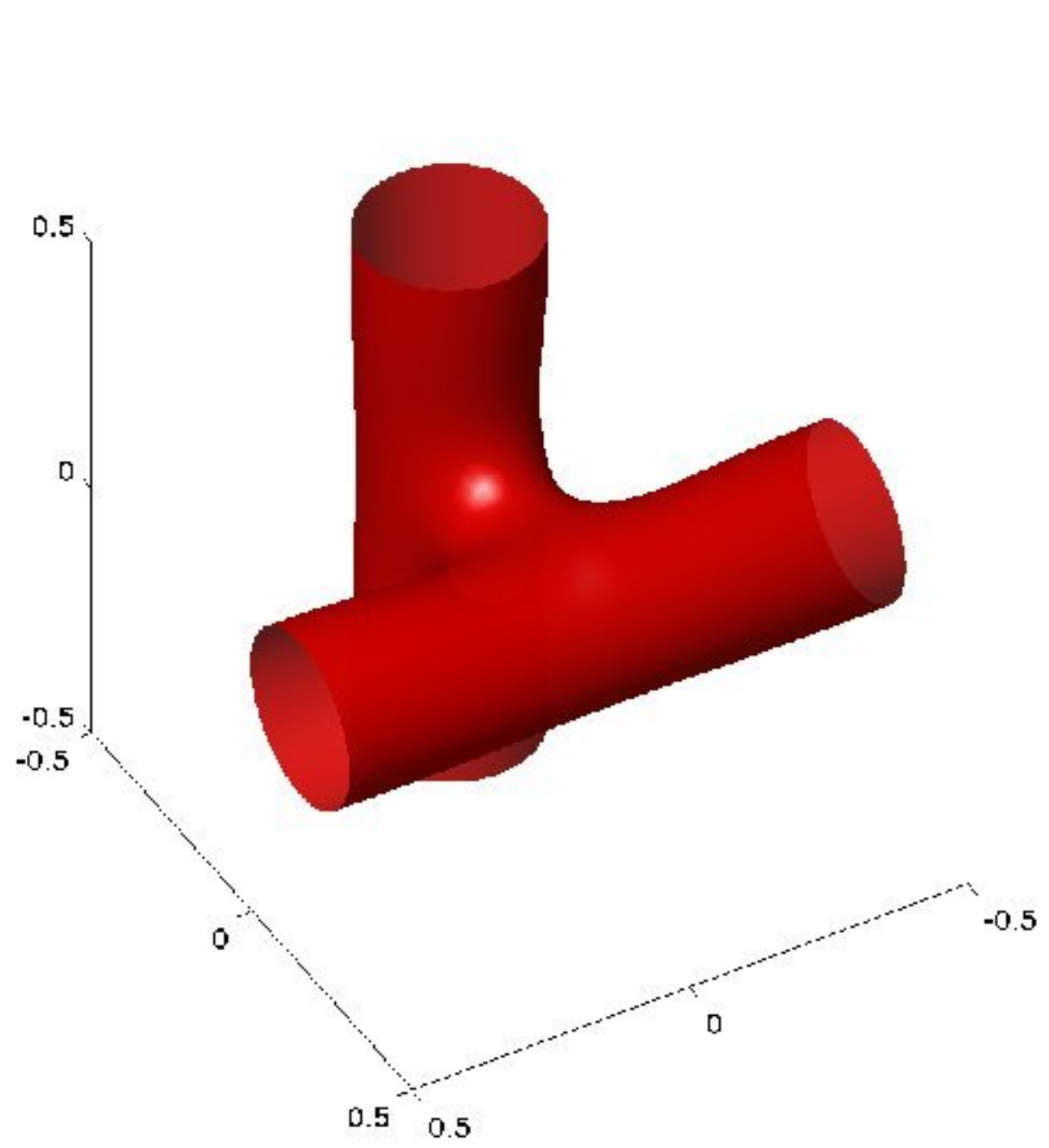} \\
\includegraphics[width=5cm]{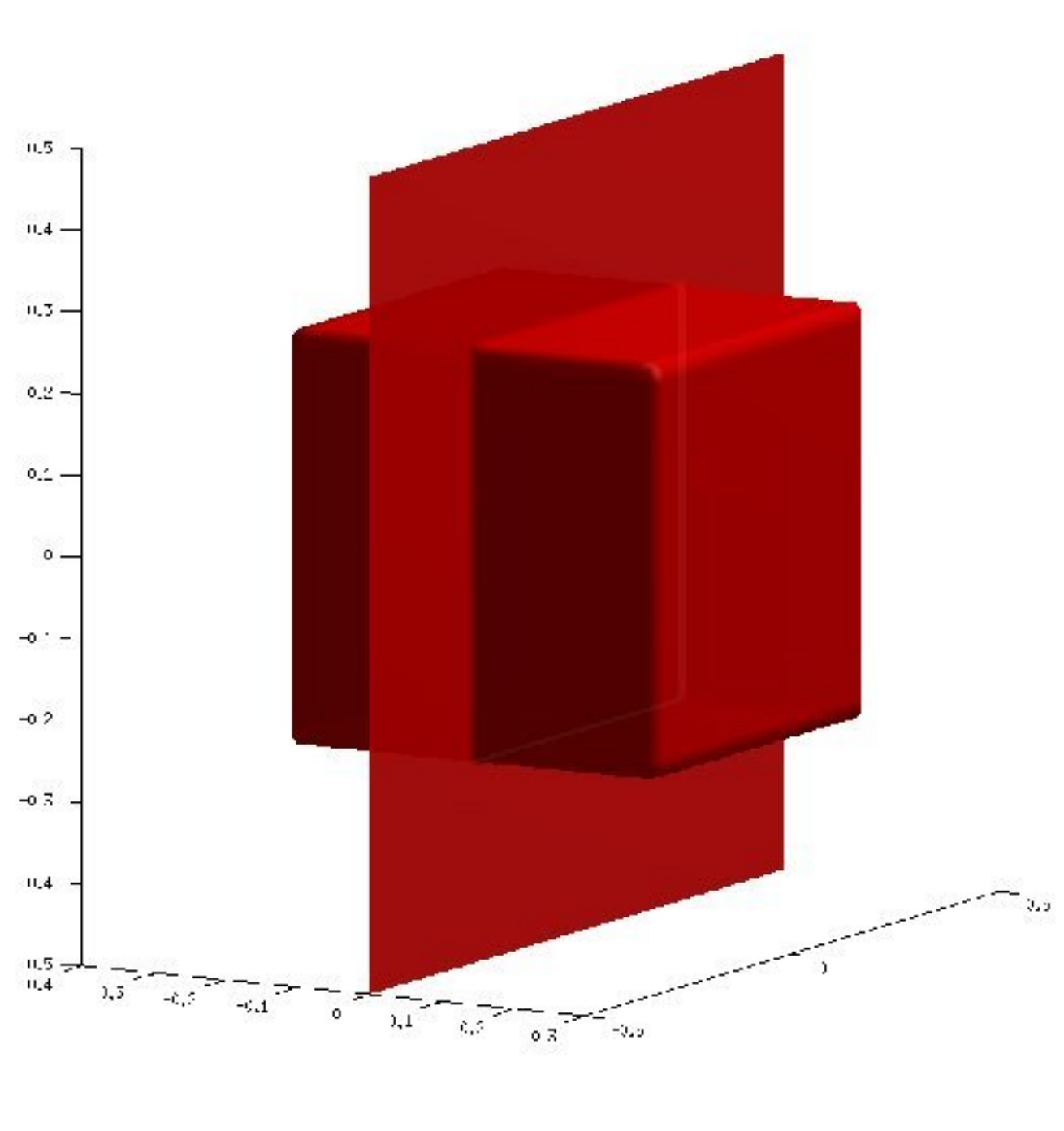}
\includegraphics[width=5cm]{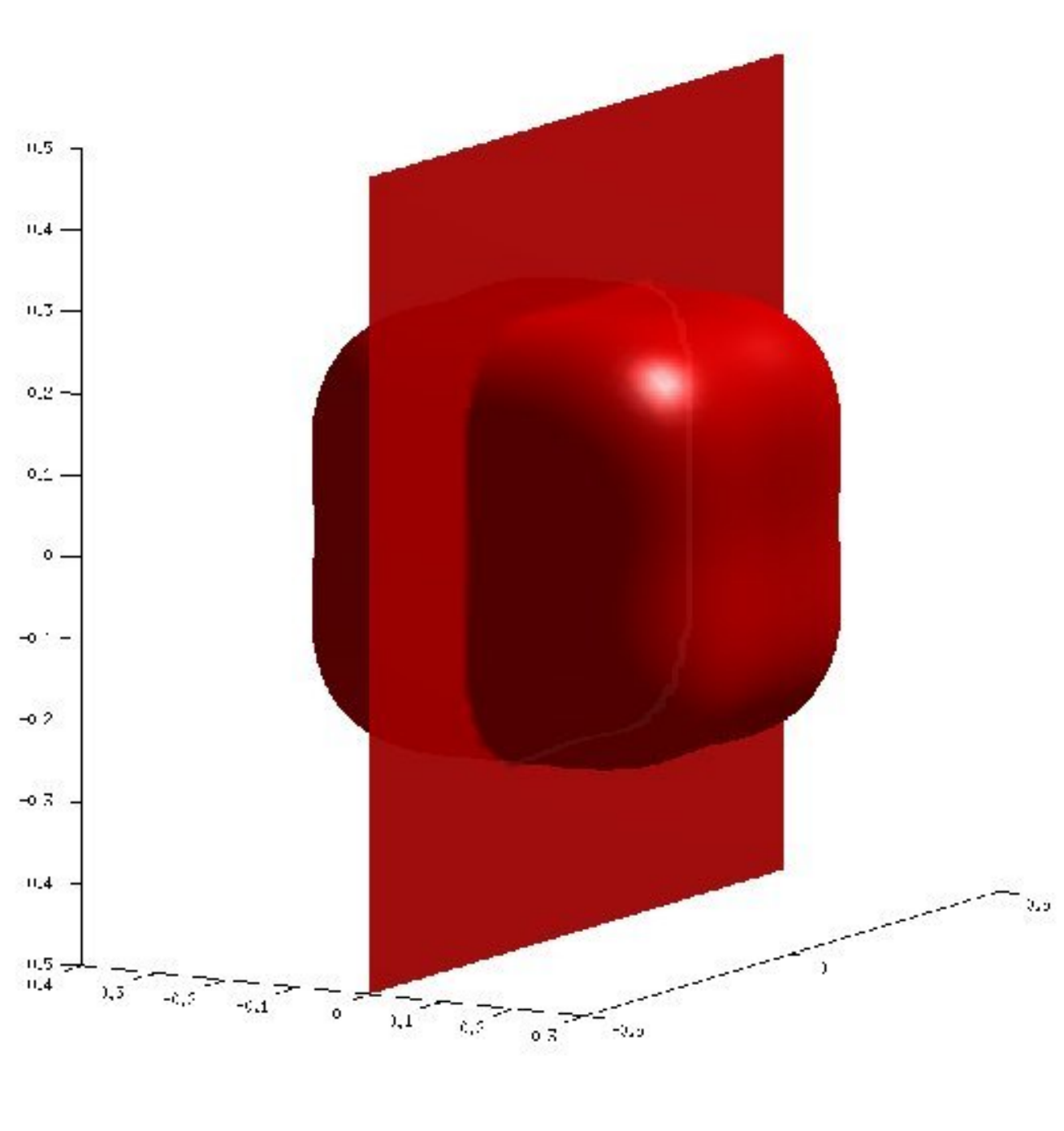}
\includegraphics[width=5cm]{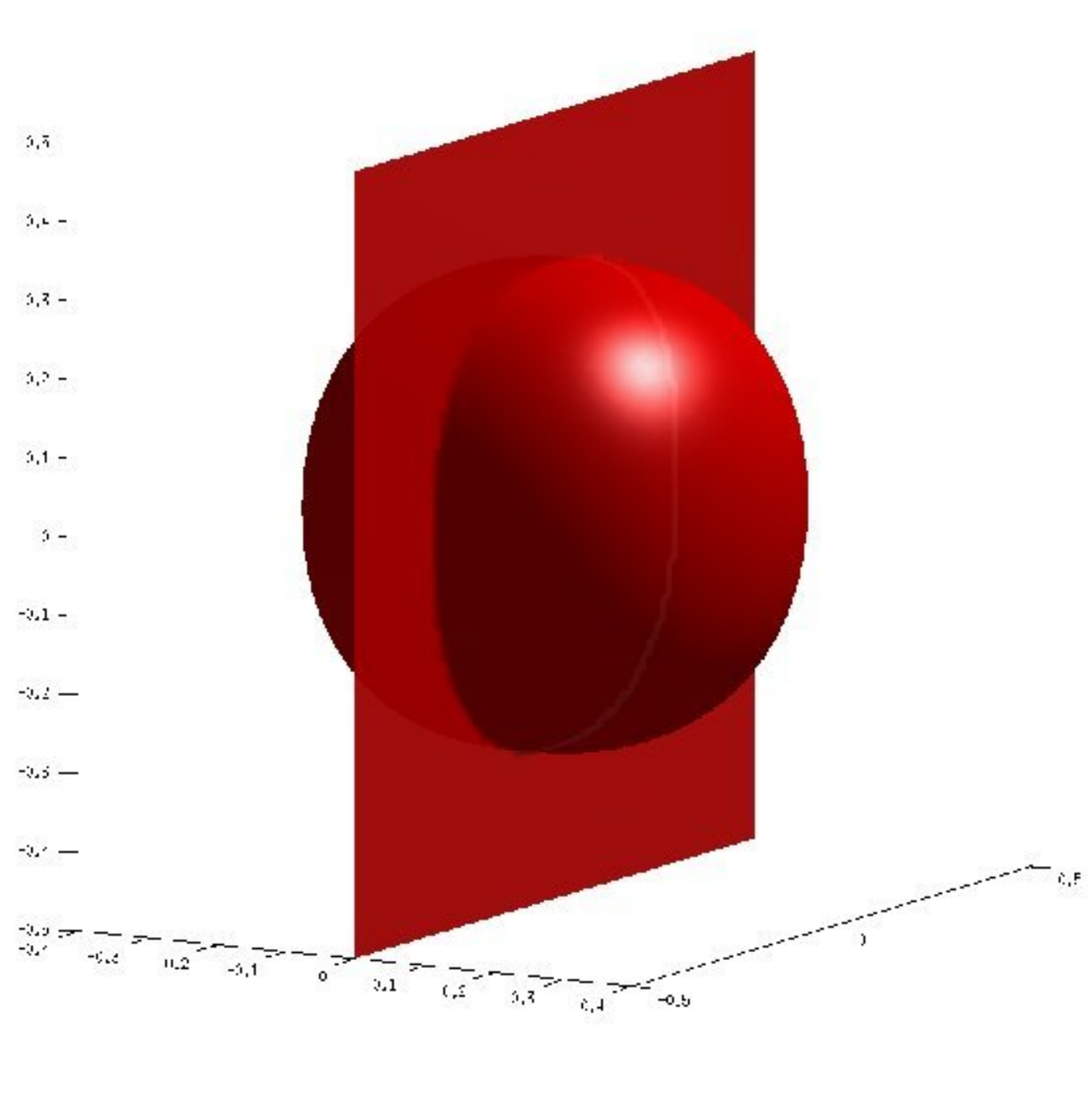}
\caption{3D-examples of evolutions by the classical diffuse flow yielding singularities.}
\label{fig:other_experiments_3D}
\end{figure}

We present additional experiments in Figure~\ref{fig:other_experiments_3D} which illustrate the formation of singularities in dimension $3$. On the first line, two spheres evolve by the classical diffuse flow. As the distance between the two spheres is about $\varepsilon$,
they merge. We take in the second experiment the initial set $\Gamma(0)$ as the union of two parallel cylinders. 
The two cylinders grow up until collision time, at which a crossing arises.
The last  example shows the evolution of a cube cut by a plane (more precisely, both the plane and the cube's boundary separate the two phases, as in the 2D situation of Figure~\ref{fig:other_experiment_2D}). The cube seems to evolve to a sphere 
without being disturbed by the presence of the plane. All these experiments show that the classical diffuse flow may yield singularities, although the comprehension of singular solutions to the Allen-Cahn equation in dimension $3$ remains uncomplete.

\paragraph{Conclusion}
In view of the above simulations, the following observations can be made on the classical diffuse approximation flow:
\begin{itemize}
 \item It is possible to simulate the crossings of more than  two interfaces;
 \item The evolution, by the classical diffuse flow, of two interfaces after crossing seems to be similar to the evolution by a smooth parametric approach. This is in favor of a varifold interpretation of the Willmore flow.
\end{itemize}

\subsection{Numerical simulations of Mugnai's flow}

We now consider the PDE system associated with Mugnai's flow in the form that we introduced in Section~\ref{sec:mugflow}: 
$$
 \begin{cases}
      \partial_t u = \Delta \mu -  \frac{1}{\varepsilon^2} W''(u) \mu  + \widetilde{\cal B}_{\sigma}(u) \\
     \mu =  \frac{1}{\varepsilon^2} W'(u) - \Delta u. \\
    \end{cases} 
 $$
where $$\widetilde{\cal B}_{\sigma}(u) = W'(u) \left[ \left( \left| \nabla \nu_{u,\sigma} \right|^2 -  \left| \mdiv \nu_{u,\sigma} \right|^2 \right) - \mrot\left(\mrot\left(  \nu_{u,\sigma}  \right) \right)\cdot \nu_{u,\sigma}\right]$$
with $\nu_{u,\sigma} =  \frac{\nabla u}{ \sqrt{|\nabla u|^2 + \sigma^2}}$. The initial conditions $u(x,0)$ and $\mu(x,0)$ have the form
$$ \begin{cases}
    u(x,0) &= \gamma \left( \frac{d(\Gamma_0)}{\varepsilon} \right) \\
    \mu(x,0) &= - \frac{1}{\varepsilon}\Delta d(\Gamma_0)  \gamma' \left( \frac{d(\Gamma_0)}{\varepsilon} \right) 
   \end{cases}
$$
We set the approximation parameter $\sigma = 10^{-3}$ and we solve numerically the system using Algorithm~2.

\begin{figure}[!ht] 
\centering
\includegraphics[width=6cm]{./test_cercle}
\includegraphics[width=6cm]{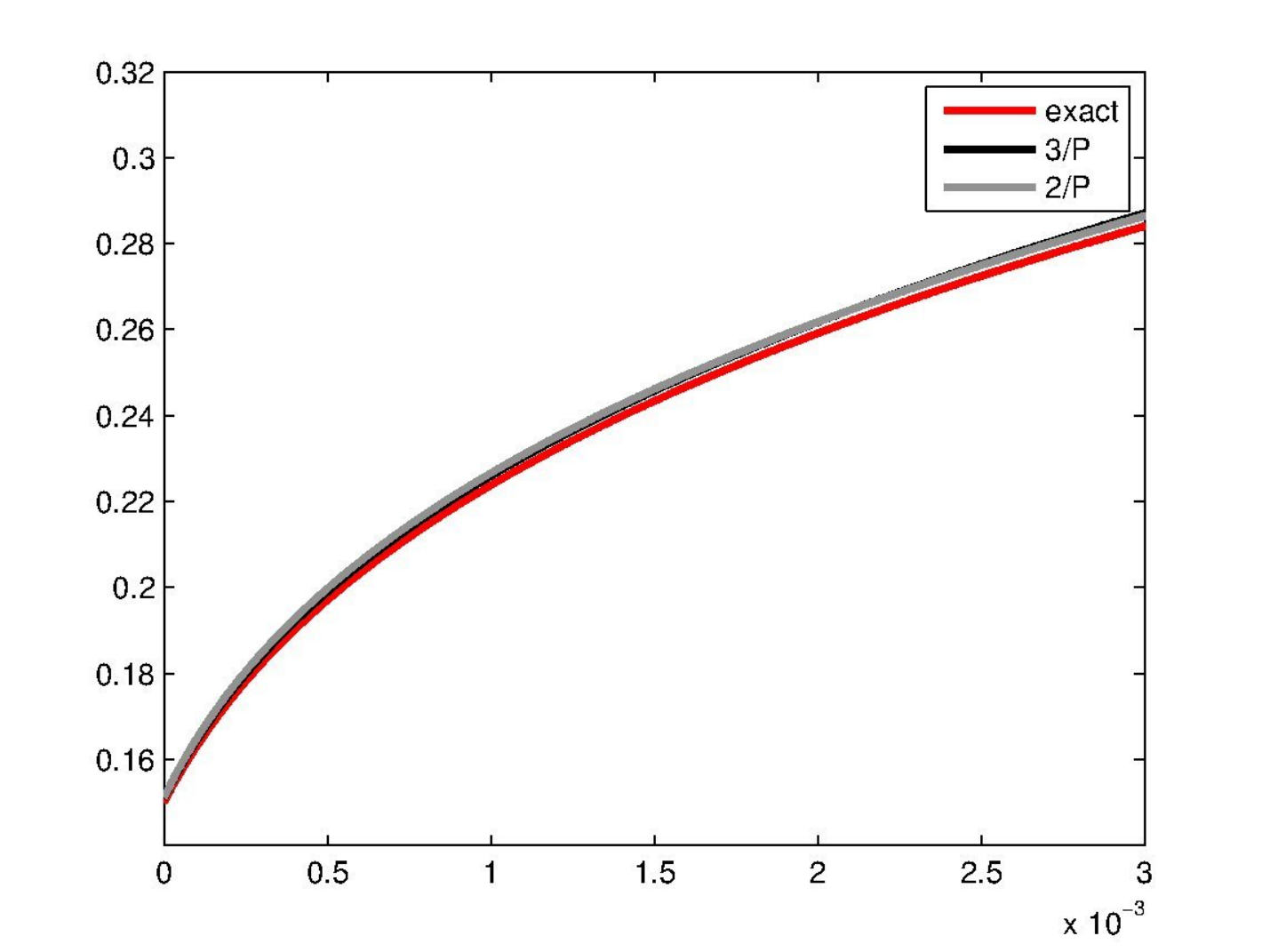}
\caption{Left: evolution of a circle by Mugnai's flow; Right : Graphs of $t\mapsto R_{\varepsilon}(t)$ for different values of $\varepsilon$}
\label{fig:error_circle_mugnai}
\end{figure}

\paragraph{Convergence of Mugnai's approximation}
The first example illustrated in Figure~\ref{fig:error_circle_mugnai} shows the evolution of a circle taken as initial set $\Gamma_0$, and the comparison with the exact solution. The numerical parameters are ${\cal P} = 2^7$, $\varepsilon = 2/{\cal P}$ or $3//{\cal P}$, and  $\delta_t = 1/2 \varepsilon^2 1/{\cal P}^2$. The smaller is $\varepsilon$, the closer the numerical flow is with respect to the continuous flow. This may indicate that the penalization term $\widetilde{\cal B}_{\sigma}(u)$ does not influence the evolution of smooth interfaces. 

\begin{figure}[!ht] 
\centering
\includegraphics[width=5cm]{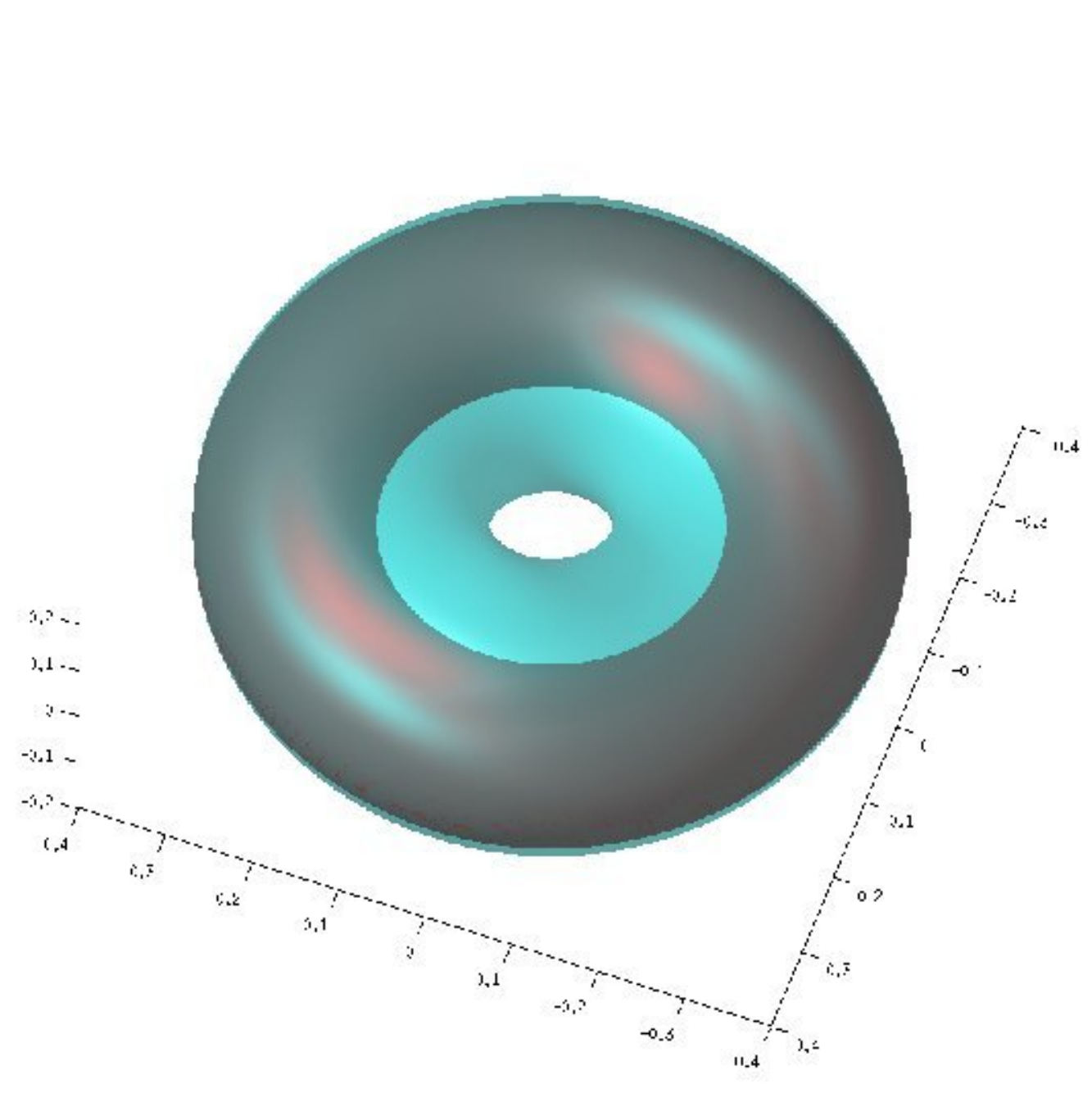}
\includegraphics[width=5cm]{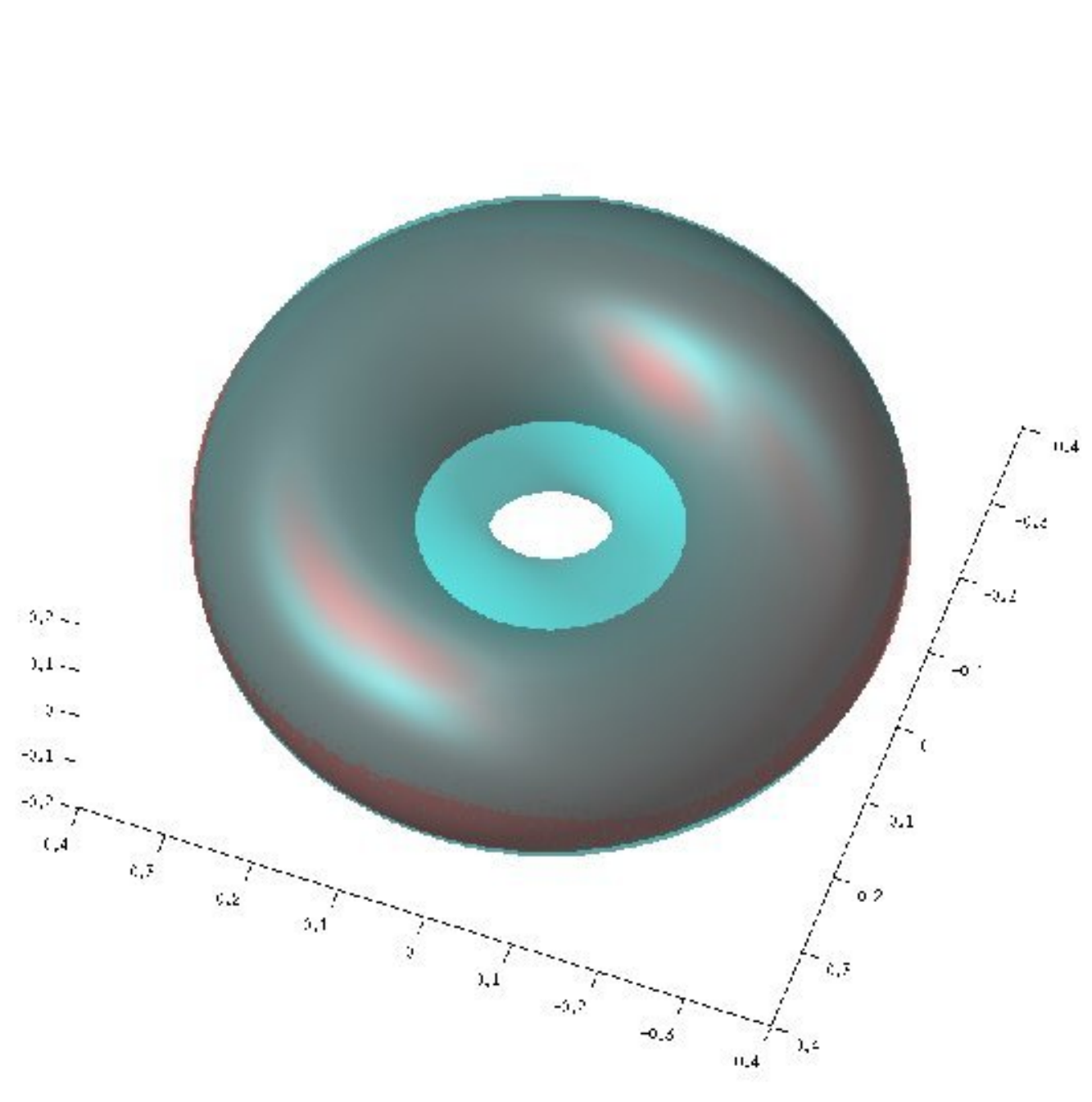}
\includegraphics[width=5cm]{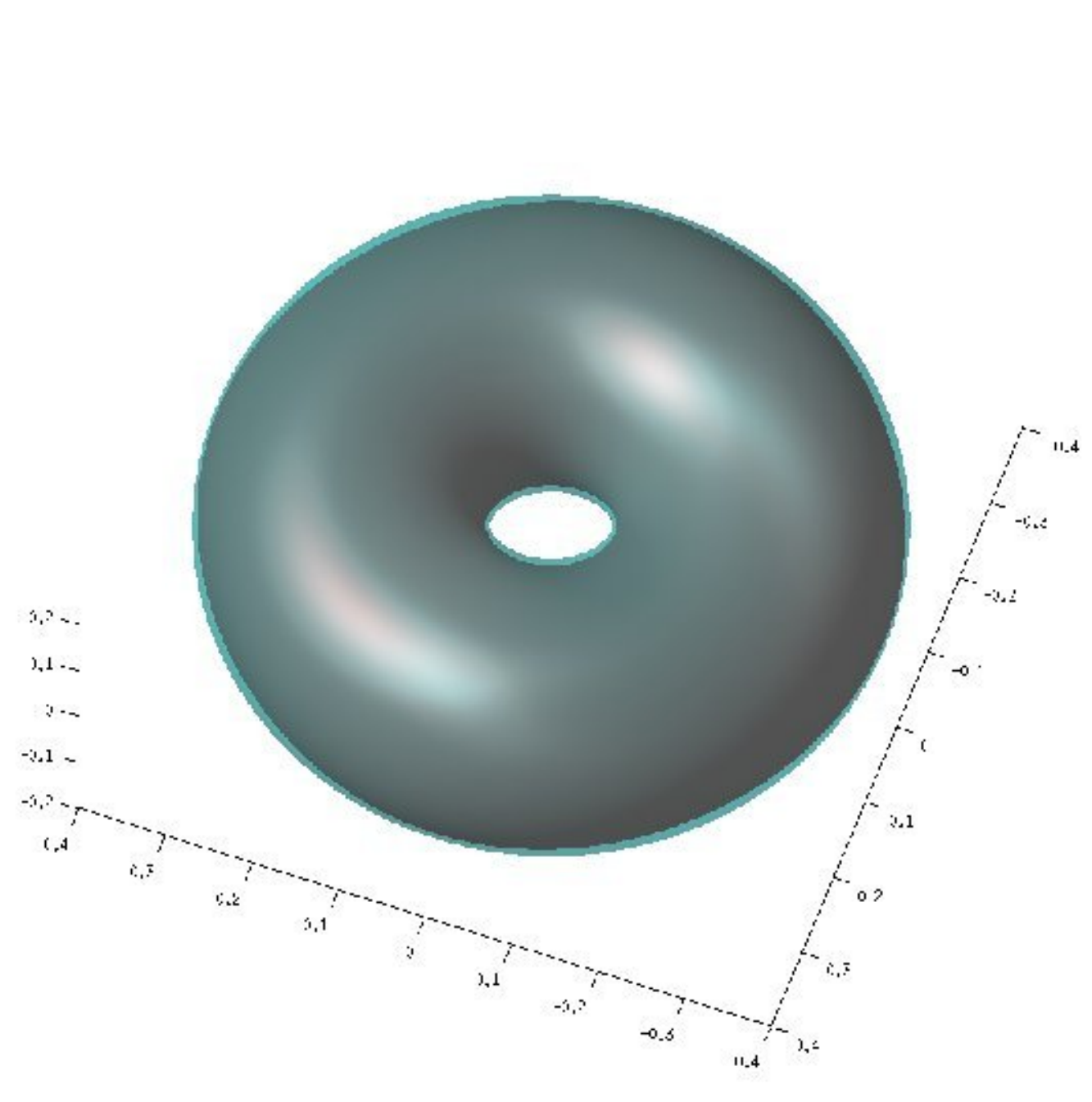} \\
\caption{Smooth evolution by Mugnai's flow of a torus in 3D. The blue torus is the target Clifford torus.}
\label{fig:tore_mugnai}
\end{figure}

This is also illustrated in Figure~\ref{fig:tore_mugnai} where an initial torus in 3D evolves to the Clifford torus.

\begin{figure}[!ht] 
\centering
\includegraphics[width=5cm]{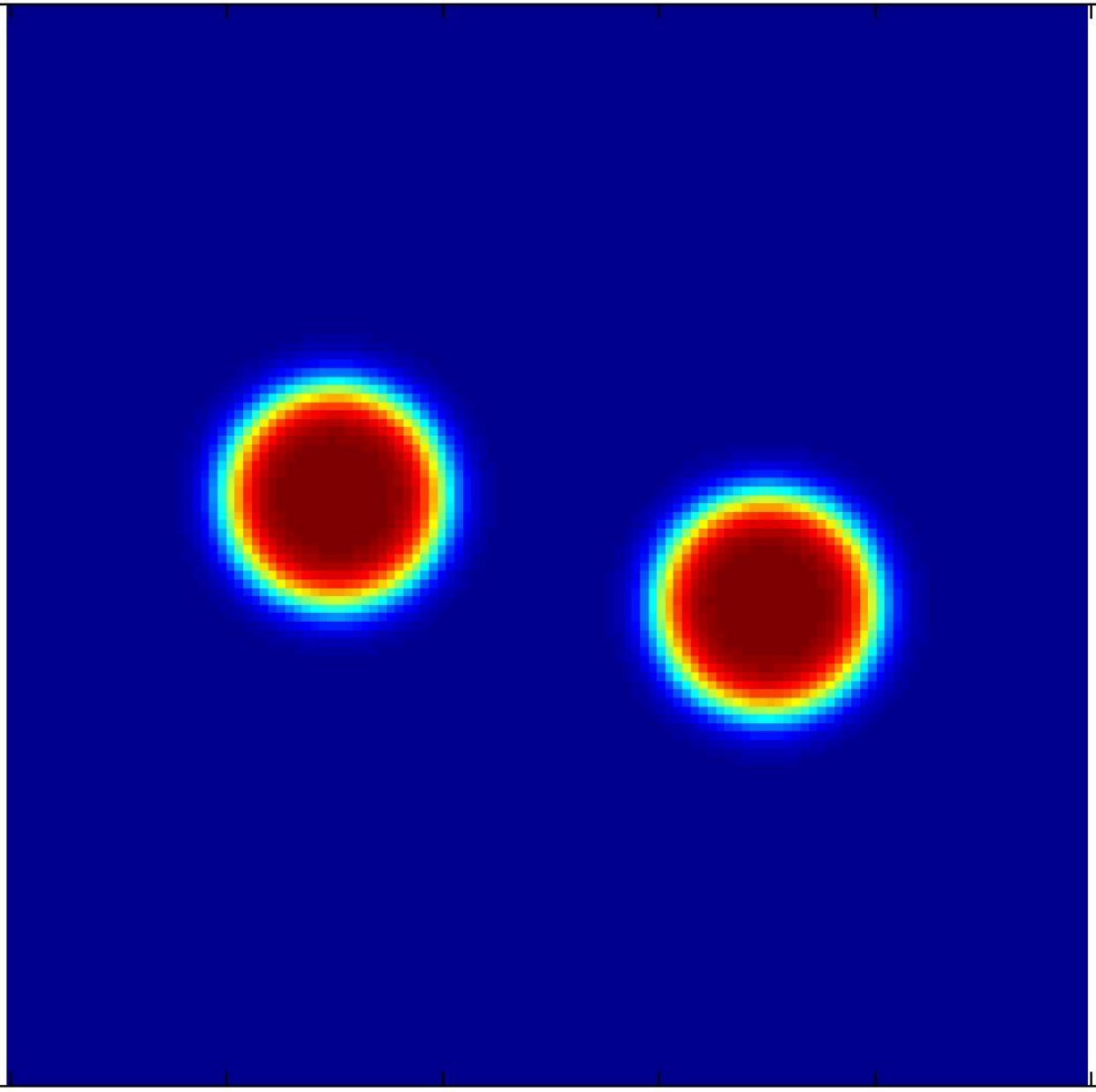}
\includegraphics[width=5cm]{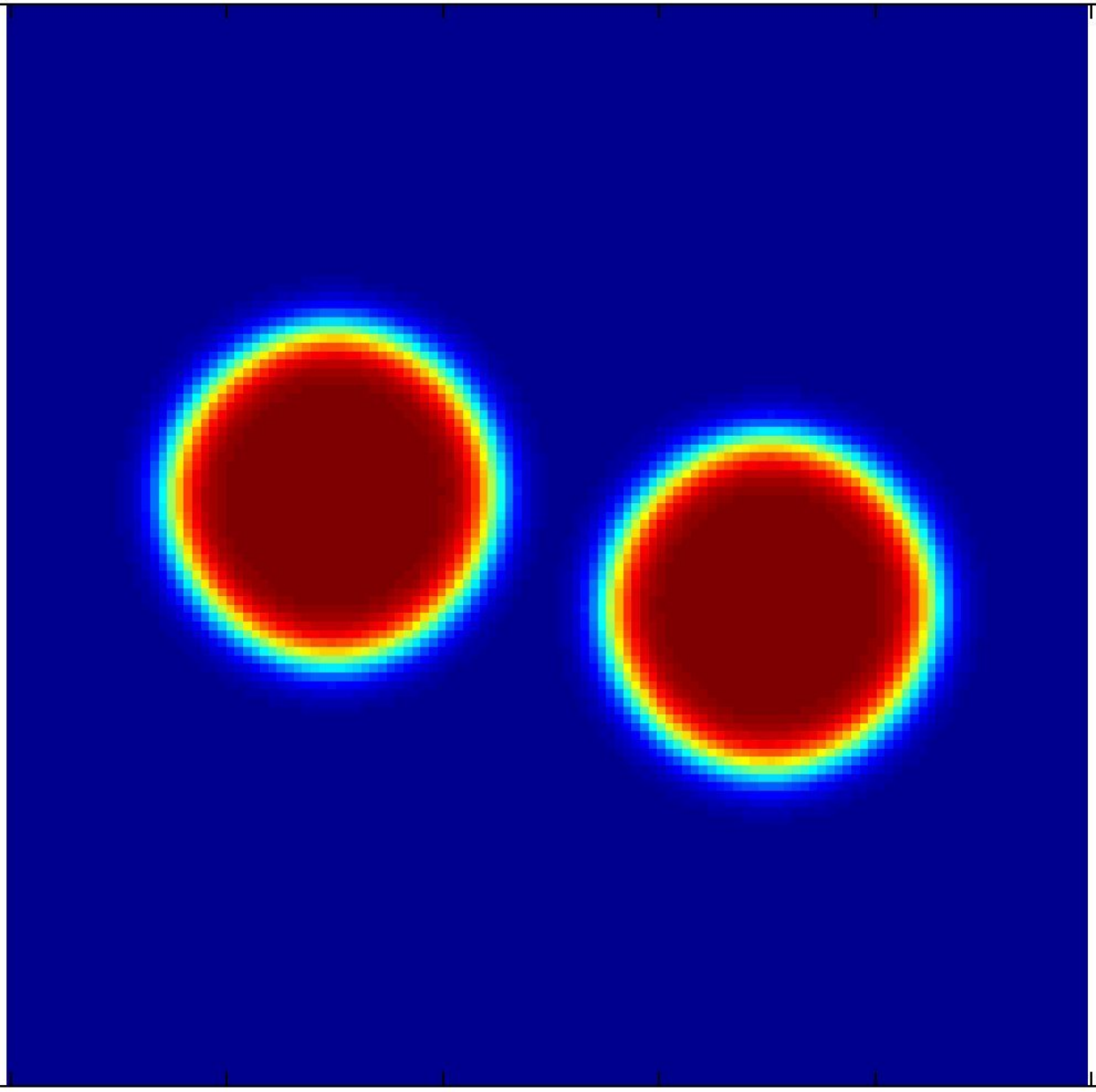}
\includegraphics[width=5cm]{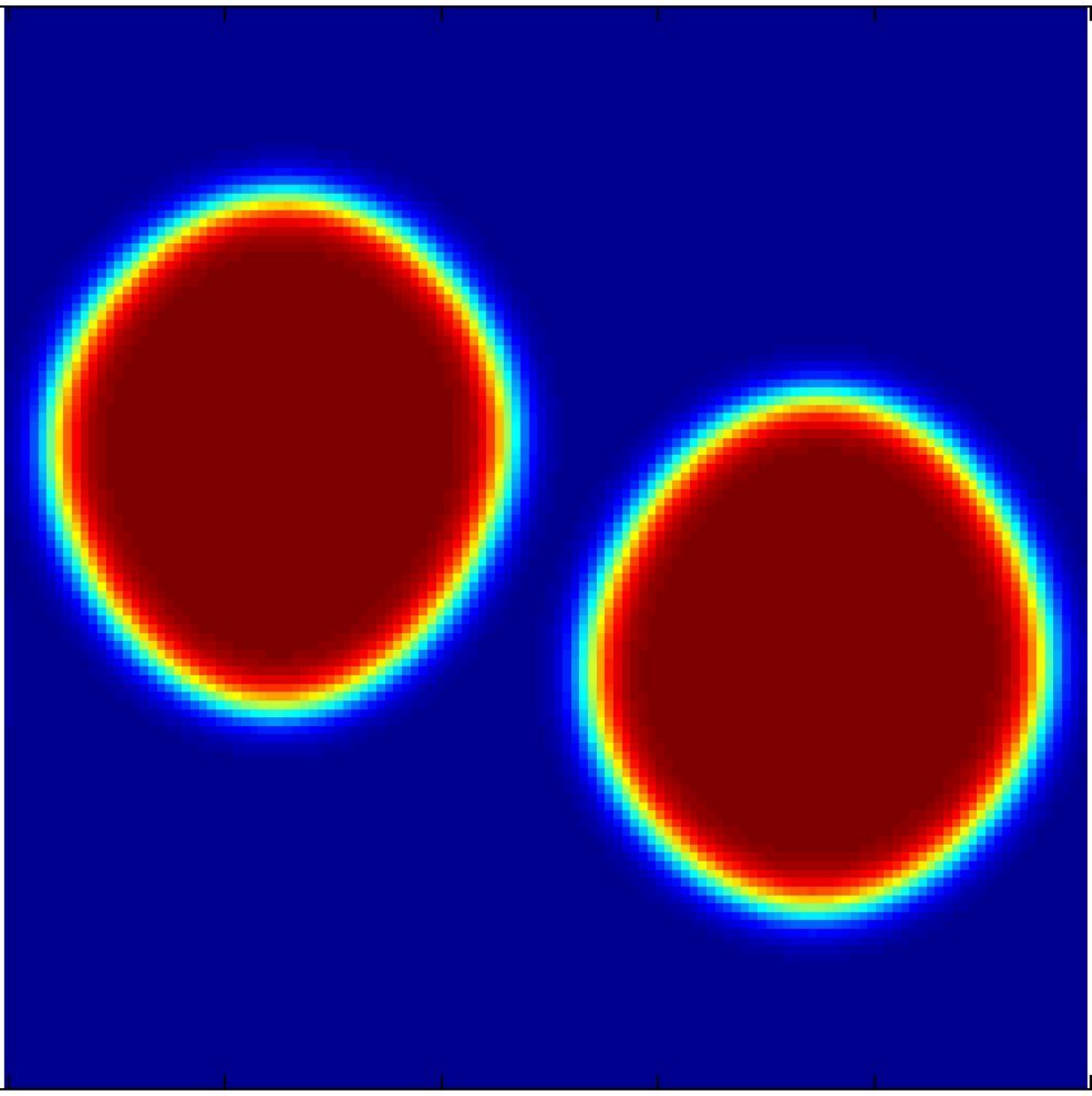} \\
\includegraphics[width=5cm]{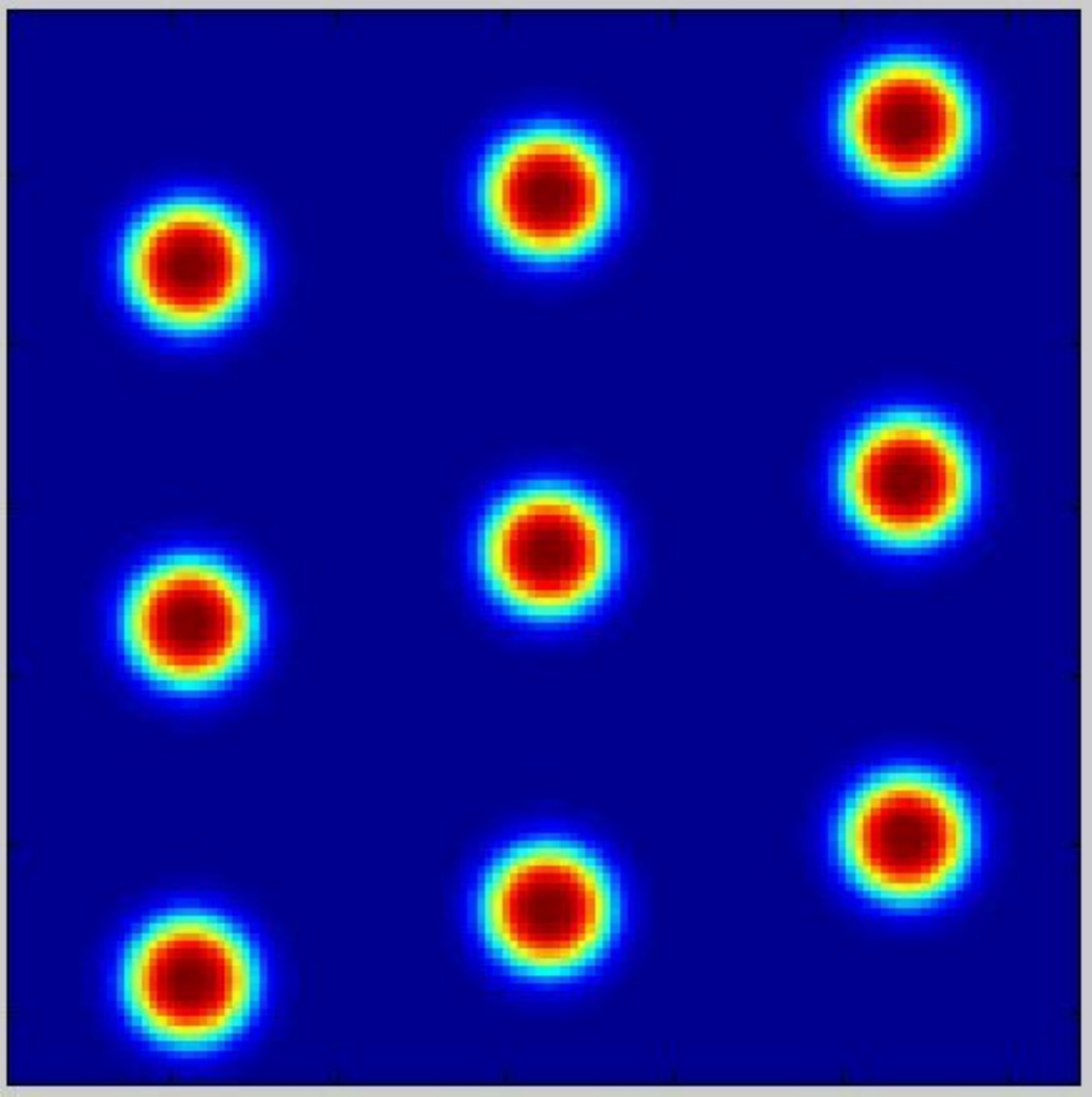}
\includegraphics[width=5cm]{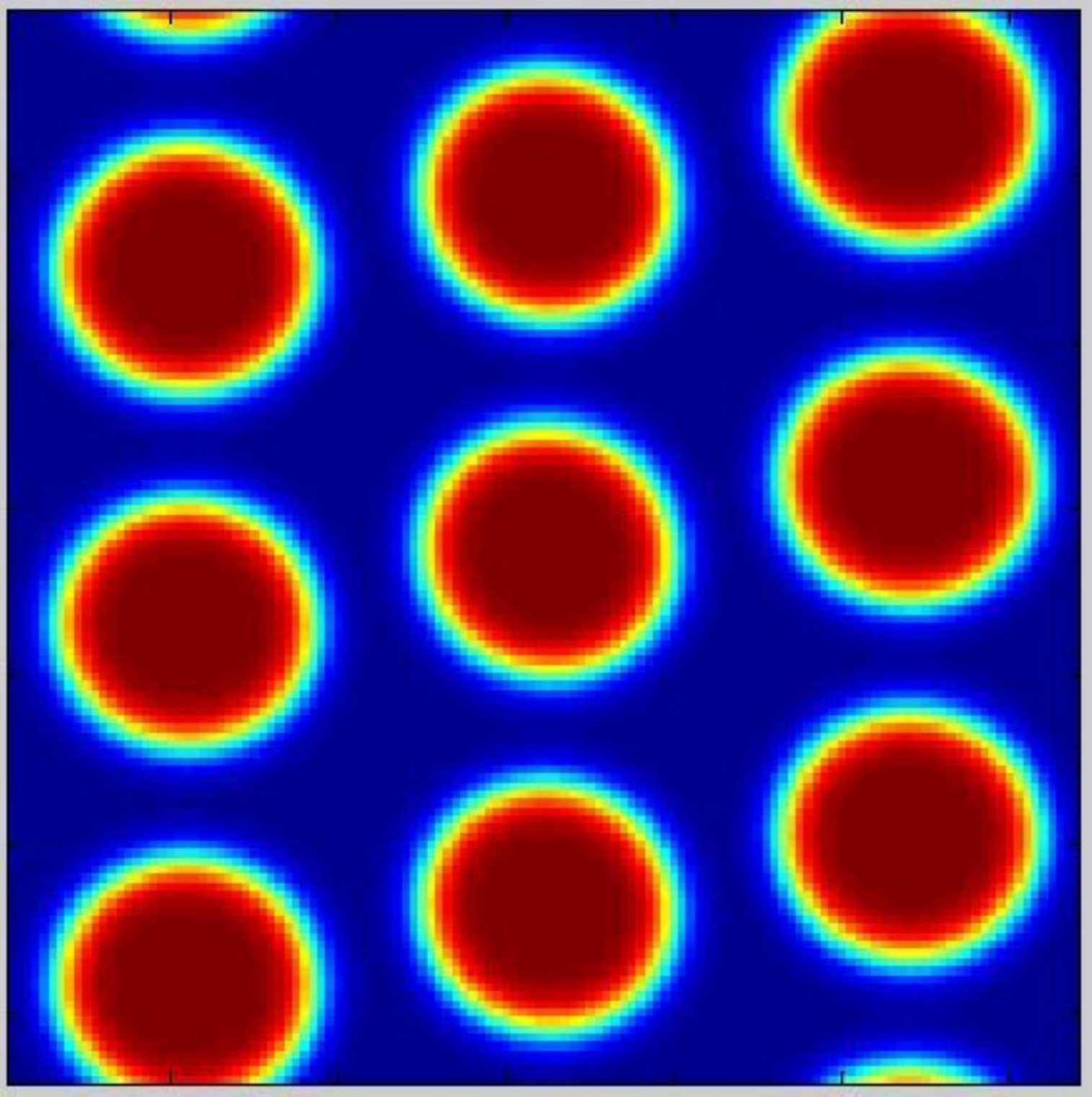}
\includegraphics[width=5cm]{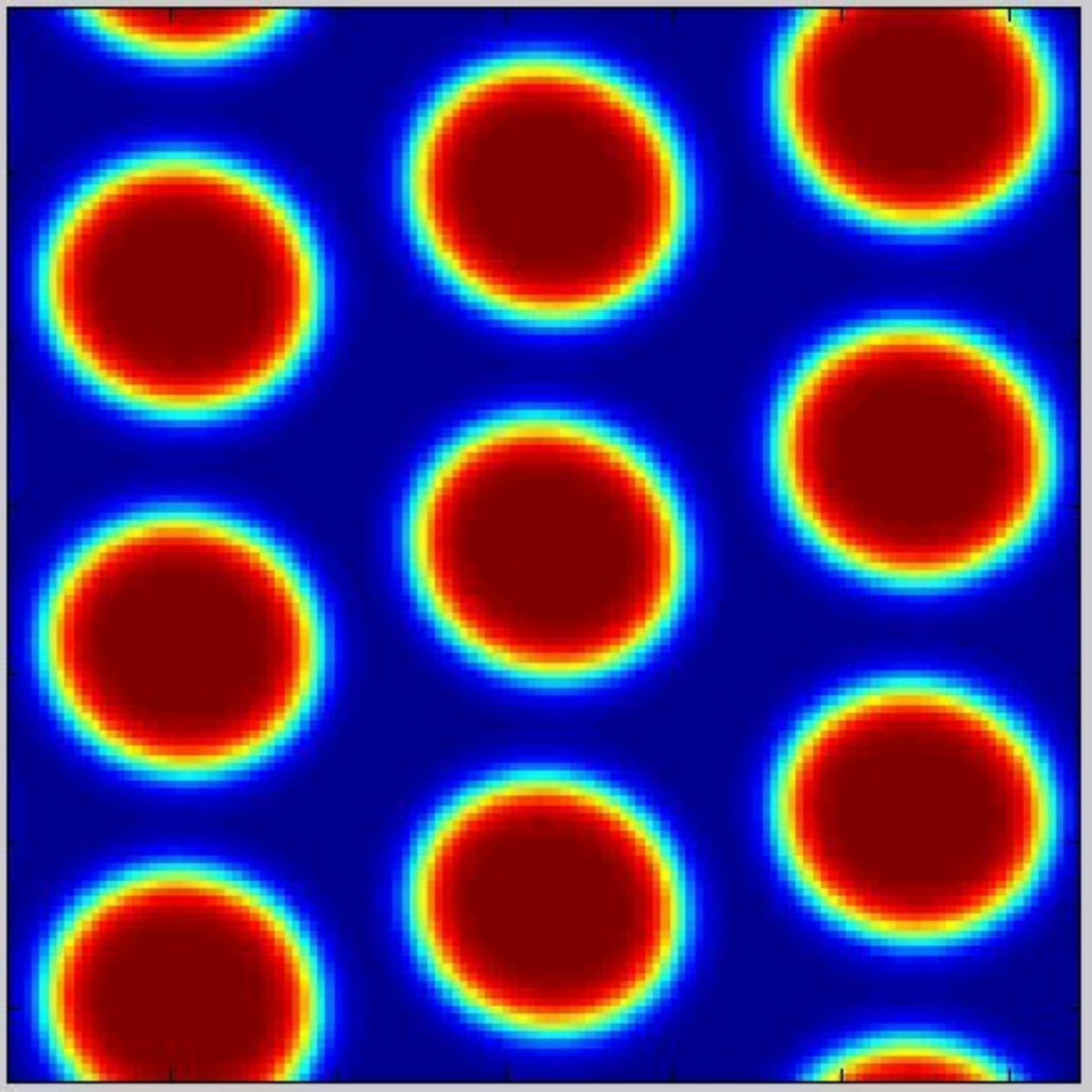} \\
\caption{Illustrations in 2D that Mugnai's flow prevents from colliding.}
\label{fig:Mugnai_other_experiment_2D}
\end{figure}

We present two experiments in Figure~\ref{fig:Mugnai_other_experiment_2D} obtained with the set of parameters ${\cal P} = 2^7$, $\varepsilon = 2/{\cal P}$ 
and $\delta_t = 1/8 \varepsilon^2 {\cal P}^{-2}$. The simulations indicate that the additional penalization term $\widetilde{\cal B}_{\sigma}(u)$ prevents the interfaces from colliding. This is coherent with what we argued in Section~\ref{sec:remarks}, i.e. that Mugnai's energy equals the classical energy plus a functional that penalizes non profile functions.

\begin{figure}[!ht] 
\centering
\includegraphics[width=5cm]{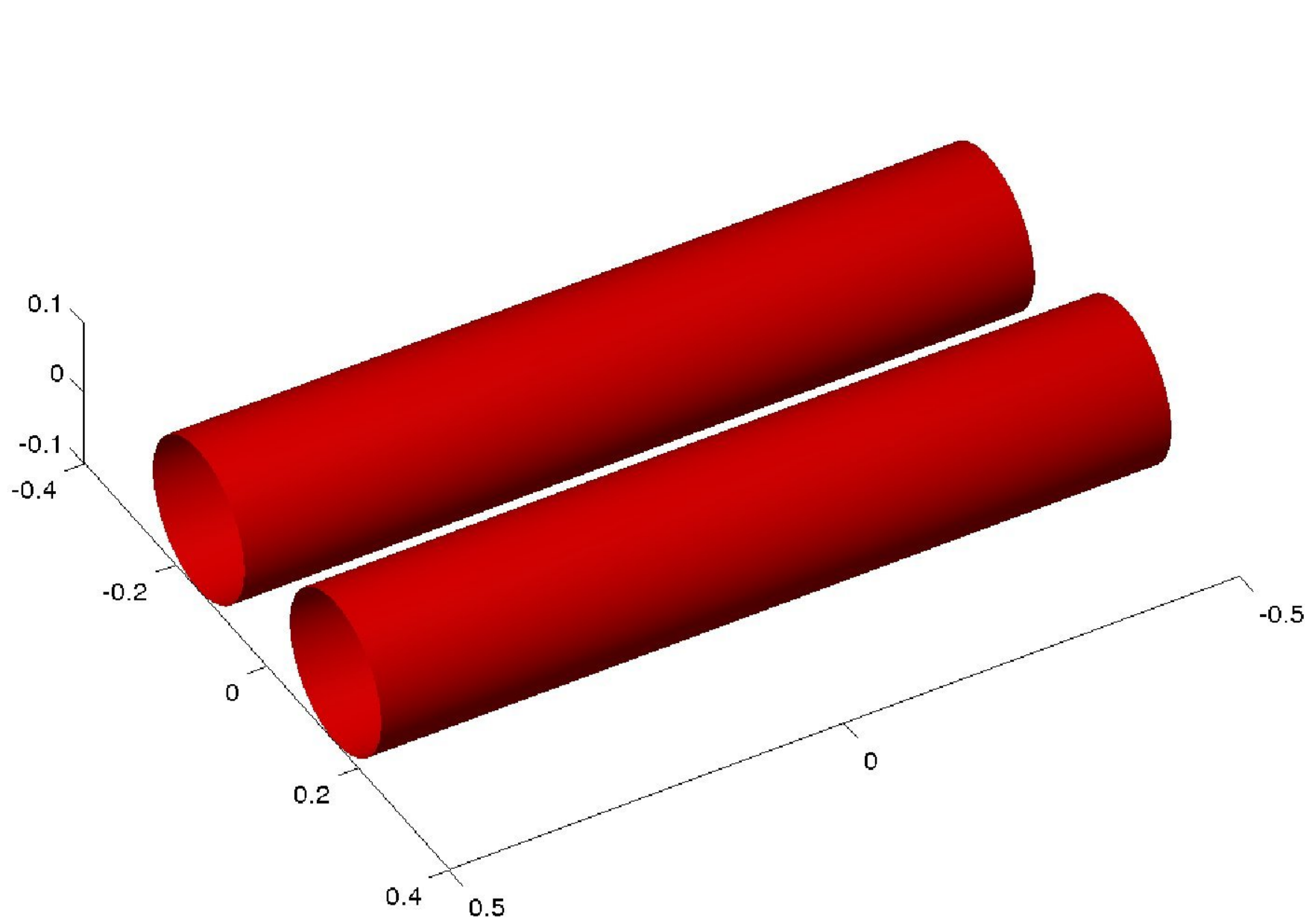}
\includegraphics[width=5cm]{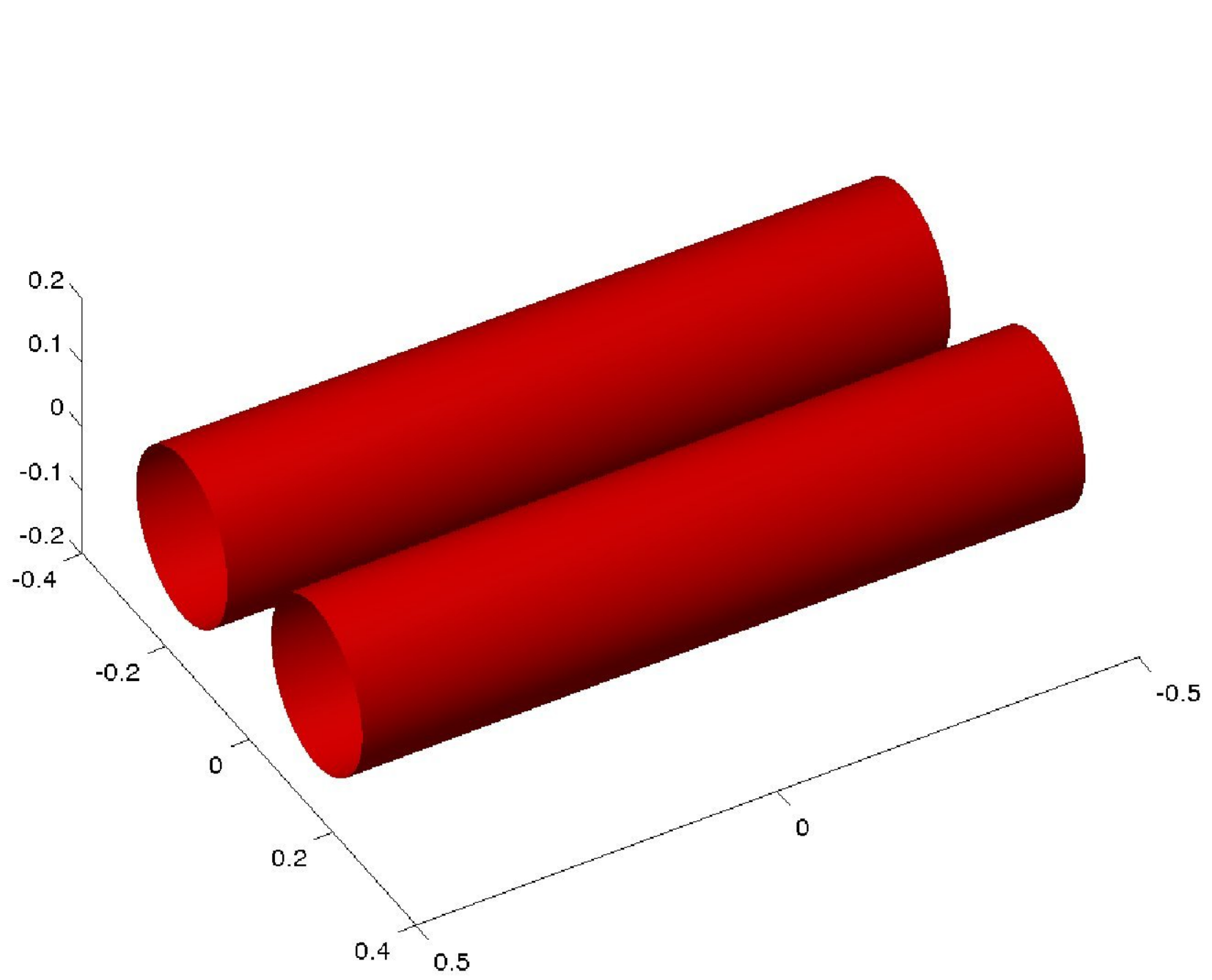}
\includegraphics[width=5cm]{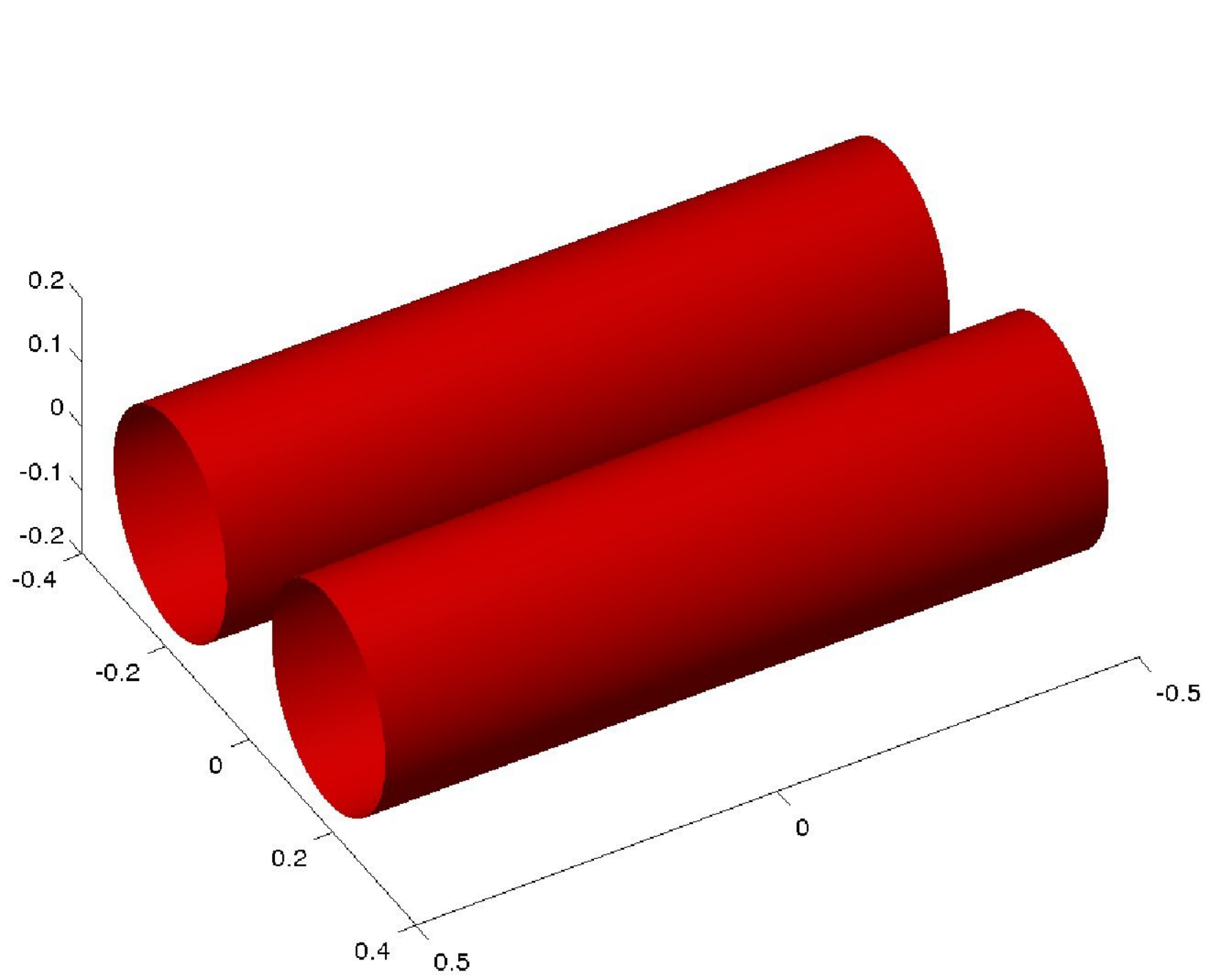} \\
\includegraphics[width=5cm]{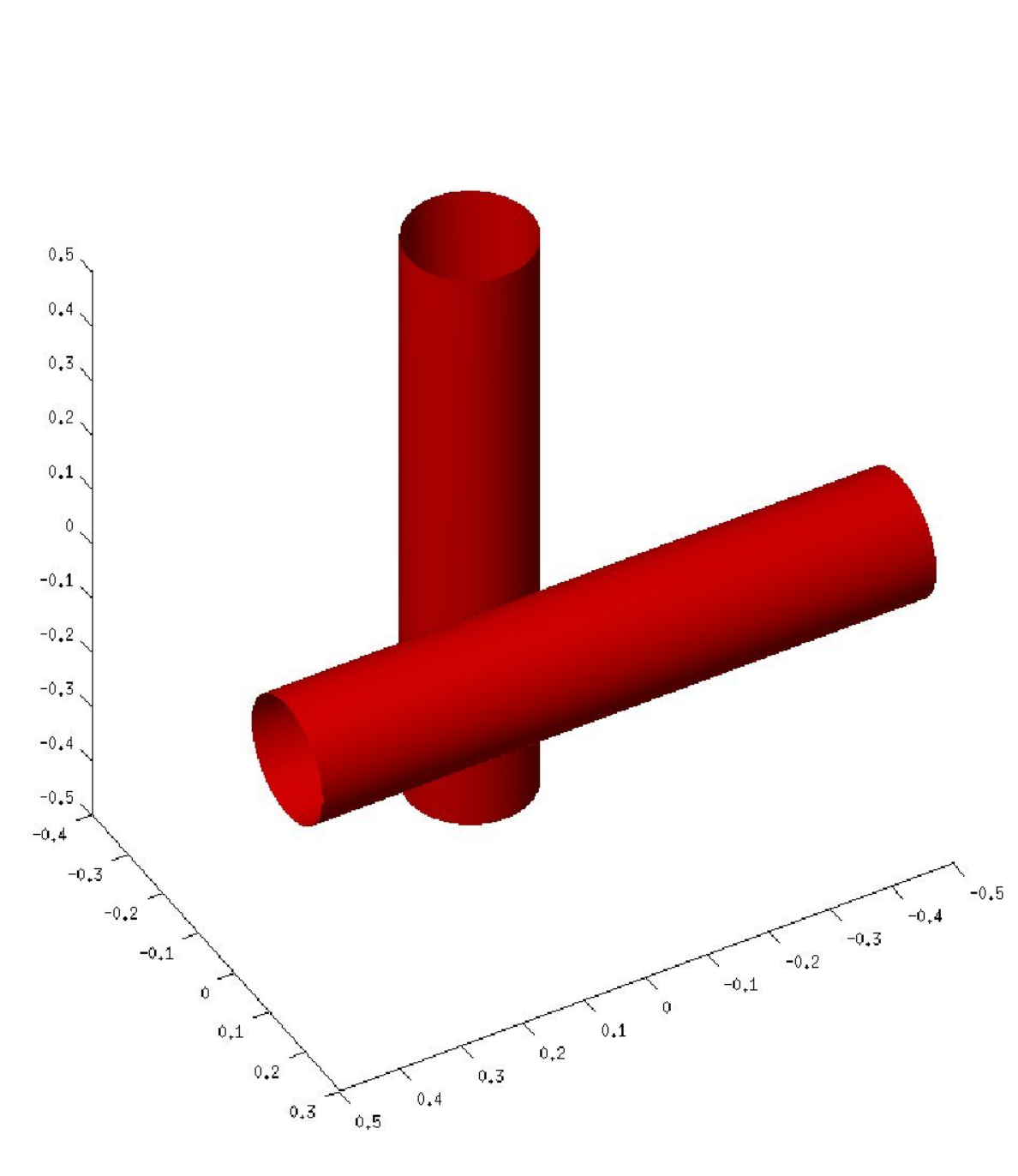}
\includegraphics[width=5cm]{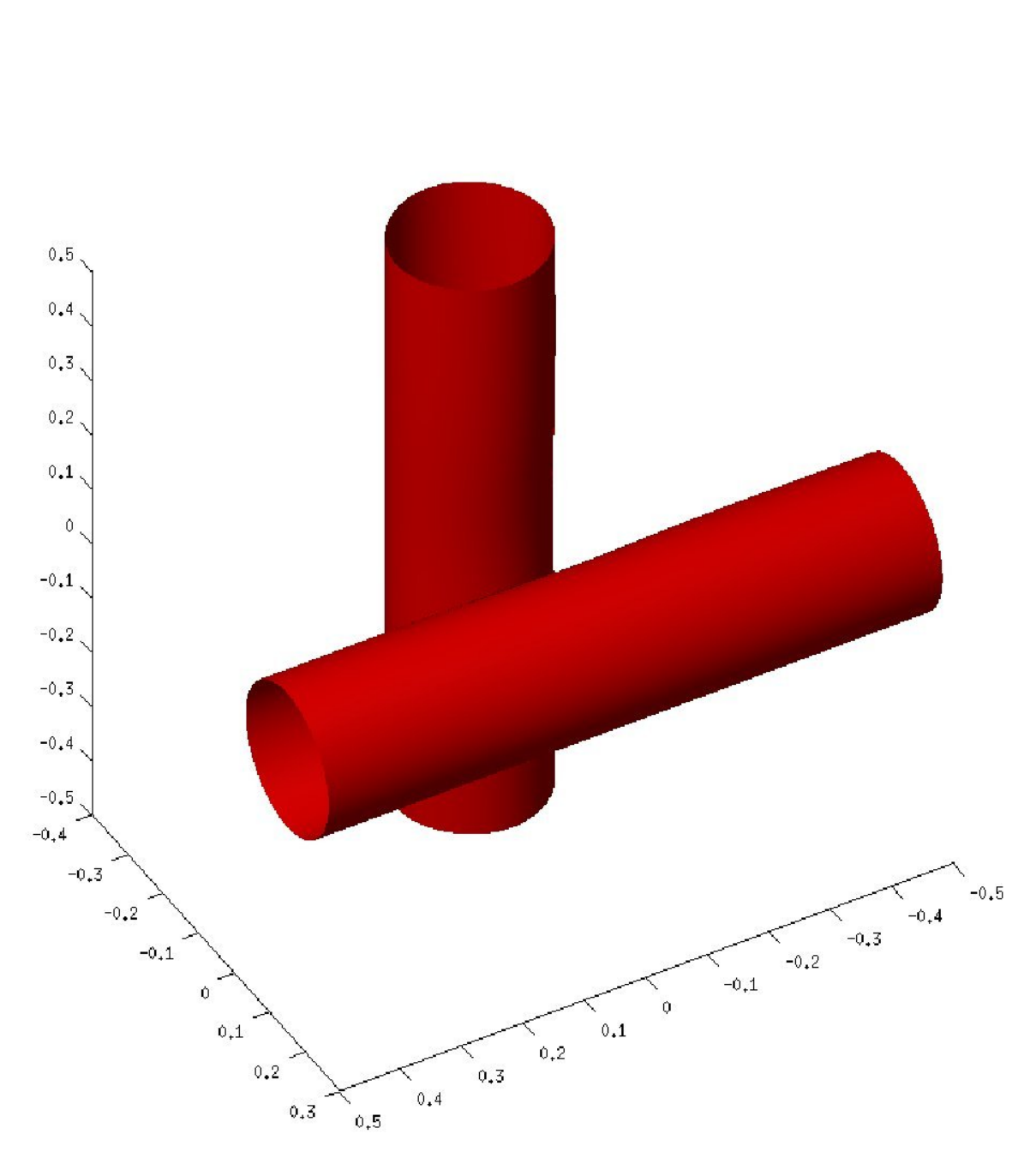}
\includegraphics[width=5cm]{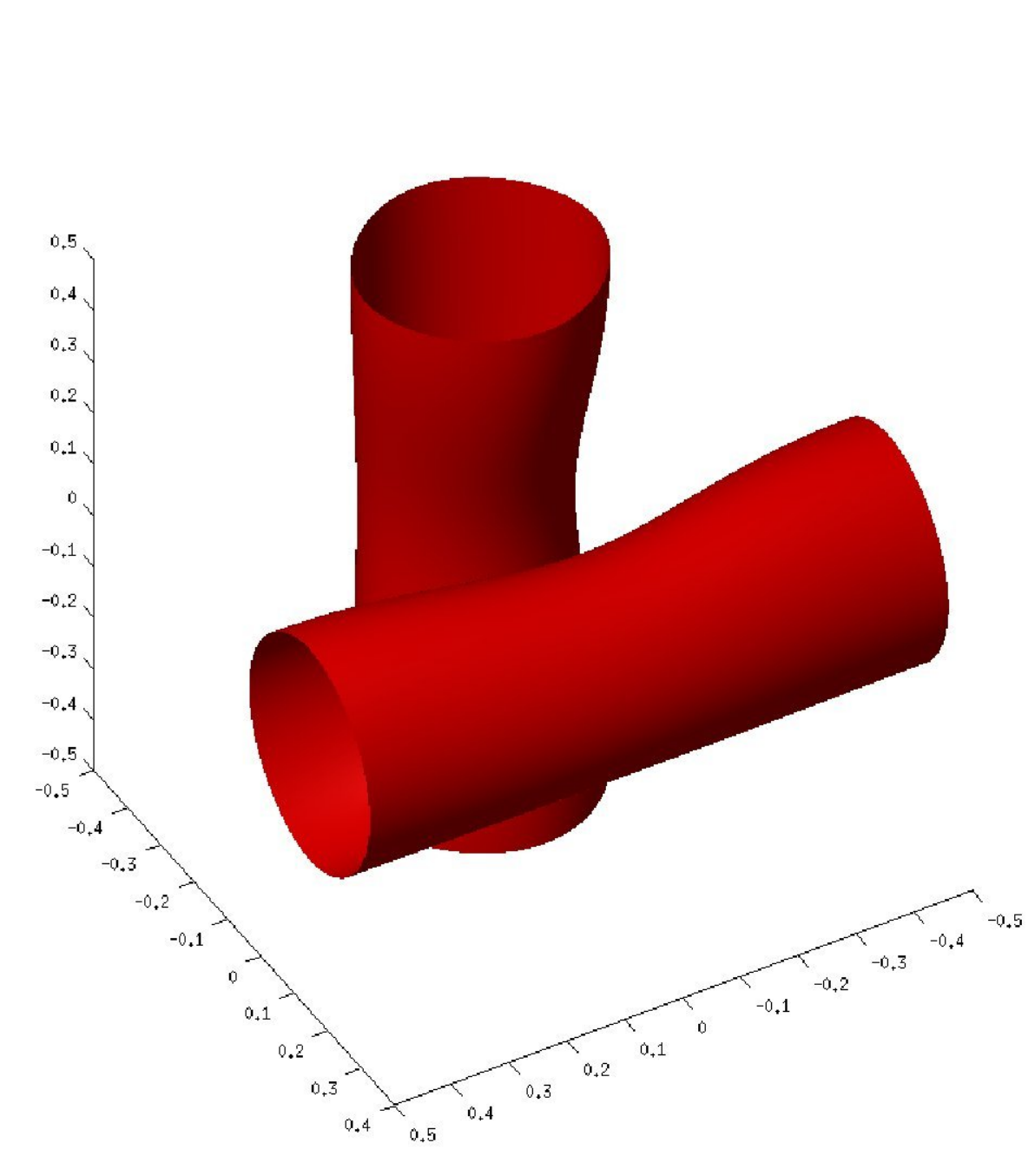} 
\caption{Illustrations in 3D that Mugnai's flow prevents from colliding. The interfaces preferably deform themselves rather than merging.}
\label{fig:Mugnai_other_experiment_3D}
\end{figure}

The same observation is illustrated in 3D on Figure~\ref{fig:Mugnai_other_experiment_3D}. Both cylinders grow up, but deform themselves rather than colliding.
  
\paragraph{Conclusion}
To conclude this experimental section on Mugnai's flow, let us observe that
\begin{itemize}
\item As long as the interfaces are smooth, Mugnai's and the classical flow behave in the same way, which was of course expected from the theoretical properties of the associated functionals. In particular, the penalization term $\widetilde{\cal B}_{\sigma}(u)$ has no critical influence on the evolution of a smooth interface, as long as the evolution remains smooth as well with the classical flow.
\item Since Mugnai's energy $\cW^{\mbox{\tiny Mu}}_\varepsilon$ $\Gamma$-converges in dimension $2$ to the relaxation of the Willmore energy, the associated flow prevents from crossing, which is confirmed by the simulations. In 3D as well, our simulations indicate that no crossing should occur. This indicates that the $\Gamma$-convergence property should also be true in 3D for Mugnai's energy, which is so far an open question that requires a better understanding of the diffuse approximation of the genus (having in mind the Gauss-Bonnet Theorem).
\end{itemize}
 
\subsection*{Acknowledgements}
The authors thank Luca Mugnai, Selim Esedoglu, Petru Mironescu, and Giovanni Bellettini for fruitful discussions.

\bibliographystyle{abbrv}
\bibliography{biblio}

\begin{thebibliography}{10}

\bibitem{Ambrosio2000}
L.~Ambrosio.
\newblock Geometric evolution problems, distance function and viscosity
  solutions.
\newblock In {\em Calculus of variations and partial differential equations
  ({P}isa, 1996)}, pages 5--93. Springer, Berlin, 2000.

\bibitem{amb-mant-96}
L.~Ambrosio and C.~Mantegazza.
\newblock Curvature and distance function from a manifold.
\newblock {\em J. Geom. Anal.}, 8(5):723--748, 1998.
\newblock Dedicated to the memory of Fred Almgren.

\bibitem{AM}
L.~Ambrosio and S.~Masnou.
\newblock A direct variational approach to a problem arising in image
  reconstruction.
\newblock {\em Interfaces and Free Boundaries}, 5:63--81, 2003.

\bibitem{Barrett:2007:PFE:1225314.1225638}
J.~W. Barrett, H.~Garcke, and R.~N\"{u}rnberg.
\newblock A parametric finite element method for fourth order geometric
  evolution equations.
\newblock {\em J. Comput. Phys.}, 222(1):441--467, Mar. 2007.

\bibitem{Barrett2008_2}
J.~W. Barrett, H.~Garcke, and R.~N\"{u}rnberg.
\newblock On the parametric finite element approximation of evolving
  hypersurfaces in {$\mathbb R^3$}.
\newblock {\em J. Comput. Phys.}, 227:4281--4307, April 2008.

\bibitem{BarrettGN08}
J.~W. Barrett, H.~Garcke, and R.~N\"{u}rnberg.
\newblock Parametric approximation of {W}illmore flow and related geometric
  evolution equations.
\newblock {\em SIAM J. Scientific Computing}, 31(1):225--253, 2008.

\bibitem{Barrett:2008}
J.~W. Barrett, H.~Garcke, and R.~N\"{u}rnberg.
\newblock A variational formulation of anisotropic geometric evolution
  equations in higher dimensions.
\newblock {\em Numer. Math.}, 109:1--44, February 2008.

\bibitem{Bellettini1997}
G.~Bellettini.
\newblock Variational approximation of functionals with curvatures and related
  properties.
\newblock {\em J. Convex Anal.}, 4(1):91--108, 1997.

\bibitem{BDP}
G.~Bellettini, G.~D. Maso, and M.~Paolini.
\newblock Semicontinuity and relaxation properties of curvature depending
  functional in 2{D}.
\newblock {\em Annali della Scuola Normale di Pisa, Classe di Scienze, $4^e$
  s\'erie}, 20(2):247--297, 1993.

\bibitem{BM1}
G.~Bellettini and L.~Mugnai.
\newblock Characterization and representation of the lower semicontinuous
  envelope of the elastica functional.
\newblock {\em Ann. Inst. H. Poincar\'e, Anal. non Lin\'eaire}, 21(6):839--880,
  2004.

\bibitem{belmug04}
G.~Bellettini and L.~Mugnai.
\newblock On the approximation of the elastica functional in radial symmetry.
\newblock {\em Calc. Var. Partial Differential Equations}, 24:1--20, 2005.

\bibitem{BM}
G.~Bellettini and L.~Mugnai.
\newblock A varifold representation of the relaxed elastica functional.
\newblock {\em Journal of Convex Analysis}, 14(3):543--564, 2007.

\bibitem{belmu10}
G.~Bellettini and L.~Mugnai.
\newblock Approximation of {H}elfrich's functional via diffuse interfaces.
\newblock {\em SIAM J. Math. Anal.}, 42(6):2402--2433, 2010.

\bibitem{belpao93}
G.~Bellettini and M.~Paolini.
\newblock Approssimazione variazionale di funzionali con curvatura.
\newblock In {\em Seminario Analisi Matematica, Univ. Bologna}, pages 87--97,
  1993.

\bibitem{belpaoqo}
G.~Bellettini and M.~Paolini.
\newblock Quasi-optimal error estimates for the mean curvature flow with a
  forcing term.
\newblock {\em Differential Integral Equations}, 8(4):735--752, 1995.

\bibitem{bellettini_1996}
G.~Bellettini and M.~Paolini.
\newblock Anisotropic motion by mean curvature in the context of {F}insler
  geometry.
\newblock {\em Hokkaido Math. J.}, 25:537--566, 1996.

\bibitem{BenceMerrimanOsher}
J.~Bence, B.~Merriman, and S.~Osher.
\newblock Diffusion generated motion by mean curvature.
\newblock {\em Computational Crystal Growers Workshop,J. Taylor ed. Selected
  Lectures in Math., Amer. Math. Soc.}, pages 73--83, 1992.

\bibitem{BobenkoSchroeder05}
A.~Bobenko and P.~Schr\"oder.
\newblock {Discrete Willmore Flow}.
\newblock In {\em Eurographics Symposium on Geometry Processing}, 2005.

\bibitem{cabre-terra}
X.~Cabr{\'e} and J.~Terra.
\newblock Saddle-shaped solutions of bistable diffusion equations in all of
  {$\mathbb R^{2m}$}.
\newblock {\em J. Eur. Math. Soc. (JEMS)}, 11(4):819--843, 2009.

\bibitem{fife}
G.~Caginalp and P.~C. Fife.
\newblock Dynamics of layered interfaces arising from phase boundaries.
\newblock {\em SIAM J. Appl. Math.}, 48(3):506--518, 1988.

\bibitem{Cahn19941045}
J.~Cahn and J.~Taylor.
\newblock Overview no. 113 -- surface motion by surface diffusion.
\newblock {\em Acta Metallurgica et Materialia}, 42(4):1045 -- 1063, 1994.

\bibitem{cahnspin}
J.~W. Cahn.
\newblock {On spinodal decomposition}.
\newblock {\em Acta Metallurgica}, 9(9):795--801, Sept. 1961.

\bibitem{CAHN:1958}
J.~W. Cahn and J.~E. Hilliard.
\newblock Free energy of a nonuniform system. {I}. {I}nterfacial free energy.
\newblock {\em The Journal of Chemical Physics}, 28(2):258--267, 1958.

\bibitem{Chen1992}
X.~Chen.
\newblock Generation and propagation of interfaces for reaction-diffusion
  equations.
\newblock {\em J. Differential Equations}, 96(1):116--141, 1992.

\bibitem{chen-96}
X.~Chen.
\newblock Global asymptotic limit of solutions of the {C}ahn-{H}illiard
  equation.
\newblock {\em J. Differ. Geom.}, 44:262--311, 1996.

\bibitem{chen_giga_goto}
Y.~G. Chen, Y.~Giga, and S.~Goto.
\newblock Uniqueness and existence of viscosity solutions of generalized mean
  curvature flow equations.
\newblock {\em Proc. Japan Acad. Ser. A Math. Sci.}, 65(7):207--210, 1989.

\bibitem{colli1}
P.~Colli and P.~Laurencot.
\newblock A phase-field approximation of the {W}illmore flow with volume
  constraint.
\newblock {\em Interfaces Free Bound}, 13:341--351, 2011.

\bibitem{colli2}
P.~Colli and P.~Laurencot.
\newblock A phase-field approximation of the {W}illmore flow with volume and
  area constraints.
\newblock {\em SIAM J. Math. Anal.}, 44:3734--3754, 2012.

\bibitem{Dang_Fife_Peletier_92}
H.~Dang, P.~C. Fife, and L.~Peletier.
\newblock Saddle solutions of the bistable diffusion equation.
\newblock {\em Zeitschrift fur angewandte Mathematik und Physik ZAMP},
  43(6):984--998, Nov. 1992.

\bibitem{degiorgi-a}
E.~De~Giorgi.
\newblock Some remarks on {$\Gamma$}-convergence and least square methods.
\newblock In G.~D. Maso and G.~Dell'Antonio, editors, {\em Composite Media and
  Homogenization Theory}, pages 135--142. Birkha{\"u}ser, Boston, 1991.

\bibitem{Deckelnick1999}
K.~Deckelnick and G.~Dziuk.
\newblock Discrete anisotropic curvature flow of graphs.
\newblock {\em M2AN Math. Model. Numer. Anal.}, 33:1203--1222, 1999.

\bibitem{Deckelnick2005}
K.~Deckelnick, G.~Dziuk, and C.~M. Elliott.
\newblock Computation of geometric partial differential equations and mean
  curvature flow.
\newblock {\em Acta Numer.}, 14:139--232, 2005.

\bibitem{delpino-etal}
M.~del Pino, M.~Kowalczyk, F.~Pacard, and J.~Wei.
\newblock Multiple-end solutions to the {A}llen-{C}ahn equation in {$\mathbb
  R^2$}.
\newblock {\em J. Funct. Anal.}, 258(2):458--503, 2010.

\bibitem{Droske04alevel}
M.~Droske and M.~Rumpf.
\newblock A level set formulation for {W}illmore flow.
\newblock {\em Interfaces and free boundaries}, pages 361--378, 2004.

\bibitem{Du_Wang_2004}
Q.~Du, C.~Liu, and X.~Wang.
\newblock A phase field approach in the numerical study of the elastic bending
  energy for vesicle membranes.
\newblock {\em J. Comput. Phys.}, 198(2):450--468, Aug. 2004.

\bibitem{Du_wang:2006}
Q.~Du, C.~Liu, and X.~Wang.
\newblock Simulating the deformation of vesicle membranes under elastic bending
  energy in three dimensions.
\newblock {\em J. Comput. Phys.}, 212(2):757--777, Mar. 2006.

\bibitem{Du_convergenceof}
Q.~Du and X.~Wang.
\newblock Convergence of numerical approximations to a phase field bending
  elasticity model of membrane deformations.
\newblock {\em Internat. J. Numer. Anal Modeling}, 4:441--459, 2007.

\bibitem{Dziuk:2008}
G.~Dziuk.
\newblock Computational parametric {W}illmore flow.
\newblock {\em Numer. Math.}, 111(1):55--80, Oct. 2008.

\bibitem{Existence_curve}
G.~Dziuk, E.~Kuwert, and R.~Sch{\"a}tzle.
\newblock Evolution of elastic curves in {${\mathbb R^n}$}: Existence and
  computation.
\newblock {\em SIAM J. Math. Anal.}, 33(5):1228--1245, 2002.

\bibitem{Esedoglu_p}
S.~{Esedoglu}.
\newblock Unpublished work, 2009.

\bibitem{Esedoglu_12}
S.~{Esedoglu}, A.~{R{\"a}tz}, and M.~{R{\"o}ger}.
\newblock {Colliding Interfaces in Old and New Diffuse-interface Approximations
  of Willmore-flow}.
\newblock {\em ArXiv e-prints}, Sept. 2012.

\bibitem{BMO_willmore}
S.~Esedoglu, S.~J. Ruuth, and R.~Tsai.
\newblock Threshold dynamics for high order geometric motions.
\newblock {\em Interfaces and Free Boundaries}, 10(3):263--282, 2008.

\bibitem{Evans_spruck}
L.~C. Evans and J.~Spruck.
\newblock Motion of level sets by mean curvature. {I}.
\newblock {\em J. Differential Geom.}, 33(3):635--681, 1991.

\bibitem{MR1772733}
P.~C. Fife.
\newblock Models for phase separation and their mathematics.
\newblock {\em Electron. J. Differential Equations}, pages No. 48, 26 pp.
  (electronic), 2000.

\bibitem{FrRuWi11}
M.~Franken, M.~Rumpf, and B.~Wirth.
\newblock A phase field based {PDE} constraint optimization approach to time
  discrete {W}illmore flow.
\newblock {\em International Journal of Numerical Analysis and Modeling}, 2011.
\newblock accepted.

\bibitem{GT}
D.~Gilbarg and N.~Trudinger.
\newblock {\em Elliptic Partial Differential Equations of Second Order}.
\newblock Springer, 1998.

\bibitem{gui}
C.~Gui.
\newblock Hamiltonian identities for elliptic partial differential equations.
\newblock {\em J. Funct. Anal.}, 254(4):904--933, 2008.

\bibitem{hart-wintn}
P.~Hartman and A.~Wintner.
\newblock On the local behavior of solutions of non-parabolic partial
  differential equations.
\newblock {\em Amer. J. Math.}, 75:449--476, 1953.

\bibitem{HsuKusnerSullivan92}
L.~Hsu, R.~Kusner, and J.~Sullivan.
\newblock Minimizing the squared mean curvature integral for surfaces in space
  forms.
\newblock {\em Experimental Mathematics}, 1(3):191--207, 1992.

\bibitem{Kowalczyk_pacard}
M.~Kowalczyk, Y.~Liu, and F.~Pacard.
\newblock Four ended solutions to the {A}llen-{C}ahn equation on the plane.
\newblock {\em Ann. Inst. H. Poincar{\'e} Anal. Non Lin{\'e}aire},
  29(5):761--781, 2012.

\bibitem{Kusner}
R.~Kusner.
\newblock Comparison surfaces for the {W}illmore problem.
\newblock {\em Pacific J. Math.}, 138(2), 1989.

\bibitem{KuScha04}
E.~Kuwert and R.~Sch{\"a}tzle.
\newblock Removability of point singularities of {W}illmore surfaces.
\newblock {\em Ann. of Math. (2)}, 160(1):315--357, 2004.

\bibitem{Existence_flow_sphere}
E.~K. Kuwert and R.~Sch{\"a}tzle.
\newblock The {W}illmore flow with small initial energy.
\newblock {\em Differential Geom}, 57(3):409--441, 2001.

\bibitem{Langer_Singer}
J.~Langer and D.~A. Singer.
\newblock {Curve straightening and a minimax argument for closed elastic
  curves}.
\newblock {\em Topology}, 24:75--88, 1985.

\bibitem{LM}
G.~Leonardi and S.~Masnou.
\newblock Locality of the mean curvature of rectifiable varifolds.
\newblock {\em Adv. Calc. of Var.}, 2(1):17--42, 2009.

\bibitem{Loreti_march}
P.~Loreti and R.~March.
\newblock Propagation of fronts in a nonlinear fourth order equation.
\newblock {\em European Journal of Applied Mathematics}, 11:203--213, 3 2000.

\bibitem{Lowengrub}
J.~Lowengrub, A.~R\"{a}tz, and A.~Voigt.
\newblock {Phase-field modeling of the dynamics of multicomponent vesicles:
  Spinodal decomposition, coarsening, budding, and fission}.
\newblock {\em Physical Review E (Statistical, Nonlinear, and Soft Matter
  Physics)}, 79(3):031926+, 2009.

\bibitem{Willmore_conjecture}
F.~C. Marques and A.~Neves.
\newblock Min-{M}ax theory and the {W}illmore conjecture.
\newblock {\em arXiv:1202.6036}, 2012.

\bibitem{MN1}
S.~Masnou and G.~Nardi.
\newblock A coarea-type formula for the relaxation of a generalized {W}illmore
  functional.
\newblock {\em J. Convex Analysis, to appear}, 2013.

\bibitem{MN2}
S.~Masnou and G.~Nardi.
\newblock Gradient {Y}oung measures, varifolds, and a generalized {W}illmore
  functional.
\newblock {\em Adv. Calc. Var., to appear}, 2013.

\bibitem{Mayer01anumerical}
U.~F. Mayer and G.~Simonett.
\newblock A numerical scheme for axisymmetric solutions of curvature driven
  free boundary problems, with applications to the {W}illmore flow.
\newblock {\em Interfaces Free Bound}, 4(1):89--109, 2002.

\bibitem{Modica1977a}
L.~Modica and S.~Mortola.
\newblock Il limite nella {$\Gamma $}-convergenza di una famiglia di funzionali
  ellittici.
\newblock {\em Boll. Un. Mat. Ital. A (5)}, 14(3):526--529, 1977.

\bibitem{Modica1977}
L.~Modica and S.~Mortola.
\newblock Un esempio di {$\Gamma-$}-convergenza.
\newblock {\em Boll. Un. Mat. Ital. B (5)}, 14(1):285--299, 1977.

\bibitem{moser}
R.~Moser.
\newblock A higher order asymptotic problem related to phase transitions.
\newblock {\em SIAM J. Math. Analysis}, 37(3):712--736, 2005.

\bibitem{mugnai2010}
L.~{Mugnai}.
\newblock {Gamma-convergence results for phase-field approximations of the
  2D-Euler elastica functional}.
\newblock {\em ArXiv e-prints}, Sept. 2010.

\bibitem{naga-tone}
Y.~Nagase and Y.~Tonegawa.
\newblock {A singular perturbation problem with integral curvature bound}.
\newblock {\em Hiroshima Mathematical Journal}, 37:455--489, 2007.

\bibitem{OlRu10}
N.~Olischl{\"{a}}ger and M.~Rumpf.
\newblock A nested variational time discretization for parametric {W}illmore
  flow.
\newblock {\em Interfaces and Free Boundaries}, 2011.
\newblock In Press.

\bibitem{Osherbook1}
S.~Osher and R.~Fedkiw.
\newblock {\em Level Set Methods and Dynamic Implicit Surfaces}.
\newblock Springer-Verlag New York, Applied Mathematical Sciences, 2002.

\bibitem{Osherbook2}
S.~Osher and N.~Paragios.
\newblock {\em Geometric Level Set Methods in Imaging, Vision and Graphics}.
\newblock Springer-Verlag, New York, 2003.

\bibitem{OsherSethian}
S.~Osher and J.~A. Sethian.
\newblock Fronts propagating with curvature-dependent speed: algorithms based
  on hamilton-jacobi formulations.
\newblock {\em J. Comput. Phys.}, 79:12--49, 1988.

\bibitem{Paolini_anisotropic}
M.~Paolini.
\newblock An efficient algorithm for computing anisotropic evolution by mean
  curvature.
\newblock In {\em Curvature flows and related topics ({L}evico, 1994)},
  volume~5 of {\em GAKUTO Internat. Ser. Math. Sci. Appl.}, pages 199--213.
  Gakk\=otosho, Tokyo, 1995.

\bibitem{pego}
R.~L. Pego.
\newblock Front migration in the nonlinear {C}ahn-{H}illiard equation.
\newblock {\em Proc. Roy. Soc. London Ser. A}, 422(1863):261--278, 1989.

\bibitem{Polthier-Hab}
K.~Polthier.
\newblock {\em Polyhedral surfaces of constant mean curvature}.
\newblock Habilitation thesis, TU Berlin, 2002.

\bibitem{Roger_schatzle_2006}
M.~R{\"o}ger and R.~Sch{\"a}tzle.
\newblock On a modified conjecture of {D}e {G}iorgi.
\newblock {\em Math. Z.}, 254(4):675--714, 2006.

\bibitem{rowlinson}
J.~Rowlinson.
\newblock Translation of {J}. {D}. van der {W}aals' "{T}he thermodynamik theory
  of capillarity under the hypothesis of a continuous variation of density".
\newblock {\em Journal of Statistical Physics}, 20(2):197--200, 1979.

\bibitem{Rusu_elastic_flow}
R.~E. Rusu.
\newblock An algorithm for the elastic flow of surfaces.
\newblock {\em Interfaces Free Bound.}, 3:229--229, 2005.

\bibitem{Serfaty}
S.~Serfaty.
\newblock Gamma-convergence of gradient flows on hilbert and metric spaces and
  applications.
\newblock {\em Disc. Cont. Dyn. Systems,}, 31(4):1427--1451, 2011.

\bibitem{tone02}
Y.~Tonegawa.
\newblock Phase field model with a variable chemical potential.
\newblock {\em Proceedings of the Royal Society of Edinburgh: Section A
  Mathematics}, 132:993--1019, 7 2002.

\bibitem{Wardetzky-et-al-07}
M.~Wardetzky, M.~Bergou, D.~Harmon, D.~Zorin, and E.~Grinspun.
\newblock Discrete quadratic curvature energies.
\newblock {\em Computer Aided Geometric Design}, 24(8--9):499--518, 2007.

\end{thebibliography}

\end{document}